\documentclass[10pt, reqno]{amsart}
\usepackage{amsfonts}
\usepackage{amsmath}
\usepackage{amsthm}
\usepackage{amssymb}  
\usepackage{mathrsfs} 
\usepackage{bm} 
\usepackage[pagebackref]{hyperref}
\hypersetup{colorlinks, linkcolor=[RGB]{30,70,100}, anchorcolor=green, citecolor=[RGB]{140,90,90}}
\usepackage{tikz-cd}

\usepackage{comment}

\usepackage[top=2.5cm, bottom=2.5cm,left=3cm,right=3cm]{geometry}

\newcommand{\bbA}{\mathbb{A}}

\newcommand{\bmzero}{\bm{\mathrm{0}}}
\newcommand{\bmone}{\bm{\mathrm{1}}}
\newcommand{\bma}{\bm{\mathrm{a}}}
\newcommand{\bmA}{\bm{\mathrm{A}}}
\newcommand{\bmb}{\bm{\mathrm{b}}}
\newcommand{\bmB}{\bm{\mathrm{B}}}
\newcommand{\bmc}{\bm{\mathrm{c}}}
\newcommand{\bmC}{\bm{\mathrm{C}}}
\newcommand{\bmd}{\bm{\mathrm{d}}}
\newcommand{\bmD}{\bm{\mathrm{D}}}
\newcommand{\bme}{\bm{\mathrm{e}}}
\newcommand{\bmE}{\bm{\mathrm{E}}}
\newcommand{\bmf}{\bm{\mathrm{f}}}
\newcommand{\bmF}{\bm{\mathrm{F}}}
\newcommand{\bmG}{\bm{\mathrm{G}}}
\newcommand{\bmH}{\bm{\mathrm{H}}}
\newcommand{\bmK}{\bm{\mathrm{K}}}
\newcommand{\bml}{\bm{\mathrm{l}}}
\newcommand{\bmL}{\bm{\mathrm{L}}}
\newcommand{\bmM}{\bm{\mathrm{M}}}
\newcommand{\bmn}{\bm{\mathrm{n}}}
\newcommand{\bmN}{\bm{\mathrm{N}}}
\newcommand{\bmO}{\bm{\mathrm{O}}}
\newcommand{\bmp}{\bm{\mathrm{p}}}
\newcommand{\bmP}{\bm{\mathrm{P}}}
\newcommand{\bmQ}{\bm{\mathrm{Q}}}
\newcommand{\bmR}{\bm{\mathrm{R}}}
\newcommand{\bms}{\bm{\mathrm{s}}}
\newcommand{\bmS}{\bm{\mathrm{S}}}
\newcommand{\bmt}{\bm{\mathrm{t}}}
\newcommand{\bmT}{\bm{\mathrm{T}}}
\newcommand{\bmu}{\bm{\mathrm{u}}}
\newcommand{\bmU}{\bm{\mathrm{U}}}
\newcommand{\bmv}{\bm{\mathrm{v}}}
\newcommand{\bmV}{\bm{\mathrm{V}}}
\newcommand{\bmw}{\bm{\mathrm{w}}}
\newcommand{\bmW}{\bm{\mathrm{W}}}
\newcommand{\bmx}{\bm{\mathrm{x}}}
\newcommand{\bmX}{\bm{\mathrm{X}}}
\newcommand{\bmy}{\bm{\mathrm{y}}}
\newcommand{\bmY}{\bm{\mathrm{Y}}}
\newcommand{\bmz}{\bm{\mathrm{z}}}
\newcommand{\bmZ}{\bm{\mathrm{Z}}}

\newcommand{\bmalpha}{\bm{\alpha}}
\newcommand{\bmbeta}{\bm{\beta}}
\newcommand{\bmgamma}{\bm{\gamma}}
\newcommand{\bmGr}{\bm{\mathrm{Gr}}}
\newcommand{\bmiota}{\bm{\iota}}
\newcommand{\bmtheta}{\bm{\theta}}

\newcommand{\calA}{\mathcal{A}}
\newcommand{\calB}{\mathcal{B}}
\newcommand{\calC}{\mathcal{C}}
\newcommand{\calD}{\mathcal{D}}
\newcommand{\calE}{\mathcal{E}}
\newcommand{\calF}{\mathcal{F}}
\newcommand{\calG}{\mathcal{G}}
\newcommand{\calH}{\mathcal{H}}
\newcommand{\calI}{\mathcal{I}}
\newcommand{\wtcalI}{\widetilde{\mathcal{I}}}
\newcommand{\calK}{\mathcal{K}}
\newcommand{\calL}{\mathcal{L}}
\newcommand{\calM}{\mathcal{M}}
\newcommand{\calN}{\mathcal{N}}
\newcommand{\calO}{\mathcal{O}}
\newcommand{\calP}{\mathcal{P}}
\newcommand{\calR}{\mathcal{R}}
\newcommand{\calS}{\mathcal{S}}
\newcommand{\calT}{\mathcal{T}}
\newcommand{\calU}{\mathcal{U}}
\newcommand{\calV}{\mathcal{V}}
\newcommand{\calW}{\mathcal{W}}
\newcommand{\calX}{\mathcal{X}}

\DeclareMathOperator{\calMat}{\mathcal{M}\bm{\mathrm{at}}}

\newcommand{\diff}{\mathrm{d}}
\newcommand{\diffalpha}{\mathrm{d\alpha}}
\newcommand{\diffbeta}{\mathrm{d\beta}}
\newcommand{\diffgamma}{\mathrm{d\gamma}}
\newcommand{\difflambda}{\mathrm{d\lambda}}
\newcommand{\difftheta}{\mathrm{d\theta}}
\newcommand{\diffa}{\mathrm{da}}
\newcommand{\diffb}{\mathrm{db}}
\newcommand{\diffc}{\mathrm{dc}}
\newcommand{\difff}{\mathrm{df}}
\newcommand{\diffs}{\mathrm{ds}}
\newcommand{\difft}{\mathrm{dt}}
\newcommand{\diffu}{\mathrm{du}}
\newcommand{\diffv}{\mathrm{dv}}
\newcommand{\diffw}{\mathrm{dw}}
\newcommand{\diffx}{\mathrm{dx}}
\newcommand{\diffy}{\mathrm{dy}}

\newcommand{\fraka}{\mathfrak{a}}
\newcommand{\frakf}{\mathfrak{f}}
\newcommand{\frakg}{\mathfrak{g}}
\newcommand{\frakh}{\mathfrak{h}}
\newcommand{\frakl}{\mathfrak{l}}
\newcommand{\frakm}{\mathfrak{m}}
\newcommand{\frakn}{\mathfrak{n}}
\newcommand{\frako}{\mathfrak{o}}
\newcommand{\frakp}{\mathfrak{p}}
\newcommand{\frakr}{\mathfrak{r}}
\newcommand{\fraks}{\mathfrak{s}}
\newcommand{\frakS}{\mathfrak{S}}
\newcommand{\frakt}{\mathfrak{t}}
\newcommand{\fraku}{\mathfrak{u}}
\newcommand{\frakU}{\mathfrak{U}}
\newcommand{\frakw}{\mathfrak{w}}
\newcommand{\frakX}{\mathfrak{X}}
\newcommand{\frakz}{\mathfrak{z}}

\newcommand{\olp}{\overline{p}}
\newcommand{\olV}{\overline{V}}
\newcommand{\olDelta}{\overline{\Delta}}
\newcommand{\olpi}{\overline{\pi}}

\newcommand{\rmA}{\mathrm{A}}
\newcommand{\rmB}{\mathrm{B}}
\newcommand{\rmC}{\mathrm{C}}
\newcommand{\rmD}{\mathrm{D}}
\newcommand{\rmf}{\mathrm{f}}
\newcommand{\rmF}{\mathrm{F}}
\newcommand{\rmG}{\mathrm{G}}
\newcommand{\rmH}{\mathrm{H}}
\newcommand{\rmK}{\mathrm{K}}
\newcommand{\rmL}{\mathrm{L}}
\newcommand{\rmm}{\mathrm{m}}
\newcommand{\whrmm}{\widehat{\mathrm{m}}}
\newcommand{\rmM}{\mathrm{M}}
\newcommand{\rmN}{\mathrm{N}}
\newcommand{\rmP}{\mathrm{P}}
\newcommand{\rmQ}{\mathrm{Q}}
\newcommand{\rmR}{\mathrm{R}}
\newcommand{\rmS}{\mathrm{S}}
\newcommand{\rmT}{\mathrm{T}}
\newcommand{\rmu}{\mathrm{u}}
\newcommand{\rmU}{\mathrm{U}}
\newcommand{\rmX}{\mathrm{X}}
\newcommand{\rmY}{\mathrm{Y}}
\newcommand{\rmZ}{\mathrm{Z}}

\newcommand{\scrA}{\mathscr{A}}
\newcommand{\scrB}{\mathscr{B}}
\newcommand{\scrC}{\mathscr{C}}
\newcommand{\scrD}{\mathscr{D}}
\newcommand{\scrF}{\mathscr{F}}
\newcommand{\scrG}{\mathscr{G}}
\newcommand{\scrI}{\mathscr{I}}
\newcommand{\scrJ}{\mathscr{J}}
\newcommand{\scrL}{\mathscr{L}}
\newcommand{\scrK}{\mathscr{K}}
\newcommand{\scrP}{\mathscr{P}}
\newcommand{\scrT}{\mathscr{T}}
\newcommand{\scrU}{\mathscr{U}}
\newcommand{\scrV}{\mathscr{V}}
\newcommand{\scrX}{\mathscr{X}}

\newcommand{\wta}{\widetilde{a}}
\newcommand{\wtb}{\widetilde{b}}
\newcommand{\wtB}{\widetilde{B}}
\newcommand{\wtrmB}{\widetilde{\mathrm{B}}}
\newcommand{\wtbmB}{\widetilde{\bm{\mathrm{B}}}}
\newcommand{\wtrmD}{\widetilde{\mathrm{D}}}
\newcommand{\wtbmD}{\widetilde{\bm{\mathrm{D}}}}
\newcommand{\wtbmF}{\widetilde{\bm{\mathrm{F}}}}
\newcommand{\wtbmG}{\widetilde{\bm{\mathrm{G}}}}
\newcommand{\wtrmG}{\widetilde{\mathrm{G}}}
\newcommand{\wtbmH}{\widetilde{\bm{\mathrm{H}}}}
\newcommand{\wtrmH}{\widetilde{\mathrm{H}}}
\newcommand{\wtcalL}{\widetilde{{\mathcal{L}}}}
\newcommand{\wts}{\widetilde{s}}
\newcommand{\wtbmS}{\widetilde{\bm{\mathrm{S}}}}
\newcommand{\wtt}{\widetilde{t}}
\newcommand{\wtw}{\widetilde{w}}
\newcommand{\wtx}{\widetilde{x}}
\newcommand{\wtX}{\widetilde{X}}
\newcommand{\wtrmX}{\widetilde{\mathrm{X}}}
\newcommand{\wtbmX}{\widetilde{\bm{\mathrm{X}}}}
\newcommand{\wty}{\widetilde{y}}
\newcommand{\wtty}{\widetilde{\widetilde{y}}}
\newcommand{\wttty}{\widetilde{\widetilde{\widetilde{y}}}}
\newcommand{\wtz}{\widetilde{z}}

\newcommand{\wtbeta}{\widetilde{\beta}}
\newcommand{\wttbeta}{\widetilde{\widetilde{\beta}}}
\newcommand{\wtttbeta}{\widetilde{\widetilde{\widetilde{\beta}}}}
\newcommand{\wtDelta}{\widetilde{\Delta}}
\newcommand{\wtGamma}{\widetilde{\Gamma}}
\newcommand{\wtlambda}{\widetilde{\lambda}}
\newcommand{\wttlambda}{\widetilde{\widetilde{\lambda}}}
\newcommand{\wtttlambda}{\widetilde{\widetilde{\widetilde{\lambda}}}}
\newcommand{\wtmu}{\widetilde{\mu}}
\newcommand{\wtnu}{\widetilde{\nu}}
\newcommand{\wtomega}{\widetilde{\omega}}
\newcommand{\wtPhi}{\widetilde{\Phi}}

\newcommand{\Z}{\mathbb{Z}}
\newcommand{\Q}{\mathbb{Q}}
\newcommand{\R}{\mathbb{R}}
\newcommand{\C}{\mathbb{C}}

\DeclareMathOperator{\SL}{\bm{\mathrm{SL}}}
\DeclareMathOperator{\GL}{\bm{\mathrm{GL}}}
\DeclareMathOperator{\SU}{\bm{\mathrm{SU}}}
\DeclareMathOperator{\SO}{\bm{\mathrm{SO}}}
\DeclareMathOperator{\Sp}{\bm{\mathrm{Sp}}}
\newcommand{\fraksl}{\mathfrak{sl}}
\newcommand{\frakgl}{\mathfrak{gl}}
\newcommand{\fraksu}{\mathfrak{su}}
\newcommand{\frakso}{\mathfrak{so}}
\newcommand{\fraksp}{\mathfrak{sp}}
\newcommand{\ep}{\varepsilon}

\newcommand{\bs}{\backslash}
\newcommand{\la}{\langle}
\newcommand{\ra}{\rangle}
\newcommand{\midd}{\;\middle\vert\;}

\newcommand{\normm}[1]{\left\lvert#1\right\rvert}
\newcommand{\norm}[1]{\left\lVert#1\right\rVert}
\newcommand{\udl}[1]{\underline{#1}}

\newcommand{\acts}{\curvearrowright} 
\newcommand{\embed}{\hookrightarrow}

\newtheorem{thm}{Theorem}[section]

\newtheorem{coro}[thm]{Corollary}
\newtheorem{defi}[thm]{Definition}
\newtheorem{exam}[thm]{Example}
\newtheorem{lem}[thm]{Lemma}
\newtheorem{prop}[thm]{Proposition}
\newtheorem{ques}[thm]{Question}
\newtheorem{rmk}[thm]{Remark}
\newtheorem{assumption}[thm]{Assumption}

\renewcommand*{\thedefinition}{\Alph{definition}} 
\newtheorem*{claim}{Claim}

\renewcommand*{\theobservation}{\Alph{observation}}

\DeclareMathOperator{\an}{\mathrm{an}}
\DeclareMathOperator{\ari}{\mathrm{ari}}
\DeclareMathOperator{\Ad}{\mathrm{Ad}}
\DeclareMathOperator{\Bary}{\mathrm{Bary}}
\DeclareMathOperator{\bdd}{\mathrm{bdd}}
\DeclareMathOperator{\Bad}{\mathrm{Bad}}
\DeclareMathOperator{\Bl}{\bm{\mathrm{Bl}}}
\DeclareMathOperator{\can}{\mathrm{can}}
\DeclareMathOperator{\conv}{\mathrm{conv}}
\DeclareMathOperator{\conn}{\mathrm{conn}}
\DeclareMathOperator{\Cone}{\mathrm{Cone}}
\DeclareMathOperator{\cor}{\mathrm{cor}}
\DeclareMathOperator{\cpt}{\mathrm{cpt}}
\DeclareMathOperator{\cvg}{\mathrm{cvg}}
\DeclareMathOperator{\CVX}{\mathrm{CVX}}
\DeclareMathOperator{\diag}{\mathrm{diag}}
\DeclareMathOperator{\divisor}{\mathrm{div}}
\DeclareMathOperator{\Extre}{\mathrm{Extre}}
\DeclareMathOperator{\Gal}{\mathrm{Gal}}
\DeclareMathOperator{\Good}{\mathrm{Good}}
\DeclareMathOperator{\Ht}{\mathrm{Ht}}
\DeclareMathOperator{\id}{\mathrm{id}}
\DeclareMathOperator{\INT}{\mathrm{INT}}
\DeclareMathOperator{\INTP}{\mathrm{INTP}}
\DeclareMathOperator{\IMTP}{\mathrm{IMTP}}
\DeclareMathOperator{\Leb}{\mathrm{Leb}}
\DeclareMathOperator{\LHS}{\mathrm{LHS}}
\DeclareMathOperator{\Lie}{\mathrm{Lie}}
\DeclareMathOperator{\loc}{\mathrm{loc}}
\DeclareMathOperator{\Map}{\mathrm{Map}}
\DeclareMathOperator{\meas}{\mathrm{meas}}
\DeclareMathOperator{\Meas}{\mathrm{Meas}}
\DeclareMathOperator{\Mor}{\mathrm{Mor}}
\DeclareMathOperator{\nc}{\mathrm{nc}}
\DeclareMathOperator{\new}{\mathrm{new}}
\DeclareMathOperator{\obs}{\mathrm{obs}}
\DeclareMathOperator{\Pic}{\mathrm{Pic}}
\DeclareMathOperator{\Plot}{\mathrm{Plot}}
\DeclareMathOperator{\Pole}{\mathrm{Pole}}
\DeclareMathOperator{\Poly}{\mathrm{Poly}}
\DeclareMathOperator{\Profile}{\mathrm{Profile}}
\DeclareMathOperator{\Prob}{\mathrm{Prob}}
\DeclareMathOperator{\Prim}{\mathrm{Prim}}
\DeclareMathOperator{\rank}{\mathrm{rank}}
\DeclareMathOperator{\rk}{\mathrm{rk}}
\DeclareMathOperator{\Rad}{\mathrm{Rad}}
\DeclareMathOperator{\red}{\mathrm{red}}
\DeclareMathOperator{\Res}{\mathrm{Res}}
\DeclareMathOperator{\semisimple}{\mathrm{ss}}
\DeclareMathOperator{\sign}{\mathrm{sign}}
\DeclareMathOperator{\slope}{\mathrm{slope}}
\DeclareMathOperator{\sma}{\mathrm{sma}}
\DeclareMathOperator{\Spec}{\mathrm{Spec}}
\DeclareMathOperator{\spl}{\mathrm{spl}}
\DeclareMathOperator{\Stab}{\mathrm{Stab}}
\DeclareMathOperator{\sta}{\mathrm{sta}}
\DeclareMathOperator{\stau}{\mathrm{stau}}
\DeclareMathOperator{\Std}{\mathrm{Std}}
\DeclareMathOperator{\supp}{\mathrm{supp}}
\DeclareMathOperator{\sys}{\mathrm{sys}}
\DeclareMathOperator{\Tr}{\mathrm{Tr}}
\DeclareMathOperator{\tr}{\mathrm{tr}}
\DeclareMathOperator{\type}{\mathrm{type}}
\DeclareMathOperator{\Val}{\mathrm{Val}}
\DeclareMathOperator{\Vol}{\mathrm{Vol}}
\DeclareMathOperator{\Zero}{\mathrm{Zero}}

\DeclareMathOperator{\Hom}{{\mathrm{Hom}}}
\DeclareMathOperator{\INC}{{\mathrm{INC}}}
\DeclareMathOperator{\Gr}{\bm{\mathrm{Gr}}}
\DeclareMathOperator{\Mat}{\bm{\mathrm{Mat}}}

\begin{document}

\title[]{Asymptotics of integral points, equivariant compactifications and equidistributions for homogeneous spaces}
\author{Runlin Zhang}

\address{College of Mathematics and Statistics, Center of Mathematics, Chongqing University, 401331, Chongqing,  China. }
\email{zhangrunlin@outlook.com}

\begin{abstract}
    Let $\bmU$ be a homogeneous variety over $\Q$ of a linear algebraic group. Let $\calU$ be an integral model and assume the existence of infinitely many integral points. Then one would like to give an asymptotic count of integral points of bounded height with the help of some height function.  In many cases, with the help of measure rigidity of unipotent flows, we reduce this problem to one on equivariant birational geometry. For instance, we show that if $\bmG$ and $\bmH$ are both connected, semisimple, simply connected and without compact factors,  then $\bmG/\bmH$ is strongly Hardy-Littlewood with respect to some height function. We also show that when $\bmH$ is ``large'' in $\bmG$ and both $\bmG$ and $\bmH$ are connected, reductive and without nontrivial $\Q$-characters, the asymptotic of integral points is the same as the volume asymptotic up to a constant for every equivariant height. Three concrete examples with explicit heights are also provided to illustrate our approach.
\end{abstract}

\maketitle

\tableofcontents
\addtocontents{toc}{\protect\setcounter{tocdepth}{1}}

\section{Introduction}

Given a set of polynomials with integral coefficients, Diophantine problems concern the existence of some integral solution or infinitely many solutions. 
Assuming there are infinitely many solutions, however, one can still ask for an asymptotic count provided there is a height function.

In this paper we consider the special case when the equations have rich symmetries: their complex solutions allow a transitive action of some linear algebraic group. 
There are many works on when an integral solution exists in this context (see \cite{BoroRudn95, Colliot-Thelene_Xu_2009_Brauer_Manin, Ellenberg_Venkatesh_2008_local_global_quadratic} for a small sample of research). 
Assuming that integral solutions do exist,
 the question remained is to count them and we are mainly interested in the case when there are infinitely many.
It turns out that the set of integral solutions has a nice structure: it is acted on by an arithmetic subgroup and decomposes into finitely many orbits.
So naturally one wants to find the asymptotic count of a single orbit first.
It was observed by \cite{DukRudSar93} that such a question is related to equidistribution of homogeneous measures on a special type of homogeneous spaces: arithmetic quotients of real points of a linear algebraic groups. Ergodic theoretic methods then kick in and solve the problem in many cases \cite{EskMcM93, EskMozSha96}.
It is also realized that the height function one use to count the solutions matters.

The purpose of the present paper is to continue the study in \cite{EskMozSha96,zhangrunlinCompositio2021} via unipotent flows, that is, the celebrated work of Ratner \cite{Rat91}, and reduces the counting problem in many cases to geometric questions about the height function.
In the remaining part of this introduction, we will present some theorems as well as problems one encounters in this process.

As far as counting integral points on homogeneous varieties is concerned, there are many existing works. Besides those already mentioned, some relevant ones are \cite{BenOh12,  ChamTschin12, Chow_2024_integral_wonderful, GolMoha14, GoroNevo12Crelle, GorOhSha09, Hassett_Tschinkel_2003_integral_points_effective_cones_moduli_stable_maps, KelKon18, Maucou07,  OhSha14, Shah00, ShaZhe18, TaklooTschinkel2013, zhangrunlinAnnalen2019}, touching on different aspects of the question.

\subsection{Hardy--Littlewood varieties}

The notion of strongly Hardy--Littlewood varieties is introduced in \cite{BoroRudn95} (compare \cite{Birch_1961_forms_in_many_variables}) for affine varieties. 
Such varieties are defined by certain nice local-to-global properties. 
To simplify discussions, we specialize to the case when $\bmU$ is isomorphic to $\bmG/\bmH$ for some semisimple linear algebraic group $\bmG$ and a semisimple $\Q$-subgroup $\bmH$. In particular, we assume that $\bmU(\Q)\neq \emptyset$.

\subsubsection{Integral points}
Let $\bmU$ be an affine homogeneous variety over $\Q$ and $\bmU(\bbA_{\Q})$ be the topological space of adelic points. Given an open compact subset $\scrK_f \subset \bmU(\bbA_{\Q,f})$  of finite adeles and a connected component $\scrK_{\infty}$ of $\bmU(\R)$, one define the set of $\scrK$-integral points $\bmU_{\scrB}(\Z)$ to be the intersection of $\scrK:=\scrK_{\infty} \times \scrK_{f}$ with the diagonal embedding to $\bmU(\Q)$ in the set of adeles.

\subsubsection{Tamagawa measures}
Let $\omega$ be a nonvanishing top-degree differential form on $\bmU$ over $\Q$. It is actually $\bmG$-invariant and induces measures $\omega_p$ on each $\bmU(\Q_p)$ ($\Q_{\infty}:=\R$) and, in fact, the products over finite $p$'s induce a measure $\omega_f$ on $\bmU(\bbA_{\Q,f})$. 
By the product formula and the triviality of global regular invertible functions on $\bmU$, the measure $\omega_f \otimes \omega_{\infty}$ is independent of the choice of $\omega$.

\subsubsection{The definition}
Let $l: \bmU(\R) \to \R_{>0}$ be a proper function and
\begin{equation*}
    B_R:= \left\{
x\in \bmU(\R) ,\;l(x )\leq R
\right\}.
\end{equation*}
We say that $\bmU$ is a strongly Hardy-Littlewood variety (with respect to $l$) if for any choice of $\scrK$ as above, 
\begin{equation*}
    \lim_{R\to \infty}
    \frac{
      \# \bmU_{\scrK}(\Z) \cap B_R 
      }
    {\omega_f\otimes\omega_{\infty} (B_f \times B_R) } =1.
\end{equation*}
In \cite{BoroRudn95}, $l$ is induced from an embedding of $\bmU$ into some affine space $\bmA^n$. 

\subsubsection{Theorems}

Now, assume further that $\bmG$ and $\bmH$ are both connected and simply connected. Also assume that the real points of every nontrivial $\Q$-factor of $\bmG$ is noncompact. Note that we have assumed $\bmU(\Q)\neq \emptyset$.
The following is \cite[Theorem 0.3]{BoroRudn95}:
\begin{thm}
    Assume that for any congruence subgroup $\Gamma \subset \bmG(\Q)$ and any rational point $x\in \bmU(\Q)$, one has
    \begin{equation}\label{equation_orbital_counting}
        \lim_{R\to \infty}
        \frac{
        \# \Gamma.x \cap B_R
        }{ \omega_{\infty}(B_R \cap \bmG(\R).x ) }
        = \frac{\Vol (\bmH_x(\R)/ \bmH_x(\R)\cap \Gamma )}{
        \Vol (\bmG(\R)/\Gamma)
        }
    \end{equation}
    where $\bmH_x$ denotes the stabilizer of $x$ in $\bmG$ and volumes are induced from compatible choices of Haar measures on $\bmG(\R),\bmH_x(\R)$ and $\bmG(\R)/\bmH_x(\R)$.
    Then $\bmU$ is strongly Hardy--Littlewood with respect to $l$.
\end{thm}

We remark that 
if $\Gamma'\subset \Gamma$ is a finite-index subgroup, then the truth of Equa.(\ref{equation_orbital_counting}) for $\Gamma'$ implies that for $\Gamma$.
Equa.(\ref{equation_orbital_counting}) has been verified in many cases. For instance, when $\bmH$ is maximal in a semisimple $\bmG$ without compact $\Q$-factors, all ``algebraic height functions'' $l$ would work.
What is proved in this paper (see Section \ref{subsection_proof_hardy_littlewood}) is that
\begin{thm}\label{theorem_Hardy_Littlewood}
    Assumptions same as above. There exists $l$ such that Equa.(\ref{equation_orbital_counting}) holds.
    In particular, $\bmU$ is always strongly Hardy-Littlewood for some $l$.
\end{thm}

Already when $\bmG=\SL_{n}$ and $\bmH= \SL_{n_1}\times...\times \SL_{n_k}$ with $\sum n_i \leq n$, this seems unknown before.
Generalizations can be obtained by  combining Theorem \ref{theorem_existence_good_height_counting}  with \cite{Wei_Xu_2016}.

\subsection{Manin conjecture for a pair}\label{subsection_introduction_manin_pair}

In the previous section, there is no emphasis on which height function $l$ one should use.

In the context of  integral points on projective varieties,  roughly speaking, Manin conjecture makes prediction on the asymptotics of integral points with respect to heights associated to the anti-canonical bundle. In the case of a pair, one should replace the anti-canonical line bundle by the log anti-canonical line bundle.

\subsubsection{Log anti-canonical line bundles}

Let $\bmX$ be a $\bmG$-equivariant compactification of $\bmU \cong \bmG/\bmH$ and $\bmD:= \bmX \setminus \bmU$ be the boundary. 
We assume that $\bmX$ is smooth and $\bmD$ is a simple normal crossing divisor. Let $\calK_{\bmX}$ be the canonical line bundle. Then $(\calK_{\bmX} \otimes  \calO_{\bmX}(\bmD) )^{\vee}$ is called the \textit{log anti-canonical} line bundle.

\subsubsection{Heights}\label{subsubsection_heights}
Without loss of generality, we only consider the height at infinity here (see Remark \ref{remark_adelic_real_height}). 
Let $\bmD=\sum_{\alpha \in \scrA} \bmD_{\alpha}$ be the decomposition of $\bmD$ into absolutely irreducible components.
Let $\bmL:= \sum \lambda_{\alpha} \bmD_{\alpha}$ with $\lambda_{\alpha} >0$ for all $\alpha \in \scrA$. View the line bundle $\calO_{\bmX}(\bmL)$ over $\R$ in the category of smooth manifolds and equip it with a smooth metric $\norm{\cdot}$. Let $\bmone_{\bmL}$ be the canonical section of $\calO_{\bmX}(\bmL)$ and define for $x\in \bmU(\R)$,
\begin{equation*}
    \Ht(x): =\norm{\bmone_{\bmL}(x)}^{-1}.
\end{equation*}
Points of height bounded by a fixed number is a compact subset of $\bmU(\R)$. 

\subsubsection{Volume}

Let $\bmD':= \sum d_{\alpha}\bmD_{\alpha}$ be a divisor supported on $\bmD$. It is not necessary to require $\bmD'$ to be effective or $d_{\alpha} \geq 0$, so $\bmone_{\bmD'}$ is only a rational section that is regular on $\bmU$.
Let $\calK_{\bmX}(\bmD'):= \calK_{\bmX} \otimes \calO_{\bmX}(\bmD')$. Assume that $\calK_{\bmX}(\bmD')$ (again, viewed in the category of smooth manifolds) is equipped with a smooth metric $\norm{\cdot}$, then a
measure, denoted as $\Vol$, on $\bmU(\R)$ can be locally defined by
\begin{equation}\label{equation_local_volume}
    \Vol = 
         \frac{
         \left\vert
         \omega
         \right\vert
         }{\norm{
         \omega \otimes \bmone_{\bmD'}  }
         }
    ,\; \text{ for any nonvanishing local volume form }\omega.
\end{equation}

\subsubsection{Haar measure}
In the homogeneous setting, a natural choice of measure is the $\bmG(\R)$-invariant Haar measure.
For this we assume the existence of a nonzero $\bmG$-invariant top-degree differential form on $\bmU$. Let $\omega_0$ be such a form and $\rmm_{\bmG/\bmH}$ be the measure on $\bmU(\R)$ induced form it, then $\rmm_{\bmG/\bmH}$ is $\bmG(\R)$ invariant. 
Let $\bmD':= -\divisor (\omega_0) $, which we shall refer to as \textit{the} anticanonical divisor, then $\bmone_{\bmD'}$ is identified with $\omega_0^{\vee}$, the unique rational section of $\calO_{\bmX}(\bmD') \cong \calK_{\bmX}^{\vee}$ (the dual of $\calK_{\bmX}$) such that the natural pairing between $\omega_0$ and $\omega_0^{\vee}$ is the constant function $1$.
Define a smooth metric on $\calK_{\bmX}(\bmD')$ by imposing $\norm{ (\omega_0 \otimes \omega_{0}^{\vee} ) (x)}=1$   for every $x$. 
Then the $\Vol$ from Equa.(\ref{equation_local_volume}) coincides with $\rmm_{\bmG/\bmH}$ here.

\subsubsection{The prediction}
Fix an integral model $\calX$ of $\bmX$ and $\calD$ of $\bmD$. Then $\calU:= \calX \setminus \calD$ is an integral model of $\bmU$. Or more generally, one can take $\calU(\Z):= \bmU_{\scrK}(\Z)$ as above.

Let $(D_{\alpha'})_{\alpha' \in \scrA'}$ be consisting of connected components of $\bmD_{\alpha}(\R)$ as $\alpha$ varies in $\scrA$. The analytic Clemens complex $\scrC^{\an}_{\R}$, elements of which consist of connected components of nonempty intersection of $D_{\alpha'}$'s, encodes the intersection patterns among boundary divisors.

Depending on the integral model, let $\scrC^{\an}_{\R,\Z}$ be the modified analytic Clemens complex consisting of faces  $F \in \scrC^{\an}_{\R}$ with 
\begin{equation*}
    \left( \bigcap_{\alpha\in F} D_{\alpha} \right)
    \bigcap 
    \left( 
    \bigcup_{x\in \calU({\Z}) } \overline{\bmG(\R)^{\circ}.x}
    \right) \neq \emptyset.
\end{equation*}
That is to say, we only consider closures of connected components of $\bmU(\R)$ where there exists at least one integral point.
Let $E_{\Q}(\bmU)$ be the abelian group of invertible functions on $\bmU$ modulo constant functions $\Q^{\times}$ and
\begin{equation*}
    b := \max\left\{
        \rk \Pic_{\Q}(\bmU) - \rk E_{\Q}(\bmU) + \# F
        \;\middle\vert \;
        F \text{ is a maximal face of } \scrC^{\an}_{\R,\Z}
    \right\}.
\end{equation*}
When $\bmU(\R)$ is connected, e.g., when $\bmG$ is semisimple and simply connected and $H^1(\R,\bmH)$ is trivial, the modified analytic Clemens complex coincide with the usual analytic Clemens complex. But in general, they could be different.

Here is a natural question inspired by Manin conjecture.
\begin{ques}\label{question_log_Manin_homogeneous_variety}
    Let $\Ht$ be a height function associated with a log anti-canonical divisor.
    Assume $\calU(\Z)\neq \emptyset$.
    Does there exist a positive number $c>0$ such that 
    \begin{equation*}
        \# \left\{
             x\in \calU(\Z)\;\middle\vert\;
             \Ht(x)\leq  R
        \right\}
        \sim 
        c R (\log R)^{b-1}
    \end{equation*}
    as $R\to +\infty$?
\end{ques}
Here is what we can do on the positive side (combine Theorem \ref{theorem_log_anti_canonical_height_general_height}, Lemma \ref{lemma_log_anti_canonical_good} and Theorem \ref{theorem_anticanonical_affine_reductive}).
\begin{thm}\label{theorem_Manin_affine_finite_volume}
Assume the following:
\begin{itemize}
    \item $\bmG$ and $\bmH$ are both connected, reductive and have no nontrivial $\Q$-characters;
    \item the identity component of the centralizer of $\bmH$ in $\bmG$ is contained in $\bmH$;
    \item  the projection of $\bmH$ to the compact $\Q$-factor of $\bmG$ is surjective.
\end{itemize}
    Then the log anti-canonical divisor is big and Question \ref{question_log_Manin_homogeneous_variety} admits a positive answer.
\end{thm}

Note that under the assumptions made in the theorem, $\rk \Pic_{\Q}(\bmU) = \rk E_{\Q}(\bmU) =0$.

\subsection{Weighted counting}

In this subsection, we assume that $\bmG$ is semisimple without compact $\Q$-simple factors and $\bmH$ is reductive and connected. However, as opposed to previous discussions, we are mainly interested in the case when $\bmH(\R)\Gamma/\Gamma$ has infinite volume.
In this case, the number of integral points of bounded height is no longer expected to be asymptotic to the volume. To remedy the situation, we consider weighted versions of the counting problem.

We also assume that $\bmZ_{\bmG}(\bmH)^{\circ}$ is contained in $\bmH$ throughout this subsection.
First we discuss an analogue of Theorem \ref{theorem_Hardy_Littlewood}. Let $\omega_{\infty}$ be an invariant measure on $\bmU(\R)$ induced from some invariant form.
\begin{thm}
    There exists a height function $\Ht$ (as in Section \ref{subsubsection_heights}) such that for every $x\in \bmU(\Q)$ and every arithmetic lattice $\Gamma$ of $\bmG$, one has
    \begin{equation*}
        \sum_{ y \in \Gamma.x \cap  B_{R,x}} 
        \bmw_y
        \sim
        \omega_{\infty}(B_{R,x})
    \end{equation*}
    where $B_{R,x}:=\left\{
    y \in \bmG(\R)^{\circ}.x ,\;\Ht(y)\leq R
    \right\}$. 
\end{thm}
See Theorem \ref{theorem_existence_good_height_1}.
We expect that for ``most'' $y\in \Gamma.x \cap  B_{R,x} $, $\bmw_y$ grows like a power of $\log R$.

For a more general height, our result is much weaker.
Let $(\bmX,\bmD)$ be a smooth $\bmG$-pair over $\Q$ such that $\bmU$ is $\bmG$-equivariantly isomorphic to $\bmX\setminus \bmD$.
Let $\Ht$ be a height function associated with some divisor whose support is equal to $\bmD$.
Let $\rmG:=\bmG(\R)^{\circ}$, $\rmH_x:= \bmH_x(\R)\cap \rmG$ and assume that $\Gamma$ is contained in $\rmG$.

Let 
\begin{equation*}
    b' := \max\left\{
         \# F
        \;\middle\vert \;
        F \text{ is a maximal face of } \scrC^{\an}_{\R,\Z}
    \right\}.
\end{equation*}
The following will be proved in Theorem \ref{theorem_log_anti_canonical_height_general_height}.
\begin{thm}\label{theorem_Manin_affine_infinite_volume}
    Let $\psi$ be a non-negative compactly supported function whose support is large enough. For $x \in \bmU(\Q)$, define
    \begin{equation*}
        \bmw_x := 
        \sum_{h\in \rmH_x /\bmH_x(\R)^{\circ}}
        \la \psi, h_* \rmm_{[\bmH_x(\R)^{\circ}]}  \ra ^{-1}.
    \end{equation*}
    Then there exists $c> 0 $ such that 
    \begin{equation*}
        \sum_{x \in \calU(\Z),\,\Ht(x)\leq R }
        \bmw_x 
        \sim c \cdot R (\log R)^{b'-1}.
    \end{equation*}
\end{thm}

\subsection{Examples}
We provide three examples of different flavours. The first two examples are spherical varieties but the third one is not.

\subsubsection{Example I}

Let $(Q_1,\Q^2)$ and $(Q_2,\Q^4)$ be two quadratic forms represented by symmetric matrices $M_{Q_i}$ for $i=1,2$.
Let $\Mat_{4,2}(\Z)$ be the set of $4$-by-$2$ matrices with integral coefficients and
\begin{equation*}
    \calU(\Z):= 
    \left\{
    M\in  \Mat_{4,2}(\Z)
    \;\middle\vert\;
    M^{\Tr} M_{Q_2} M = M_{Q_1}
    \right\}.
\end{equation*}
satisfying:
\begin{itemize}
    \item both $(Q_1,\R^2)$ and $(Q_2,\R^4)$ are split. And $(Q_1,\Q^2)$ is $\Q$-anisotropic;
    \item $\calU(\Q)$ is nonempty.
\end{itemize}
Hence we can find $(Q_1',\Q^2)$ such that $(Q_2,\Q^4)$ is isomorphic to $(Q_1,\Q^2)\oplus (Q_1',\Q^2)$. We further assume that
\begin{itemize}
    \item the orthogonal group of $(Q_1,\Q^2)$ is not isomorphic to the orthogonal group of $(Q_1',\Q^2)$ as linear algebraic groups over $\Q$.
\end{itemize}

Let $\norm{\cdot}$ denote the Euclidean norm of a matrix. 

\begin{thm}
     As $R$ tends to infinity,
    there exists some constant $c>0$ such that 
    \begin{equation*}
         \# \left\{
             M\in \calU(\Z) ,\; \norm{M}\leq R
        \right\}  \sim c R^2 \log R.
    \end{equation*}
\end{thm}

See Theorem \ref{theorem_main_example_1} for more details.

\begin{rmk}
    As suggested to us by Rudnick, such a result should be known, but we fail to identify a reference. Also, the many restrictions put here should be unnecessary. From the point view of the current paper, a discussion of the general case (especially, lifting the dimension restrictions on quadratic forms) seems possible if the explicit construction of resolution of singularities is understood. But this is too complicated to be discussed here.
\end{rmk}

\subsubsection{Example II}

Let $\Lambda(n,1)$ denote the set of splittings of $\Z^{n+1}$ as a direct sum of two subgroups, one of which has rank $n$. Each element $A\oplus B$ of $\Lambda(n,1)$ can be associated with $(v,M)$ where $v\in \Z^{n+1}$ is an integral vector and $M$ is an element in $\wedge^{n}\Z^{n+1}$ where $A=\Z .v$, $M=v_1\wedge ...\wedge v_n$ with $B=\oplus_{i=1}^n \Z.v_i$. Moreover, $v,M$ are unique up to sign. Below $\norm{\cdot}$ denotes the natural Euclidean norm.

\begin{thm}
    For two positive integers $\lambda_1,\lambda_2$, we have 
    \begin{equation*}
        \#\{[(v,M)]\in\Lambda(n,1), \;
        \norm{v} ^{\lambda_1} \norm{M}^{\lambda_2} \leq R \} \sim 
        \begin{cases}
            c_{\lambda_1,\lambda_2} \cdot R^{\frac{n}{\lambda_1}} \log (R) ,\quad
            & \lambda_1=\lambda_2
            \\
            c_{\lambda_1,\lambda_2} \cdot R^{\frac{n}{\min\{ \lambda_1,\lambda_2 \} } },\quad
            & \lambda_1 \neq \lambda_2
        \end{cases}
    \end{equation*}
    for some $ c_{\lambda_1,\lambda_2} > 0$.
\end{thm}

See Theorem \ref{theorem_main_example_2} below for details.

\subsubsection{Example III}\label{subsection_introduction_example_III}

The third example concerns with the space of triangles. 
Let $\calM_3(\Z)$ denote the set of three ordered linearly independent lines in $\Q^3$. For an element $(\bml_1,\bml_2,\bml_3) \in \calM_3(\Z)$, we take $\bmv_i$  to be a nonzero vector on $\bml_i$.
Here are two ways of measuring the ``complexity'' of such an element (of course, it is independent of the choice of $\bmv_i$):
\begin{equation*}
    \begin{aligned}
        &\Ht_1(\bml_1,\bml_2,\bml_3):=
        \frac{
        \prod_i \norm{\bmv_i}^2
        }{
        \norm{\bmv_1\wedge \bmv_2}
        \norm{\bmv_1\wedge \bmv_3}
        \norm{\bmv_2\wedge \bmv_3}
        }\\
        &\Ht_2(\bml_1,\bml_2,\bml_3):=
        \frac{
        \norm{\bmv_1\wedge \bmv_2}^2
         \norm{\bmv_1\wedge \bmv_3}^2
         \norm{\bmv_2\wedge \bmv_3}^2
        }
        {
        (\prod_i \norm{\bmv_i}) \cdot \norm{\bmv_1\wedge \bmv_2 \wedge\bmv_3}^3
        }.
    \end{aligned}
\end{equation*}

\begin{thm}
    Let $\kappa_1,\kappa_2 >0$ and 
    $\Ht(\bmx):= \Ht_{1}(\bmx)^{\kappa_1} \Ht_{2}(\bmx)^{ \kappa_2 }$ for $\bmx \in \bmM_3(\R)$.
    Then for some constant $c_{\kappa_1,\kappa_2}>0$,
    \begin{equation*}
        \sum_{
         \left\{
          \bmx \in \calM_3(\Z) \;\middle\vert\;
          \Ht(x) \leq R
        \right\} }  \bmw_x
        \sim 
        \begin{cases}
             c_{\kappa_1,\kappa_2} \cdot R^{\frac{8}{3} \kappa_1^{-1} } \cdot \log(R)
            \quad & \kappa_1 =\kappa_2  \\
             c_{\kappa_1,\kappa_2} \cdot R^{\frac{8}{3} \max\{\kappa_1^{-1},\kappa_2^{-1} \} } 
            \quad & \kappa_1 \neq \kappa_2.
        \end{cases}
    \end{equation*}
\end{thm}

See Theorem \ref{theorem_main_example_3} below for more details, especially the definition of $\bmw_x$.

\subsection{Outline of the proof and future directions}

\subsubsection{Outline}
Here is an outline of the proof as well as this paper.

First, it is known that if $\bmU$ admits transitive action of a linear algebraic group $\bmG$, then its set of integral points admits an action of some arithmetic subgroup $\Gamma$ and furthermore, there are only finitely many orbits (see Theorem \ref{theorem_finiteness_orbits}). 
Fix such an orbit $\Gamma.x$.
One wishes to compare the asymptotics of number of points of bounded height on $\Gamma.x$ with volume of points of bounded height on $\bmG(\R)^{\circ}.x$. 
The basics of height function are discussed in Section \ref{subsection_arithmetic_height}.
We explain what ``volume'' is from the point view of invariant gauge forms and discuss its boundary behaviours in Section \ref{subsection_invariant_gauge_forms}. Particularly, Theorem \ref{theorem_anticanonical_affine_reductive} and \ref{theorem_anticanonical_general_arithmetic} are important when one applies results from \cite{ChamTschin12} in later sections.

How to make this comparison? 
 It is discovered in \cite{DukRudSar93} that the comparison can be made if one solves certain equidistribution problem on $\bmG(\R)^{\circ}/\Gamma$.
The equidistribution problem alluded is about limiting behaviour of certain average of translates of a fixed homogeneous measure, which looks like:
\begin{equation*}
    \frac{1}{\Vol(B_{R,x})}
    \int_{[g]\in B_{R,x}} g_* \rmm_{\bmH_x(\R)\Gamma/\Gamma}
    \, \Vol([g]).
\end{equation*}
Here $\bmH_x$ is the stabilizer of $x$ in $\bmG$, $\Vol$ refers to certain natural measure on $\bmG(\R)/\bmH_x(\R)$ and $B_{R,x}$ are those cosets $g\bmH_x(\R)$ satisfying $\Ht(g.x)\leq R$.  We elaborate on this in Section \ref{section_equidistribution_orbital_counting}.

The study of the limiting behaviors of the integrand in the above expression is the subject of \cite{EskMozSha96} and \cite{zhangrunlinCompositio2021}. Perhaps the most crucial input behind these works is the classification of unipotent-invariant ergodic probability measures due to Ratner \cite{Rat91}.
We will review the relevant equidistribution and nondivergence results in Section \ref{section_equidistribution_nondivergence}.

Theorems from \cite{EskMozSha96, zhangrunlinCompositio2021}
are most useful when $\bmH$  is large in $\bmG$. When this is not the case, one may do the centralizer trick by considering, for instance, $\bmG \times \bmZ/ \bmH \cdot \Delta(\bmZ)$ where $\bmZ$ is contained in the centralizer of $\bmH$ in $\bmG$.
This does not change the underlying variety, but by considering a larger automorphism group, we implicitly put restrictions on the heights under consideration.
The trick works better when $\bmH$ is connected and reductive.

Anyway, assuming $\bmH$ is large in $\bmG$ in some sense, one can try to analyze the possible limiting measures of $(g_*\rmm_{\bmH_x(\R)\Gamma/\Gamma})$ as $g$ varies.
The best thing is that for ``most'' $g$ tending to infinity, the limit is the $\bmG(\R)^{\circ}$-invariant Haar measure on $\bmG(\R)^{\circ}/\Gamma$. When this fails, we say that \textit{focusing} happens. The crucial thing is to understand whether focusing\footnote{We will rarely use the term focusing in the main body of the paper. The precise definition can be found in \cite{EskMozSha96}.} happens or not. 

Preferably, focusing can be detected by finitely many vectors in finitely many linear representations. 
If so, then we will manage to interpret these conditions in the language of algebraic varieties and their closed subvarieties and furthermore, by resolution of singularities for a pair, into language of smooth projective varieties and their simple normal crossing divisors.
Once this is done, the equidistribution theorem of \cite{ChamTschin12} becomes very useful.
We will discuss, in a rather abstract way, the relation between this and the counting problem in Section \ref{section_good_heights}.
Then in Section \ref{section_equidistribution_focusing}, we will provide concrete group theoretical conditions under which abstract setup in Section \ref{section_good_heights} can be actually realized.
These will further be illustrated by three concrete examples in the last three sections.

\subsubsection{What remains to be done?}

We list some questions that are related to the discussion of the current paper. Of course, one natural task is to verify the log Manin conjecture (which should correspond to the ``non-focusing'' case because of the ``rigid'' condition) recently proposed by \cite{Santens_Manin_integral_points_2023_arxiv}. Here we want to be a bit more specific.

Let $\bmU$ be a homogeneous variety over $\Q$ isomorphic to $\bmG/\bmH$.  Assume $\bmU(\Q)\neq \emptyset$. For $x\in \bmU(\Q)$, we let $\bmH_x$ be the stabilizer of $x$ in $\bmG$.
Let $\rmG:= \bmG(\R)^{\circ}$, $\rmH_x:= \bmH_x(\R)\cap \rmG$ and $\Gamma$ be an arithmetic subgroup contained in $\rmG$.
Let $\rmm_{[\rmH_x]}$, $\rmm_{[\rmG]}$ and  $\rmm_{\rmG/\rmH_x}$ denote invariant measures (so we need to assume that they exist) on $\rmH_x\Gamma/\Gamma$, $\rmG/\Gamma$ and $\rmG/\rmH_x$ respectively. We require these measures to be compatible. Namely, if $\rmm_{\rmH_x}$ and $\rmm_{\rmG}$ are the corresponding invariant measures on $\rmH_x$ and $\rmG$ respectively that induce 
$\rmm_{[\rmH_x]}$ and $\rmm_{[\rmG]}$, then for every compactly supported function $f$ on $\rmG$, one has
\begin{equation*}
    \int 
    \left(  \int f(gh) \, \rmm_{\rmH_x}(h) \right)
    \rmm_{\rmG/\rmH_x}([g])
    = \int f(x) \, \rmm_{\rmG}(x).
\end{equation*}

\subsubsection{Finite-volume case}

Assume 
\begin{itemize}
    \item $\bmG$ and $\bmH$ are both connected;
    \item all the nontrivial $\Q$-factors of $\bmG$ are noncompact over $\R$;
    \item $\bmG$ and $\bmH$ have no nontrivial $\Q$-characters. Equivalently, $\bmG(\R)/\Gamma$ and $\bmH(\R)\Gamma/\Gamma$ have finite volume.
\end{itemize}

We expect that there exists a class of good heights (or, good boundary divisors) such that 
\begin{equation}\label{equation_expectation_finite_volume}
     \# \left\{
       y \in \Gamma.x 
       \;\middle\vert\;
       \Ht(y)\leq R
    \right\} \sim 
    \frac{
    \normm{  \rmm_{ [\rmH_x] } }
    }{ 
    \normm{
    \rmm_{[\rmG]}
    }
    }
    \cdot \rmm_{\rmG/\rmH_x}\left(
     \left\{
       y \in \rmG.x 
       \;\middle\vert\;
       \Ht(y)\leq R
    \right\}
    \right)
\end{equation}
for every $x\in \bmU(\Q)$ and every arithmetic lattice $\Gamma$ (of course, it is sufficient to verify this for some finite-index subgroup of $\Gamma$).

\begin{ques}
    Assume that $\bmH$ is large in $\bmG$ and all the invariant measures exist.
    Let $(\bmX,\bmD)$ be a smooth $\rmG$-pair over $\Q$ with $\bmX\setminus \bmD$ equivariantly isomorphic to $\bmU$.
    Is Equa.(\ref{equation_expectation_finite_volume}) true for the log anti-canonical height on $(\bmX,\bmD)$?
\end{ques}
Being large does not have a precise definition. For instance, one may say that $\bmH$ is large in $\bmG$ if one of the following is true:
\begin{itemize}
    \item[1.] $\bmG$ is semisimple and $\bmG/\bmH$ is a spherical variety;
    \item[2.] $\bmG$ and $\bmH$ are reductive and the identity component of the centralizer of $\bmH$ in $\bmG$ is contained in $\bmH$.
\end{itemize}

\subsubsection{Infinite-volume case}

Assume that
\begin{itemize}
    \item $\bmG$ is connected, semisimple and all $\Q$-simple factors are noncompact;
    \item  $\bmH$ is connected, reductive and has nontrivial $\Q$-characters.
\end{itemize}
In this case, we no longer expect the asymptotic of integral points is the same as the volume asymptotic on $\rmG/\rmH$. However, we do have a replacement of the volume asymptotic.

Recall $\scrP^{\max}_{\bmH}$ is the collection of maximal parabolic $\Q$-subgroups containing $\bmH$. For each $\bmP \in \scrP^{\max}_{\bmH}$, we fix an integral vector\footnote{Implicitly, an integral structure $\frakg_{\Z}$ of $\frakg$ is fixed.} $\bmv_{\bmP} \in \wedge^{\dim\bmP} \frakg$ that lifts the Lie algebra of $\bmP$. 
Let $\Ht$ be a height function on $\bmU(\R)$ and let $x\in \bmU(\Q)$.
Let $\bmH'_x$ be the largest connected $\Q$-subgroup of $\bmH_x$ where all the $\Q$-characters of $\bmH_x$ vanish.
For $\eta, R>0$, define the following region
\begin{equation*}
    B_{\eta,R ,x}:=
    \left\{
      [g] \in \rmG/\rmH_x' \;\middle\vert\;
      \Ht(g.x)\leq R,\;
      \norm{\Ad(g) \bmv_{\bmP} }\geq \eta,\;\forall\,
      \bmP\in \scrP^{\max}_{\bmH_x}
    \right\}.
\end{equation*}
In favorable situations, we expect that the asymptotic of $ \rmm_{\rmG/\rmH'_x} \left(
     B_{\eta,R ,x}
    \right)$ as $R \to +\infty$ is independent of $\eta$.

\begin{ques}
    Fix $\eta>0$.
    For what height functions does
    \begin{equation*}
         \# \left\{
       y \in \Gamma.x 
       \;\middle\vert\;
       \Ht(y)\leq R
    \right\} \sim 
    \frac{1}{ 
    \normm{
    \rmm_{[\rmG]}
    }
    }
    \cdot \rmm_{\rmG/\rmH'_x} \left(
     B_{\eta,R ,x}
    \right)\;
    \text{ as }R\to +\infty
    \end{equation*}
     hold for every $x\in \bmU(\Q)$ and every arithmetic subgroup $\Gamma$?
\end{ques}

\begin{ques}
    Given $\eta>0$, what is the asymptotic of $\rmm_{\rmG/\rmH'_x} \left(
     B_{\eta,R ,x}
    \right)$?
\end{ques}

In the case of Example III (Section \ref{subsection_introduction_example_III}) discussed above, it seems that $\rmm_{\rmG/\rmH'_x} \left(
     B_{\eta,R ,x}
    \right)$ can be expressed in terms of the norms of sections of certain line bundles.

\section*{Notations}

For a linear algebraic group $\bmG$ over a field $k$, let the gothic letter $\frakg$ denote its Lie algebra over $k$.
When the notation looks cumbersome, we also use $\Lie(\bmG)$.
Let $\frakX^*(\bmG)$ (resp. $\frakX^*_k(\bmG)$) be the group of characters (resp. $k$-characters) of $\bmG$.
Similarly $\frakX_*^{k}(\bmG)$  is the group of cocharacters over $k$.
Following \cite{BorSer73}, when $k=\Q$, we let 
\begin{equation*}
    {}^{\circ}\bmG:=\bigcap_{\alpha \in \frakX^*_{\Q}(\bmG)} \ker \alpha^2.
\end{equation*}
Let $\bmS_{\bmG}$ be the quotient $\Q$-split torus $\bmG/ {}^{\circ}\bmG $ and
$p^{\spl}$ be the natural quotient morphism $\bmG \to \bmS_{\bmG}$.

Assume $\bmG$ to be connected and
let $\bmR_{\bmu}(\bmG)$ be its unipotent radical, 
then $\bmG$ is equal to
\begin{equation*}
    \bmG^{\red} \ltimes \bmR_{\bmu}(\bmG)
    =
    (\bmG^{\semisimple} \cdot \bmZ(\bmG)) \ltimes \bmR_{\bmu}(\bmG)
    =
    \left( \bmG^{\cpt}\cdot \bmG^{\nc} \cdot \bmZ(\bmG)^{\an} \cdot \bmZ(\bmG)^{\spl} \right) \ltimes \bmR_{\bmu}(\bmG)
\end{equation*}
where $\bmG^{\red}$ is a connected reductive $\Q$-subgroup (called a Levi subgroup) lifting $\bmG/\bmR_{\bmu}(\bmG)$, $\bmG^{\semisimple}= [\bmG^{\red},\bmG^{\red}]$ is a connected semisimple $\Q$-subgroup, $\bmZ(\bmG)$ is the identity component of the center of $\bmG^{\red}$ (called the central torus), $\bmG^{\cpt}$ is the product of $\Q$-simple factors of $\bmG^{\semisimple}$ that are $\R$-anisotropic, $\bmG^{\nc}$
is the product of $\Q$-simple factors of $\bmG^{\semisimple}$ that are $\R$-isotropic.
A $\Q$-torus $\bmT$ can be written uniquely as an almost direct product of $\bmT^{\an}\cdot \bmT^{\spl}$, a $\Q$-anisotropic subtorus and a $\Q$-split subtorus. This explains the notation $\bmZ(\bmG)^{\an} $ and $ \bmZ(\bmG)^{\spl}$ above.

Let $\bbA_{\Q}$ denote the topological adelic ring of the field of rational numbers $\Q$. $\Gal_{\Q}$ denotes the Galois group of the extension $\overline{\Q}/\Q$.

We say that $(\bmX,\bmD)$ is a $\bmG$-pair over $\Q$ if $\bmX$ is a projective variety over $\Q$, $\bmD$ is a closed subvariety over $\Q$ and $\bmX$ is equipped with the action of $\bmG$ such that $\bmD$ is preserved under this action. 
A pair $(\bmX,\bmD)$ is said to be a smooth pair over $\Q$ if $\bmX$ is a smooth projective variety over $\Q$ and  $\bmD$ is a simple/strict normal crossing divisor over $\Q$.

For locally finite measures $(\nu_n)$ and $\nu$ on a locally compact second countable Hausdorff space $X$, following \cite{ShaZhe18}, we say that $[\nu_n]$ converges to $[\nu]$ and write $\lim_{n\to \infty} [\nu_n] = [\nu]$ if there exists $(a_n)\subset \R_{>0}$ such that $\lim_{n\to \infty} a_n \nu_n =\nu$ under the weak-$*$ topology.

\section{Arithmetic and Geometry of Equivariant Compactifications}

For this section we further assume the following notations:
\begin{itemize}
    \item  $\bmX$ is a smooth projective variety over $\Q$ of dimension $d$;
    \item  $\bmD \subset \bmX$ is a simple normal crossing divisor over $\Q$ and  $\bmU$ is the complement of $\bmD$ in $\bmX$;
    \item $(\bmD_{\alpha})_{\alpha\in \scrA}$ are the irreducible components of $\bmD$ over the algebraic closure $\overline{\Q}$. This set and its index set $\scrA$ are equipped with an action of $\Gal_{\Q}$;
    \item  $\bmG$ is a linear algebraic group over $\Q$ and $\bmH$ is a $\Q$-subgroup of $\bmG$.
\end{itemize}

\subsection{Integral points, adelic heights and orbits of arithmetic groups}\label{subsection_arithmetic_height}

\subsubsection{Integral points}
Let $\Val$ be the set of valuations on $\Q$ up to equivalence, identified with the set of prime numbers and the Archimedean place $\{\infty\}$. By convention $\Q_{\infty}:=\R$.
Let $\Val_f$ denote only the set of finite primes. 
For a finite set $S\subset \Val_f$, let $\Z_S$  denote the set of rational numbers that are integral outside $S$.

Consider the space of adelic points given by certain restricted product
\begin{equation*}
    \bmU(\bbA_{\Q}) := {\prod}^{\prime}_{p\in \Val} \bmU(\Q_p)  \times \bmU(\R)
\end{equation*}
equipped with adelic topology (see \cite{Conrad_adelic_points_2012} for details).
Let $\scrK_f$ be a nonempty open compact subset of $\prod_{p\in \Val_f} \bmU(\Q_p)$. More concretely, there exists a finite set $S\subset \Val_f$ and a smooth model $(\calX,\calD)$ of $(\bmX,\bmD) $ over $\Spec \Z_S$ such that 
$\scrK_f =  K_S \times \prod_{p \notin S} K_p $ for some open compact subset $K_S \subset \prod_{p\in S} \bmU(\Q_p)$ and 
\begin{equation}\label{equation_S_K_p}
    K_p = (\calX \setminus \calD) (\Z_p) ,\quad \forall \, p\notin S.
\end{equation}
For $\scrK= \scrK_f \times K_{\infty}$ where $K_{\infty}$ is a union of some connected components of $\bmU(\R)$, define
\begin{equation*}
    \bmU_{\scrK}(\Z):= \bmX(\Q ) \bigcap \scrK_f \times K_{\infty}.
\end{equation*}
This adelic point of view contains the classical case of taking the integral points of some integral model as a special case. 

\subsubsection{Adelic heights}\label{subsection_heights}
As a reference, see \cite{Chamber-Loir_Tschinkel_2010_Igusa_integral}.

Let $(L,s)$ be a line bundle over $\bmX$ together with a global section $s\in \Gamma(\bmX, L)$. 
An \textit{adelic metric} $\left(\norm{\cdot}_{p}\right)_{p\in \Val}$ is the following:
\begin{itemize}
    \item[1.]  for each $p\in \Val$, the $\norm{\cdot}_p$ is a smooth metric on the analytic line bundle $L(\Q_p)$ over the analytic manifold $\bmX(\Q_p)$;
    \item[2.] there exist a finite set $S\subset \Val_f$ and a smooth model $\calL \to \calX$ over $\Spec \Z_S$  such that $\norm{\cdot }_p$ is  defined by
    \begin{equation}\label{equation_S_adelic_metric}
        \norm{s(x)}_p \leq 1 \iff
        s(x) \in  \wtx^* \calL
    \end{equation}
    where $x\in \bmX(\Q_p)$ and $\wtx \in \calX(\Z_p)$ is a lift of $x$.
\end{itemize}
By saying that a metric $\norm{\cdot}_p$ is \textit{smooth}, we mean that for every nonzero local section $s$, the map $x \mapsto \norm{s(x)}_p $ is a smooth function on the locus where $s$ does not vanish.
For simplicity, $L$ together with this adelic metric is called an \textit{adelic line bundle}.

If $L= \calO_{\bmX}(\bmL)$ for some effective divisor $\bmL$ and $\bms_{\bmL}$ is the canonical section, then for some finite set $S$,  
\begin{equation}\label{equation_S_metric_equal_to_one}
    \norm{\bms_{\bmL}(x) }_{p} = 1  ,\quad
    \forall \, p\notin S,  \;   x \in K_p.
\end{equation}
Indeed, since $\bms_{\bmL}$ is a local generator of $L$ over $\Q$, it is also a local generator of $\calL$ over $\bmZ_{S}$ for some  finite $S$.
Take $S$ larger such that Equa.(\ref{equation_S_adelic_metric}) also holds.
Then Equa.(\ref{equation_S_metric_equal_to_one}) follows from Equa.(\ref{equation_S_adelic_metric}) above.

Assume that the support of $\bmL$  is contained in $\bmD$. We define an \textit{adelic height function} on $\bmU(\bbA_{\Q})$ by
\begin{equation*}
    \Ht \left(  (x_p)
    \right) := \left( 
    \prod_{p\in\Val} \norm{\bms_{\bmL}(x_p)}_p
    \right) ^{-1}
\end{equation*}
Now we choose a finite $S\subset\Val_f$ large enough such that Equa.(\ref{equation_S_K_p}, \ref{equation_S_adelic_metric}, \ref{equation_S_metric_equal_to_one}) hold and moreover, the integral models $\calX$ appearing in Equa.(\ref{equation_S_K_p}) and  (\ref{equation_S_adelic_metric}) coincide over $\Z_S$.
Then
\begin{equation*}
    \Ht \left(  (x_p)
    \right)
    = 
    \left( \norm{\bms_{\bmL}(x_{\infty})}_{\infty} \cdot 
    \prod_{p\in S} \norm{\bms_{\bmL}(x_p)}_p
    \right ) ^{-1},\;\forall \,(x_p)\in \scrK.
\end{equation*}
A great advantage of adelic height is that, by the product formula, this height is independent of the choice of sections.

\subsubsection{Orbits of arithmetic subgroups}

Let $S\subset \Val_f$, $\scrK$ and  an adelic line bundle $(L, (\norm{\cdot}_p)_{p\in \Val})$ be the same as  last subsection.
Assume that $\bmX$ is a  $\bmG$-equivariant compactification of $\bmG/\bmH$ and $\bmD$ is the complement of $\bmG/\bmH$.
By enlarging $S$, we assume further that we have a smooth model $\calG$ of $\bmG$ over $\Z_S$ and that the action $\bmG\times \bmU \to \bmU$ over $\Q$ extends to $\calG \times \calU \to \calU$ over $\Z_{S}$.
Note the following:
\begin{lem}
    For $p\in \Val_f$, the subgroup
    \begin{equation*}
         \Stab_{\bmG}(\norm{\cdot}_p) :=
         \left\{
         g\in \bmG(\Q_p) \;\middle\vert\;
         \norm{\bms_{\bmL}}_{g.x} = \norm{\bms_{\bmL}}_{x} ,\;\forall \, x\in K_p
        \right\}
    \end{equation*}
   contains an open and compact subgroup of $\bmG(\Q_p)$.
\end{lem}

\begin{proof}
    Indeed, the function $x \mapsto \norm{s_{\bmL}(x)}_{p}$ is smooth and hence locally constant on $\calU(\Q_p)$ and in particular on $K_p$. By compactness and the totally disconnectedness of p-adic topology, we can cover $K_p$ by disjoint open sets such that $\norm{s_{\bmL}(x)}_{p}$ is constant on each open set. By choosing a small enough open compact subgroup $K \subset \bmG(\Q_p)$, we can guarantee that each open subset is $K$-stable. Such a $K$ is thus contained in $ \Stab_{\bmG}(\norm{\cdot}_p)$, proving the claim.
\end{proof}

Let 
\begin{equation*}
\Gamma_p :=\begin{cases}
        \text{an open, compact subgroup of }\Stab_{\bmG}(\norm{\cdot}_p),\quad p\in S
        \\
        \calG(\Z_p) ,\quad p\notin S,
    \end{cases}
\end{equation*}
and
\begin{equation*}
    \Gamma:= \bmG(\Q)\cap  \prod_{p\in \Val_f} \Gamma_p.
\end{equation*}
Discussions above prove the following:
\begin{lem}\label{lemma_orbit_arithmetic_group_same_local_height}
    The $\Gamma$ defined above is an arithmetic subgroup of $\bmG(\Q)$ that preserves $\bmU_{\scrK}(\Z)$. Moreover, 
\begin{equation*}
    \Ht_p(x)= \Ht_p(\gamma.x) ,\;
    \forall\, p\in \Val_f,\; x= (x_p) \in \bmU_{\scrK}(\Z),\;\gamma\in \Gamma
\end{equation*}
where $\Ht_p(x):= \norm{\bms_{\bmL}(x_p)}_p^{-1}$.
\end{lem}

The theorem below follows from \cite[Proposition 7.13, 7.14]{Gille_Moret-Bailly_2013}. Special cases can be proved using reduction theory (see \cite[Remark 6.4]{Bor2019}).

\begin{thm}\label{theorem_finiteness_orbits}
    There are only finitely many $\Gamma$-orbits on $\bmU_{\scrK}(\Z)$.
\end{thm}

 We sketch a proof following \cite{Wei_Xu_2016} assuming $\bmH$ is connected for the sake of completeness.

\begin{proof}[Sketch of proof assuming $\bmH$ to be connected]
The conclusion is going to follow from a spreading-out argument, finiteness of Shafarevich--Tate groups and finiteness of class numbers of linear algebraic groups.
We define
\begin{equation*}
    x \sim_{\loc} y \iff \Gamma_p.x = \Gamma_p.y,\;\forall p\in \Val_f
\end{equation*}
 for two elements $x=(x_p)$ and $y=(y_p)$ in $\bmU_{\scrK}(\Z)$. Let $[x]_{\loc}$ denote the equivalence class containing $x$.
As a first step, we note that $\{[x]_{\loc}\}$ is finite (compare \cite[Lemma 1.6.4]{BoroRudn95}).
Indeed, for $p\notin S$, there always exists $\gamma_p \in \Gamma_p$ such that $y_p =\gamma_p x_p$ by Hensel's lemma and Lang's theorem on triviality of the Galois cohomology group $H^1(F,\bmH)$ for a finite field $F$ and connected $\bmH$. So it remains to show that  for each $p\in S$
, $K_p \cap \bmU(\Q)$ is finite modulo the equivalence relation defined by $x\sim y$ iff $\Gamma_p .x = \Gamma_p .y$, which follows from the finiteness of $K_p/\Gamma_p$.

We define an equivalence relation $\sim_{\Q}$ on $\bmU_{\scrK}(\Z)$ by
\begin{equation*}
    x\sim_{\Q} y \iff \bmG(\Q).x =\bmG(\Q).y.
\end{equation*}
Then we show that for each $x$, the number of equivalence classes $[x]_{\loc}/\sim_{\Q}$ is finite. 
Indeed, by inspecting the following commutative diagram
\begin{equation*}
    \begin{tikzcd}
         1 \arrow[r,""]& \bmH(\Q) \arrow[r,""]\arrow[d,""]& \bmG(\Q) \arrow[r,""]\arrow[d,""] & \bmU(\Q) \arrow[r,""]\arrow[d,""] &
         H^1(\Q,\bmH)
         \arrow[d]\\
        1\arrow[r,""]& \bmH(\bbA_{\Q})\arrow[r,""]& \bmG(\bbA_{\Q}) \arrow[r,""]& \bmU(\bbA_{\Q}) \arrow[r,""]& \prod_{p\in\Val_f} H^1(\Q_p,\bmH),
    \end{tikzcd}
\end{equation*}
we have an injection 
\begin{equation*}
    [x]_{\loc}/\sim_{\Q} \embed \ker\left(H^1(\Q,\bmH) \to \prod_p H^1(\Q_p,\bmH)\right).
\end{equation*}
But the latter is finite (See \cite[Theorem 7.1]{BorelSerre1964_finiteness_Galois_cohomology}).

For $y\in [x]_{\loc}$, let $[y]_{\loc,\Q}:= [x]_{\loc} \cap [y]_{\Q}$ denote the equivalence class containing $y$. As the final step, we show that $\Gamma \backslash [y]_{\loc,\Q}$ is finite for every $x\in \bmU_{\scrK}(\Z)$ and $y\in [x]_{\loc}$. For each $z \in [y]_{\loc,\Q}$, we find $q_z \in \bmG(\Q)$ and $\gamma_z \in \prod \Gamma_p$ such that $z =q_z y =\gamma_z y$. Then the map
\begin{equation*}
    \begin{aligned}
        \Gamma \backslash \left[ y \right]_{\loc ,\Q} &\to  \bmH_y(\Q) \bs \bmH_y(\bbA_{\Q}) / \bmH_{y}(\bbA_{\Q})\cap \prod \Gamma_p \\
        z &\mapsto q_z^{-1}\gamma_z
    \end{aligned}
\end{equation*}
is well-defined and injective. As the right hand side is finite (see \cite[Theorem 5.8]{PlatonovRapinchuk2023Vol1} or \cite[Theorem 5.1]{Borel1963_finiteness_adele_groups}), we are done.

\end{proof}

\begin{rmk}\label{remark_adelic_real_height}
    In light of Lemma \ref{lemma_orbit_arithmetic_group_same_local_height} and Theorem \ref{theorem_finiteness_orbits}, it makes no difference to count with respect to the adelic height or the height at infinity.
\end{rmk}

\subsection{Invariant gauge forms}\label{subsection_invariant_gauge_forms}

Assume that $(\bmX,\bmD)$ is a smooth $\bmG$-pair over $\Q$ with $\bmU=\bmX \setminus \bmD$  isomorphic to $\bmG/\bmH$ $\bmG$-equivariantly. Let $o \in \bmU(\Q)$ denote the identity coset.
Let $\calT_{\bmX}$ be the tangent bundle of $\bmX$, $\calT_{\bmX}^*$ be its dual and $\calK_{\bmX}:= \det (\calT_{\bmX}^*) $ be the canonical line bundle.
Let $\la \cdot,\cdot \ra$ denote the natural pairing between $\calT_{\bmX}$ and $\calT_{\bmX}^*$ and their wedge products.

An element $g\in \bmG(\C)$ acts on $\bmX(\C)$ and hence on tangent vectors by $(g,v)\mapsto g_* v $ and differential forms by $(g,\omega)\mapsto g^*\omega$. From the definition, one has that for every analytically open subset $\calO \subset \bmX(\C)$, $\omega \in \Gamma(\calO,\calT^*_{\bmX})$ and $\partial \in \Gamma(\calO, \calT_{\bmX})$,
\begin{equation*}
    \la \omega, \partial   \ra _{g.x} 
    = \la g^*\omega, (g^{-1})_* \partial \ra _{x},\;
    \forall \, x\in g^{-1}\calO,  \; g\in \bmG(\C).
\end{equation*}

Identify the fibre $\calT_o(\bmX)$ (resp. $\calT^*_o(\bmX)$) of the tangent (resp. cotangent) bundle at $o$ with $\frakg/\frakh$ (resp. $(\frakg/\frakh)^{*}$).  The natural pairing on
$( (\frakg/\frakh)^* , \frakg/\frakh)$ is then identified with that on $(\calT^*_o(\bmX), \calT_o(\bmX))$. Let $\Delta_{\bmG}$ (resp. $\Delta_{\bmH}$) denote the determinant character of the adjoint action of $\bmG$ (resp. $\bmH$) on $\frakg$ (resp. $\frakh$), that is, $\Delta_{\bmG}(g)= \det(\Ad(g),\frakg)$ and $\Delta_{\bmH}(h) = \det(\Ad(h),\frakh)$.

Since $\bmH$ stabilizes $o$, its actions on $\calT_{\bmX}$ and $\calT^*_{\bmX}$ induce actions on $\calT_o(\bmX)$ and $\calT_o^*(\bmX)$. These actions are identified with the adjoint action of $\bmH$ on $\frakg/\frakh$ and $(\frakg/\frakh)^*$ respectively.
Therefore,
\begin{equation*}
    \begin{aligned}
        \det (h_*, \calT_o(\bmX))
        = \det ( \Ad(h), \frakg/\frakh )
        = \frac{ \det (\Ad(h), \frakg) 
        }
        {
        \det (\Ad(h), \frakh) 
        } = \frac{\Delta_{\bmG}(h)}{\Delta_{\bmH}(h)}.
    \end{aligned}
\end{equation*}
and
\begin{equation*}
    \det( (h^{-1})^* , \calT^*_o(\bmX) ) = 
    \det( h_*, \calT_o(\bmX) ) ^{-1} = \frac{\Delta_{\bmH}(h)}{\Delta_{\bmG}(h)}.
\end{equation*}
Define $(\Delta_{\bmG/\bmH})^{-1} :=  \Delta_{\bmG}^{-1}\Delta_{\bmH} $, an element of $\mathfrak{X}^{*}_{\Q}(\bmH)$.

For a $\Q$-character $\alpha$ of $\bmH$, let $ \calL_{\alpha}$ be the GIT quotient of $\bmG \times \Spec\Q[x]$ by the right action of $\bmH$ defined by
\begin{equation*}
    (g,\lambda)\cdot h := ( gh, \alpha(h)^{-1} \lambda  ).
\end{equation*}
It is a $\bmG$-linearized line bundle over $\bmG/\bmH$. On the other hand, given a $\bmG$-linearized line bundle $\calL$ over $\bmG/\bmH$, the (left) action of $\bmH$ on the fibre of $o$ gives a $\Q$-character $\alpha_{\calL}$ of $\bmH$. One can check the following:

\begin{lem}
   Let $\calL$ be a $\bmG$-linearized line bundle, then $\calL \cong \calL_{\alpha_{\calL}}$. In particular, $\calK_{\bmG/\bmH} \cong \calL_{\Delta_{\bmG/\bmH}^{-1}}$.
\end{lem}

\begin{lem}
    Let $\chi$ be a $\Q$-character of $\bmH$,
    \begin{equation*}
        \Gamma(\bmG/\bmH, \calL_{\chi})
        \cong \left\{
        f \in \Q[\bmG]\;\middle\vert\;
        f(gh) = \chi^{-1}(h) f(g),\;\forall (g,h)\in \bmG\times\bmH
        \right\}.
    \end{equation*}
\end{lem}

By  \cite[2.2, 2.3]{Knop_Kraft_Vust_Picard_G_variety}, the nonvanishing global sections of $\calK_{\bmG/\bmH}$ consist of $\bmG$-eigenvectors. Combined with the lemma above, we have
\begin{lem} \label{lemma_global_nonvanishing_section_canonical_line_bundle}
     The nonvanishing global sections of $\calK_{\bmG/\bmH}$ consist of 
    \begin{equation*}
        \Gamma(\bmG/\bmH, \calK_{\bmG/\bmH} )^{\times}
         \cong \left\{
        f \in \Q[\bmG]   \;\middle\vert\; f= \lambda \chi ,\;
        \exists\,
        \lambda \in\Q^{\times},\; \chi \in \frakX_{\Q}^*(\bmG)\text{ with }
        \chi\vert_{\bmH} = \Delta_{\bmG/\bmH}
        \right\}.
    \end{equation*}
    Hence we may label, up to a scalar, nonvanishing global sections of $\calK_{\bmG/\bmH}$ as $\omega_{\bmG/\bmH}^{\chi}$ with $\chi$ ranging over extensions of $\Delta_{\bmG/\bmH}$ to $\bmG$.
    In particular, $\calK_{\bmG/\bmH}$ is trivial iff $\Delta_{\bmH}$ extends to a $\Q$-character on $\bmG$.
\end{lem}

If both $\Delta_{\bmG/\bmH}$ and $\chi$ are trivial characters, we abbreviate $\omega_{\bmG/\bmH}:= \omega_{\bmG/\bmH}^{\chi}$ and call it an \textit{invariant gauge form}.
One may wish to compare the discussions above with those in 
\cite[Chapter 1]{Ragh72} or \cite[Chapter 2]{Weil_1982_adele_algebraic_groups}.

\begin{lem}
    If $\calK_{\bmG/\bmH}$ is trivial, then $\bmH$ is observable in $\bmG$.
\end{lem}

\begin{proof}
    By Lemma \ref{lemma_global_nonvanishing_section_canonical_line_bundle}, $\Delta_{\bmH}^{-1}$ extends to a character $\beta$ on $\bmG$ and take some nonzero vector $\bmv_{\beta}$ in the one dimensional representation where $\bmG$ acts by $\beta$.
    Let $\bmN_{\bmG}(\bmH)^{\beta}$ be the stabilizer of $\bmv_{\frakh}\otimes \bmv_{\beta}$, where $\bmv_{\frakh}\in \wedge^{\dim \bmH} \frakg$ is a lift of $\frakh$.
    
    Then $\bmN_{\bmG}(\bmH)^{\beta}$ is observable in $\bmG$ by definition.
    On the other hand, $\bmH$ is contained in $\bmN_{\bmG}(\bmH)^{\beta}$ as a normal subgroup and is hence observable in $\bmN_{\bmG}(\bmH)^{\beta}$. By the transitivity of observability, $\bmH$ is observable in $\bmG$.
\end{proof}

\begin{rmk}
    The converse is not true. For instance, if $\bmG:=\SL_3$ and 
    $\bmH$ is generated by elements of the form
    \begin{equation*}
    \left[
        \begin{array}{ccc}
          t   &  &  \\
             & t & \\
             && t^{-2}
        \end{array} \right]
        \cdot 
        \left[
        \begin{array}{ccc}
          1  &  &  \\
             & 1 & s \\
             && 1
        \end{array} \right].
    \end{equation*}
    Then $\bmH$ is observable in $\bmG$ yet $\Delta_{\bmH}$  is not extendable to $\bmG$.
\end{rmk}

\subsection{Boundary components of anti-canonical divisors}

Keep assumptions in last subsection and assume that $\calK_{\bmG/\bmH}$ is a trivial line bundle.
Here we are interested in the anti-canonical divisor $-\divisor(\omega_{\bmG/\bmH}^{\chi})$ on $\bmX$ where $\chi \in \mathfrak{X}^{*}_{\Q}(\bmG)$ is a fixed extension of $\Delta_{\bmG/\bmH}$.

Note that $\frakg$ naturally maps to $\Gamma(\bmX,\calT_{\bmX})$, inducing 
\begin{equation*}
    \wedge^d \frakg \to 
\Gamma(\bmX, \wedge^d \calT_{\bmX}) \cong \Gamma(\bmX,(\calK_{\bmX})^*).
\end{equation*}
 For $\bmv=v_1\wedge...\wedge v_d \in \wedge^d \frakg$, let $\partial_{\bmv}$ denote its image.
Then for each $\bmv \in \wedge^d \frakg$, 
\begin{equation*}
    \varphi_{\bmv}(x):= \la \omega_{\bmG/\bmH}^{\chi} , \partial_{\bmv} \ra_{x}
\end{equation*}
defines a rational function on $\bmX$ that is regular on $\bmU$.
Thus, to understand the boundary behaviour of  $\omega^{
\chi}_{\bmG/\bmH}$, it suffices to understand that of $\varphi_{\bmv}$ and $\partial_{\bmv}$.

Take a nonzero element $\bml_{\frakh} \in \wedge^d (\frakg/\frakh)^* \embed \wedge^d \frakg^{*} $. 
So $\left (\omega^{\chi}_{\bmG/\bmH} \right)_o =c_{\chi} \bml_{\frakh} $ for some $c_{\chi}\neq 0$.

\begin{lem}\label{lemma_computation_varphi_1}
     For $g\in \bmG(\C)$,
    $\varphi_{\bmv}(g.o) = c_{\chi} \chi(g) \cdot 
    \la \, \bml_{\frakh}, \Ad(g)^{-1} \bmv \,\ra
    $.
\end{lem}

\begin{proof}
    This follows from a direct computation:
    \begin{equation*}
        \begin{aligned}
            \la \omega_{\bmG/\bmH}^{\chi} , \partial_{\bmv}  \ra_{g.o}
           &= \chi(g) \cdot  \la ( g^{-1})^* \omega_{\bmG/\bmH}^{\chi}, \partial_{\bmv} \ra_{g.o} \\
           &=  \chi(g)  \cdot 
           \la \omega_{\bmG/\bmH}^{\chi} , ( g^{-1})_* \partial_{\bmv} \ra_{o} \\
           &=  \chi(g) \cdot 
           \la \omega_{\bmG/\bmH}^{\chi} ,  \partial_{\Ad(g)^{-1}\bmv} \ra_{o} \\
           &=  c_{\chi}\chi(g)  \cdot 
           \la \,\bml_{\frakh}, \Ad(g)^{-1}\bmv \ra.
        \end{aligned}
    \end{equation*}
\end{proof}

\subsubsection{Dual}
To simplify the conclusion of Lemma \ref{lemma_computation_varphi_1}, we use a dual operation.

Fix $\omega_0 \in \wedge^{\dim \bmG} \frakg^*  $ and $\omega_0^* \in \wedge^{\dim \bmG} \frakg $  with $\la \omega_0^*, \omega_0 \ra =1$.
For two positive integers $l+l'= \dim \bmG$, we have a linear isomorphism:
\begin{equation*}
    \begin{aligned}
        \wedge^{l} \frakg 
        &\to
        \wedge^{l'} \frakg^* \cong  (\wedge^{l'} \frakg) ^*
        \\
        \bmv  &\mapsto 
        \bmv^{\vee}
    \end{aligned}
\end{equation*}
defined by
\begin{equation*}
    \bmv \wedge \bmw = \la \bmv^{\vee}, \bmw \ra \omega_0^*,\; \forall\, \bmw \in \wedge^{l'}\frakg.
\end{equation*}
Note that $\la \bmv^{\vee}, \bmw \ra = \la \bmv, \bmw^{\vee} \ra$.
Furthermore,
\begin{equation*}
\begin{aligned}
    \la (\Ad(g)\bmv)^{\vee} , \bmw \ra \omega_0^*
         &= \Ad(g) \bmv \wedge \bmw 
         = \Ad(g) (\bmv \wedge \Ad(g)^{-1}\bmw)\\
         &= \Ad(g) \left(  \la \bmv^{\vee} , \Ad(g^{-1})\bmw \ra \omega_0^*
         \right)  
         \\
         &=\Delta_{\bmG}(g)    \la \bmv^{\vee} , \Ad(g^{-1})\bmw  \ra \omega_0^*.
\end{aligned}
\end{equation*}
Thus, $(\Ad(g)\bmv)^{\vee} = \Delta_{\bmG} (g) \Ad^*(g) \bmv^{\vee}$.
So we see that $(\cdot)^{\vee}$ intertwines $\Ad $ and $\Ad^* \otimes \Delta_{\bmG}$, or equivalently, $\Ad \otimes \Delta_{\bmG}^{-1}$ and $\Ad^* $.

Now take a $\Q$-vector $\bmv_{\frakh}$ lifting $\frakh$ such that $\bmv_{\frakh}^{\vee}= \bml_{\frakh}$.
Then we have
\begin{equation*}
    \begin{aligned}
        \la \, \bml_{\frakh}, \Ad(g)^{-1} \bmv \ra 
        &= 
        \la  \Ad^*(g)\bmv^{\vee}_{\frakh}, \bmv
        \ra
        = \Delta_{\bmG}(g)^{-1} \la \left( \Ad(g)\bmv_{\frakh} \right)^{\vee}, \bmv \ra
        \\
        &=
         \Delta_{\bmG}(g)^{-1}  \la \Ad(g) \bmv_{\frakh}, \bmv^{\vee} \ra.
    \end{aligned}
\end{equation*}
Combined with Lemma \ref{lemma_computation_varphi_1} above, we have
\begin{lem}\label{lemma_computation_varphi_2}
For $g\in \bmG(\C)$,
    $
    \varphi_{\bmv}{(g.o)}= 
    c_{\chi}\chi(g)\Delta_{\bmG}(g)^{-1}
    \la \Ad(g) \bmv_{\frakh}, \bmv^{\vee} \ra.
    $
\end{lem}

\subsubsection{The conclusion}

For $\bmv \in \wedge^d\frakg$ and $\alpha \in \scrA$,
we define integers $d_{\alpha}$, $d_{\alpha}^{\bmv}$ and $d_{\alpha}^{\varphi, \bmv}$ by
\begin{equation*}
    \begin{aligned}
        -\divisor(\omega_{\bmG/\bmH}^{\chi}) := \sum_{\alpha \in \scrA} d_{\alpha} \bmD_{\alpha} ,\;
        \divisor(\partial_{\bmv}):= \sum_{\alpha \in \scrA} d^{\bmv}_{\alpha} \bmD_{\alpha},\; 
        -\divisor(\varphi_{\bmv}) := \sum_{\alpha \in \scrA} d^{\varphi,\bmv}_{\alpha} \bmD_{\alpha}.
    \end{aligned}
\end{equation*}

By the definition of $\varphi_{\bmv}$, we have
\begin{equation}
    -\divisor(\omega_{\bmG/\bmH}^{\chi}) =
     \divisor \partial_{\bmv}- \divisor(\varphi_{\bmv}) 
     ,\; \text{or }
    \;
    d_{\alpha}= d^{\bmv}_{\alpha} + d^{\varphi,\bmv}_{\alpha} ,\;\forall\,\alpha \in \scrA.
\end{equation}

Each $\bmD_{\alpha}$ has dimension strictly smaller than $d$, one sees that $d^{\bmv}_{\alpha}\geq 1$ for every $\bmv$ and $\alpha$. 
As a corollary of Lemma \ref{lemma_computation_varphi_2} we have
\begin{thm}\label{theorem_pole_anti_canonical_divisor}
    Let $\alpha \in \scrA$, $x\in \bmU(\C)$ and $\bmH_x$ be the stabilizer of $x$ in $\bmG$. Let $(g_n)$ be a sequence in $\bmG(\C)$ such that
    \begin{itemize}
        \item[(1)] $\lim_{n \to \infty} g_n.x \in \bmD_{\alpha}(\C)^{\circ}$ and $\normm{\chi(g_n)\Delta_{\bmG}(g_n)^{-1}}=1$;
        \item[(2)] $\left( \Ad(g_n).\bmv_{\frakh_x}  \right)$ is bounded away from $\bmzero$.
    \end{itemize}
    Then $d_{\alpha}\geq 1$. If additionally
    \begin{itemize}
        \item[(3)] $\left( \Ad(g_n).\bmv_{\frakh_x}  \right)$ diverges to $\infty$,
    \end{itemize}
    then $d_{\alpha}\geq 2$.
\end{thm}
Here $\bmD_{\alpha}(\C)^{\circ}$ denotes the points in $\bmD_{\alpha}(\C)$ that do not lie in any other $\bmD_{\beta}(\C)$.

\begin{proof}
    Without loss of generality assume $x=o$ and $\bmH_x=\bmH$. 
    Note that after passing to a subsequence,
    \begin{equation*}
        \begin{aligned}
            \left( \Ad(g_n).\bmv_{\frakh}  \right) \text{ bounded away from } \bmzero
            &\implies l \left( \Ad(g_n).\bmv_{\frakh}  \right) > c,\, \exists \, l\in (\wedge^d \frakg)^*,\, c >0;\\
            \left( \Ad(g_n).\bmv_{\frakh}  \right) \text{ diverges } 
            &\implies l \left( \Ad(g_n).\bmv_{\frakh}  \right) \text{ diverges},\, \exists \, l\in (\wedge^d \frakg)^*.
        \end{aligned}
    \end{equation*}
     By Lemma \ref{lemma_computation_varphi_2} and our assumption, the first line implies the existence of  $\bmv \in \wedge^d \frakg$ such that $ d_{\alpha}^{\varphi,\bmv} \geq 0$, implying $d_{\alpha} \geq 1+  d_{\alpha}^{\varphi,\bmv} \geq 1$. Similarly, the second line implies $d_{\alpha}\geq 2$.
\end{proof}

There are two special cases that we would like to mention.

\begin{thm}\label{theorem_anticanonical_affine_reductive}
    Assume that  $\bmH$ is reductive and $\chi\cdot \Delta_{\bmG}^{-1}$ is trivial. Then $d_{\alpha} \geq 1$  for every  $\alpha \in \scrA$. If additionally, $\bmG$ is reductive and $\bmZ_{\bmG}(\bmH^{\circ})^{\circ} \subset \bmH$, then $d_{\alpha} \geq 2$ for every  $\alpha \in \scrA$. 
\end{thm}

Note that when $\bmH$ is reductive, $\Delta_{\bmH}$ is trivial and hence $\calK_{\bmG/\bmH}$ is trivial. And if $\chi$ corresponds to the trivial extension of $\Delta_{\bmH}$, then $\chi\cdot \Delta_{\bmG}^{-1}$ is trivial.

\begin{proof}
If $\bmG$ is also reductive, then $\bmG.\bmv_{\frakh}$ is closed. 
Indeed, for $\bmH':= \bmN_{\bmG}(\bmH)^{\bmone}$, the stabilizer of $\bmv_{\frakh}$, the identity component of its centralizer in $\bmG$ is in $\bmH'$, thus the closedness follows from the result of \cite{Kemp78}.
For every $\alpha \in \scrA$ and $x\in \bmD_{\alpha}(\C)^{\circ}$, find $(g_n)\subset \bmG(\C)$ with $\lim g_n.o =x$. Then Theorem \ref{theorem_pole_anti_canonical_divisor} applies and completes the proof. 

Now assume moreover that $\bmZ_{\bmG}(\bmH^{\circ})^{\circ} \subset \bmH$. Let $(g_n)$ and $x$ be the same as above. Thus $(g_n)$ is unbounded modulo $\bmH$, and under our assumption, this in turn implies that $(g_n)$ is also unbounded modulo the stabilizer of $\bmv_{\frakh}$. Note that identity component of the normalizer of a connected reductive group (taken to be $\bmH^{\circ}$ here) is the same as the identity connected component of  the group generated by its centralizer and itself. So $(\Ad(g_n).\bmv_{\frakh})$ diverges and we conclude by Theorem \ref{theorem_pole_anti_canonical_divisor}.

When $\bmG$ may not be reductive, we still claim that $\bmzero \notin \overline{\bmG(\C).\bmv_{\frakh}}$, which is sufficient by discussion above.
Find a Levi decomposition over $\Q$: $\frakg= \frakl \ltimes \fraku$ such that $\frakh \subset \frakl$. Here $\fraku$ denotes the Lie algebra of the unipotent radical $\bmR_{\bmu}(\bmG)$ of $\bmG$ and $\frakl$ is the Lie algebra of a maximal reductive subgroup $\bmL$ of $\bmG$.

Decompose $\wedge^{\dim \bmH}\frakg = \wedge^{\dim \bmH}\frakl \oplus W$ where
\begin{equation*}
    W:=\left\{ v\in \wedge^{\dim \bmH}\frakg\;\middle\vert\;
         v \wedge \bmu= \bmzero ,\; \forall \, \bmu \in \wedge^{\dim\fraku} \fraku
    \right\}.
\end{equation*}
For an element $g\in \bmG(\C)$, written as $g= u\cdot l $ with $u \in \bmR_{\bmu}(\C)$ and $l\in \bmL(\C)$, we have
\begin{equation*}
   \Ad(g)\bmv_{\frakh} = \Ad(ul) \bmv_{\frakh} = \Ad(l) \bmv_{\frakh} + w
\end{equation*}
for some $w\in W$. 
Hence $\norm{\Ad(ul)\bmv_{\frakh}}$ is bounded from below by some positive multiple of $\norm{\Ad(l)\bmv_{\frakh}}$.
But $\Ad(\bmL(\C))\bmv_{\frakh}$ is closed and can never approach $0$.  So our claim follows.
\end{proof}

Therefore, if $\bmG$ and $\bmH$ are both semisimple and the centralizer of $\bmH$ in $\bmG$ is finite, the log anti-canonical bundle of $\bmX$ is big by \cite[Proposition 5.1]{HassettTanimotoTschinkel2015_balanced_line_bundle}.
The reductivity assumption on $\bmH$ can not be dropped (see Example \ref{example_SL2_Unipotent} below).  Also, in the second part, the assumption $\bmZ_{\bmG}(\bmH^{\circ})^{\circ}\subset \bmH$ is necessary. 
Nevertheless, we have 
\begin{thm}\label{theorem_anticanonical_general_arithmetic}
    Let $\Gamma_{\C} \subset \bmG(\Q(i))$ be an arithmetic subgroup of $\Res_{\Q/\Q(i)} (\bmG) $ and $\alpha \in \scrA$. Assume that for some $x\in \bmU(\Q(i))$ and sequence $(\gamma_n)$ in $\Gamma_{\C}$ one has $\lim_{n \to \infty} \gamma_n.x \in \bmD_{\alpha}(\C)^{\circ}$,  then $d_{\alpha}\geq 1$.
    Assume additional that $(\gamma_n)$ can be chosen to be unbounded modulo $\bmN_{\bmG}(\bmH)^{\bmone}(\C)\cap \Gamma_{\C}$, then $d_{\alpha}\geq  2$.
\end{thm}

Here $\Res_{\Q/\Q(i)}(\bmG)$ denotes the restriction of scalar of a linear algebraic group $\bmG$ (see \cite[11.4,12.4]{Spr98}), and $\bmN_{\bmG}(\bmH)^{\bmone}$ is the stabilizer of $\bmv_{\frakh}$ under the adjoint action. 

Related results are obtained in \cite[Theorem 2.7]{Hassett_Tschinkel_1999_equivariant_compactification_G_a} for bi-equivariant compactifications of the additive groups $\bmG_{a}^{n}$.

\begin{proof}

 Take $x\in \bmU(\Q(i))$ and let $\bmH_x$ be the stabilizer of $x$ in $\bmG$, which is a subgroup defined over $\Q(i)$. Let $\bmv_{\frakh_x} \in (\wedge^{\dim \bmH} \frakg) \otimes_{\Q} \Q(i)$ be a $\Q(i)$-vector lifting $\frakh_x$.
Define $\varphi^x_{\bmv}$ by modifying the definition of $\varphi_{\bmv}$ with $o$ replaced by $x$ and $\bmH$ replaced by $\bmH_x$. 
Let $\Gamma_{\C}$ be an arithmetic subgroup of $\Res_{\Q(i)/\Q}(\bmG)$ viewed as a subgroup of $\bmG(\C)$.
Note that $\chi \Delta_{\bmG}^{-1}$ is defined over $\Q(i)$ and thus $(\chi \Delta_{\bmG}^{-1})(\Gamma_{\C})$ is discrete in $\C$, implying that $(\chi \Delta_{\bmG}^{-1})(\Gamma_{\C})$ lies in the unit circle of $\C$, that is, $ \normm{\chi(\gamma)\Delta_{\bmG}(\gamma)^{-1}}=1$ for $\gamma \in \Gamma_{\C}$.
As $\Ad(\Gamma_{\C}) \bmv_{\frakh_x}$ is discrete in $\wedge^{\dim \bmH}\frakg \otimes {\C}$, by similar arguments as those used in the previous theorem, we  complete the proof.
\end{proof}

\begin{exam}\label{example_SL2_Unipotent}
Let $\bmG:=\SL_2$ and $\bmH$ be the upper triangular unipotent subgroup. Thus, $\bmG/\bmH$ has an invariant gauge form and can be identified with $\bmA^2_\Q \setminus \{0,0\}$.  Hence $\bmG/\bmH$ embeds in $\bmP_{\Q}^2$ as an open and dense subvariety.
The action of $\bmG$ on $\bmG/\bmH$ extends to $\bmP_{\Q}^2$.

    Let $\bmX$ be the blow-up of $\bmP^2_{\Q}$ at $\{[0:0:1]\}$, then $\bmX= \bmU \sqcup \bmP \sqcup \bmE$ where $\bmP$ is the proper transform of a hyperplane section and $\bmE$ is the exceptional divisor.
    By explicit computation, $- \divisor(\omega_{\bmG/\bmH}) = - \bmE + 3 \bmP$. Therefore, not all $d_{\alpha}$'s are positive here.
    Nevertheless, one can verify that $-\divisor (\omega_{\bmG/\bmH})$ is still an effective divisor. Indeed, $-\bmE$ is linearly equivalent to the proper transform of $\{y=0\}$.
    It is actually big, see \cite{Fu_Zhang_2013_compact_tori_automorphism} for a general statement.
\end{exam}

\subsection{Analytic Clemens complex}\label{subsection_analytic_Clemens_complex}
Keep the assumptions from last section. 
Also assume that  $\max\{d_{\alpha}\}>0$.
As a reference, see \cite[Section 3.1]{Chamber-Loir_Tschinkel_2010_Igusa_integral}. Since we are assuming $\bmD$ to be a simple normal crossing divisor rather than just normal crossing, the discussions are simplified. For instance, if $\sigma(\bmD_{\alpha}) \neq \bmD_{\alpha}$ for $\sigma_{\neq \id} \in \Gal(\C/\R)$, then $\bmD_{\alpha}(\C) \cap \bmX(\R) = \emptyset$.

Let the analytic irreducible components of $\bmD(\R)$ be indexed as $ (D_{\alpha} )$ for  $\alpha \in \scrA_{\R}^{\an}$.
There is a natural map from $\scrA_{\R}^{\an}$ to $\scrA$.
For a subset $I \subset \scrA^{\an}_{\R}$, let $D_I:= \cap_{\alpha \in I} D_{\alpha}$.
Since $\rmG:= \bmG(\R)^{\circ}$ is connected, each connected component of $D_{I}$ is $\rmG$-invariant.

We define a partially ordered set by
\begin{equation}
    \begin{aligned}
        &\scrC^{\an}_{\R}:=
        \left\{
         (I,Z) 
        \;\middle\vert \;
        I \subset \scrA^{\an}_{\R},\,
        Z \text{ is a connected component of }D_I \neq \emptyset
        \right\},\\
        & (I,Z)  \prec (I',Z') \iff I \subset I',\; Z'\subset Z.
    \end{aligned}
\end{equation}

For $x\in \bmU(\R)$, let 
\begin{equation*}
    \scrA_{\R,x}^{\an}:= \left\{
        \alpha \in \scrA^{\an}_{\R}\;\middle\vert\;
        D_{\alpha} \cap \overline{\rmG.x} \neq \emptyset
    \right\}.
\end{equation*}
And $\scrC^{\an}_{\R,x}$ is defined similarly as $\scrC^{\an}_{\R}$.
Note that $\rmG.x$ is exactly the connected component of $\bmU(\R)$ containing $x$.

Take $\bmL:= \sum_{\alpha \in \scrA} \lambda_{\alpha} \bmD_{\alpha}$ to be another divisor defined over $\R$ such that $ \lambda_{\alpha } \geq 0$ for all $\alpha \in \scrA$ and
\begin{equation*}
    d_{\alpha} >0 \implies \lambda_{\alpha}>0.
\end{equation*}
Thus points of bounded heights (associated with $\bmL$) have finite volume. For every $x\in \bmU(\R)$ such that  $\lambda_{\alpha} >0$ for some $\alpha\in \scrA_{\R,x}^{\an}$, let
\begin{equation*}
    a_x:= \max_{
       \alpha \in \scrA^{\an}_{\R,x},\,\lambda_{\alpha}\neq 0
    } \left\{
       \frac{d_{\alpha}-1}{\lambda_{\alpha}}
    \right\},
\end{equation*}
and 
\begin{equation*}
    \scrC^{\an}_{\R,x}(\bmL):= \left\{
            (I,Z)\in \scrC_{\R,x}^{\an}
            \;\middle\vert\;
             \frac{d_{\alpha}-1}{\lambda_{\alpha}} = a_x,\,\forall \, \alpha\in I
    \right\}.
\end{equation*}
Also let $b_{x}:= \dim \scrC^{\an}_{\R,x}(\bmL)$.

Let $\Vol$ denote the measure on $\bmU(\R)$ induced from $\omega_{\bmG/\bmH}^{\chi}$.
Let $\Ht_{\bmL}$ be the height at $\infty$ associated with certain smooth metric on $\calO_{\bmX}(\bmL)$ together with its canonical section.
For $R >0$, let $B_{R,x}:= \left\{ y \in \rmG.x, \;\Ht(y) < R \right\}$.
The following is proved in 
\cite[Theorem 4.7, Corollary 4.8]{Chamber-Loir_Tschinkel_2010_Igusa_integral}:
\begin{thm}\label{theorem_equidistribution_Chamber-Loir_Tschinkel}
Take $x\in \bmU(\R)$ such that  $\lambda_{\alpha} >0$ for some $\alpha\in \scrA_{\R,x}^{\an}$.
Assume $a_x \geq 0$ and $b_{x} \geq 1$.
 As $R$ tends to infinity,  there exists $c_x>0$ such that 
  \begin{equation*}
  \Vol (B_{R,x}) \sim
  \begin{cases}
        c_x R^{a_x} \log(R)^{b_x-1} \quad & a_x \neq 0\\
         c_x \log(R)^{b_x} \quad & a_x =0.
  \end{cases}
  \end{equation*}
  The family of probability measures on $\bmX(\R)$
 \begin{equation*}
     \nu_{R,x}:= \frac{
      \bmone_{B_{R,x}  } \cdot \Vol
     }{
     \Vol(B_{R,x})
     }
 \end{equation*}
  has a limit $\nu$ under the weak-$*$ topology. 
  Furthermore, $\nu$ is a sum of smooth measures on  $Z$ as $(I,Z)$ varies over faces of $\scrC^{\an}_{\R,x}(\bmL)$ of dimension $b_x$.
\end{thm}

\section{Equidistribution and Orbital Counting}\label{section_equidistribution_orbital_counting}

In this section, we explain the relationship between equidistribution and orbital counting. Compared to \cite{DukRudSar93,EskMcM93}, the possibility of ``focusing'' is also considered. On the other hand, we put more restrictions on the height functions, which do not apply to those coming from Riemannian metrics on symmetric spaces.

\subsection{Orbital counting follows from equidistribution}
Let $G$ be a locally compact second countable topological group and $H$ be a closed subgroup. Let $\Gamma$ be a discrete subgroup of $G$, $\rmm_{G}$ (resp. $\rmm_{H}$) be a left $G$-invariant locally finite measure on $G$ (resp. on $H$) and $o\in G/H$ denote the identity coset.

\begin{assumption}\label{assumption_invariant_measures}
We assume that 
\begin{itemize}
    \item[(1)] $G/\Gamma$ has a $G$-invariant locally finite measure $\rmm_{[G]}$;
    \item[(2)] $H\Gamma/\Gamma$ is closed in $G/\Gamma$ and supports an $H$-invariant locally finite measure $\rmm_{[H]}$;
    \item[(3)] there exists a $G$-invariant locally finite measure $\rmm_{G/H}$ on $G/H$ and that the triple $(\rmm_{G},\rmm_{H},\rmm_{G/H})$ is \textit{compatible}, that is, 
    \begin{equation}
         \int f(x)\, \rmm_{G}(x)
    =\int \int f(gh)\,\rmm_{H}(h) \,\rmm_{G/H}([g])
    \end{equation}
    for every compactly supported continuous function $f$ on $G$.
\end{itemize}
\end{assumption}

Item ($1$) always holds in the arithmetic setting.
By \cite[Lemma 1.4]{Ragh72}, item ($3$) above is equivalent to the coincidence of the modular functions $\Delta_G$  and  $\Delta_H$ on $H$. For instance, if $G$ is a semisimple Lie group, then $\Delta_G$ is trivial and such a $\rmm_{G/H}$ exists iff $\Delta_H$ is trivial.

Let $l:G/H \to \R_{> 0}$ be a continuous function and for a positive number $R$, define
\begin{equation*}
\begin{aligned}
        B_R:=&\left\{
    x \in G/H  \;\vert\; l(x)\leq R
    \right\}.
\end{aligned}
\end{equation*}
If we wish to emphasize the role of $l$, an upper index will be added $B_R^l$.
Assume that $\rmm_{G/H}(B_R)$ is finite for each $R$.

A family of $(B_R)_{R\geq 1}$ is said to be \textit{well-rounded} (see \cite[Proposition 1.3]{EskMcM93}) iff
for every $0<\ep<1$, there exists an open neighborhood $\calO_{\ep}$ of $\id\in G$ such that 
\begin{equation}\label{equation_well_rounded}
    (1-\ep) \rmm_{G/H} \left( 
    \bigcup_{g\in \calO_{\ep}} gB_R
    \right)
    \leq \rmm_{G/H}(B_R) 
    \leq 
     (1+\ep) \rmm_{G/H} \left( 
    \bigcap_{g\in \calO_{\ep}} gB_R
    \right).
\end{equation}

\begin{rmk}
This is different from the notion of (two-sided) well-roundedness as introduced in \cite{GoroNevo12Crelle} which does not seem to hold in many examples considered in this paper.
\end{rmk}

Here we assume that $\calO_{\ep}$ can be arranged such that additionally the following is true for all $R>1$:
\begin{equation}\label{equation_strong_well_rounded}
     B_{(1-\ep)R} \subset 
     \bigcap_{g\in \calO_{\ep}} gB_R 
     \subset 
     \bigcup_{g\in \calO_{\ep}} gB_R \subset B_{(1+\ep)R}.
\end{equation}
This property allows us to prove counting results in the presence of focusing and it is satisfied by all the heights considered in this paper.

A family $(C_R)_{R\geq 1}$ of positive real numbers is said to have a \textit{polynomial asymptotic} if there exist $c_1>0$, $c_2,c_3\geq 0$ such that $\left(
c_1R^{c_2}(\log{R})^{c_3}
\right)/C_R \to 1$ as $R \to +\infty$.

For the application when $\rmm_{[H]}$ is an infinite measure, a weighted variant is also considered. Assume that $\bmw: G/H \to \R_{>0}$  satisfies the following:
For every $\ep>0$, there exists an open neighborhood $\calO_{\ep}$ of $\id\in G$ such that
\begin{equation}\label{equation_weight_function}
         (1-\ep) \bmw_{gx} \leq \bmw_{x} \leq 
         (1+\ep) \bmw_{gx},\quad
         \forall x\in G/H,\; \forall g\in \calO_{\ep}.
\end{equation}
The weighted orbital counting problem is to study the asymptotic of the function $\Phi^{\bmw}_{R}: G/\Gamma \to \R$ defined by
\begin{equation*}
    \Phi^{\bmw}_{R}([g]) := 
    \frac{1}{C_R} \sum_{
       x \in g\Gamma.o \cap  B_R
    } \bmw_x
\end{equation*}
as $R$ tends to $+\infty$.

\begin{thm}\label{theorem_equidistribution_imply_counting}
    Assume that Assumption \ref{assumption_invariant_measures} is met, Equa.(\ref{equation_well_rounded},\ref{equation_strong_well_rounded},\ref{equation_weight_function}) hold and $(C_R)$ has polynomial asymptotic.
    Assume further that there exists a nonzero locally finite measure $\mu_{\infty}$ on $G/\Gamma$ such that in the weak-$*$ topology,
    \begin{equation}\label{Equation_equidistribution_average_weighted_homogeneous_known}
        \lim_{R\to \infty}  \frac{1}{C_R}
        \int_{[g]\in B_R}  \left( \bmw_{g.o} g_* \rmm_{[H]} \right) \,\rmm_{G/H}([g]) = \mu_{\infty}.
    \end{equation}
    Then there exists a positive continuous function $f_{\infty}$  such that $\mu_{\infty} = f_{\infty}\cdot \rmm_{[G]}$ and $\lim_{R\to \infty} \Phi^{\bmw}_{R}([g]) = f_{\infty}([g])$ for every $[g]\in G/\Gamma$.
    Actually, if Equa.(\ref{Equation_equidistribution_average_weighted_homogeneous_known}) holds only against compactly supported continuous functions supported in some fixed open neighborhood of the identity coset, then we still have that $\mu_{\infty} = f_{\infty}\cdot \rmm_{[G]}$ when restricted to this neighborhood with $f_{\infty}$ continuous and  $\lim_{R\to \infty} \Phi^{\bmw}_{R}([g]) = f_{\infty}([g])$ in this neighborhood.
\end{thm}

\subsection{Proof of the theorem}

The proof follows similar lines as \cite[Section 5]{EskMcM93} or \cite[Section 2]{DukRudSar93} with some modifications.

\subsubsection{Step 1, limit exists and is continuous.}

We first show that the limit $\lim_{R\to \infty} \Phi^{\bmw}_{R}([g_0])$ exists and is continuous.
Fix $[g_0]\in G/\Gamma$, 

For $0<\ep<1$, choose $\calO_{\ep}$ satisfying  
Equa.(\ref{equation_well_rounded}, \ref{equation_strong_well_rounded}, \ref{equation_weight_function}). By shrinking to a smaller one we assume that $\calO_{\ep}=\calO_{\ep}^{-1}$.

Then we choose $V'_{\ep}\subset V_{\ep}\subset \calO_{\ep}$ to be two families of open neighborhoods of identity such that the closure of $V'_{\ep}$ is contained in $V_{\ep}$ and
     \begin{equation*}
        \lim_{\ep\to 0} 
        \frac{\mu_{\infty}(V_{\ep}'[g_0])}
        {\mu_{\infty}(V_{\ep}[g_0])} =
        \lim_{\ep\to 0} 
        \frac{\rmm_{[G]}(V_{\ep}'[g_0])}
        {\rmm_{[G]}(V_{\ep}[g_0])}
        =1.
    \end{equation*}
Also choose a continuous function $f_{\ep}$ with 
$1_{V'_{\ep}[g_0]}\leq f_{\ep} \leq 1_{V_{\ep}[g_0]}$.

By an unfolding and folding argument as in \cite{DukRudSar93, EskMcM93},   we have
\begin{equation*}\label{equation_folding_unfolding}
\begin{aligned}
    \la \Phi^{\bmw}_R, f_{\ep} \ra_{\rmm_{[G]}}
    := &
    \int_{x\in G/\Gamma} \Phi^{\bmw}_R(x) f_{\ep}(x)\, \rmm_{[G]}(x)
    \\
    =& \frac{1}{C_R} \int_{[g]\in B_R} \bmw_{g.o}
   \left( \int_{x\in G/\Gamma} f_{\ep}(gx) \rmm_{[H]}(x) \right)
   \rmm_{G/H}([g])
   \\
   =& \int_{x\in G/\Gamma} f_{\ep}(x) \mu_R(x)
\end{aligned}
\end{equation*}
where 
\begin{equation*}
    \mu_R:= \frac{1}{C_R}
        \int_{[g]\in B_R} \left(\bmw_{g.o} g_* \rmm_{[H]} \right) \rmm_{G/H}([g]) .
\end{equation*}
Therefore,
    \begin{equation}\label{EquaProofProp2.2}
    \begin{aligned}
         \limsup_{R\to \infty}  \la \Phi_R^{\bmw}, 1_{V_{\ep}'[g_0]} \ra_{\rmm_{[G]}}  
        &\leq 
        \lim_{R\to \infty}  \la \Phi_R^{\bmw}, f_{\ep} \ra_{\rmm_{[G]}}  \\
        &=  \int f_{\ep}(x) \mu_{\infty}(x)
        \leq \liminf_{R\to \infty} \la \Phi_R^{\bmw}, 1_{V_{\ep}[g_0]} \ra_{\rmm_{[G]}}.
    \end{aligned}
    \end{equation}
    Also note that
    \begin{equation}\label{equation_estimate_f_epsilon}
        \mu_{\infty} (V_{\ep}'[g_0])
        \leq 
        \int f_{\ep}(x) \mu_{\infty}(x)
        \leq 
        \mu_{\infty} (V_{\ep}[g_0]).
    \end{equation}
     Assume $\ep$ to be small enough so that the natural map $\calO_{\ep}\to \calO_{\ep}[g_0]$ is a homeomorphism. Then $\lambda_0\rmm_{G}\vert_{\calO_\ep}$ is identified with $\rmm_{[G]}\vert_{\calO_{\ep}[g_0]}$ under this homeomorphism for some $\lambda_0>0$.

     Now we can start to estimate:
    \begin{equation*}
        \begin{aligned}
                \la \Phi_R^{\bmw}, 1_{V_{\ep}'[g_0]} \ra_{\rmm_{[G]}} =& 
                \frac{1}{C_R} \int _{[g]\in V_{\ep}'[g_0]} 
                \left(\sum_{x\in g\Gamma. o \cap B_R} \bmw_{x}\right)
                \,\rmm_{[G]}([g])
                \\
                =&
                \frac{1}{C_R} \int _{g\in V_{\ep}'} 
                 \left(\sum_{x\in g_0\Gamma. o \cap g^{-1} B_R} \bmw_{gx}  \right)
                \,\lambda_0\rmm_{G}(g)
                \\
                (\text{by Equa.}(\ref{equation_weight_function}))\;
                \geq&
               \frac{ \rmm_{[G]}(V_{\ep}'[g_0])
               }
               {(1+\ep) C_R} 
               \left(\sum_{x\in g_0\Gamma. o \cap (\cap_{g\in V'_{\ep}} g^{-1} B_R) } \bmw_{x}  \right)
                \\
                (\text{by Equa.}(\ref{equation_strong_well_rounded}))\;
                \geq& 
                \frac{ \rmm_{[G]}(V_{\ep}'[g_0])
               }
               {(1+\ep) C_R} 
               \left(\sum_{x\in g_0\Gamma. o \cap  B_{(1-\ep)R} } \bmw_{x}  \right)\\
               (\text{by definition of }\Phi^{\bmw}_R)\;
               =&
               \frac{ \rmm_{[G]}(V_{\ep}'[g_0])}{1+\ep}
               \frac{C_{(1-\ep)R}}{C_{R}}
                 \Phi^{\bmw}_{(1-\ep)R}([g_0]).
        \end{aligned}
    \end{equation*}
    By taking the $\limsup_{R \to \infty}$ and making use of Equa.(\ref{EquaProofProp2.2}, \ref{equation_estimate_f_epsilon}), we get 
    \begin{equation*}
    \begin{aligned}
         & \rmm_{[G]}(V_{\ep}'[g_0])\limsup_{R\to \infty}
         \frac{C_{(1-\ep)R}}{C_{R}}
        \Phi^{\bmw}_{(1-\ep)R}([g_0])
        \\
        &\leq (1+\ep)\limsup_{R\to \infty} \la \Phi_R^{\bmw}, 1_{V_{\ep}'[g_0]} \ra_{\rmm_{[G]} }
        \leq   (1+\ep)\mu_{\infty} (V_{\ep}[g_0]).
    \end{aligned}
    \end{equation*}
    Similarly we have 
    \begin{equation*}
         \rmm_{[G]}(V_{\ep}[g_0])\liminf_{R\to \infty}
          \frac{C_{(1+\ep)R}}{C_{R}}
        \Phi^{\bmw}_{(1+\ep)R}([g_0])
        \geq   (1-\ep) \mu_{\infty} (V'_{\ep}[g_0]).
    \end{equation*}
    Combining these two while replacing $(1-\ep)R$ by $R$ in the first inequality and $(1+\ep)R$ by $R$ in the second inequality, we get
    \begin{equation*}
        \begin{aligned}
            &
            (1-\ep) \liminf_{R \to \infty} \frac{C_R}{C_{R(1+\ep)^{-1}}}
            \frac{
            \mu_{\infty}(V'_{\ep}[g_0])
            }{
            \rmm_{[G]}(V_{\ep}[g_0])
            } 
            \leq 
             \liminf_{R\to \infty} \Phi^{\bmw}_R([g_0])  \\
             &\leq  \limsup_{R\to \infty} \Phi^{\bmw}_R([g_0]) \leq
              (1+\ep) \limsup_{R\to \infty} \frac{C_{R(1-\ep)^{-1}}}{C_{R}}
            \frac{
            \mu_{\infty}(V_{\ep}[g_0])
            }{
            \rmm_{[G]}(V'_{\ep}[g_0])
            } 
        \end{aligned}
    \end{equation*}
   By our assumptions, the $\limsup_{\ep \to 0}$ of the left end coincides with the $\liminf_{\ep \to 0}$ of the right end. Hence we must have
   \begin{equation*}
       \lim_{R\to\infty} \Phi_{R}^{\bmw}([g_0]) = \lim_{\ep \to 0}
       \frac{
            \mu_{\infty}(V_{\ep}[g_0])
            }{
            \rmm_{[G]}(V_{\ep}[g_0])
            } .
   \end{equation*}
    Call this limit $f([g_0])$. It is continuous because for $u\in \calO_{\ep}$,
    \begin{equation*}
    \begin{aligned}
        f(u[g_0])= &
        \lim_{R\to \infty} \frac{1}{C_R}
        \sum_{x \in 
        ug_0\Gamma.o\cap B_R 
        } \bmw_{x}
        = \lim_{R\to \infty} \frac{1}{C_R}
        \sum_{x \in 
        g_0\Gamma.o\cap u^{-1}B_R 
        } \bmw_{ux}
        \\
        \geq &
        \lim_{R\to \infty} \frac{1}{C_R}
        \sum_{x \in 
        g_0\Gamma.o\cap B_{(1-\ep)R}
        } \bmw_{x} (1+\ep)^{-1}
        \geq 
        (1+\ep)^{-1}
        \liminf_{R\to \infty} \frac{C_{(1-\ep)R}}{C_R}
        f([g_0]),
    \end{aligned}
    \end{equation*}
   which converges to $f ([g_0] )$ as $\ep \to 0$. A similar argument shows that the limit is no larger than $f ([g_0])$. Thus $\lim_{u \to \id}f(u[g_0]) =f([g_0])$, proving the continuity.
    
    By similar arguments, it is not hard to see that for any two $g_1,g_0$ in $G$, there exists a positive constant $c_{g_1,g_0}$ such that $f([g_1g_0])\geq c_{g_1,g_0}f([g_0])$. Therefore $f$ is either constantly equal to $0$ (which will be excluded below) or strictly positive.

    \subsubsection{Step 2, absolute continuity.}
    Now we switch our attention to $\mu_{\infty}$ and show that it is absolutely continuous with respect to $\rmm_{[G]}$.
    That is to say, we need to show $\mu_{\infty}(E)=0$ whenever $\rmm_{[G]}(E)=0$.

    Recall that $\rmm_{G}$ is the Haar measure on $G$ whose restriction to a fundamental set induces some constant multiple, which we may just assume to be $1$ for simplicity, of $\rmm_{[G]}$.
     By a Fubini-type argument, we first observe that if $\rmm_{[G]}(E)=0$, then for any probability measure $\lambda$ on $G/\Gamma$ and any bounded nonempty open set $U$ in $G$, we have
    \begin{equation*}
        \left(\int_{u\in U} u_*\lambda \,\rmm_{G}(u)
        \right)
        (E)=0.
    \end{equation*}
    
     Now take a bounded measurable set $E$ with $\rmm_{[G]}(E)=0$, there exists a family of shrinking bounded measurable sets $(E_i)_{i \in \Z^+}$ such that $E=\cap E_i$ and $\mu_{\infty}(\partial E_i)=0$. For $u\in \calO_{\ep}$,
    \begin{equation*}
    \begin{aligned}
          u_*\mu_{\infty}(E_i) 
          & =
          \lim_{R\to \infty} \frac{1}{C_R}\int_{[g]\in u\cdot B_R}g_*\rmm_{\rmH} (E_i)
          \bmw_{u^{-1}[g]}
          \,\rmm_{G/H}([g]) \\
          &\leq 
          (1-\ep)^{-1}
          \lim_{R\to \infty} \frac{1}{C_R}
          \int_{[g]\in B_{(1+\ep)R} }g_*\rmm_{\rmH} (E_i)
          \bmw_{ [g] }
          \,\rmm_{G/H} ([g]) \\
          &\leq 
           (1-\ep)^{-1}
         \left( \limsup_{R\to \infty} \frac{C_{(1+\ep)R}}{C_R} \right)
          \mu_{\infty}(E_i).
    \end{aligned}
    \end{equation*}
     A lower bound for $u_*\mu_{\infty}(E_i) $ can be similarly obtained. Thus we have positive number $\ep'$, which converges to $0$ as $\ep$ does so, such that for every $u\in \calO_{\ep}$,
     \begin{equation*}
        (1-\ep') \mu_{\infty}(E_i) \leq 
         u_*\mu_{\infty}(E_i)  \leq (1+\ep') \mu_{\infty}(E_i).
     \end{equation*}
     Hence
     \begin{equation*}
    \begin{aligned}
        &\left| \mu_{\infty}(E_i)- \frac{1}{\rmm_{G}(\calO_{\ep})}\int_{\calO_{\ep}} u_*\mu_{\infty}(E_i) \,\rmm_{G}(u) \right|
        \leq \ep' \mu_{\infty}(E_i).
    \end{aligned}
    \end{equation*}
    Letting $i$ tend to infinity, we get
    \begin{equation*}
         \left|\mu_{\infty}(E)- \frac{1}{\rmm_{G}(\calO_{\ep})}\int_{u\in \calO_{\ep}} u_*\mu_{\infty}(E) \rmm_{G}(u) \right| \leq \ep'\mu_{\infty}(E).
    \end{equation*}
    However, by our observation, the measure of $E$ with respect to ``Haar average'' of any probability measure is $0$. So we are left with
    \begin{equation*}
         \mu_{\infty}(E) \leq \ep'\mu_{\infty}(E_i)
    \end{equation*}
    for $\ep'$ arbitrarily small, which forces $\mu_{\infty}(E) =0$. This ends the proof of absolute continuity.
    
    \subsubsection{Step 3, completing the proof.}
    
    Write $\mu_{\infty}= \psi \cdot \rmm_{[G]}$ for some non-negative measurable function $\psi$.
    By what has been shown in Step 1, for any $[g_0]\in G/\Gamma$,
    \begin{equation*}
        f([g_0])=
        \lim_{\ep \to 0}
        \frac{\mu_{\infty}(V_{\ep}[g_0])}{\rmm_{[G]}(V_{\ep}[g_0])}
        =
        \lim_{\ep \to 0}
        \frac{1}{\rmm_{G}(V_{\ep})} \int_{u\in V_{\ep}} \psi(u[g_0]) \,\rmm_{G}(u).
    \end{equation*}
    Therefore $f=\psi$ almost surely. As $\mu_{\infty}=\psi\cdot  \rmm_{[G]}$  is nonzero, we see that $f$ can not be the $0$ function. As mentioned towards the end of step 1, this implies that $f$ is strictly positive.

\section{Equidistribution and Nondivergence}\label{section_equidistribution_nondivergence}

In this section we collect some results on equidistribution and nondivergence of homogeneous measures on a special type of homogeneous spaces: arithmetic quotients of real points of linear algebraic groups over $\Q$.
The structure of linear algebraic groups over $\Q$ and reduction theory for arithmetic subgroups are well-understood (see \cite{BorHar62, Spr98, Bor2019}).
For instance, arithmetic quotients have finite volume iff the algebraic group has no non-torsion $\Q$-characters.
The study of dynamics of subgroup action on such homogeneous spaces comes later. The full classification of unipotent-invariant ergodic measures is obtained by Ratner \cite{Rat91} with a shorter proof given by \cite{MarTom94}. The analysis of ergodic components of a unipotent-invariant measure is also possible by the linearization method developed in \cite{DanMar91}.
The results presented below (mainly taken from \cite{EskMozSha96, zhangrunlinCompositio2021, zhangrunlinDcds2022}) 
 rely on the above work as well as the theory of $(C,\alpha)$-functions of \cite{KleMar98} (see \cite{EskMozSha97} for a slightly different approach).

\subsection{Notations}\label{subsection_standing_assumptions_equidistribution_nondivergence}

Throughout this section we adopt the following notations:
\begin{itemize}
    \item For a linear algebraic group $\bmG$ over $\R$, write $\rmG:=\bmG(\R)^{\circ}$.  When it is understood that there is some ambient algebraic group $\bmG$ over $\R$,
    let $\rmL:=\bmL(\R)\cap \rmG$ for an $\R$-subgroup $\bmL$. If $\bmH$ is another $\R$-subgroup contained in $\bmL$, let $\rmL^{\rmH}:=\rmL^{\circ}\cdot \rmH$.
    \item  Let $\Gamma \subset \rmG \cap \bmG(\Q)$ be an arithmetic lattice;
    \item For an observable $\Q$-subgroup $\bmL$ of $\bmG$, let $\rmm_{[\rmL]}$ denote a locally finite $\rmL$-invariant measure supported on $\rmL\Gamma/\Gamma$. For other $\rmL^{\circ}\subset L \subset \rmL$, we similarly define $\rmm_{[L]}$, supported on $L\Gamma/\Gamma$.
    \item For an algebraic subgroup $\bmL$ of $\bmG$, let $\bmZ_{\bmG}(\bmL)$ be the centralizer of $\bmL$ in $\bmG$.
\end{itemize}
 We fix an integral structure (that is, a lattice) $\frakg_{\Z}$ on $\frakg$ and hence on its exterior powers. Whenever an arithmetic lattice is fixed, we further require that $\frakg_{\Z}$ is preserved by $\Ad(\Gamma)$.
 Also, Euclidean metrics on these spaces are fixed.
 For every parabolic $\Q$-subgroup $\bmP$, take $\bmv_{\bmP}$ to be the unique (up to $\pm 1$) primitive integral vector in $\wedge^{\dim \bmP} \frakg_{\Z}$ that lifts its Lie algebra.
Also define
 \begin{equation*}
     \begin{aligned}
             \scrP_{\bmH}:= &
            \left\{
              \text{proper parabolic } \Q\text{-subgroups of }\bmG \text{ containing } \bmH
            \right\}. \\
            \scrP^{\max}_{\bmH}:= &
            \left\{
              \text{maximal proper parabolic } \Q\text{-subgroups of }\bmG \text{ containing } \bmH
            \right\}.
    \end{aligned}
\end{equation*}

Recall that a subgroup $\bmL$ of a linear algebraic group $\bmG$ is said to be \textit{observable} iff there exists a finite-dimensional representation of $\bmG$ and a vector whose stabilizer subgroup is equal to $\bmL$. If $\bmL$ is assumed to be a $\Q$-subgroup, then $\rmL\Gamma$ is closed iff $\bmL$ is observable in $\bmG$ (\cite[Corollary 7]{Weiss98}).
Reductive groups and unipotent groups are automatically observable in any ambient group.
A more comprehensive treatment of observable subgroups can be found in \cite{Grosshans97}.
If one is only interested in the case when $\rmm_{[\rmH]}$ is finite, which implies the observability, then one can ignore the word observable or replace it by having finite invariant measure on the quotient by an arithmetic subgroup.

Related to observable subgroups is the notion of \textit{epimorphic subgroups}. A subgroup $\bmL$ of $\bmG$ is said to be epimorphic iff the smallest observable subgroup containing $\bmL$ is $\bmG$. Alternatively, for every $v$ in a finite-dimensional representation of $\bmG$, $v$ fixed by $\bmL$ implies that $v$ is fixed by $\bmG$.

\subsection{Subgroup convergence and equidistribution}

\begin{defi}\label{definition_convergence_subgroups}
    Let $\bmL$ be a connected linear algebraic group over $\Q$ and $(\bmM_n)$ be a sequence of connected observable $\Q$-subgroups of $\bmL$, we say that $(\bmM_n)$ \textit{converges} to $\bmL$ iff there is no proper connected observable subgroup of $\bmL$ that contains $\bmM_n$ for infinitely many $n$'s. 
\end{defi}

Note that $(\bmM_n)$ always converges after passing to a subsequence.
It is not true that if $(\bmM_n)$ converges to $\bmL$, then the associated homogeneous measures also converge. Below is a special case when this does happen.

\begin{thm}\label{theorem_equidistribution_EMS}
Let $\bmH$ and $\bmL$ be two connected observable $\Q$-subgroups of a connected linear algebraic group $\bmG$ over $\Q$.
Let $(\gamma_n)$ be a sequence in $\Gamma$ such that $(\gamma_n \bmH \gamma_n^{-1})$ converges to $\bmL$.
Then there exists a sequence of positive numbers $(a_n)$ such that
$\lim_{n \to \infty} \frac{1}{a_n}(\gamma_n)_* \rmm_{[\rmH^{\circ}]} = \rmm_{[\rmL^{\circ}]}$ under the weak-$*$ topology.
\end{thm}

This is proved in \cite[Theorem 1.3]{zhangrunlinCompositio2021} building on previous work of \cite{EskMozSha96} (see also 
\cite{OhSha14, ShaZhe18, zhangrunlinAnnalen2019}).

The asymptotic of $(a_n)$ is unique. Let $\psi$ be a compactly supported non-negative continuous function on $\rmG/\Gamma$ with $\la \psi, \rmm_{[\rmL^{\circ}]} \ra \neq 0$, the conclusion of last theorem unwraps to the following:
\begin{equation*}
    \lim_{n \to \infty}
    \frac{
     \la f, (\gamma_n)_*\rmm_{[\rmH^{\circ}]} \ra
    }{
     \la \psi, (\gamma_n)_*\rmm_{[\rmH^{\circ}]} \ra
    } = 
    \frac{
     \la f, \rmm_{[\rmL^{\circ}]} \ra
    }{
     \la \psi,\rmm_{[\rmL^{\circ}]} \ra
    }.
\end{equation*}

\subsubsection{Definition of $p^{\sta}$}

For a connected linear algebraic group $\bmG$, let $\bmG^{\ari}:=\left( \overline{\overline{\Gamma}} \right)^{\circ}$, the identity component of the Zariski closure of $\Gamma$.
Let $p^{\sta}$ be the natural quotient map $\bmG \to \bmG/ \bmG^{\ari}$. Actually  $\bmG^{\ari}$ is independent of the choice of the arithmetic subgroup $\Gamma$.
More explicitly, recall 
\begin{equation*}
    \bmG = \left( \bmG^{\semisimple} \cdot \bmZ(\bmG) \right) \ltimes \bmR_{\bmu}(\bmG)
    =  \left( (\bmG^{\cpt}\cdot \bmG^{\nc})
     \cdot \bmZ(\bmG)  \right)
     \ltimes \bmR_{\bmu}(\bmG).
\end{equation*}
Also, $\bmZ(\bmG)$ is an almost direct product $\bmZ(\bmG)^{\cpt} \cdot \bmZ(\bmG)^{\spl} \cdot \bmZ(\bmG)^{\ari}$ where $\bmZ(\bmG)^{\spl}$ is a $\Q$-split torus,  $\bmZ(\bmG)^{\cpt}$ is an $\R$-anisotropic $\Q$-torus and $\bmZ(\bmG)^{\ari}$ is the identity component of the Zariski closure of some/any arithmetic lattice of $\bmZ(\bmG)$. 
In terms of the character groups, these tori can be described as follows:

 Let $\Lambda_0 \leq \frakX^*(\bmZ(\bmG))$ be the primitive subgroup consisting of $\Q$-characters, i.e. those fixed by $\Gal_{\Q}$. It admits a unique $\Gal_{\Q}$-stable complementary primitive subgroup $\Lambda_1$ (in the sense that $\Lambda_0 \cap \Lambda_1 =\{0\}$ and $\Lambda_0+ \Lambda_1\leq \frakX^*(\bmZ(\bmG))$ has finite index). Then $\bmZ(\bmG)^{\spl}$ is the common kernel of characters in $\Lambda_1$. Similarly, let $\Lambda_{2}$ be the primitive subgroup of $\frakX^*(\bmZ(\bmG))$ fixed by the smaller subgroup $\Gal_{\overline{\Q} \cap \R }$ and $\Lambda_{3}$ be the $\Gal_{\overline{\Q} \cap \R }$-stable complement of $\Lambda_{2}$ in $\frakX^*(\bmZ(\bmG))$. As $\Gal_{\overline{\Q} \cap \R }$ is normal in $\Gal_{\Q}$, $\Lambda_2$ and $\Lambda_3$ are $\Gal_{\Q}$-stable. Then $\bmZ(\bmG)^{\cpt}$ corresponds to $\Lambda_2$ and $\bmZ^{\ari}$ corresponds to the primitive subgroup generated by $\Lambda_0+\Lambda_3$.

With these notations, $\bmG^{\ari}$ is the normal subgroup $(\bmG^{\nc}\cdot \bmZ(\bmG)^{\ari})\ltimes \bmR_{\bmu}(\bmG)$ by Borel density lemma.

\subsubsection{Criteria on convergence to the full group}

\begin{lem}\label{lemma_existence_convergence_to_full'}
Let $\bmG$ be a connected linear algebraic group over $\Q$, $\bmN$ be a connected normal $\Q$-subgroup of $\bmG$ and $\bmH$ be an observable $\Q$-subgroup contained in $\bmN$. Let $\Gamma$ be an arithmetic lattice of $\bmG$.
The followings are equivalent:
\begin{itemize}
        \item[(1)] There exists $(g_n)\subset \rmG$ such that 
        $\left[ (g_n)_* \rmm_{[\rmH^{\circ}]}  \right] \to \left[  \rmm_{[\rmN^{\circ}]} \right]$;
        \item[(2)] There exists $(\gamma_n)\subset \Gamma$ such that 
        $\left[ (\gamma_n)_* \rmm_{[\rmH^{\circ}]}  \right] \to \left[  \rmm_{[\rmN^{\circ}]} \right]$;
        \item[(3)] $p^{\sta}(\bmH)=p^{\sta}(\bmN)$ and any proper $\Q$-subgroup that is normalized by $\bmN \cdot \bmG^{\ari}$ and  
                    contains $\bmH^{\circ}$ must also contain $  \bmN$.
\end{itemize}
\end{lem}

\begin{proof}
Part (1) and (2) are easily seen to be equivalent. That  (2) implies (3) is quite direct.
It only remains to explain why (3) implies (2).

By \cite[Lemma 4.14]{zhangrunlinCompositio2021}, if (2) were not true, then we can find an observable $\Q$-subgroup $\bmL$ that is properly contained in $\bmN$, contains $\bmH^{\circ}$ and is normalized by some arithmetic subgroup $\Gamma'$. Hence $\bmL$ is normalized by $\bmG^{\ari}$. But $p^{\sta}(\bmH^{\circ})=p^{\sta}(\bmH)= p^{\sta}(\bmN)$, thus   $\bmL$ is normalized by $\bmN \bmG^{\ari}$, which can not be true by (3). 
\end{proof}

In the special case of $\bmN= \bmG$ we get
\begin{lem}\label{lemma_existence_convergence_to_full}
Let $\bmG$ be a connected linear algebraic group over $\Q$ and $\bmH$ be an observable $\Q$-subgroup. Let $\Gamma$ be an arithmetic lattice.
The followings are equivalent:
\begin{itemize}
        \item[(1)] There exists $(g_n)\subset \rmG$ such that 
        $\left[ (g_n)_* \rmm_{[\rmH^{\circ}]}  \right] \to \left[  \rmm_{[\rmG]} \right]$;
        \item[(2)] There exists $(\gamma_n)\subset \Gamma$ such that 
        $\left[ (\gamma_n)_* \rmm_{[\rmH^{\circ}]}  \right] \to \left[  \rmm_{[\rmG]} \right]$;
        \item[(3)] $p^{\sta}$ restricted to $\bmH$ is surjective and $\bmH^{\circ}$ is not contained in any proper normal $\Q$-subgroup of $\bmG$.
\end{itemize}
\end{lem}

\subsection{Nondivergence criteria}

\subsubsection{Nondivergence of homogeneous closed subsets}

Recall that a $\Q$-group $\bmL$  is  said to be $\Q$-anisotropic iff the set of $\Q$-cocharacters of $\bmL$ is finite, which is equivalent to $\bmL(\R)/\Gamma$ being compact for some/any arithmetic subgroup $\Gamma$ of $\bmL$. The following is from  \cite[Theorem 1.7]{zhangrunlinCompositio2021}.

\begin{thm}\label{theorem_nondivergence}
    Assume that $\bmH$ and $\bmG$ are reductive linear algebraic groups over $\Q$ and $\Gamma$ is an arithmetic subgroup of $\bmG$ contained in $\rmG$.
    If $\bmZ_{\bmG}(\bmH^{\circ})/\bmZ(\bmH)$ is $\Q$-anisotropic, then there exists a bounded subset $\calB_{\Gamma}\subset \rmG$ such that
    \begin{equation*}
        \rmG = \calB_{\Gamma} \cdot \Gamma \cdot \rmH^{\circ}.
    \end{equation*}
    Moreover, if $\bmL$ is an observable $\Q$-subgroup containing $\bmH$, then $\bmL$ is reductive.
\end{thm}

Here is a more general version, which is restated and proved in the appendix (see Theorem \ref{theorem_divergence_observable}). 
 \begin{thm}\label{theorem_divergence_observable_1}
 Let $\bmG$ be a connected linear algebraic group over $\Q$, $\bmH$ be a connected observable $\Q$-subgroup and
 $\Gamma$ be an arithmetic subgroup of $\bmG$ contained in $\rmG$.
 Fix a maximal reductive subgroup $\bmG^{\red}$ of $\bmG$ and write $\bmG=\bmG^{\red} \ltimes \bmR_{\bmu}(\bmG)$. Also fix a Cartan involution and hence a maximal compact subgroup $\rmK$ of $\bmG^{\red}(\R)$.
 Given a sequence $(g_n) \subset \rmG$, after passing to a subsequence, there exist a sequence $(h_n)\subset \rmH^{\circ}$, $(\gamma_n)\subset \Gamma$ and a parabolic $\Q$-subgroup $\bmP$ such that the following holds. Write $g_n h_n \gamma_n^{-1} = k_n a_n p_n$ using horospherical coordinates of  $(\bmP,\rmK)$.
 Then
 \begin{itemize}
     \item[(1)]  $(p_n)$ is bounded;
     \item[(2)]  $\alpha(a_n) \to 0$ for every $\alpha\in \Delta^{\red}(\bmA_{\bmP,\rmK},\bmP)$;
     \item[(3)]  if $(a_n)$ is unbounded, there exist a $\Q$-representation $\bmV$ of $\bmG$ factoring through $\bmG/\bmR_{\bmu}(\bmG)$ and $\bmv\in \bmV(\Q)$ such that the line spanned by $\bmv$ is preserved by $\bmP$, $\bmv$ is fixed by  $\gamma_n \bmH \gamma_n^{-1}$ for all $n$ and $\lim_{n\to \infty} a_n. \bmv = \bmzero$;
     \item[(4)] $(\gamma_n \bmH \gamma_n^{-1})$ strongly converges to some observable subgroup of $\bmG$.
 \end{itemize}
 \end{thm}

Strongly convergence is stronger than convergence from Definition \ref{definition_convergence_subgroups}.
\begin{defi}\label{definition_strong_convergence_subgroups}
    Let $(\bmA_n)$ be a sequence of connected $\Q$-subgroups of a connected linear algebraic group $\bmC$ over $\Q$ and $\bmE$ be another connected $\Q$-subgroup of $\bmC$. We say that $(\bmA_n)$ strongly converges to $\bmE$ if
    \begin{itemize}
        \item[1.] every $\bmA_n$ is contained in $\bmE$;
        \item[2.] for every subsequence $(\bmA_{n_k})$, $\bmE$ is the smallest $\Q$-subgroup containing all $\bmA_{n_k}$'s.
    \end{itemize}
\end{defi}

\begin{coro}
Let $\bmG,\bmH,\Gamma$ be the same as in last theorem. Let $(g_n)$ be a sequence in $\rmG$ such that $\left( g_n \rmH^{\circ}\Gamma/\Gamma \right)$ diverges topologically, then there exist a $\Q$-representation $\rho: \bmG \to \SL(\bmV)$, an integral structure $\bmV(\Z)$ on $\bmV$ and a sequence of $\rho(\bmH)$-fixed nonzero vectors $(v_n)\subset \bmV(\Z)$ such that $\rho(g_n).v_n \to \bmzero$.
\end{coro}

\subsubsection{Nondivergence of a bounded piece}

Let $\bmG$ be a connected linear algebraic group over $\Q$ and $\bmH$ be a $\Q$-subgroup.
For a nonempty open bounded subset $\calO_{\rmH}$ of $\rmH^{\circ}\Gamma/\Gamma$, let $\rmm_{\calO_{\rmH}}$ be the restriction of $\rmm_{[\rmH]}$ to $\calO_{\rmH}$ and $\rmm_{\calO_{\rmH}}^{\bmone}$ be the unique probability measure proportional to $\rmm_{\calO_{\rmH}}$.
Recall it follows from the definition that ${}^{\circ}\rmG= \left\{ g \in \rmG \midd \chi(g) =1 ,\; \forall \chi \in \frakX^*_{\Q}(\bmG)  \right\}$.

\begin{thm}\label{theorem_nondivergence_bounded_piece}
    Fix a nonempty open bounded subset $\calO_{\rmH}$ of $\rmH^{\circ}\Gamma/\Gamma$. Then
    \begin{itemize}
        \item[(1)] for every bounded subset $\calB \subset \rmG/\Gamma$, there exists $\eta_{\calB}>0$ such that for every $g\in {}^{\circ}\rmG$,
        \begin{equation*}
            g\calO_{\rmH}\cap \calB \neq \emptyset \implies
            \norm{\Ad(g)\bmv_{\bmP}} > \eta_{\calB},\;\forall \, \bmP\in \scrP^{\max}_{\bmH^{\circ}}.
        \end{equation*}
         \item[(2)] for every $\eta>0$, there exists a bounded subset $\calB_{\eta} \subset \rmG/\Gamma$ such that for every $g\in {}^{\circ}\rmG$,
        \begin{equation*}
            \norm{\Ad(g)\bmv_{\bmP}} > \eta,\;\forall \, \bmP\in \scrP^{\max}_{\bmH^{\circ}}
             \implies g\calO_{\rmH}\cap \calB_{\eta} \neq \emptyset.
        \end{equation*}
        For every $\ep\in (0,1)$, by enlarging $\calB_{\eta}$ if possible, one can arrange that
        for every $g\in {}^{\circ}\rmG$,
        \begin{equation*}
            \norm{\Ad(g)\bmv_{\bmP}} > \eta,\;\forall \, \bmP\in \scrP^{\max}_{\bmH^{\circ}}
             \implies (g_* \rmm^{\bmone}_{\calO_{\rmH}}) \left( \calB_{\eta} \right) > 1-\ep.
        \end{equation*}
    \end{itemize}
\end{thm}

 Part ($1$) follows directly. It remains to prove part ($2$). Due to the $(C,\alpha)$-good property, it suffices to prove the first claim of part ($2$). The special case when $\bmG$ is semisimple is proved in \cite[Theorem 1.3]{zhangrunlinDcds2022}. In order to reduce the general case to this, we need a little preparation.

\subsubsection{Parabolic subgroups and radicals}

Let $r$ denote the dimension of $\bmR(\bmG)$, the radical of $\bmG$ (i.e., the maximal connected normal solvable subgroup of $\bmG$).  
Let $p^{\semisimple}: \bmG \to \bmG/\bmR(\bmG)$ be the natural quotient morphism. For a non-negative integer $l$, let $\wedge^l \diff p^{\semisimple}$ be the induced morphism between exterior powers of Lie algebras $\wedge^l \frakg \to \wedge^l \frakg/\frakr(\frakg) $.
Fix an integral vector $\bmv_{\frakr} \in \wedge^{r} \frakg$ that lifts the Lie algebra $\frakr({\frakg})$ of $\bmR(\bmG)$. 

For every $\omega \in \wedge^l \frakg/\frakr(\frakg)$, choose some $\wtomega \in  
(\diff p^{\semisimple})^{-1} (\omega) \subset 
\wedge^l\frakg$. Define a map $\varphi$ from $\wedge^l \frakg/\frakr(\frakg)$ to $\wedge^{l+r} \frakg$ by
\begin{equation*}
    \omega \mapsto \varphi(\omega):= \wtomega \wedge \bmv_{\frakr}.
\end{equation*}
One can check that $\varphi(\omega)$ is independent of the choice of $\wtomega$ and $\varphi$ is a linear injection. Thus, there exists $C>1$ such that
\begin{equation}\label{equation_norm_bounded_quotient_by_radical}
    C^{-1} \norm{\varphi(\omega)}
    \leq \norm{ \omega }
    \leq C \norm{ \varphi(\omega) },\quad
    \forall\,
    \omega \in \wedge^l \frakg/\frakr(\frakg).
\end{equation}
Moreover,  $\varphi$ is ${}^{\circ}\rmG$-equivariant since ${}^{\circ}\rmG$ fixes the vector $\bmv_{\frakr}$.

In addition to this, we need another simple fact:
\begin{lem}
     Every parabolic subgroup of $\bmG$ contains  $\bmR(\bmG)$. Consequently, the map $\bmP \mapsto \bmP/\bmR(\bmG)$ defines a bijection between parabolic $\Q$-subgroups of $\bmG$ and those of $\bmG/\bmR(\bmG)$. It also induces a bijection between $\scrP^{\max}_{\bmH^{\circ}}$ and $\scrP^{\max}_{p^{\semisimple}(\bmH^{\circ})}$.
\end{lem}
By \cite[Theorem 6.2.7]{Spr98}, a parabolic subgroup contains some Borel subgroup and hence $\bmR(\bmG)$. The rest of the claim follows from this.

For $\bmP \in \scrP^{\max}_{p^{\semisimple}(\bmH^{\circ})}$, let $\bmQ:=  (p^{\semisimple})^{-1}(\bmP) \in \scrP^{\max}_{\bmH^{\circ}}$. Note that $\bmv_{\bmQ}=c_{\bmP}\varphi(\bmv_{\bmP})$ for some positive constant $c_{\bmP}$ bounded away from $0$ and $+\infty$ as $\bmP$ varies.
Therefore, by  Equa.(\ref{equation_norm_bounded_quotient_by_radical}), we find some $C_1>1$ such that for every $\bmP \in \scrP^{\max}_{p^{\semisimple}(\bmH^{\circ})}$ and $g\in {}^{\circ}\rmG$,
\begin{equation}\label{equation_compare_parabolic_semisimple_quotient}
    C_1^{-1}\norm{\Ad(g)\bmv_{\bmQ}} \leq 
    \norm{
    \Ad(p^{\semisimple}(g) ) \bmv_{\bmP}
    }
    \leq C_1 \norm{\Ad(g)\bmv_{\bmQ}}.
\end{equation}

\subsubsection{Proof of Theorem \ref{theorem_nondivergence_bounded_piece}}
For simplicity, write $\bmG':= \bmG/\bmR(\bmG)$ and $\rmG':= \bmG'(\R)^{\circ}$.

The morphism $p^{\semisimple}$ induces $\overline{p}^{\semisimple}: \rmG/\Gamma \to  \rmG'/p^{\semisimple}(\Gamma)$.
By \cite[Theorem 1.3]{zhangrunlinDcds2022}, there exists a bounded subset $ \calB_{\eta}' \subset \rmG'/p^{\semisimple}(\Gamma)$ such that
for every $g\in \rmG'$,
\begin{equation}\label{equation_nondivergence_semisimple}
    \norm{
    \Ad(g) \bmv_{\bmP}
    } > \eta, \;\forall\, \bmP\in \scrP^{\max}_{p^{\semisimple}(\bmH^{\circ})}
    \implies g \overline{p}^{\semisimple}(\calO_{\rmH}) \cap \calB_{\eta}' \neq \emptyset.
\end{equation}
On the other hand, take a bounded open subset $\calB^{\spl}\subset \bmS_{\bmG}(\R)^{\circ}$ containing $\overline{p}^{\spl}(\calO_{\rmH})$.
Let $\overline{p}^{\semisimple}_{\calB}$ be the restriction of $\overline{p}^{\semisimple}$ to the closure of $(\overline{p}^{\spl})^{-1 }\left(
\calB^{\spl}
\right)$. 
Then $\overline{p}^{\semisimple}_{\calB}$ is a proper map.
 Define a bounded subset of $\rmG/\Gamma$ by
\begin{equation*}
    \calB_{\eta}:=
    (\overline{p}^{\semisimple})^{-1} \left(
           \calB'_{C_1^{-1} \eta}
    \right) \bigcap
    (\overline{p}^{\spl})^{-1}
    \left(
    \calB^{\spl}
    \right).
\end{equation*}

Now fix some $g\in {}^{\circ}\rmG$ such that $\norm{\Ad(g)\bmv_{\bmP}}> \eta$ for every $\bmP\in \scrP^{\max}_{\bmH^{\circ}}$. 
Since $\overline{p}^{\spl}(g\calO_{\rmH}) = \overline{p}^{\spl}(\calO_{\rmH}) $, the containment $g\calO_{\rmH} \subset (\overline{p}^{\spl})^{-1}
    \left(
    \calB^{\spl}
    \right)$  is assured.
For the other part, by Equa.(\ref{equation_compare_parabolic_semisimple_quotient}),
\begin{equation*}
    \norm{
    \Ad(p^{\semisimple}(g)) \bmv_{\bmP} 
    }> C_1^{-1}\eta
     ,\quad
    \forall \, \bmP \in \scrP^{\max}_{p^{\semisimple}(\bmH^{\circ})}.
\end{equation*}
Hence we are done by Equa.(\ref{equation_nondivergence_semisimple}).

\subsection{Refined versions}

In this subsection, $\bmG$ is a connected linear algebraic group over $\Q$ and $\bmH$ is a $\Q$-subgroup. Let $\bmH^{\bmone}:={}^{\circ}\bmG \cap \bmH$. Assume that the natural injective morphism $\bmH/\bmH^{\bmone} \to \bmS_{\bmG}$ is an isomorphism. 

\subsubsection{Explicit form of $(a_n)$} \label{subsection_explicit_a_n}

We choose a lift $\wtbmS^{\bmH}$, a $\Q$-split torus, of $\bmS_{\bmG}$ in $\bmH$.
Then the natural product map
\begin{equation*}
    \wtbmS^{\bmH}(\R)^{\circ} \times \bmH^{\bmone} (\R)^{\circ} \to \bmH(\R)^{\circ}
\end{equation*}
is a homeomorphism. From this we deduce that  the natural map
\begin{equation}\label{equation_decomposition_H_Gamma_lift_S_G}
    \wtbmS^{\bmH}(\R)^{\circ} \times \bmH^{\bmone} (\R)^{\circ}/ \bmH^{\bmone} (\R)^{\circ}\cap \Gamma 
    \to \bmH(\R)^{\circ}/\bmH(\R)^{\circ}\cap \Gamma
\end{equation}
is a homeomorphism.
Based on this, for $g\in {}^{\circ}\rmG $ and $\eta>0$, define
\begin{equation}\label{equation_decomposition_HGamma}
    \Omega_{\bmH,g,\eta}^{\bmone} := 
    \left\{
       [h]\in \bmH^{\bmone}(\R)^{\circ}/\bmH^{\bmone}(\R)^{\circ}\cap \Gamma \;\middle\vert\;
       \norm{
       \Ad(gh)\bmv_{\bmP}
       } >\eta,\;\forall\,\bmP \in \scrP^{\max}_{\bmH^{\circ}}
    \right\}.
\end{equation}

\begin{thm}\label{theorem_equidistribution_explicit_a_n}
Assume additionally that $\bmH$ is observable.
    Let $(g_n)$ be a sequence in ${}^{\circ}\rmG$ such that $\lim_{n \to \infty} \left[ (g_n)_* \rmm_{[\rmH^{\circ}]} \right] = \left[ \rmm_{[\rmG]} \right]$. Then for every $\eta>0$, 
    \begin{equation}\label{equation_convergence_to_full_Haar_explicit_a_n}
        \lim_{n\to \infty}
        \frac{1}{
        \rmm_{[(\rmH^{\bmone})^{\circ}]} \left(
                \Omega^{\bmone}_{\bmH,g_n,\eta}
        \right)
        }
        (g_n)_* \rmm_{[\rmH^{\circ}]}  = \rmm_{[\rmG]}.
    \end{equation}
\end{thm}

In Equa.(\ref{equation_convergence_to_full_Haar_explicit_a_n}), the Haar measures are chosen as follows:
Choose some Haar measure $\rmm_{\rmS}$ on $\wtbmS^{\bmH}(\R)^{\circ}$ and $\rmm_{[(\rmH^1)^{\circ}]}$ on $\bmH^{\bmone}(\R)^{\circ}/\bmH^{\bmone}(\R)^{\circ}\cap \Gamma$. By Equa.(\ref{equation_decomposition_H_Gamma_lift_S_G}), define $\rmm_{[\rmH^{\circ}]}:= \rmm_{\rmS} \otimes \rmm_{[(\rmH^1)^{\circ}]}$. Similar to Equa.(\ref{equation_decomposition_H_Gamma_lift_S_G}), one has 
\begin{equation*}
    \wtbmS^{\bmH}(\R)^{\circ} \times {}^{\circ}\rmG/  \Gamma 
    \cong \rmG/ \Gamma,
\end{equation*}
and we define $\rmm_{[\rmG]}:= \rmm_{\rmS} \otimes \rmm_{[{}^{\circ}\rmG]}^{\bmone} $.
The truth of Equa.(\ref{equation_convergence_to_full_Haar_explicit_a_n}) is then independent of the choice of $\rmm_{\rmS}$ and  $\rmm_{[(\rmH^1)^{\circ}]}$.

\subsubsection{Two lemmas}

The proof of Theorem \ref{theorem_equidistribution_explicit_a_n} is essentially contained in \cite{zhangrunlinCompositio2021} by using Theorem \ref{theorem_nondivergence_bounded_piece}. We outline a proof here for the sake of completeness.
Theorem \ref{theorem_equidistribution_explicit_a_n} can be directly deduced from the following two lemmas:

\begin{lem}\label{lemma_asymptotic_volume_polytope}
    Assumption as in Theorem \ref{theorem_equidistribution_explicit_a_n}. There exists a sequence of positive real numbers $(\eta_n)$ diverging to $+\infty$ such that for any $\eta>0$,
    \begin{equation*}
        \lim_{n\to \infty}
        \frac{
        \rmm_{[(\rmH^{\bmone})^{\circ}]} \left(
                \Omega^{\bmone}_{\bmH,g_n,\eta_n}
        \right)
        }{
        \rmm_{[(\rmH^{\bmone})^{\circ}]} \left(
                \Omega^{\bmone}_{\bmH,g_n,\eta}
        \right)
        } =1.
    \end{equation*}
\end{lem}

\begin{lem}\label{lemma_deep_in_polytope_equidistribute}
     Assumption as in Theorem \ref{theorem_equidistribution_explicit_a_n}.
     Let $(\eta_n)$ be a sequence of real numbers diverging to $+\infty$. Then for every sequence $(h_n)\subset \Omega^{\bmone}_{\bmH,g_n,\eta_n}$ and nonempty open bounded subset $\calO_{\rmH}\subset \rmH^{\circ}\Gamma/\Gamma$, 
     \begin{equation*}
         \lim_{n\to \infty}
         (g_nh_n)_* \rmm^{\bmone}_{\calO_{\rmH}} =
         \rmm^{\bmone}_{\calO_{\rmH}^{\spl}}
     \end{equation*}
     where 
     \begin{equation*}
         \calO_{\rmH}^{\spl}:= (\overline{p}^{\spl})^{-1}\left(
            \overline{p}^{\spl}(\calO_{\rmH})
         \right).
     \end{equation*}
\end{lem}

In the above, $\rmm_{\calO_{\rmH}^{\spl}}$ denotes the restriction of $\rmm_{[\rmG]}$ to  ${\calO_{\rmH}^{\spl}}$. The upper index $(\cdot )^{\bmone}$ is used to denote the unique probability measure  proportional to $(\cdot)$.

\subsubsection{Proof of Lemma \ref{lemma_deep_in_polytope_equidistribute}}

By Theorem \ref{theorem_nondivergence_bounded_piece}, $g_nh_n \omega_n = \delta_n \gamma_n $ for some bounded sequence $(\omega_n) \subset \rmH$ lifting elements in $\calO_{\rmH}$, some bounded sequence $(\delta_n)$ in $\rmG$ and $(\gamma_n)\subset \Gamma$.
Write $\delta_n = \delta_n^1 \delta_n^2$ where $\delta_n^1 \in \wtbmS_{\bmH}(\R)^{\circ}$ and $\delta_n^2 \in {}^{\circ}\rmG$. Both $(\delta_n^1)$ and $(\delta_n^2)$ are bounded.
Passing to a subsequence, may assume $( \omega_n )$ converges to $\omega_{\infty}$ and $(\delta_n^1)$ converges to $\delta^1_{\infty}$. 
Note that $\delta_{\infty}^1 \omega_{\infty}^{-1} \in {}^{\circ}\rmG$.
Thus it is sufficient to establish the claim with $g_n$ replaced by $\gamma_n$.

By the work of \cite{EskMozSha96} (see \cite[Section 4]{zhangrunlinCompositio2021} for details), it is sufficient to exclude the possibility that $(\gamma_n \bmH^{\circ} \gamma_n^{-1})$ is contained in some proper $\Q$-subgroup of $\bmG$ for infinitely many $n$'s.
So assume this possibility and let $\bmL$ be such a $\Q$-subgroup. If $\bmL$ is observable, this would violate the assumption that $ \lim \left[(g_n)_*\rmm_{[\rmH^{\circ}]}\right] = \left[ \rmm_{[\rmG]} \right]$.
So let us assume that $\bmL$ is epimorphic in $\bmG$.
By Lemma \ref{lemma_epimorphic_subgroups_parapbolic} below, we assume that $\bmL=\bmP$ is a maximal proper parabolic $\Q$-subgroup.

Find an infinite subsequence $(n_k)$ such that $\gamma_{n_k} \bmH^{\circ} \gamma_{n_k}^{-1}\subset \bmP$. Let $\bmP_k:= \gamma_{n_k}^{-1} \bmP \gamma_{n_k}$, then by assumption 
\begin{equation*}
    \norm{
    \Ad(g_{n_k} h_{n_k} )\bmv_{\bmP_{k}}
    } \to + \infty.
\end{equation*}
 On the other hand, there exists $C>1$ such that 
 \begin{equation*}
      C^{-1}\norm{
    \Ad(\gamma_{n_k} )\bmv_{\bmP_{k}}
    } 
      \leq 
      \norm{
    \Ad(g_{n_k}h_{n_k} )\bmv_{\bmP_{k}}
    }  \leq 
     C\norm{
    \Ad(\gamma_{n_k} )\bmv_{\bmP_{k}}
    } .
 \end{equation*}
 Therefore, $\norm{
    \Ad(\gamma_{n_k} )\bmv_{\bmP_{k}}
    }  \to +\infty$.
  On the other hand, $\bmP= \gamma_{n_k}\bmP_{k}  \gamma_{n_k}^{-1} $, thus for some $C>1$,
  \begin{equation*}
      C^{-1} \norm{\bmv_{\bmP}} \leq
  \norm{
    \Ad(\gamma_{n_k} )\bmv_{\bmP_{k}}
    } \leq C \norm{\bmv_{\bmP}},
  \end{equation*}
   which are bounded and unbounded at the same time, leading to a contradiction.

\subsubsection{Epimorphic subgroups}

\begin{lem}\label{lemma_epimorphic_subgroups_parapbolic}
   Let $\bmA$ be a linear algebraic group over $\Q$.
    Every proper epimorphic $\Q$-subgroup of $\bmA$ is contained in a proper parabolic $\Q$-subgroup.
\end{lem}
\begin{proof}
    Let $\bmC$ be epimorphic in $\bmA$, we need to show that $\bmC$ is contained in some proper parabolic $\Q$-subgroup.
    Let $\bmU$  be the nontrivial unipotent radical of $\bmC$ and $\bmv_{\bmU} \in \wedge^{\dim \bmU}\fraka$ be a vector lifting the Lie algebra of $\bmU$. 
    
    When $\bmA$ is reductive, $\bmv_{\bmU}$ is an unstable $\bmA$-vector. Thus, the stabilizer of the line spanned by $\bmv_{\bmU}$ is contained in a proper parabolic $\Q$-subgroup $\bmP$ of $\bmA$ (see \cite[Corollary 5.1]{Kemp78}). In particular, $\bmC$ is contained in $\bmP$.

    In general, consider the morphism $p^{\red}: \bmA \to \bmA /\bmR_{\bmu}(\bmA)$. Then $p^{\red}(\bmC)$ is epimorphic in $\bmA / \bmR_{\bmu}(\bmA)$ and hence is either contained in a proper maximal parabolic $\Q$-subgroup of $\bmA/\bmR_{\bmu}(\bmA)$ or equal the full $\bmA /\bmR_{\bmu}(\bmA)$. In the former case, we are done by taking the preimage. In the latter case, we see that $\bmU$ is contained in $\bmR_{\bmu}(\bmA)$. Hence $\bmA/\bmC$ is affine by \cite[Theorem 7.1]{Grosshans97}. Thus $\bmC$ is observable, a contradiction.
\end{proof}

\subsubsection{Proof of Lemma \ref{lemma_asymptotic_volume_polytope}}
\label{subsection_unbounded_polytopes}

Thanks to the properties of convex polytopes (See \cite[Section 4]{ShaZhe18}), it suffices to show that there exists $(\eta_n)$ diverging to $+\infty$ such that $\Omega^{\bmone}_{\bmH,g_n,\eta_n}\neq \emptyset$.

By assumption, there exist $(h_n)\subset \rmH^{\circ}$, $(b_n)\subset \rmG$ bounded and $(\gamma_n)\subset \Gamma$ such that 
\begin{equation*}
    g_nh_n= b_n \gamma_n\; \text{for all }n.
\end{equation*}
Since $p^{\spl}(h_n)= p^{\spl}(g_n^{-1}b_n\gamma_n) = p^{\spl}(b_n)$ remains bounded, we can find $(\omega_n)\subset \rmH^{\circ}$ bounded and $(h_n')\subset(\rmH^{\bmone})^{\circ}$ such that $h_n=h_n'\omega_n$. Therefore, it suffices to find diverging $(\eta_n)$ such that $\Omega^{\bmone}_{\bmH,\gamma_n,\eta_n}\neq \emptyset$.

For simplicity, let $\rmH^{\bmone 0}:= \bmH^{\bmone}(\R)^{\circ}$, $\rmS_{\rmH^{\bmone 0}}:= \rmH^{\bmone 0}/ \rmH^{\bmone 0}\cap {}^{\circ}(\bmH^{\bmone})(\R)$ and $p: \rmH^{\bmone 0} \to \rmS_{\rmH^{\bmone 0}}$ be the natural quotient map. One sees that $p$ factors through $\overline{p}: \rmH^{\bmone 0}/ \rmH^{\bmone 0}\cap \Gamma \to \rmS_{\rmH^{\bmone 0}}$. Then for any $g\in {}^{\circ}\rmG$,  $\Omega^{\bmone}_{\bmH,g,\eta}$ descends along $\overline{p}$, i.e, $\Omega^{\bmone}_{\bmH,g,\eta}= \overline{p}^{-1} \left(  \exp\left(  \calP^{\bmone}_{\bmH,g,\eta} \right) \right)$ where
\begin{equation*}
    \calP^{\bmone}_{\bmH,g,\eta}:=
    \left\{ x \in \Lie(\rmS_{\rmH^{\bmone 0}}) \;\middle\vert\;
         \norm{\Ad(g \wtt ) \bmv_{\bmP} } > \eta,\;\forall\,\bmP\in \scrP^{\max}_{\bmH^{\circ}},\; \wtt\in p^{-1}(\exp(x))
    \right\}.
\end{equation*}
For $\bmP \in \scrP^{\max}_{\bmH^{\circ}}$, let $l_{\bmP}$ be the induced linear functional on $\Lie(\rmS_{\rmH^{\bmone 0}})$ such that 
\begin{equation*}
    \Ad(\wtt) \bmv_{\bmP} = \exp ( l_{\bmP}(x) )\cdot \bmv_{\bmP},\; \forall\,\wtt\in p^{-1}(\exp(x)). 
\end{equation*}
Let $\Phi$ denote the collection of linear functionals $l$ on $\Lie(\rmS_{\rmH^{\bmone 0}})$ such that $l= l_{\bmP}$ for some $\bmP\in \scrP^{\max}_{\bmH^{\circ}}$.
Then
\begin{equation*}
    \calP^{\bmone}_{\bmH,g,\eta} =
    \left\{ x \in \Lie(\rmS_{\rmH^{\bmone 0}}) \;\middle\vert\;
         l(x) > \log(\eta)-
         \log \inf_{\bmP\in \scrP^{\max}_{\bmH^{\circ}},\, l_{\bmP}=l }
         \norm{\Ad(g) \bmv_{\bmP}}, \;\forall\,l\in \Phi
    \right\}.
\end{equation*}
After passing to a subsequence, we assume that for each $l\in \Phi$, 
\begin{equation*}
    \inf_{\bmP\in \scrP^{\max}_{\bmH^{\circ}},\, l_{\bmP}=l }
         \norm{\Ad(\gamma_n) \bmv_{\bmP}}
    \;\text{ either remains constant or diverges to infinity.}
\end{equation*}
Let $\Phi_{\infty}$ consist of those $l\in \Phi$ with $\left( \inf
\norm{\Ad(\gamma_n) \bmv_{\bmP}} \right)$ diverging to infinity and let $\wtPhi$ be its complement. Let $\Phi_0$ be consisting of $l_0 \in \wtPhi$ such that there exists $(a_{l})_{l\in \wtPhi} \subset \Z_{\geq 0}$ such that $\sum a_{l}l=0$ and $a_{l_0}\neq 0$. Let $\Phi_1$ be the complement of $\Phi_0$ in $\wtPhi$. By \cite[Lemma 3.4]{zhangrunlinCompositio2021}, if $\Phi_0=\emptyset$, then $\calP^{\bmone}_{\bmH,\gamma_n,\eta_n} $ is nonempty for some $(\eta_n)$ diverging to $+\infty$ and hence the proof is complete. 

Note that by assumption, $\gamma_n \bmH^{\circ} \gamma_n^{-1}$ can not be contained in any proper observable $\Q$-subgroup of $\bmG$ for infinitely many $n$'s. Let us assume that $\Phi_0$ is nonempty and find a  proper observable $\Q$-subgroup containing $\gamma_n \bmH^{\circ} \gamma_n^{-1}$ for infinitely many $n$, which is a contradiction.

So take $(l_1,...,l_k) \subset \wtPhi$ and $(a_1,...,a_k) \subset \Z^+$ with $\sum a_i l_i =0$. 
As $\left( \inf \norm{\Ad(\gamma_n)\bmv_{\bmP}} \right)$ is bounded and $(\Ad(\gamma_n)\bmv_{\bmP})$, as $n$ and $\bmP$ vary, is discrete, after passing to a subsequence, for $i=1,...,k$, we find $(\bmP_n^{i})\subset \scrP^{\max}_{\bmH^{\circ}}$ such that 
\begin{equation*}
    \Ad(\gamma_n) \bmv_{\bmP_n^i}  = \Ad(\gamma_1) \bmv_{\bmP_1^i}
    \quad \forall \, n\in \Z^+.
\end{equation*}

Let $v_0 := \Ad(\gamma_1) \bigotimes_{i=1}^k \bmv_{\bmP_1^i}^{\otimes a_i} $, then $v_0$ lies in a $\bmG$-representation space and is a $\gamma_n \bmH^{\circ} \gamma_n^{-1}$-weight vector. Since $\sum a_i l_i = 0$, $v_0$ is fixed by $\gamma_n \rmH^{\bmone 0} \gamma_n^{-1}$ for all $n$. Replacing $v_0$ by $v_0^{\otimes m}$ for some positive integer $m$, we assume that $\gamma_n \bmH^{\bmone} \gamma_n^{-1}$ fixes $v_0$ for all $n$.

Let $\beta_n \in \frakX^{*}(\bmH^{\circ})$ be such that $\gamma_n h \gamma_n^{-1} v_0 = \beta_n(h) v_0$ for every $h\in \bmH^{\circ}$. 
By rigidity of diagonalizable groups (see \cite[Proposition 3.2.8]{Spr98}),  there are only finitely many possibilities for $\beta_n$. By passing to a subsequence, assume $\beta_n=\beta_1$ for all $n$. As $\beta_1$ trivializes on $\bmH^{\bmone}$, it factors through some $\overline{\beta}\in \frakX^*(\bmS_{\bmG})$. Let $W_{\beta}$ be a one-dimensional representation of $\bmG$, factoring through $\bmS_{\bmG}$, defined by
\begin{equation*}
    g\cdot w_{\beta}:= \overline{\beta}^{-1}\left(
      p^{\spl}(g)
    \right) w_{\beta},\quad\forall\, w_{\beta}\in W_{\beta}.
\end{equation*}
Fix some nonzero $w_{\beta}\in W_{\beta}$.
Then $\gamma_n \bmH^{\circ} \gamma_n^{-1}$ is contained in the stabilizer subgroup of $v_0\otimes w_{\beta}$ for all $n$.
But the stabilizer of $v_0 \otimes w_{\beta}$ is not the full $\bmG$ since it is contained in $ \cap_{i=1}^k \gamma_1 \bmP_1^i\gamma_1^{-1}$.
So we are done.

\section{Good Height Functions}\label{section_good_heights}

In this section, we reduce the counting problem to various properties of equivariant compactifications. Those actually satisfying these properties will be constructed in the next section under suitable assumptions.

Notations from Section \ref{subsection_standing_assumptions_equidistribution_nondivergence} are inherited except that $\bmH^{\bmone}$ (used to be $\bmH \cap {}^{\circ}\bmG$) is replaced by $\bmH \cap ({}^{\circ}\bmG)^{\circ}$  here.
Also, for a smooth $\bmG$-pair $(\bmX, \bmD)$ over $\Q$ such  that $\bmX \setminus \bmD$ is $\bmG$-equivariantly isomorphic to $\bmU$, notations/definitions from Section \ref{subsection_analytic_Clemens_complex} are kept.

\subsection{Notations and assumptions}\label{subsection_standing_assumptions_good_heights}

We start with some notations and assumptions.

\begin{itemize}
    \item $\bmG$ is a connected linear algebraic group over $\Q$ and let $\bmG^{\bmone}:= ({}^{\circ}\bmG)^{\circ}$;
     \item $\bmU$ is a variety over $\Q$ that is homogeneous under some action of $\bmG$ and assume that $\bmU(\Q)\neq \emptyset$;
    \item Assume that the canonical line bundle $\calK_{\bmU}$ on $\bmU$ is trivial and $\omega_{\bmU}$ is an invariant gauge form on $\bmU$;
    \item For every $x\in \bmU(\Q)$, let $\bmH_x$ (resp.  $\bmH^{\bmone}_x$) be the stabilizer of $x $ in $\bmG$ (resp. $\bmG^{\bmone}$). We assume that $\bmH_x$'s are connected.
\end{itemize}

Recall $\rmG:= \bmG(\R)^{\circ}$ and for every $\R$-subgroup $\bmH$ of $\bmG$, $\rmH:= \bmH(\R)\cap \rmG$. 
Note that $\rmG^{\bmone}= \rmG \cap \bmG^{\bmone}(\R)$  is connected but in general, $\rmH$ could have different connected components.
For $x\in \bmU(\Q)$, choose an invariant measure $\rmm_{\rmG/\rmH_x}$ on $\rmG/\rmH_x$ that is identified with the measure induced from $\omega_{\bmU,\infty}$ under the orbit map. We choose invariant measures $\rmm_{\rmG}$ on $\rmG$ and $\rmm_{\rmH_x}$ on $\rmH_x$ such that the triple $(\rmm_{\rmG}, \rmm_{\rmH_x}, \rmm_{\rmG/\rmH_x})$ is compatible. 
We similarly ask ($\rmm_{\rmG^{\bmone}}, \rmm_{\rmH^{\bmone}_x}, \rmm_{\rmG^{\bmone}/\rmH^{\bmone}_x}$ )  to be compatible.
When $\bmG/\bmH_x$ is isomorphic to $\bmG^{\bmone}/\bmH^{\bmone}_x$, we further ask $\rmm_{\rmG^{\bmone}/\rmH^{\bmone}_x}$ to be identified with $\omega_{\bmU,\infty}$. 
We also require that $\rmm_{[\rmH_x]}$  is compatible with $\rmm_{\rmH_x}$, namely the triple consisting of $\rmm_{[\rmH_x]}$, $\rmm_{\rmH_x}$ and the counting measure on $\Gamma \cap \rmH_x$ is compatible.

\subsection{Definition of good and ok heights}

Let $B_R^{\Ht}:= \left\{ x\in \bmU(\R),\; \Ht(x)\leq R \right\}$ for a function $\Ht: \bmU(\R) \to \R_{\geq 0}$.
And if $x\in \bmU(\R)$, let $B_{R,x}^{\Ht}:= \left\{ y\in \rmG.x ,\; \Ht(y)\leq R \right\}$.
The upper index $\Ht$ may be dropped if the function is clear from the context.
 Note that any arithmetic lattice $\Gamma$ of $\bmG$ has a finite-index subgroup contained in $\rmG$ and hence in ${}^{\circ}\rmG = \rmG^{\bmone}$.

\subsubsection{Good and ok heights}

    If $\Pic(\bmU)$ is torsion, we say that a function $\Ht: \bmU(\R) \to \R_{\geq 0}$ is \textbf{good} if for every $x\in \bmU(\Q)$ and every arithmetic subgroup $\Gamma$ of $\bmG$ contained in $\rmG^{\bmone}$, one has that as $R$ tends to infinity,
    \begin{equation*}
        \# \Gamma.x \cap B_{R}^{\Ht} \sim 
        \frac{
        \normm{  \rmm_{[\rmH_x^{\bmone}]}  }
        }{
         \normm{   \rmm_{[\rmG^{\bmone}]}  }
        }
        \rmm_{\rmG^{\bmone}/\rmH^{\bmone}_x} \left(
         B_{R,x}^{\Ht}
        \right).
    \end{equation*}
 The function $\Ht$ is said to be \textbf{ok} if the above asymptotic holds up to some constant $c>0$ possibly dependent on $x$ and $\Gamma$.

For $\eta>0$, $g\in \rmG$ and $x\in \bmU(\Q)$,
let 
\begin{equation*}
    \wtB^{\Ht}_{R,x,\eta}:= 
    \left\{
       [g] \in  \rmG^{\bmone}/ \rmH_x \cap \Gamma
       \;\middle\vert\;
       \Ht(g.x)\leq R,\;
       \norm{\Ad(g)\bmv_{\bmP}} > \eta,\; \forall \, \bmP \in \scrP^{\max}_{\bmH^{\circ}_x}
    \right\}
\end{equation*}
where $\scrP^{\max}_{\bmH^{\circ}_x}$, defined in last section, denotes the set of maximal proper parabolic $\Q$-subgroups containing $\bmH^{\circ}_x$. 

 Without assuming $\Pic(\bmU)$ to be torsion, we say that a function $\Ht: \bmU(\R) \to \R_{\geq 0}$ is \textbf{good} if for every $x\in \bmU(\Q)$, $\eta>0$ and every arithmetic subgroup $\Gamma$ of $\bmG$ contained in $\rmG^{\bmone}$, one has 
    \begin{equation*}
        \# \Gamma.x \cap B_R^{\Ht} \sim 
        \frac{
        1 
        }{
         \normm{   \rmm_{[\rmG^{\bmone}]}  }
        }
        \rmm_{[\rmG^{\bmone}]_x} \left(
         \wtB^{\Ht}_{R,x,\eta}
        \right)
    \end{equation*}
    where $\rmm_{[\rmG^{\bmone}]_x}$ denotes the $\rmG^{\bmone}$-invariant measure on $\rmG^{\bmone}/ \rmH_x \cap \Gamma$ compatible with $\rmm_{\rmG^{\bmone}}$ and the counting measure on $\rmH_x \cap \Gamma$.
     The function $\Ht$ is said to be \textbf{ok} if the above asymptotic holds up to some constant $c>0$ possibly dependent on $x$ and $\Gamma$.
    One can check that this is compatible with the above when $\Pic(\bmU)$ is torsion (that is, when $\rmm_{[\rmH_{x}^{\bmone}]}$ is a finite measure).

\subsubsection{Weighted versions}

    Unfortunately, in this paper we are not able to exhibit any good or ok height when $\Pic(\bmU)$ is not torsion. Nevertheless, here is a weighted version.  
    For $\eta>0$ and $y \in \bmU(\Q)$,  let
    \begin{equation*}
        \Omega_{y,\eta}:= \left\{
         [h] \in \rmH^{\bmone}_y/\rmH^{\bmone}_y \cap \Gamma
         \;\middle\vert\;
         \norm{\Ad(h) \bmv_{\bmP}} \geq \eta,
         \;\forall\,
         \bmP\in \scrP^{\max}_{\bmH^{\circ}_y}
        \right\}.
    \end{equation*}

    We say that $\Ht: \bmU(\R) \to \R_{\geq 0}$ is \textbf{good with weights} if for every $x\in \bmU(\Q)$ and arithmetic subgroup $\Gamma$ contained in $\rmG^{\bmone}$, there exists $\eta>0$ such that for every $c>0$, as $R\to +\infty$,
        \begin{equation}\label{equation_good_with_weights}
            \sum_{y \in \Gamma.x \cap B_R^{\Ht}} \bmw_y
            \sim
            \frac{ 1 }{
         \normm{   \rmm_{[\rmG^{\bmone}]}  }
            }
            \rmm_{\rmG^{\bmone}/\rmH^{\bmone}_x} \left(
             \rmG^{\bmone}.x \cap B_R^{\Ht}  
            \right)
        \end{equation}
        where 
            $\bmw_y^{-1} := \max \left\{ \rmm_{[\rmH^{\bmone}_y]}\left(
            \Omega_{y,\eta}
            \right),c \right\}$.
    Here we require that $ \rmm_{[\rmH^{\bmone}_y]}$ is identified with $ \rmm_{[\rmH^{\bmone}_x]}$ when conjugating by $\gamma$.
    
    The naive definition of ``ok with weights'' does not seem natural for us. Here is a slightly different one.
    We say that $\Ht$ is \textbf{ok with weights}  if for every $x\in \bmU(\Q)$ and  arithmetic subgroup $\Gamma$ contained in $\rmG^{\bmone}$, there exists a compactly supported non-negative function $\psi: \rmG^{\bmone}/\Gamma \to \R $ such that
    \begin{itemize}
        \item[(1)] $\la \psi, (\gamma h)_*\rmm_{[(\rmH^1_x)^{\circ}]} \ra \neq 0$
        for all $\gamma \in \Gamma$ and $h\in \rmH^1_x$;
        \item[(2)] for some constant $c>0$,
        \begin{equation*}
            \sum_{y \in \Gamma.x \cap B_R^{\Ht}} \bmw^{\psi}_y
            \sim
            c\cdot \frac{ 1 }{
         \normm{   \rmm_{[\rmG^{\bmone}]}  }
            }
            \rmm_{\rmG^{\bmone}/\rmH^{\bmone}_x} \left(
             \rmG^{\bmone}.x \cap B_R^{\Ht}  
            \right)
        \end{equation*}
        where for $y=\gamma.x$, we set 
        \[
        \bmw_y^{\psi}:= \sum_{h\in \rmH^1_x/(\rmH^1_x)^{\circ}} 
        \la \psi, (\gamma h)_*\rmm_{  [ (\rmH^{\bmone}_x)^{\circ} ]  } \ra ^{-1}
        = \sum_{h\in \rmH^1_y/(\rmH^1_y)^{\circ}} 
        \la \psi, h_*\rmm_{ [\rmH_y^{\circ}]} \ra ^{-1}
        .
        \]
    \end{itemize}
    We usually take $\bmG^{\bmone}=\bmG$ in the discussion of heights that are ok with weights.
    
\subsection{Conditions implying good heights}

Here are some additional conditions on the pair $(\bmX,\bmD)$ that could imply good heights.
 We always assume that $\Ht$ is as constructed in Section \ref{subsection_heights}. So implicitly here is some $\bmL:=\sum \lambda_{\alpha}\bmD_{\alpha}$, a divisor supported on $\bmD$ satisfying:
\begin{itemize}
    \item For every $\alpha \in \scrA^{\an}_{\R,\Q}$, $\lambda_{\alpha} \geq 0$ and if $d_{\alpha}-1 \geq 0$, then $\lambda_{\alpha} >0 $.
\end{itemize}

Let $\bmB:=\sum_{\alpha \in \scrB} \bmD_{\alpha}$ be a closed $\Q$-subvariety of $\bmD$ that is a union of irreducible divisors indexed by some subset $\scrB $ of $\scrA$. We let $\scrB_{\R}^{\an}$  be the subset of $\scrA_{\R}^{\an}$ corresponding to those contained in $\bmB(\R)$. $\scrB_{\R,x}^{\an}$ and $\scrB_{\R,\Q}^{\an}$ are defined in a similar way.

\subsubsection{Conditions for good heights}

We say that $(\bmX,\bmD,\bmB)$ satisfies condition $(\rmB 1)$ if
\begin{itemize}
     \item[($\rmB 1$)]\label{condition_B_1} For every $x\in \bmU(\Q)$ and $(g_n)\subset \rmG$ such that $\lim_{n \to \infty} g_n.x \in \bmB(\R)$, one has 
        \begin{equation*}
            \lim_{n \to \infty} \left[
            (g_n)_* \rmm_{[\rmH^{\circ}_x]}
              \right] = \left[
            \rmm_{[\rmG]}   \right].
        \end{equation*}
\end{itemize}
For $x\in \bmU(\Q)$, $R>0$ and $\Ht: \bmU(\R) \to \R_{\geq 0}$, let 
\begin{equation*}
    \nu_{R,x}:= \frac{1}{\omega_{\bmU,\infty}\left(  B_{R,x}^{\Ht}  \right) } 
    \cdot \omega_{\bmU,\infty}\vert_{ B_{R,x}^{\Ht} }
\end{equation*}
be a probability measure on $\bmX(\R)$.

We say that $(\bmX,\bmD,\bmB, \Ht)$ satisfies condition $(\rmB \rmH 1)$ if
\begin{itemize}
    \item[($\rmB \rmH 1$)]\label{condition_BH_1} For every $x\in \bmU(\Q)$, $\nu_{\infty,x}:=\lim_{R\to \infty} \nu_{R,x}$ exists and $\supp(\nu_{\infty,x}) \subset \bmB(\R)$.
\end{itemize}

\begin{lem}\label{lemma_condition_good_height}
    Assume that \hyperref[condition_B_1]{$(\rmB 1)$} and \hyperref[condition_BH_1]{$(\rmB \rmH 1)$} hold.  If $\Pic(\bmU)$ is torsion, then $\Ht$ is good. In general, $\Ht$ is good with weights.
\end{lem}

\subsubsection{Proof of Lemma \ref{lemma_condition_good_height}}

For $x \in \bmU(\Q)$, by assumption, there exists some sequence $(g_n)\subset {}^{\circ}\rmG$ such that $ \lim \left[(g_n)_* \rmm_{[\rmH^{\circ}_x]} \right] = \left[ \rmm_{[\rmG]}   \right]$. Hence $p^{\spl}(\bmH_x)=\bmS_{\bmG}$ by Lemma \ref{lemma_existence_convergence_to_full} and Theorem \ref{theorem_equidistribution_explicit_a_n} is applicable. For simplicity, we write $\wtbmS_x$ for $\wtbmS^{\bmH_x}$ and $\Omega^{\bmone}_{x,g,\eta}$ for $\Omega^{\bmone}_{\bmH_x,g,\eta}$. Thus for every $\eta>0$,
\begin{equation*}
    \lim_{n\to \infty}
    \frac{
    (g_n)_* \rmm_{[\rmH_x^{\circ}]}
    }
    {
    \rmm_{
    [(\rmH_x^{\bmone})^{\circ}]
    } \left( 
       \Omega^{\bmone}_{x,g_{n},\eta}
    \right)
    } = \rmm_{[\rmG]}
    .
\end{equation*}
In the following, for $g\in \rmG$, we let $s_x\in \rmH_x$ be such that $p^{\spl}(g)=p^{\spl}(s_x)$ and $g_x:= gs_x^{-1}$. 
In other words, we choose $g_x \in {}^{\circ}\rmG \cap g \rmH_x$. 
The coset space $g_x \rmH_x^{\bmone}$ is uniquely determined and $\Omega^{\bmone}_{x,g_x,\eta}$ is independent of the choice of the coset representatives.

For a subset $\calB \subset \rmG/\Gamma$, let 
\begin{equation*}
    \calB_x:= \left(
         p^{\spl}\vert_{\wtbmS_x(\R)^{\circ}} 
    \right)^{-1} \left(
     p^{\spl}(\calB)
    \right).
\end{equation*}
Assume that $\calB$ is nonempty, open, bounded and contains the identity coset, then there exists $\eta_{\calB}>0$ such that for every $g\in \rmG$,
\begin{equation}\label{equation_nondivergence_contained_in_polytope}
    \left\{
      [h] \in \rmH_x^{\circ}\Gamma/\Gamma \;\middle\vert\;
      g_x[h]\in \calB
    \right\}
    \subset \calB_x \times \Omega^{\bmone}_{x,g_x,\eta_{\calB}}.
\end{equation}
under the natural homeomorphism as in Equa.(\ref{equation_decomposition_H_Gamma_lift_S_G}).
In particular,  ${g_*\rmm_{[\rmH^{\circ}_x]}}  \vert_{\calB} =0 $  if $\Omega^{\bmone}_{x,g_x,\eta_{\calB}} = \emptyset $.
Since $\rmH_x/\rmH_x^{\bmone}$ is connected, $\pi_0(\rmH_x^{\bmone})$ surjects on $\pi_0(\rmH_x)$. 
So actually ${g_*\rmm_{[\rmH_x]}}  \vert_{\calB} =0 $ when $\Omega^{\bmone}_{x,g_x,\eta_{\calB}} = \emptyset $.
Fix some $c,\eta>0$, let 
\begin{equation*}
    \bmw_{g.x}:= \left(
        \max\{
        \normm{
         \rmH_x^{\bmone}/ (\rmH_x^{\bmone})^{\circ}(\rmH_x^{\bmone}\cap \Gamma)
        }
        \cdot 
        \rmm_{[(\rmH_x^{\bmone})^{\circ}]} \left(
       \Omega^{\bmone}_{x,g_x,\eta} 
         \right) , c
        \}
    \right)^{-1}.
\end{equation*}
Note that $  \normm{ \rmH_x/\rmH_x^{\circ}(\rmH_x\cap \Gamma) }=
 \normm{ \rmH^{\bmone}_x/ (\rmH^{\bmone}_x)^{\circ}(\rmH^{\bmone}_x\cap \Gamma) }
$ since $\rmH_x=\rmH_x^{\circ} \cdot \rmH_x^{\bmone}$, $\rmH_x^{\circ}\cap \rmH_x^{\bmone} = (\rmH_x^{\bmone})^{\circ}$ and $\rmH_x\cap \Gamma=  \rmH^{\bmone}_x\cap \Gamma$. We use $C_2$ to denote this number.

When $g = \gamma \in \Gamma$ and $y=\gamma.x$, $\bmw_{\gamma.x}$ coincides with the $\bmw_y$ as defined following Equa.(\ref{equation_good_with_weights}).
Let us note immediately that there exists $C_3(\eta,c)> 1$ such that 
\begin{equation*}
    \begin{aligned}
             &C_2 \rmm_{[(\rmH^{\bmone}_x)^{\circ}]}\left(  \Omega_{x,g_x,\eta} \right) <  c \implies 
            C_2 \rmm_{[(\rmH^{\bmone}_x)^{\circ}]}\left(  \Omega_{x,g_x,\eta_{\calB}} \right) < C_3(\eta,c) c;\\
             &
             C_2 \rmm_{[(\rmH^{\bmone}_x)^{\circ}]}\left(  \Omega_{x,g_x,\eta} \right) \geq  c \implies 
             C_2\rmm_{[(\rmH^{\bmone}_x)^{\circ}]}\left(  \Omega_{x,g_x,\eta_{\calB}} \right) < C_3(\eta,c) \rmm_{[(\rmH^{\bmone}_x)^{\circ}]}\left(  \Omega_{x,g_x,\eta} \right).
    \end{aligned}
\end{equation*}
So we have 
\begin{equation}\label{equation_compare_polytope_different_size}
     \bmw_{g.x} \cdot  \rmm_{[(\rmH^{\bmone}_x)^{\circ}]}\left(  \Omega_{x,g_x,\eta_{\calB}} \right) < 
     C_2^{-1} C_3(\eta,c).
\end{equation}

In virtue of Theorem \ref{theorem_equidistribution_imply_counting}\footnote{The remaining hypothesis either follows from Theorem \ref{theorem_equidistribution_Chamber-Loir_Tschinkel} or can be verified directly.}, it only remains to prove that (we are using the normalization of Haar measures from Theorem  \ref{theorem_equidistribution_explicit_a_n}. In particular, one normalizes $\rmm_{[\rmG^{\bmone}]}$ to be a probability measure.)
\begin{equation*}
 \frac{1}{\rmm_{\rmG/\rmH_x}  (B^{\Ht}_{R,x})  }
 \int_{g.x\in B^{\Ht}_{R,x}}
 \bmw_{g.x}\cdot \left(
  g_* \rmm_{[\rmH_x]} 
 \right)\,\rmm_{\rmG/\rmH_x} \text{ converges to } \rmm_{[\rmG]}
\end{equation*}
against all test functions $\psi \in C_c(\calB)$.

Fix  $\psi \in C_c(\calB)$ and $\ep>0$.
By assumption and Theorem  \ref{theorem_equidistribution_explicit_a_n}, Lemma \ref{lemma_asymptotic_volume_polytope}, for $R$ sufficiently large, we find $\Good \sqcup \Bad =B_{R,x}^{\Ht}$ such that 
\begin{itemize}
    \item[1.] 
    \begin{equation*}
        \frac{
        {\rmm_{\rmG/\rmH_x}}(\Bad)
        }
        {
        {\rmm_{\rmG/\rmH_x} }(B_{R,x}^{\Ht})
        } \leq 
        \ep \left( \norm{\psi}_{\sup}\cdot  C_3(\eta,c) \cdot 
        \max\{\rmm_{\wtbmS_x(\R)^{\circ}}(\calB_x), \rmm_{[\rmG]}(\calB)\}  \right) ^{-1}
    \end{equation*}
    \item[2.]
    \begin{equation*}
        \left\vert
            \la \psi, 
             \frac{
                 g_* \rmm_{[\rmH_x]}
             }{
             C_2 \rmm_{[ (\rmH^{\bmone}_x)^{\circ} ]} \left(
              \Omega_{x,g_{x},\eta}^{\bmone}   
              \right)
              }  \ra - 
              \la \psi, \rmm_{[\rmG]} \ra
        \right\vert \leq \ep,\; \forall\, g.x \in \Good.
            \end{equation*}
     \item[3.] $ C_2\rmm_{[ (\rmH^{\bmone}_x)^{\circ} ]} \left(
              \Omega_{x,g_{x},\eta}^{\bmone}   
              \right) > c$ for every $g.x \in \Good$.
\end{itemize}
For every $g\in \rmG$, by Equa.(\ref{equation_compare_polytope_different_size}) and the discussion following Equa.(\ref{equation_nondivergence_contained_in_polytope}),
\begin{equation*}
\begin{aligned}
       \bmw_{g.x} \left\vert \la \psi, g_*\rmm_{[\rmH_x]} \ra \right\vert
       &\leq 
       C_2 \cdot \bmw_{g.x}\norm{\psi}_{\sup}  \cdot 
       \rmm_{\wtbmS_x(\R)^{\circ}} \otimes \rmm_{[(\rmH_x^1)^{\circ}]}(\calB_x \times \Omega^{\bmone}_{x,g_x,\eta_{\calB}})
        \\
    &\leq C_3(\eta,c)\cdot \norm{\psi}_{\sup} \rmm_{\wtbmS_x(\R)^{\circ}}(\calB_x)   .
\end{aligned}
\end{equation*}
Therefore, item $1$ above implies that 
\begin{equation*}
\begin{aligned}
     & \left\vert  \frac{1
    }{\rmm_{\rmG/\rmH_x}  (B^{\Ht}_{R,x})  }
   \int_{[g]\in \Bad}
   \bmw_{g.x}
   \la \psi,  
    g_* \rmm_{[\rmH_x]} 
    \ra 
    - \la \psi, \rmm_{[\rmG]} \ra 
    \,\rmm_{\rmG/\rmH_x}
    \right\vert   \\
    \leq & \,
    \frac{
    \rmm_{\rmG/\rmH_x}  (\Bad) 
    }{
    \rmm_{\rmG/\rmH_x}  (B^{\Ht}_{R,x}) 
    } \left(
       C_3(\eta,c) \norm{\psi}_{\sup} \rmm_{\wtbmS_x(\R)^{\circ}}(\calB_x) +    \norm{\psi}_{\sup} \rmm_{[\rmG]}(\calB)
    \right)\\
    \leq & \, 2\ep.
\end{aligned}
\end{equation*}
On the other hand, $\bmw_{g.x}^{-1} = C_2\rmm_{[ (\rmH^{\bmone}_x)^{\circ} ]} \left(
              \Omega_{x,g_{x},\eta}^{\bmone}   
              \right) $  when $g.x \in \Good$. So,
\begin{equation*}
\begin{aligned}
     & \left\vert  \frac{1
    }{\rmm_{\rmG/\rmH_x}  (B^{\Ht}_{R,x})  }
   \int_{[g]\in \Good}
   \bmw_{g.x}
   \la \psi,  
    g_* \rmm_{[\rmH_x]} 
    \ra
    - \la \psi, \rmm_{[\rmG]} \ra 
     \,\rmm_{\rmG/\rmH_x}
    \right\vert   \\
    \leq &\,
     \frac{   \rmm_{\rmG/\rmH_x}(\Good)
     }{
     \rmm_{\rmG/\rmH_x}(B^{\Ht}_{R,x})
     } \ep \leq \ep.
\end{aligned}
\end{equation*}
The proof is thus complete.

\subsubsection{When is the log anti-canonical height good?}

\begin{itemize}
     \item[($\rmB 2$)]\label{condition_B_2} For every $x\in \bmU(\Q)$ and $y \in \overline{\rmG.x}$, one has $\overline{\rmG.y}\cap \bmB(\R)\neq \emptyset$.
\end{itemize}

\begin{itemize}
    \item[($\rmK 1$)]\label{condition_K_1} For every $\alpha \in \scrA_{\R,\Q}^{\an}$, $d_{\alpha} -1 >0$.
\end{itemize}

\begin{lem}\label{lemma_log_anti_canonical_good}
    If condition \hyperref[condition_B_2]{$(\rmB 2)$} and \hyperref[condition_K_1]{$(\rmK 1)$} hold, then condition \hyperref[condition_BH_1]{$(\rmB \rmH 1)$} holds for the log anti-canonical height 
    $\Ht= \Ht_{-(\bmK_{\bmX}+\bmD)}$.
\end{lem}

\begin{proof}

Fix $x\in \bmU(\Q)$ and take $\bmL= -(\bmK_{\bmX}+\bmD)$.  By condition \hyperref[condition_K_1]{$(\rmK 1)$},  $\lambda_{\alpha}=d_{\alpha}-1>0$ and hence $\frac{d_{\alpha}-1}{\lambda_{\alpha}}=1$ for every $\alpha \in \scrA^{\an}_{\R,x}$.
Thus, $\scrC^{\an}_{\R,x} = \scrC^{\an}_{\R,x}(\bmL)$.
By Theorem \ref{theorem_equidistribution_Chamber-Loir_Tschinkel}, $\nu_{\infty,x}$ exists and is supported on $D_{F}$ as $F$ varies over maximal faces of $\scrC^{\an}_{\R,x}$. It remains to argue that $D_{F} \subset \bmB(\R)$ for every such $F$.

Indeed, since $D_{F}$ is $\rmG$-invariant and closed, $D_{F}\cap D_{\beta} \neq \emptyset$ for some $\beta \in \scrB^{\an}_{\R,x}$ by condition \hyperref[condition_B_2]{$(\rmB 2)$}. By maximality of $F$, we must have $\beta \in F$ and therefore, $D_{F} \subset D_{\beta} \subset \bmB(\R)$.

\end{proof}

\subsubsection{When do good heights exist?}

\begin{itemize}
     \item[$(\rmB 4)$]\label{condition_B_4} For every $x\in \bmU(\Q)$, $d_{\alpha} -1 >0$ for some $\alpha \in \scrB_{\R,x}^{\an}$. In particular, $\scrB_{\R,x}^{\an}$ is nonempty.
\end{itemize}

\begin{lem}\label{lemma_exists_good_heights}
    If condition \hyperref[condition_B_4]{$(\rmB 4)$} holds, then there exists an effective divisor $\bmL$ supported on $\bmD$ such that \hyperref[condition_BH_1]{$(\rmB \rmH 1)$} holds for $\Ht:= \Ht_{\bmL}$.
\end{lem}

\begin{proof}

For every $x\in \bmU(\Q)$,
choose $\beta_x \in \scrB_{\R, x}$ with $d_{\beta_x} -1 >0$. Choose $(\lambda_{\alpha})_{\alpha \in \scrA}$ such that 
\begin{equation*}
    \min_{x\in \bmU(\Q)} \left\{
    \frac{d_{\beta_x}-1}{\lambda_{\beta_x}}
    \right\}
    > \max_{\alpha \in \scrA\setminus \{\beta_x\} }
    \left\{
     \frac{d_{\alpha}-1}{\lambda_{\alpha}}
     \right\}.
\end{equation*}
For instance, one can pick
\begin{equation*}
\lambda_{\alpha} :=
    \begin{cases}
        d_{\alpha}-1  & \text{ if } \alpha \in \{\beta_x\} \\
        2(d_{\alpha}-1)  & \text{ if } d_{\alpha}-1>0 \,\text{ and }\,\alpha\notin \{ \beta_x\} \\
        1  &\text{ otherwise}
    \end{cases}
\end{equation*}
Then by Theorem \ref{theorem_equidistribution_Chamber-Loir_Tschinkel}, for every $x\in \bmU(\Q)$, $\supp(\nu_{\infty,x} ) \subset \bmB(\R)$.
    
\end{proof}

\subsubsection{Ok heights}\label{subsubsection_ok_heights}

Assume that $\bmH_x^{\circ}$ has no nontrivial $\Q$-characters.
For $x\in \bmU(\Q)$ and $\Gamma \subset \bmG(\Q)\cap \rmG$, let $\Psi_x: \rmG.x \to \Prob(\rmG/\Gamma)$ be defined by $\Psi_x(g.x):= g_* \rmm^{\bmone}_{[\rmH_x]}$.
Let $ \mu_{R,x}:= \Bary\left( (\Psi_x)_*(\nu_{R,x}) \right) \in \Prob(\rmG/\Gamma)$.
Here $\Bary: \Prob(\Prob(\rmG/\Gamma)) \to \Prob(\rmG/\Gamma)$ is defined by $\nu \mapsto \int_{\phi \in \Prob(\rmG/\Gamma)} \phi \; \nu(\phi)$.
\begin{itemize}
    \item[$(\rmH 1)$]\label{condition_H_1} For every $x\in \bmU(\Q)$ and arithmetic subgroup $\Gamma \subset \bmG(\Q)\cap \rmG$, $\bmH_x^{\circ}$ has no nontrivial $\Q$-characters and $\lim_{R\to \infty} \mu_{R,x} $ exists in $\Prob(\rmG/\Gamma)$.
\end{itemize}

The following is a direct corollary of Theorem \ref{theorem_equidistribution_imply_counting} and \ref{theorem_equidistribution_Chamber-Loir_Tschinkel}.
\begin{lem}\label{lemma_H_1_imply_ok}
    If \hyperref[condition_H_1]{$(\rmH 1)$} holds, then $\Ht$ is ok.
\end{lem}

\begin{itemize}
    \item[$(\rmD 1)$] \label{condition_D_1}
    For every $x\in \bmU(\Q)$ and arithmetic subgroup $\Gamma \subset \bmG(\Q)\cap \rmG$, $\bmH_x^{\circ}$ has no nontrivial $\Q$-characters and $\Psi_x$
    extends continuously to $\overline{\rmG.x}  \to  \Prob(\rmG/\Gamma) \cup \{0\}$ where the closure is taken in $\rmX^{\cor}$, the manifold with corners associated with $(\bmX,\bmD)$ (see Section \ref{subsection_manifolds_corners}).
    \item[$(\rmS 1)$]\label{condition_S_1} For every $x\in \bmU(\Q)$ and arithmetic subgroup $\Gamma$ in $\bmG(\Q)\cap \rmG$,  there exists a bounded subset $B$ of $\rmG/\Gamma$ such that $ g\rmH_x^{\circ}\Gamma/\Gamma$ intersects with $B$ for every $g\in \rmG$.
\end{itemize}

\begin{lem}\label{lemma_ok_heights}
    If \hyperref[condition_D_1]{$(\rmD 1)$} and \hyperref[condition_S_1]{$(\rmS 1)$} hold, then  \hyperref[condition_H_1]{$(\rmH 1)$} holds for  $\Ht= \Ht_{\bmL}$ for every $\bmL$ from Section \ref{subsection_analytic_Clemens_complex}, namely, $\lambda_{\alpha}\geq 0$ and $d_{\alpha}>0 \implies \lambda_{\alpha}>0$ where $\bmL=\sum \lambda_{\alpha}\bmD_{\alpha}$.
\end{lem}

\begin{proof}
    Indeed, $\Psi_x$ extends continuously to $\overline{\rmG.x} \to \Prob(\rmG/\Gamma)$.
    By Theorem \ref{theorem_variant_equidistribution_Chamber-Loir_Tschinkel}, $\lim \nu_{R,x}$ exists in $\Prob(X^{\cor})$. Applying $\Bary\circ (\Psi_x)_*$, we get that $\lim \mu_{R,x}$ exists in $\Prob(\rmG/\Gamma)$.
\end{proof}

With more care, \hyperref[condition_S_1]{$(\rmS 1)$} can be weakened as
\begin{itemize}
    \item[$(\rmD \rmS 1)$]\label{condition_DS_1}
    There exists a closed subset $D$ of $\bmB(\R)$ such that for every $(g_n)$ and every $x\in \bmU(\Q)$ with $\lim g_n.x $ not in $D$,
    there exists a bounded subset $B$ of $\rmG/\Gamma$ such that $g_n\rmH_x^{\circ}\Gamma/\Gamma$ intersects with $B$ for every $n$.
    Moreover, for every $F\subset \scrC^{\an}_{\R,\Q}$, if $\Leb_F$ denotes a smooth measure on $D_F$, then $\Leb_F(D)=0$.
\end{itemize}

With a similar proof, noting that $ \mu_{\infty,x} \left((\Psi_x)^{-1}(\{\bmzero\}) \right) =0$, one has:
\begin{lem}\label{lemma_ok_heights_possible_diverges}
    If \hyperref[condition_D_1]{$(\rmD 1)$} and \hyperref[condition_DS_1]{$(\rmD \rmS 1)$}  hold, then  \hyperref[condition_H_1]{$(\rmH 1)$} holds for $\Ht= \Ht_{\bmL}$ for every $\bmL$ from Section \ref{subsection_analytic_Clemens_complex}.
\end{lem}

\subsubsection{Ok with weights}

Assume that $\bmG$ has no nontrivial $\Q$-characters.
Take $x\in \bmU(\Q)$ and an arithmetic subgroup $\Gamma$.
We allow $\bmH_x^{\circ}$  to have nontrivial $\Q$-characters. Take a non-negative compactly supported function $\psi$ on $\rmG/\Gamma$ such that 
$\la \psi, g_*\rmm_{[\rmH_x^{\circ}]} \ra \neq 0$
for all $g \in \rmG$. Such a $\psi$ exists if \hyperref[condition_S_1]{$(\rmS 1)$} holds.
Define 
\begin{equation*}
\begin{aligned}
      \Psi_x^{\psi}: \rmG.x &\to \Prob^{\psi}(\rmG/\Gamma)  \\
      g.x & \mapsto \alpha^{\psi}_g \left( \rmm^{\psi}_{[\rmH_x]}  \right).
\end{aligned}
\end{equation*}
$\Prob^{\psi}(\rmG/\Gamma)$ collects all locally finite measures such that $\la \psi, \nu \ra=1$. And
if $ \rmm_{[\rmH_x]}  = \sum_{h\in \rmH_x/ \rmH_x^{\circ}(\rmH_x\cap \Gamma) } h_* \rmm_{  [\rmH_x^{\circ}]  }$, then 
\[
 \alpha^{\psi}_g \left( \rmm^{\psi}_{[\rmH_x]}  \right):= 
  \normm{ \rmH_x/ \rmH_x^{\circ}(\rmH_x\cap \Gamma) }^{-1}
   \sum_{h\in \rmH_x/ \rmH_x^{\circ}(\rmH_x\cap \Gamma) }  \frac{
   (gh)_* \rmm_{ [\rmH_x^{\circ}] }
   } 
   {\la \psi,   (gh)_* \rmm_{ [\rmH_x^{\circ}] }
    \ra
   }  
\]
See Section \ref{subsection_measure_compactification_infinite} for more.
Let 
\begin{equation*}
      \mu_{R,x}^{\psi}: = \Bary\left( (\Psi_x^{\psi})_*(\nu_{R,x})
       \right) \in \Prob^{\psi}(\rmG/\Gamma).
\end{equation*}

\begin{itemize}
    \item[$(\rmH 2)$]\label{condition_H_2} For every $x\in \bmU(\Q)$ and arithmetic subgroup $\Gamma \subset \bmG(\Q)\cap \rmG$,  there exists a non-negative compactly supported continuous function $\psi$ on $\rmG/\Gamma$ such that  $\la \psi, g_*\rmm_{[\rmH_x^{\circ}]} \ra \neq 0$
    for all $g \in \rmG$
    and $\lim_{R\to \infty} \mu^{\psi}_{R,x} $ exists in $\Prob^{\psi}(\rmG/\Gamma)$.
\end{itemize}

\begin{lem}\label{lemma_H_2_imply_ok_weights}
   Assume that $\bmG$ has no nontrivial $\Q$-characters. If \hyperref[condition_H_2]{$(\rmH 2)$} holds, then $\Ht$ is ok with weights.
\end{lem}

The proof is analogous to that in Section \ref{subsubsection_ok_heights}, except that there are weights and Theorem \ref{theorem_equidistribution_imply_counting} should  be applied to $H= \rmH_x^{\circ}$ rather than $H= \rmH_x$.

\begin{itemize}
    \item[$(\rmD 2)$]\label{condition_D_2} For every $x\in \bmU(\Q)$ and arithmetic subgroup $\Gamma \subset \bmG(\Q)\cap \rmG$, there exists $\psi$ satisfying the paragraph above and $\Psi^{\psi}_x$
    extends continuously to $\overline{\rmG.x}  \to  \Prob^{\psi}(\rmG/\Gamma)$ where the closure is taken in $\rmX^{\cor}$, the manifold with corners associated with $(\bmX,\bmD)$.
\end{itemize}

\begin{lem}\label{lemma_ok_heights_weights}
   Assume that $\bmG$ has no nontrivial $\Q$-characters. If \hyperref[condition_D_2]{$(\rmD 2)$} and \hyperref[condition_S_1]{$(\rmS 1)$} hold, then \hyperref[condition_H_2]{$(\rmH 2)$} holds for every $\Ht=\Ht_{\bmL}$.
\end{lem}

The proof is analogous to that of Lemma \ref{lemma_ok_heights} with $\Prob^{\psi}(\rmG/\Gamma)$ in the place of $\Prob(\rmG/\Gamma)$.

\section{Equivariant Compactifications and Focusing}\label{section_equidistribution_focusing}

Notations from Section \ref{subsection_standing_assumptions_equidistribution_nondivergence} and 
\ref{subsection_standing_assumptions_good_heights} are inherited.
Furthermore, $\bmG$ and $\bmH$ are  connected linear algebraic groups over $\Q$ and $\bmH$ is assumed to be observable in $\bmG$.
For two $\Q$-subgroups $\bmA,\bmB$ of $\bmG$, define a $\Q$-subvariety
\[
\bmZ(\bmA,\bmB):= \left\{  g\in \bmG \midd g\bmA g^{-1} \subset \bmB \right\}.
\]

We will construct compactifications and heights of $\bmG/\bmH$ satisfying various properties from Section \ref{section_good_heights}. Two guiding examples are 
\begin{itemize}
     \item[(1)]  $\bmG$ and $\bmH$ are connected, reductive and $\bmZ_{\bmG}(\bmH)^{\circ} \subset \bmH$;
    \item[(2)] $\bmG$ is $\Q$-split, semisimple and $\bmH$ contains a maximal unipotent subgroup of $\bmG$.
\end{itemize}

\subsection{Finiteness of intermediate subgroups}
Let
\begin{equation*}
\begin{aligned}
    \INT(\bmH,\bmG)
     &:= \left\{
            \bmH \subset  \bmL \subsetneq \bmG \;\middle\vert\;
            \bmL \text{ is connected}
        \right\}
\end{aligned}
\end{equation*}
be the set of intermediate connected closed subgroups between $\bmH$ and $\bmG$. 
Let $\INT_{\Q}(\bmH,\bmG)$ (resp. $\INT^{\obs}(\bmH,\bmG)$) collect elements in $\INT(\bmH,\bmG)$ that are defined over $\Q$ (resp. that are observable). Let $\INT_{\Gamma}(\bmH,\bmG)$ denote the even smaller subset consisting of $\bmL\in \INT_{\Q}^{\obs}(\bmH,\bmG)$ such that  $(\gamma_n\bmH \gamma_n^{-1})$ converges to $\bmL$ (see Definition \ref{definition_convergence_subgroups}) for some sequence $(\gamma_n) \subset \Gamma$.
We say that $(\bmG,\bmH)$ satisfies condition (F1) if
\begin{itemize}
    \item[(F1)]\label{condition_F_1} $\INT_{\Gamma}(\bmH,\bmG)$ is finite.
\end{itemize}

\begin{lem}\label{lemma_F_1}
    Assume that $\bmG$ is reductive and one of the following is true:
     \begin{itemize}
         \item[(1)] $\bmH$ is reductive and $\bmZ_{\bmG}(\bmH)^{\circ}$ is contained in $\bmH$;
         \item[(2)] $\bmH$ contains a maximal unipotent subgroup.
     \end{itemize}
     Then \hyperref[condition_F_1]{$(\rmF 1)$} holds.
\end{lem}

From the proof, $\INT(\bmH,\bmG)$ is finite in case (1), which is no longer true in case (2).

\subsubsection{Proof of Lemma \ref{lemma_F_1}, part 1}

Under our assumptions, the central torus $\bmZ(\bmG)$ of $\bmG$ is contained in $\bmH$, so it suffices to show this for $\bmH/\bmZ(\bmG)\subset \bmG/\bmZ(\bmG)$.
When $\bmH/\bmZ(\bmG)$ is semisimple, this already follows from \cite[Lemma 3.1]{EinMarVen09}.
 
In general, every intermediate observable subgroup is reductive by Theorem \ref{theorem_nondivergence}.
Then, up to $\bmG(\overline{\Q}) $-conjugacy, there are only finitely many $\bmL\in \INT^{\obs}(\bmH,\bmG)$ (see for instance, \cite[Lemma A.1]{EinMarVen09} and the beginning part of \cite[Section 4]{zhanghan_zhangrunlin_nondivergence_1_arxiv}).
Fix finitely many representatives $\{\bmL_1,...,\bmL_k\} \subset  \INT^{\obs}(\bmH,\bmG)$. By Lemma \ref{lemma_conjugacy_subgroups}, for each $i=1,...,k$, find $(f^i_j)_{j=1,...,l_i} \subset \bmZ(\bmH,\bmL_i)(\overline{\Q}) $ such that 
\begin{equation*}
    \bmZ(\bmH,\bmL_i)( \overline{\Q})  = \bigcup_{j=1}^{l_k} \bmL_i( \overline{\Q} )f^i_j \bmH( \overline{\Q} )
\end{equation*}
where we used the fact that $\bmH$ is of finite-index in $\bmN_{\bmG}(\bmH)$ under our assumptions. Now we claim that
\begin{equation*}
    \INT^{\obs}(\bmH,\bmG)= \left\{  (f^i_j)^{-1} \bmL_i f^i_j  
    \;\middle\vert\; i=1,...,k;\;j=1,...,l_i
    \right\}.
\end{equation*}
Indeed, for $\bmL \in \INT^{\obs}(\bmH,\bmG)$, we find $g\in \bmG(\overline{\Q})$ and $i\in\{1,...,k\}$ such that $g\bmL g^{-1}=\bmL_i$. So $g \in \bmZ(\bmH,\bmL_i)( \overline{\Q})$. Find $l\in \bmL_i( \overline{\Q} )$, $h\in \bmH( \overline{\Q} )$ and $j\in \{1,...,l_i\}$ such that $g=l f^i_j h$. Then,
\begin{equation*}
    \bmL = g^{-1} \bmL_i g
    = (f^i_j h)^{-1} \bmL_i f^i_j h = (f^i_j)^{-1} \bmL_i f^i_j. 
\end{equation*}
So we are done.

 \subsubsection{Structure of $\bmL$ containing $\bmU_{\max}$}\label{section_structure_containing_maximal_unipotent}

  Fix a Borel subgroup $\bmB$ containing a maximal unipotent subgroup $\bmU_{\max}$ and a maximal torus $\bmT$ in $\bmB$. Let $\Phi$ be the set of (nonzero) weights. Thus $\frakg= \left( \oplus_{\alpha \in \Phi} \frakg_{\alpha} \right) \oplus \frakt$. Let $\bmL$ be a connected closed subgroup containing $\bmU_{\max}$, then one can show that its Lie algebra is
  \begin{equation}\label{equation_horospherical_structure}
      \frakl= \left(
      \oplus_{\alpha\in \Phi_{\frakl}} \frakg_{\alpha}
      \right) \oplus (\frakl \cap \frakt),\quad
      \exists \, \Phi_{\frakl} \subset \Phi.
  \end{equation}
  Let $\frakm_{\frakl}$ be the Lie subalgebra generated by $\oplus_{\alpha\in \Phi_{\frakl}} \frakg_{\alpha}$ and $\bmM_{\bmL}$ the corresponding connected algebraic subgroup. Then there exists a torus $\bmT_{\bmL} $ in $\bmT$ such that $\bmL = \bmM_{\bmL}\cdot \bmT_{\bmL}$ and $\frakl=\frakm_{\frakl}\oplus \frakt_{\frakl}$.
  So we get an injection
  \begin{equation*}
  \begin{aligned}
        \INT(\bmU,\bmG) &\embed 2^{\Phi} \times \{ \text{subtori of }\bmT \} \\
        \bmL &\mapsto (\Phi_{\frakl}, \bmT_{\bmL}).
  \end{aligned}
  \end{equation*}

\subsubsection{Proof of Lemma \ref{lemma_F_1}, part 2}

Assume that $\bmH$ contains a maximal unipotent subgroup and $\bmL\in \INT_{\Gamma}(\bmH,\bmG)$. By Section \ref{section_structure_containing_maximal_unipotent}, $\bmL=\bmM_{\bmL}\cdot \bmT_{\bmL}$ with $\bmM_{\bmL}$ normalized by the torus $\bmT_{\bmL}$ and $\frakl= \frakm_{\frakl}\oplus \frakt_{\frakl}$. When $\bmL$ is defined over $\Q$,  both $\bmM_{\bmL}$ and $\bmT_{\bmL}$ can be arranged to be over $\Q$. 
Let $\pi: \bmL \to \bmL/\bmM_{\bmL}$ be the natural quotient morphism.

Let $(\gamma_n)\subset \Gamma$ be such that $(\gamma_n \bmH \gamma_n^{-1})$ converges to $\bmL$. Then $\left(\pi(\gamma_n \bmH \gamma_n^{-1}) \right)= (\pi(\bmH))$ also converges to $\bmL/\bmM_{\bmL}$, which is impossible unless $\pi(\bmH)= \bmL/\bmM_{\bmL}$. Hence $\bmL= \bmM_{\bmL}\cdot \bmH$. This implies that $\bmL$ is determined by $\bmM_{\bmL}$. So there are only finitely many possibilities.

\subsection{Compactifications using intermediate subgroups}\label{section_algebraic_compactification_intermediate_groups}

\subsubsection{Condition (F2)}

To detect nondivergence we introduce
\begin{itemize}
    \item[(F2)]\label{condition_F_2} there are only finitely many parabolic  $\Q$-subgroups containing $\bmH$, i.e., $\scrP_{\bmH}$ is finite.
\end{itemize}
For simplicity we abbreviate $\INT:=  \INT_{\Gamma}(\bmH,\bmG) $ and $\INTP:= \INT_{\Gamma}(\bmH,\bmG) \cup \scrP_{\bmH}$ in this subsection. 
Clearly $\INTP$ is finite if $\INT_{\Q}^{\obs}(\bmH,\bmG)$ is, but the converse does not hold. For instance $\bmG= \SL_n$ ($n\geq 3$) and $\bmH$ is a maximal unipotent $\Q$-subgroup. Another case is when $\bmG= \SL_n \times \bmD$ ($n\geq 3$) with $\bmD$ embeded in $\SL_n$ as a $\Q$-anisotropic maximal $\Q$-torus and $\bmH$ is the diagonal embedding of $\bmD$ in $\bmG$.

For a parabolic $\Q$-subgroup $\bmP$, we let $\scrV(\bmP)$ consist of all triples $(\rho,\bmV_{\rho},\bmv_{\rho})$ where
\begin{itemize}
    \item[1.] $\bmV_{\rho}$ is a finite dimensional  $\Q$-vector space;
    \item[2.] $\rho: \bmG \to \bmG/\bmR_{\bmu}(\bmG) \to \GL(\bmV_{\rho})$ is a $\Q$-representation of $\bmG$ factoring through $\bmG/\bmR_{\bmu}(\bmG)$;
    \item[3.] $\bmv_{\rho}\in \bmV_{\rho}(\Q)$ is not fixed by $\rho(\bmP) $ but the line spanned by $\bmv_{\rho}$ is preserved by $\rho(\bmP)$.
\end{itemize}
By abuse of notation, we simply write an element of $\scrV(\bmP)$ as $\rho$ or $\bmv_{\rho}$ when no confusion might arise. 
And we let $\alpha_{\bmv_\rho}$ denote the $\Q$-character of $\bmP$ attached to $\bmv_{\rho}$.
Let $\pi'_{\bmP}: \bmP \to \bmS'_{\bmP}$ be a $\Q$-split quotient torus of $\bmP$ such that every $\Q$-character of $\bmP$ factors through $\pi'_{\bmP}$. By rigidity of diagonalizable groups, 
\[
    \left\{
        \pi_{\bmP}'(g\bmH g^{-1}) \midd g\in \bmZ(\bmH,\bmP) 
    \right\}
    =\left\{   \bmS_1,...,\bmS_k \right\}
\]
is finite. For each $\bmS_i$, let 
\[
     \scrV(\bmP)^{\bmS_i, \bmH  }:=\left\{
        \rho \in \scrV(\bmP) \midd \bmv_{\rho} \text{ is fixed by }\bmS_i \text{ and } \bmH
     \right\}.
\]
We find a finite subset $\scrV_i \subset \scrV(\bmP)^{\bmS_i,\bmH}$ such that $\left\{ \alpha_{\bmv_{\rho}} \midd \rho \in \scrV(\bmP)^{\bmS_i,\bmH} \right\} $ is contained in the cone spanned by $\left\{ \alpha_{\bmv_{\rho}} \midd \rho \in \scrV_i \right\}$. Let $\scrV(\bmH,\bmP):= \bigcup \scrV_i$ and $\scrV(\bmH):= \cup_{\bmP \in \scrP_{\bmH}} \scrV(\bmH,\bmP)$. 


\subsubsection{Definition of the compactification}\label{subsection_definition_compactification_intermediate_groups}
For a finite dimensional $\Q$-vector space $\bmV$, let $\bmP(\bmV)$ be the associated projective variety over $\Q$.
For $\bmL \in \INT$, let $\bmv_{\frakl}$ be a vector in $\bmV_{\frakl}:= \wedge^{\dim \bmL}\frakg$ that lifts $\frakl$.
In order that $\bmv_{\frakl}$ is fixed by $\bmH$ we need to assume that the modular character of $\bmL$ is trivial on $\bmH$. If $\bmH$ is assumed to have no nontrivial $\Q$-characters, then this automatically holds.
To get an embedding of $\bmG/\bmH$, we fix a rational $\bmG$-representation $\bmV_{\bmH}$  and $\bmv_{\bmH}\in \bmV_{\bmH}(\Q)$ such that the stabilizer of $\bmv_{\bmH}$  in $\bmG$ is exactly $\bmH$.

\begin{defi}\label{defiComptfyX_1}
Assume both $\bmG,\bmH$ are connected and \hyperref[condition_F_1]{$(\rmF 1)$}, \hyperref[condition_F_2]{$(\rmF 2)$} hold. 
Moreover, assume that the modular character of $\bmL$ is trivial on $\bmH$ for every $\bmL \in \INT$.
Let
\begin{equation*}
    \begin{aligned}
            \iota_{\bmH}^{\INTP}:\, \bmG/\bmH &\to 
            \bmV_{\bmH} \oplus
            \bigoplus_{\bmL \in \INT}  \bmV_{\frakl} \oplus \bigoplus_{\rho \in \scrV(\bmH)}  \bmV_{\rho}
            \\
            & \to \bmP(\bmV_{\bmH}\oplus \Q) \oplus \bigoplus_{\bmL \in \INT}  \bmP(\bmV_{\frakl} \oplus \Q)    \oplus \bigoplus_{\rho \in \scrV(\bmH)}  \bmP(\bmV_{\rho} \oplus \Q)    
    \end{aligned}
\end{equation*}
where the first arrow sends $ g\bmH$ to $g.\bmv_{\bmH} \oplus (\oplus g. \bmv_{\frakl}) \oplus (\oplus g.\bmv_{\rho})$.
Let $\bmX^{\INTP}_{\bmH}$ be the Zariski closure of $\iota_{\bmH}^{\INTP}(\bmG/\bmH)$ and $\rmX^{\INTP}_{\rmH}$ be the analytic closure of $\rmG.o$ in $\bmX^{\INTP}_{\bmH}(\R)$ where $o$ denotes the image of the identity coset.
\end{defi}

By extending every representation $\rho$ here to $(\rho \oplus 1, \bmV \oplus \Q )$, we see that $\iota_{\bmH}^{\INTP}$ is $\bmG$-equivariant.

\subsubsection{Nondivergence}

Every element $\bmv$ in $\rmX^{\INTP}_{\rmH}$ can be written as $ \oplus [\bmv^{\bmL}:t^{\bmL}]$ with $t^{\bmL} \in \{0,1\}$. 
Let $(g_n)$ be a sequence in $\rmG$ and $\bmv = \lim_{n\to \infty} g_n.o$. If $\bmv^{\rho}=\bmzero$ for some $\rho \in \scrV(\bmH) $, then $\left( g_n \rmH \Gamma/\Gamma \right)$ diverges topologically in $\rmG/\Gamma$ by Mahler's criterion. The converse is also true by Theorem \ref{theorem_divergence_observable_1}:
\begin{lem}\label{lemma_compactification_nondivergence}
    Assume \hyperref[condition_F_1]{$(\rmF 1)$} and \hyperref[condition_F_2]{$(\rmF 2)$} hold.
    If $(g_n)$ is a sequence in $\rmG$ such that $\bmv = \lim_{n\to \infty} g_n.o$ exists in $\rmX^{\INTP}_{\rmH}$ and $\bmv^{\rho} \neq \bmzero$ for every  $\rho \in \scrV(\bmH) $, then $g_n \rmH^{\circ}\Gamma/\Gamma$ intersects with some bounded subset in $\rmG/\Gamma$ for all $n$.
\end{lem}

\begin{proof}
By Theorem \ref{theorem_divergence_observable_1}, if the conclusion were not true, then there exist $(\gamma_n)\subset \Gamma$, $(h_n) \subset \rmH^{\circ}$, a parabolic $\Q$-subgroup $\bmP$ containing all $\gamma_n \bmH \gamma_n^{-1}$ and a nonzero rational $\bmP$-eigenvector $\bmv$ that is fixed by $\gamma_n\bmH \gamma_n^{-1}$ and lives in a $\bmG$-representation that factors through $\bmG/\bmR_{\bmu}(\bmG)$ such that $g_n h_n \gamma_n^{-1}.\bmv \to \bmzero$. 
So $g_n \gamma_n^{-1}.\bmv \to \bmzero$.
 Since $\gamma_n^{-1} \bmP \gamma_n \in \scrP_{\bmH}$ for all $n$ and  $\scrP_{\bmH}$ is assumed to be finite, $\gamma_n^{-1} \bmP \gamma_n =\gamma_1^{-1} \bmP \gamma_1$ for all $n$ after passing to a subsequence. 

Let $\bmP':= \gamma_1^{-1} \bmP \gamma_1 \in \scrP_{\bmH}$ and $\bmv' := \gamma_1^{-1}.\bmv$. Then $\bmv'$ is a $\bmP'$-eigenvector fixed by $\gamma_1^{-1}\gamma_n \bmH \gamma_n^{-1}\gamma_1$ for all $n$. Passing to a subsequence, we find $i_0$ such that $\pi'_{\bmP'}(  \gamma_1^{-1}\gamma_n \bmH \gamma_n^{-1}\gamma_1  ) = \bmS_{i_0} $ for all $n$. 
So $\diff \alpha_{\bmv'}$ is a $\Q_{\geq 0}$-linear combination of $\{\diff \alpha_{\bmv_{\rho}} ,\; \rho \in \scrV_{i_0} \}$. In particular, passing to a further subsequence, we can find a fixed $\rho_0\in  \scrV_{i_0} \subset \scrV(\bmH,\bmP')$ such that (note that $\gamma_n^{-1} \gamma_1 \in \bmP'\cap\Gamma$)
\[
    g_n.\bmv' = \pm g_n \gamma_n^{-1} \gamma_1 .\bmv'  \to \bmzero
    \implies 
    g_n.\bmv_{\rho_0} \to \bmzero.
\]
In particular, $\bmv^{\rho_0}= \bmzero$ and this is a contradiction.
\end{proof}

\begin{defi}\label{definition_B_INTP}
    Assume condition \hyperref[condition_F_1]{$(\rmF 1)$} and \hyperref[condition_F_2]{$(\rmF 2)$} hold and that the modular character of $\bmL$ is trivial on $\bmH$ for every $\bmL \in \INT$.
    Let  $\bmB^{\INTP}_{\bmH}$ be the $\bmG$-invariant closed $\Q$-subvariety of  $\bmX^{\INTP}_{\bmH}$ consisting of  elements of the form $\oplus_{\bmL \in \{\bmH\} \sqcup \INT \sqcup \scrV(\bmH)} [\bmv^{\bmL}:t^{\bmL}] $ with $t^{\bmL}=0$ for every $\bmL\in \INT$.
\end{defi}

\begin{lem}\label{lemma_sufficient_condition_B_1_one_orbit}
    Assume that $\bmG$ is reductive, or more generally, the unipotent radical of $\bmG$ commutes with $\bmG$.
    Also assume that Condition \hyperref[condition_F_1]{$(\rmF 1)$}, \hyperref[condition_F_2]{$(\rmF 2)$} hold and that the modular character of $\bmL$ is trivial on $\bmH$ for every $\bmL \in \INT$. 
    Let $(g_n)$ be a sequence in $\rmG$ such that $(g_n.o)$ converges to $\bmv \in \rmX^{\INTP}_{\rmH}$, then
    \begin{equation*}
        \lim_{n\to \infty} \left[ (g_n)_* \rmm_{[\rmH]} \right]  = \left[ \rmm_{[\rmG]} \right] \iff
        \bmv \in \bmB^{\INTP}_{\rmH}(\R).
    \end{equation*}
\end{lem}

In particular, assuming the equivalent conditions in Lemma \ref{lemma_existence_convergence_to_full},  $\bmB^{\INTP}_{\bmH}$ is nonempty.

\begin{proof}

   The nontrivial direction is ``$\Longleftarrow$''.
   So we assume $g_n.\bmv_{\frakl} \to \infty$ for all $\bmL \in  \INT$.
   
   By Theorem  \ref{theorem_divergence_observable_1}, we can find $(h_n)\subset \rmH^{\circ}$, $(\gamma_n)\subset \Gamma$, a parabolic $\Q$-subgroup $\bmP$ and an observable $\Q$-subgroup $\bmO$ such that  $\gamma_n \bmH \gamma_n^{-1}$ is contained in $\bmP$ and strongly converges to $\bmO$. 
   Replacing $\gamma_n$ by $\gamma_1^{-1}\gamma_n$ if necessary, we assume $\bmH \subset \bmO$ and hence $\bmO \in \INT$.
   Thus, $\gamma_n^{-1}\bmO \gamma_n \in \INT$ for all $n$.
   But $\INT$ is finite, by passing to a subsequence, we assume that $\gamma_n^{-1}\bmO \gamma_n = \gamma_1^{-1}\bmO \gamma_1=: \bmO'$ for all $n$. 
   So $\gamma_n^{-1}\gamma_1$ normalizes $\bmO'$,
    implying $\gamma_n^{-1}\gamma_1. \bmv_{\frako'}= \pm \bmv_{\frako'}$. 
    Also, the modular character of $\bmO'$ is assumed to vanish on $\bmH$.
    Thus $g_n.\bmv_{\frako'} = \pm g_n h_n\gamma_n^{-1} \gamma_1 .\bmv_{\frako'}$.
    
   Let $g_nh_n \gamma_n^{-1}\gamma_1 = k_n a_n p_n$ 
   using the horospherical coordinate attached to $\bmP'$ (and some choice of maximal reductive subgroup of $\bmG$ and Cartan involution), then Theorem  \ref{theorem_divergence_observable_1}, together with our assumption on $\bmG$, asserts that $\alpha(a_n) $ is bounded for all character appearing in the Lie algebra of  $\bmP'$ and in particular, of $\bmO'$. Since $(p_n)$ is bounded, we have that $g_n h_n\gamma_n^{-1} \gamma_1. \bmv_{\frako'}= k_n a_n p_n.\bmv_{\frako'}$ is bounded. This is a contradiction against our assumption that $g_n \bmv_{\frako'}$ should diverge to infinity.
\end{proof}

\subsection{Conjugacy of subgroups}

Recall
$
     \bmZ(\bmH,\bmL):=\left\{
          g\in \bmG\;\middle\vert\;
          g\bmH g^{-1} \subset \bmL
     \right\}
$,
for two $\Q$-subgroups $\bmH,\bmL$ of $\bmG$, admits an action of $\bmL$ by multiplying from the left and an action of $\bmN_{\bmG}(\bmH)$ from the right. One naturally wonders whether the action has finitely many orbits.
We formulate two conditions related to this. 
We say that $(\bmG,\bmH)$ satisfies (C1) or (N1) if
\begin{itemize}
    \item[(C1)]\label{condition_C_1} $\bmZ(\bmH,\bmL)$ decomposes into finitely many orbits under the action of $\bmL \times \bmN_{\bmG}(\bmH)$ for every $\bmL \in \INT(\bmH,\bmG)$;
    \item[(N1)]\label{condition_N_1} $\Gamma \cap \bmH$ is a finite index subgroup of $\Gamma \cap \bmN_{\bmG}(\bmH)$ for an(y) arithmetic subgroup $\Gamma$.
\end{itemize}

 \begin{lem}\label{lemma_conjugacy_subgroups}
     Assume that one of the following is true:
     \begin{itemize}
         \item[(1)]  $\bmG$ and $\bmH$ are connected, reductive and $\bmZ_{\bmG}(\bmH)^{\circ} \subset \bmH$;
         \item[(2)]  $\bmG$ is a $\Q$-split semisimple group and
         $\bmH$ contains a maximal unipotent subgroup of $\bmG$.
     \end{itemize}
     Then \hyperref[condition_C_1]{$(\rmC 1)$} and \hyperref[condition_N_1]{$(\rmN 1)$} hold.
 \end{lem}

 Item $1$ is true by \cite[Lemma 5.2]{EskMozSha96} (cf. \cite[Theorem 8.1]{Richardson_Conjugacy_Annals_1967}).
 Item $2$ will be proved in the following subsections. Note that ($2$) might fail if $\bmG$ is not $\Q$-split, for instance, when $\bmG= \Res_{\Q(\sqrt{2})/\Q}(\SL_2)$ (but the discussion in last subsection still applies).
 We first note the following consequence of these conditions plus Theorem \ref{theorem_finiteness_orbits}:
 \begin{prop}\label{proposition_focusing}
     Assume condition \hyperref[condition_C_1]{$(\rmC 1)$} holds.
     Let $\Gamma$ be an arithmetic subgroup  of $\bmG$, then there exists a finite-index subgroup $\Gamma'$ with the following property:
     The action of $\bmL \cap \Gamma' \times \bmN_{\bmG}(\bmH)\cap \Gamma'$ on $\bmZ(\bmH,\bmL) \cap \Gamma'$ has finitely many orbits for every $\bmL$. 
     In particular, if \hyperref[condition_N_1]{$(\rmN 1)$} also holds, then the action of $\bmL \cap \Gamma' \times \bmH \cap \Gamma'$ has finitely many orbits.
 \end{prop}

 Indeed, $\Gamma'$ can be  taken to be the intersection of $\Gamma$ with any congruence subgroup.

\begin{rmk}
   The statement \cite[Lemma 5.2]{EskMozSha96} seems to suggest that \hyperref[condition_C_1]{$(\rmC 1)$} holds whenever $\bmG$ is reductive. This is wrong. One can find unipotent subgroups $\bmH$ and $\bmL$ of $\bmG:=\SL_6$ such that \hyperref[condition_C_1]{$(\rmC 1)$} fails (even after replacing the action of $\bmL$ by that of $\bmN_{\bmG}(\bmL)$).
\end{rmk}

  \subsubsection{ Proof of Lemma \ref{lemma_conjugacy_subgroups}:
  verify \hyperref[condition_N_1]{(N1)}}

  We explain that $\bmN_{\bmG}(\bmH) \cap \Gamma$  is virtually contained in $\bmH$ under the assumption that $\bmG$ is $\Q$-split and $\bmH$ contains a maximal unipotent subgroup $\bmU_{\max}$.
   It suffices to show that $\bmN_{\bmG}(\bmH)^{\circ}/\bmH^{\circ}$ is a $\Q$-split torus (without assuming $\bmG$ to be $\Q$-split, $\bmN_{\bmG}(\bmH)^{\circ}/\bmH^{\circ}$ is only proved to be a $\Q$-torus). As $\bmG$ is $\Q$-split, we assume that $\bmU_{\max}$ and $\bmT$ (from Section \ref{section_structure_containing_maximal_unipotent}) are defined over $\Q$.

  By Equa.(\ref{equation_horospherical_structure}), $\bmH$ is normalized by $\bmT$. Hence $\bmN_{\bmG}(\bmH)$ contains a Borel subgroup, implying that it is a parabolic subgroup $\bmP$. In particular, it is connected. By assumption, $\bmH$ contains $\bmR_{\bmu}(\bmP)$ and so $\bmH/\bmR_{\bmu}(\bmP)$ is a normal subgroup of $\bmP/\bmR_{\bmu}(\bmP)$. But $\bmH/\bmR_{\bmu}(\bmP)$ contains a maximal unipotent subgroup of $\bmP/\bmR_{\bmu}(\bmP)$, so $\bmH/\bmR_{\bmu}(\bmP)$  contains every semisimple factor of $\bmP/\bmR_{\bmu}(\bmP)$. Hence $\bmP/\bmH$ is a quotient of the central torus of $\bmP/\bmR_{\bmu}(\bmP)$, which is $\Q$-split since $\bmG$ is $\Q$-split.
  
  \subsubsection{Proof of Lemma \ref{lemma_conjugacy_subgroups}:
  verify \hyperref[condition_C_1]{(C1)}
  }

  Condition \hyperref[condition_C_1]{$(\rmC 1)$} is a geometric one and holds without assuming $\bmG$ to be $\Q$-split. 

  Take $g\in \bmZ(\bmH,\bmL)(\overline{\Q}) $, so $g\bmU_{\max}g^{-1}\subset g\bmH g^{-1} \subset \bmL$. As both $g\bmU_{\max}g^{-1}$ and $\bmU_{\max}$ are maximal unipotent subgroups of $\bmL$, there exists $l_g\in \bmL(\overline{\Q}) $ such that 
  \begin{equation*}
      l_g g \bmU_{\max} (l_g g)^{-1} =\bmU_{\max}
      \subset \bmH_g:= l_g g \bmH (l_g g)^{-1} \subset \bmL.
  \end{equation*}
  Thus $l_g g$ belongs to the normalizer of $\bmU_{\max}$, which is exactly $\bmB$. By the structure of subgroups containing $\bmU_{\max}$, $\bmH$ is normalized by $\bmT$ and hence $\bmB$. So we conclude that $l_g g \in \bmN_{\bmG}(\bmH)$ and $ \bmZ(\bmH,\bmL)(\overline{\Q})  = \bmL(\overline{\Q})\cdot \bmN_{\bmG}(\bmH)(\overline{\Q})$.

\subsubsection{Another compactification by intermediate groups}

For each $\bmL \in \INT$, 
we fix a rational $\bmG$-representation $\bmV_{\bmL}$  and $\bmv_{\bmL}\in \bmV_{\bmL}(\Q)$ such that the stabilizer of $\bmv_{\bmL}$  in $\bmG$ is exactly $\bmL$. Here we do not need to assume that the modular character of $\bmL$ is trivial on $\bmH$.

\begin{defi}\label{defiComptfyX_2}
Assume both $\bmG,\bmH$ are connected and \hyperref[condition_F_1]{$(\rmF 1)$}, \hyperref[condition_F_2]{$(\rmF 2)$} hold. 
Let
\begin{equation*}
    \begin{aligned}
            \iota_{\bmH}^{\IMTP}:\, \bmG/\bmH &\to 
            \bigoplus_{\bmL \in \INT}  \bmV_{\bmL} \oplus \bigoplus_{\rho \in \scrV(\bmH)}  \bmV_{\rho}
            \to  \bigoplus_{\bmL \in \INT \sqcup \scrV(\bmH)}  \bmP(\bmV_{\bmL} \oplus \Q) \\
             g\bmH &\mapsto (\oplus g. \bmv_{\bmL}) \oplus (\oplus g.\bmv_{\rho})
    \end{aligned}
\end{equation*}
Let $\bmX^{\IMTP}_{\bmH}$ be the Zariski closure of $\iota_{\bmH}^{\IMTP}(\bmG/\bmH)$ and $\rmX^{\IMTP}_{\rmH}$ be the analytic closure of $\rmG.o$ in $\bmX^{\IMTP}_{\bmH}(\R)$ where $o$ denotes the image of the identity coset. Define $\bmB_{\bmH}^{\IMTP}$ in the same way as Definition \ref{definition_B_INTP}.
\end{defi}

Statement analogous to Lemma \ref{lemma_compactification_nondivergence} holds with the same proof.
The analogue of Lemma \ref{lemma_sufficient_condition_B_1_one_orbit} also holds assuming additionally \hyperref[condition_C_1]{$(\rmC 1)$} and \hyperref[condition_N_1]{$(\rmN 1)$}.

\subsubsection{Focusing and intermediate subgroups}

For a sequence $(g_n)$ in $\rmG$, assuming $g_n \rmH^{\circ}\Gamma/\Gamma$ intersects with some bounded subset of $\rmG/\Gamma$ for all $n$, then $g_n h_n =\delta_n \gamma_n$ for some $(h_n)\subset \rmH^{\circ}$, $(\gamma_n)\subset \Gamma$ and bounded $(\delta_n)\subset \rmG$. So the limiting behaviour of $g_n \rmH^{\circ}\Gamma/\Gamma$ is essentially the same as $\gamma_n \rmH^{\circ}\Gamma/\Gamma$. A sequence  $(\gamma_n)$ in $\Gamma$ is said to be \textbf{clean} if for every $\bmL \in \INT$, one of the following holds
\begin{itemize}
    \item[1.] $(\gamma_n.\bmv_{\bmL})$ is bounded, or equivalently, $(\gamma_n)$ is bounded modulo $\bmL$;
    \item[2.] $(\gamma_n.\bmv_{\bmL})$ diverges, or equivalently, $(\gamma_n)$ diverges modulo $\bmL$.
\end{itemize}
For a clean sequence $(\gamma_n)$, let 
\begin{equation*}
    \scrL((\gamma_n)):= \left\{
    \bmL \in \INT \;\middle\vert\;
    (\gamma_n.\bmv_{\bmL}) \text{ is bounded }
    \right\}.
\end{equation*}

\begin{lem}\label{lemma_minimal_subgroup_bounded}
    Let $(\gamma_n)$ be a clean sequence in $\Gamma$. Further assume condition \hyperref[condition_C_1]{$(\rmC 1)$} and \hyperref[condition_N_1]{$(\rmN 1)$} hold. Then for every $\bmL_1,\bmL_2 \in \scrL((\gamma_n))$, there exists $\bmL' \in \scrL((\gamma_n))$ that is contained in $\bmL_1\cap \bmL_2$. In particular, $\scrL((\gamma_n))$ has a unique minimum element, denoted as $\bmL((\gamma_n))$.
\end{lem}

\begin{proof}
    Since $(\gamma_n \bmv_{\bmL_1})$ and $(\gamma_n \bmv_{\bmL_2})$ are bounded, we must have that $(\gamma_n (\bmv_{\bmL_1}\oplus \bmv_{\bmL_2}))$ is bounded, that is, $(\gamma_n)$ is bounded modulo $\bmL_1 \cap \bmL_2$. 
    Without loss of generality we assume $(\gamma_n)$ is contained in $\bmL_3:=(\bmL_1\cap \bmL_2)^{\circ}$.
    But we do not know whether $\bmL_3$ belongs to $ \INT$.
    If not, then by Lemma \ref{lemma_existence_convergence_to_full}, at least one of the following is true
    \begin{itemize}
        \item[1.] $\bmH$ is contained in a proper normal $\Q$-subgroup of $\bmL_3$;
        \item[2.] $\bmH\bmL_3^{\ari}$ is a proper subgroup of $\bmL_3$.
    \end{itemize}
    Let $\bmL_4$ be the identity component of the proper normal $\Q$-subgroup in case $1$ and be $\bmH\bmL_3^{\ari}$  in case $2$. Note that $\bmL_4$ is observable in $\bmG$. 
    In the second case, $(\gamma_n)$ is clearly bounded modulo $\bmL_4$ by the definition of $\bmL_3^{\ari}$.
     For the first case, we invoke the condition \hyperref[condition_C_1]{$(\rmC 1)$}, \hyperref[condition_N_1]{$(\rmN 1)$} and Proposition \ref{proposition_focusing} to see that $\gamma_n = l^4_{\gamma_n} f_{i_n} h_{\gamma_n}$ for some $l^4_{\gamma_n} \in \bmL_4 \cap \Gamma$, $h_{\gamma_n} \subset \bmH \cap \Gamma$ and $f_{i_n}$ belongs to certain finite subset of $\bmZ(\bmH,\bmL_4)$. 
    Since $(\gamma_n)$ belongs to $\bmL_3$,  $(f_{i_n})$ is contained in $\bmZ(\bmH,\bmL_4)\cap \bmL_3$.
    And since $\bmL_4$ is normal in $\bmL_3$, we find that $\gamma_n
    = f_{i_n}  (f_{i_n}^{-1} l^4_{\gamma_n} f_{i_n}) h_{\gamma_n}
    $ is bounded modulo $\bmL_4$.

    Continuing this process one will end with the smallest element in $\scrL((\gamma_n))$.
\end{proof}

\subsection{Resolution of singularities and manifolds with corners}

Assume $\rmm_{[\rmH]}$ is finite for the moment.
It is not true that in general, the natural map\footnote{This was called $\Psi_o$ in Section \ref{subsubsection_ok_heights}.}
\begin{equation*}
\begin{aligned}
        \Psi_{\bmH}: \rmG.o  &\to \Prob(\rmG/\Gamma) \\
        g.o &\mapsto g_* \rmm_{[\rmH]}^{\bmone}
\end{aligned}
\end{equation*}
extends to a continuous map from $\rmX^{\IMTP}_{\rmH}$ to $\Prob(\rmG/\Gamma)  \sqcup \{0\}$.
This is because the limiting measure of $(\gamma_n)_* \rmm_{[\rmH]}^{\bmone} $ is of the form $\rmm^{\bmone}_{[\rmL^{\circ}\rmH]}$ but
the stabilizer of $\rmm_{[\rmL^{\circ}\rmH]}$ may not contain $\rmL$ even though it contains a Zariski-dense subset. 

The issue would disappear if one had assumed that $\rmL^{\circ}\rmH = \rmL$ for all $\bmL\in\INT_{\Gamma}(\bmH,\bmG)$. This is indeed the case, for instance, when $\bmG$ is semisimple and $\bmH$ contains a maximal $\R$-split torus (see \cite{Matsu64} or \cite[Theorem 14.4]{BorelTits65}). As we are not assuming this, one is forced to pass to the nonalgebraic world. It will be shown (see Theorem \ref{theorem_almost_algebraic_compactification_dominate_measure_compactification_finite_volume} and \ref{theorem_almost_algebraic_compactification_dominate_measure_compactification_infinite_volume} below) that under some assumptions $\Psi_{\bmH}$ extends to the closure in a manifold with corners.

\subsubsection{Resolution of singularities}

As a first step, one applies the resolution of singularity due to Hironaka \cite{Hironaka_resolution_singularities}. The desired statement can be found in \cite[Theorem 1.0.2]{Wlodarczyk2005} and the lifting of algebraic group actions is explained in \cite[Proposition 3.9.1]{Kollar07}.
\begin{thm}\label{theorem_resolution_of_singularity}
    Let $(\bmX,\bmD)$ be a $\bmG$-pair such that $\bmX \setminus \bmD$ is homogeneous.
    There exists a smooth $\bmG$-pair $(\wtbmX, \wtbmD)$  over $\Q$ and a $\bmG$-equivariant morphism $\pi: \wtbmX \to \bmX$ such that $\pi$ is an isomorphism restricted to $\pi^{-1}(\bmX \setminus \bmD)= \wtbmX\setminus \wtbmD$.
\end{thm}

Assume \hyperref[condition_F_1]{$(\rmF 1)$}, \hyperref[condition_F_2]{$(\rmF 2)$} hold and $\bmB^{\IMTP}_{\bmH}$ is nonempty. Apply Theorem \ref{theorem_resolution_of_singularity} to the pair $(\bmX^{\IMTP}_{\bmH}, \bmB^{\IMTP}_{\bmH})$.
We get $\pi_1: (\bmX_1,\bmB_1) \to (\bmX^{\IMTP}_{\bmH}, \bmB^{\IMTP}_{\bmH})$ where $ (\bmX_1,\bmB_1) $ is a smooth pair over $\Q$. Let $\bmD_1$ be the preimage of $\bmD_{\bmH}^{\IMTP}$ under $\pi_1$ where $\bmD_{\bmH}^{\IMTP}:= \bmX_{\bmH}^{\IMTP}\setminus \iota_{\bmH}^{\IMTP}(\bmG/\bmH)$.
Apply Theorem \ref{theorem_resolution_of_singularity} again to the pair $(\bmX_1,\bmD_1)$ to get $\pi_2: (\bmX_2,\bmD_2) \to (\bmX_1,\bmD_1)$ where $ (\bmX_2,\bmD_2) $ is a smooth pair. 
If $\bmB_2$ denotes the preimage of the (Cartier) divisor $\bmB_1$, then $\bmB_2$ is a union of certain irreducible components of $\bmD_2$.
Define $\wtbmX^{\IMTP}_{\bmH}:= \bmX_2$, $\wtbmD^{\IMTP}_{\bmH}:= \bmD_2$ and $\wtbmB^{\IMTP}_{\bmH}:= \bmB_2$.
If one is only interested in whether $(g_n)_*\rmm_{[\rmH]}$ converges to the full Haar measure, then it is not necessary to pass to the nonalgebraic world.

\subsubsection{Manifolds with corners}\label{subsection_manifolds_corners}

Let $(\bmX,\bmD)$ be a smooth $\bmG$-pair over $\Q$ with $\bmX\setminus \bmD$ being homogeneous under the action of $\bmG$. We take $o\in (\bmX\setminus\bmD)(\Q)$. Let $\rmX$ be the analytic closure of $\rmG.o$ and $\rmD$ be its intersection with $\bmD(\R)$. We define in this section the associated ``manifold with corners'', denoted as $(\rmX^{\cor},\rmD^{\cor})$, together with an equivariant continuous map $\pi^{\cor}$ onto $(\rmX,\rmD)$.

For $i=1,2,...,d$, let $\bme_i(x_1,...,x_d):= x_i$ be the standard coordinate functions on $\R^d$.
Then $\rmX$ can be glued by local charts $(\calO_i,\varphi_i)_{i\in \scrI}$, where $(\calO_i)$ is an open covering of $\rmX$, such that for each $i\in\scrI$, 
\begin{itemize}
    \item[1.] $\varphi_i: \calO_i \to \varphi_i(\calO_i)$ is a homeomorphism onto a connected open neighborhood of $\bmzero \in \R^d$;
    \item[2.] there exists $0\leq r_i\leq d$ such that
    \[
       \varphi_i( \rmD \cap \calO_i ) = \varphi(\calO_i) \cap \left( \bigcup_{j=1,...,r_i} \{\bme_j = 0\} \right).
    \]
\end{itemize}
By shrinking $\calO_i$'s if necessary, we further assume the following:
For every $(i,j)$ such that $\calO_{ij}:=\calO_i \cap \calO_j \neq \emptyset$, there exist $I_{ij}\subset \{1,...,r_i\}$, $I_{ji}\subset \{1,...,r_j\}$ and a bijection $\tau_{ij}:I_{ij} \to I_{ji}$ such that
\begin{itemize}
    \item[1.] $\varphi_i(\calO_{ij}) \cap \{ \bme_k=0,\;\forall\, k\in I_{ij} \} \neq \emptyset$,
    $\varphi_j(\calO_{ji}) \cap \{ \bme_k=0,\;\forall\, k\in I_{ji} \} \neq \emptyset$;
    \item[2.] for every $k\in \{1,...,r_i\} \setminus I_{ij}$, 
    $\varphi_i(\calO_{ij}) \cap \{ \bme_k=0 \} = \emptyset$, every $k\in \{1,...,r_j\} \setminus I_{ji}$, 
    $\varphi_j(\calO_{ji}) \cap \{ \bme_k=0 \} = \emptyset$;
    \item[3.] for every $k\in I_{ji}$, there is a number $\epsilon^{ij}_k \in \{-1, 1\}$ such that 
    the homeomorphism 
    \[
    \varphi_j \circ \varphi_i^{-1} : \varphi_i(\calO_{ij}) \to \varphi_j(\calO_{ij})
    \]
    sends $\{\bme_{\tau_{ij}^{-1}(k) }>0\}$ to $ \{\bme_{k}\epsilon^{ij}_k>0\}$.
\end{itemize}
We define a new topological space $\rmX^{\cor}$ by ``cutting along the boundaries''. 
For $\kappa = \kappa(\cdot ) \in \{-1,1\}^{r_i}$ and $i\in \scrI$, let 
\begin{equation*}
\begin{aligned}
        &\calO^{\kappa}_i:= \left\{
         x\in \calO_i,\;
         \kappa(k) \bme_k(x)\geq 0, \;\forall\, k = 1,...,r_i
         \right\};\\
        &\calO^{\kappa,\circ}_i:= \left\{
         x\in \calO_i,\;
         \kappa(k) \bme_k(x) > 0, \;\forall\, k = 1,...,r_i
         \right\}.
\end{aligned}
\end{equation*}
To distinguish $\calO_i^{\kappa}$ (resp. $\calO_i$) as an independent topological space from a subset of $\rmX$, let $\iota_i^{\kappa} : \calO_i^{\kappa} \to \rmX$ (resp. $\iota_i: \calO_i \to \rmX$) be the natural inclusion map. Their inverses are well defined and continuous on the image.

We define an equivalence relation $\sim$ on $\bigsqcup_{i\in \scrI,\kappa\in \{-1,1\}^{r_i}} \calO_i^{\kappa}$ as follows:
$x \sim y$ iff there exist $i,j \in \scrI$ and $\kappa_i \in \{-1,1\}^{r_i}, \kappa_j \in \{-1,1\}^{r_j} $ such that
\[
    x \in \calO^{\kappa_i}_i,\;   y \in \calO^{\kappa_j}_j,\;
    \iota_i^{\kappa_i}( \calO^{\kappa_i,\circ}_i ) \cap  \iota_j^{\kappa_j}(\calO^{\kappa_j,\circ}_j) \neq \emptyset      
    \text{ and }\iota_i^{\kappa_i}(x)  = \iota_j^{\kappa_j}(y).
\]
Let 
$\rmX^{\cor}:= \bigsqcup_{i,\kappa} \calO_i^{\kappa} /\sim$ be equipped with the quotient topology. 

The natural map from $\calO_i^{\kappa}$ to $\rmX^{\cor}$ is an open embedding.
Also, the natural continuous map $\bigsqcup_{i,\kappa} \calO_i^{\kappa} \to \rmX$ factors through some map $\pi^{\cor}: \rmX^{\cor} \to \rmX$, which is continuous by the definition of quotient topology. Let $\rmD^{\cor}$ be the preimage of $\rmD$ under $\pi^{\cor}$. From the construction we see that

\begin{lem}
    The restriction of $\pi^{\cor} $ to $\rmX^{\cor} \setminus \rmD^{\cor}$, or to each $\calO_i^{\kappa}$, is a homeomorphism onto its image.
\end{lem}

As the structures are compatible with the group action, one can also check that the $\rmG$-action on $\rmX$ lifts to a continuous $\rmG$-action on $\rmX^{\cor}$.

Let us also point out that the proof of Theorem \ref{theorem_equidistribution_Chamber-Loir_Tschinkel} in \cite{Chamber-Loir_Tschinkel_2010_Igusa_integral} also yields similar equidistribution  statements on the manifolds with corners:
\begin{thm}\label{theorem_variant_equidistribution_Chamber-Loir_Tschinkel}
    Under the weak-$*$ topology, the family of probability measures on $\rmX^{\cor}$
 \begin{equation*}
     \mu_{R,x}:= \frac{
      \bmone_{B_{R,x}  } \cdot \Vol
     }{
     \Vol(B_{R,x})
     }
 \end{equation*}
  has a limit $\nu$ in $\Prob(\rmX^{\cor})$ as $R \to +\infty$.
  Moreover, the support of $\nu$ is equal to the union of $Z^{\cor}:= (\pi^{\cor})^{-1}(Z)$ as $(I,Z)$ varies over faces of $\scrC^{\an}_{\R,x}(\bmL)$ of dimension $b_x$.
\end{thm}

\subsubsection{Measure compactifications, finite volume cases}\label{subsection_measure_compactification_finite}

Whereas Theorem \ref{theorem_equidistribution_EMS} gives a satisfactory description for the limiting behaviour of $(g_n)_* \rmm_{[\rmH^{\circ}]}$, what is really needed for the counting problem is the limiting behaviour of $(g_n)_* \rmm_{[\rmH]}$. For this assume $\bmH$ to be connected and apply the following corollary to the weak approximation property \cite[Theorem 7.7]{PlaRap94}:
\begin{lem}\label{lemma_Q_points_dense_in_real}
    Let $\bmA$ be a connected linear algebraic group.
    Under the analytic topology, $\bmA(\Q)$ is dense in $\bmA(\R)$.
\end{lem}

By Lemma \ref{lemma_Q_points_dense_in_real}, we choose a finite subset $(c_i^{\bmH})_{i \in I} \subset \bmH(\Q)$ such that 
\begin{equation*}
    {\rmH} = \bigsqcup_{i\in I} c_i^{\bmH} \rmH^{\circ},
\end{equation*}
and let
\begin{equation*}
    \Gamma^{\sma}:= \bigcap_{i \in I} c_i^{\bmH} \Gamma (c_i^{\bmH})^{-1}.
\end{equation*}
Taking some $c_i^{\bmH}$ to be the identity, $ \Gamma^{\sma}$ is a finite index subgroup of $\Gamma$.
We define two equivalence relations on $I$. We say that $i\sim_1 j$ iff the image of $c_i^{\bmH}$ in  $\rmH/\rmH^{\circ} (\rmH \cap \Gamma)$ equals to that of $c_j^{\bmH}$.
For $\bmL \in \INT_{\Gamma}(\bmH,\bmG)$, we say that $i \sim_{\bmL} j$ iff the image of $c_i^{\bmH}$ in  $\rmH/\rmH \cap (\rmL^{\circ}\cdot (\rmL \cap \Gamma))$ equals to that of $c_j^{\bmH}$.
One can check that
\begin{equation*}
    \rmH\Gamma = \bigsqcup_{i\in I/\sim_1} c_i^{\bmH} \rmH^{\circ} \Gamma,\quad
     \rmL^{\rmH}\Gamma = \bigsqcup_{i\in I/\sim_{\bmL}} c_i^{\bmH} \rmL^{\circ}\Gamma
\end{equation*}
where $\rmL^{\rmH}:= \rmL^{\circ}\rmH$. Hence,
\begin{equation*}
\begin{aligned}
     & \rmm^{\bmone}_{[\rmH]} = \normm{I/\sim_1}^{-1}
    \sum_{ [i] \in I/\sim_1}  (c_i^{\bmH})_* \rmm^{\bmone}_{[\rmH^{\circ}]}  ,
    \\
    &  \rmm^{\bmone}_{[\rmL^{\rmH}]} = \normm{I/\sim_{\bmL}}^{-1}
    \sum_{ [i]  \in I/\sim_{\bmL}} (c_i^{\bmH})_* \rmm^{\bmone}_{[\rmL^{\circ}]}  .
\end{aligned}
\end{equation*}

\begin{lem}\label{lemma_limit_fat_homogeneous_measures_finite}
    Assume $\bmH$ is connected and $\rmm_{[\rmH]}$ is finite.
    Let $(\gamma_n)$ be a sequence in $\Gamma^{\sma}$ such that $(\gamma_n \bmH \gamma_n^{-1})$ converges to $\bmL$, 
    then $\lim_{n \to \infty} ({\gamma_n})_{*} \rmm^{\bmone}_{[\rmH]}  =  \rmm^{\bmone}_{[\rmL^{\rmH}]}$.
\end{lem}

To avoid redundancy, the reader is referred to Lemma \ref{lemma_limit_fat_homogeneous_measures_infinite} below for the proof in a similar case.

\subsubsection{The measure compactification is dominated by manifolds with corners}

For the sake of simplicity, write $\rmX:= \wtrmX_{\rmH}^{\IMTP,\cor}$ and $\rmD:= \wtrmD_{\rmH}^{\IMTP,\cor}$ in this subsection.

\begin{thm}\label{theorem_almost_algebraic_compactification_dominate_measure_compactification_finite_volume}
    Assume \hyperref[condition_F_1]{$(\rmF 1)$}, \hyperref[condition_F_2]{$(\rmF 2)$}, \hyperref[condition_C_1]{$(\rmC 1)$} and \hyperref[condition_N_1]{$(\rmN 1)$} hold and that $\bmH$ is connected with no nontrivial $\Q$-characters. Let $(x_n)$ be a sequence in $\rmX \setminus \rmD$ converging to some $x_{\infty} \in \rmX$. Then $\left( \Psi_{\bmH}(x_n) \right)$ is convergent in $ \Prob(\rmG/\Gamma)\sqcup \{ \bmzero \} $. Consequently, $\Psi_{\bmH}$ extends to a continuous map $\overline{\Psi}_{\bmH} : \rmX \to \Prob(\rmG/\Gamma)\sqcup \{ \bmzero \} $.
\end{thm}

   Such an extension $\overline{\Psi}_{\bmH}$ is necessarily $\rmG$-equivariant.

\begin{proof}
Let $\pi$ denote the natural projection $\rmX \to \rmX^{\IMTP}_{\rmH}$. 
Write 
    \begin{equation*}
        \pi(\bmx_{\infty}) =: \bmv_{\infty} = \oplus_{\bmL\in \INT \cup \scrV(\bmH)}[\bmv^{\bmL}_{\infty} : t^{\bmL}_{\infty}].
    \end{equation*}
If $\bmv^{\rho}_{\infty}= \bmzero$ for some $\rho \in \scrV(\bmH)$, then $ \left( \Psi_{\bmH}(x_n) \right)$ converges to $\bmzero$. Assume otherwise, by Lemma \ref{lemma_compactification_nondivergence}, we find $(\gamma_n)\subset \Gamma^{\sma}$ and a bounded sequence $(\delta_n)\subset \rmG$ such that $x_n = \delta_n \gamma_n .o$. 
Let $\bmL_{\infty}:= \bmL((\gamma_n))$ by Lemma \ref{lemma_minimal_subgroup_bounded}, noting that $(\gamma_n)$ is a clean sequence since $(\delta_n \gamma_n[\bmv_{\bmL}:1])$ converges. Let $\delta_{\infty} \in \rmG$ be such that $\lim \delta_n \gamma_n \bmv_{\bmL_{\infty}} = \delta_{\infty} \bmv_{\bmL_{\infty}}$.
As we are going to apply Proposition \ref{proposition_focusing} (to $\Gamma^{\sma}$ in the place of $\Gamma$), we may replace $\Gamma^{\sma}$ by a finite-index subgroup  from the beginning and assume that the conclusion of Proposition \ref{proposition_focusing} holds with $\Gamma'=\Gamma^{\sma}$.

By assumption, $(\gamma_n.\bmv_{\bmL_{\infty}})$ is bounded. Take a finite set $\calF \subset \Gamma^{\sma}$ such that for every $n$, there exists $\gamma\in \calF$ such that
$    \gamma_n . \bmv_{\bmL_{\infty}} = \gamma . \bmv_{\bmL_{\infty}}$.
By ignoring finitely many terms and replace $\calF$ by a smaller subset, for each $\gamma\in \calF$,  the number of $n$'s such that $\gamma_{n} .\bmv_{\bmL_{\infty}} = \gamma .\bmv_{\bmL_{\infty}}$ is infinite. 
Fix $\gamma \in \calF$ and $(\gamma_{n_k})$ with $\gamma_{n_k}. \bmv_{\bmL_{\infty}} = \gamma . \bmv_{\bmL_{\infty}}$.
We claim that $ \left( \gamma^{-1}\gamma_{n_k} \bmH \gamma_{n_k}^{-1} \gamma \right)$ converges to $ \bmL_{\infty} $.  

Write $l_k := \gamma^{-1}\gamma_{n_k} \in \bmL_{\infty}\cap \Gamma^{\sma}$.
If the claim were not true,  by passing to a further subsequence, we may assume that $l_k \bmH l_k^{-1} $'s are all contained in some $\bmL \in \INT_{\Gamma}(\bmH,\bmG)$, a proper subgroup of $\bmL_{\infty}$.
By Proposition \ref{proposition_focusing} and passing to a further subsequence, we assume $l_k=l_k'f h_k$ where $(l_k')\subset \bmL \cap \Gamma^{\sma} $, $ (h_k) \subset \bmH \cap \Gamma^{\sma}$ and $f \in \bmZ(\bmH,\bmL)  \cap \Gamma^{\sma}$. Note that actually $f$ is contained in $\bmL_{\infty}$. Write $\bmL_{f}:= f^{-1} \bmL f$, then it is also an element of $\INT_{\Gamma}(\bmH,\bmG)$ contained in $\bmL_{\infty}$. But $(l_k)$, and hence $(\gamma_{n_k})$, is bounded modulo $\bmL_f$. This is a contradiction against the minimality of $\bmL_{\infty}$ and the claim is proved.

 Then $(\delta_{n_k})$ converges to $\delta_{\infty} \gamma^{-1}$ modulo $\gamma \bmL_{\infty} \gamma^{-1}$ or equivalently, modulo $\gamma \rmL_{\infty} \gamma^{-1}$. Let $\calF_{\gamma} \subset \gamma \rmL_{\infty} \gamma^{-1}$ be a finite subset such that every limit point of $(\delta_n)$ in $\rmG/\gamma \rmL_{\infty}^{\rmH} \gamma^{-1}$ is equal to the image of $\delta_{\infty} c$ for some $c\in \calF_{\gamma}$.
By Theorem \ref{lemma_limit_fat_homogeneous_measures_finite}, we have 
\begin{equation*}
    \lim_{k \to \infty} \delta_{n_k}\gamma_{n_k} \rmm^{\bmone}_{[\rmH]}
    \in\left\{ (\delta_{\infty} c \gamma )_*\rmm^{\bmone}_{[\rmL_{\infty}^{\rmH}]} \midd c\in \calF_{\gamma}
    \right\}.
\end{equation*}

Up to now, we have shown that whenever $(x_n)$ converges to $x_{\infty}$, all limit points of $\Psi_{\bmH}(x_n)$ are contained in the following finite set:
\begin{equation*}
    F:= \left\{
    (\delta_{\infty} c \gamma )_*\rmm^{\bmone}_{[\rmL_{\infty}^{\rmH}]}
    \;\middle\vert\;
    \gamma \in \calF,\; c\in \calF_{\gamma}
    \right\}.
\end{equation*}

 There exist neighborhoods $\calN_{\nu} \subset \Prob(\rmG/\Gamma)$ of $\nu\in F$ such that $\calN_{\nu}\cap \calN_{\mu} =\emptyset$ for $\nu\neq \mu \in F$.
By discussion above, there exists an open neighborhood $\calN_{x_{\infty}}$  of $x_{\infty}$ in $\rmX^{\cor}$ such that
\begin{equation*}
        \Psi_{\bmH} (\calN_{x_{\infty}} \setminus \rmD ) \subset 
        \bigsqcup_{\nu\in F} \calN_{\nu}.
\end{equation*}    
By the construction of manifolds with corners, there exists a smaller neighborhood $\calN_{x_{\infty}}'$  of $x_{\infty}$ such that $\calN_{x_{\infty}}' \setminus \rmD$ is connected (this is the only place where manifolds with corners are needed!).
Thus $\Psi_{\bmH}(\calN'_{x_{\infty}} \setminus \rmD ) \subset 
\calN_{\nu}$ for a unique $\nu \in F$. This shows that $\lim_{n\to \infty} \Psi_{\bmH}(x_n)$ exists and is equal to $\nu$  for all sequences $(x_n)$ converging to $x_{\infty}$.
And the proof is now complete.
\end{proof}

\subsubsection{Measure compactifications, infinite volume cases}\label{subsection_measure_compactification_infinite}

For this subsection, we assume that $\bmG$ and $\bmH$ are connected, reductive and $\bmZ_{\bmG}(\bmH)^{\circ} \subset \bmH$. 
Here we are mainly interested in the case when $\rmm_{[\rmH]}$ is infinite.

By Theorem \ref{theorem_nondivergence}, for any arithmetic subgroup $\Gamma \subset \rmG\cap \bmG(\Q)$, there exists an open bounded subset $\calB_{\Gamma} \subset \rmG$ such that $\rmG = \calB_{\Gamma} \cdot \Gamma \cdot \rmH$. Therefore, for a  positive function $\psi$ whose support is large enough, one has $\la \psi, \nu \ra >0$ for all $\nu$ of the form $g_* \rmm_{[\rmL^{\circ}]}$ with $g\in \rmG$ and $\bmL\in \INT_{\Gamma}(\bmH,\bmG)$.
Fix such a $\psi$, define
\begin{equation}\label{definition_prob_psi}
    \Prob^{\psi}(\rmG/\Gamma):=\left\{
         \nu \in \Meas(\rmG/\Gamma) \;\middle\vert\;
         \la \psi, \nu \ra =1
    \right\},
\end{equation}
equipped with the weak-$*$ topology where $\Meas(\rmG/\Gamma)$ denotes the set of locally finite Borel measures on $\rmG/\Gamma$.
We will be concerned with measures $ \nu \in \Meas(\rmG/\Gamma) $ such that $\nu= \sum_{i=1}^k \nu_i$ where
\begin{itemize}
    \item each $\supp \nu_i$ is connected and $\supp \nu_i \cap \supp \nu_j = \emptyset$ for $i\neq j$;
    \item $ \la \psi, \nu_i \ra \neq 0$ for all $i$.
\end{itemize}
For such a measure $\nu$, let $\nu^{\psi}: = 
k^{-1} \sum_{i=1}^k \nu_i /\la \psi, \nu_i \ra =k^{-1} \sum_{i=1}^k \nu_i^{\psi} \in \Prob^{\psi}(\rmG/\Gamma)$. 

\begin{rmk}
  The definition is made so that Lemma \ref{lemma_limit_fat_homogeneous_measures_infinite}  below holds.
  If $\bmH(\R) \to \bmG(\R)/\bmG(\R)^{\circ}$ is surjective, then we could have defined $\nu^{\psi}:= \frac{\nu}{\la\psi, \nu \ra}$.
\end{rmk}

For $g\in \rmG$, if both $\nu $ and $g_* \nu$ enjoy the above properties, we define $\alpha^{\psi}_g(\nu^{\psi}):= (g_*\nu)^{\psi} \in \Prob^{\psi}(\rmG/\Gamma)$. 
Let 
\begin{equation*}
    X^{\meas, \psi}_{\bmH}: = \overline{
 \alpha^{\psi}_{\rmG} \left(  
 \rmm^{\psi}_{[\rmH]}
 \right) 
 }  \subset \Prob^{\psi}(\rmG/\Gamma)
 .
\end{equation*}
Similar to Section \ref{subsection_measure_compactification_finite}, one has 
\begin{equation*}
\begin{aligned}
     & \rmm^{\psi}_{[\rmH]} = \normm{I/\sim_1}^{-1}
    \sum_{ [i] \in I/\sim_1} ( (c_i^{\bmH})_* \rmm_{[\rmH^{\circ}]}  )^{\psi},
    \\
    &  \rmm^{\psi}_{[\rmL^{\rmH}]} = \normm{I/\sim_{\bmL}}^{-1}
    \sum_{ [i]  \in I/\sim_{\bmL}} ( (c_i^{\bmH})_* \rmm_{[\rmL^{\circ}]}  )^{\psi}.
\end{aligned}
\end{equation*}

\begin{lem}\label{lemma_limit_fat_homogeneous_measures_infinite}
    Let $(\gamma_n)\subset \Gamma^{\sma}$ be a sequence such that $(\gamma_n \bmH \gamma_n^{-1})$ converges to $\bmL \in \INT(\bmH,\bmG)$, then $\lim_{n \to \infty} \alpha^{\psi}_{\gamma_n} \left(\rmm^{\psi}_{[\rmH]} \right)  =  \rmm^{\psi}_{[\rmL^{\rmH}]}$.
\end{lem}

\begin{proof}
    For every $n,i$, let $\gamma_n^i:= (c_i^{\bmH})^{-1} \gamma_n c_i^{\bmH} \in \Gamma$. Then
    \begin{equation*}
        (\gamma_n c_i^{\bmH})_* \rmm_{[\rmH^{\circ}]}
    = (c_i^{\bmH} \gamma^i_n )_* \rmm_{[\rmH^{\circ}]}.
    \end{equation*}
     Also, $\left( \gamma_n^i \bmH (\gamma_n^i)^{-1} \right)$ converges to 
    $(c_i^{\bmH})^{-1} \bmL c_i^{\bmH}$, which is equal to $\bmL$ as $\bmH \subset \bmL$.
    By Theorem \ref{theorem_equidistribution_EMS}, 
    $\lim \alpha^{\psi}_{\gamma_n^i} 
    \left(  \rmm^{\psi}_{[\rmH^{\circ}]}  \right) = \rmm^{\psi}_{[\rmL^{\circ}]}$ and so 
    $ \lim \alpha^{\psi}_{\gamma_n c_i^{\bmH}} \left(   \rmm^{\psi}_{[\rmH^{\circ}]}  \right) = \left(  (c_i^{\bmH})_*\rmm_{[\rmL^{\circ}]}  \right)^{\psi}
    $.
    Consequently,
    \begin{equation*}
        \begin{aligned}
            \lim_{n \to \infty} \alpha^{\psi}_{\gamma_n}
            \left( \rmm^{\psi}_{[\rmH]} \right)
             &=
             \normm{I/\sim_1}^{-1}
             \sum_{ [i] \in I/\sim_1} 
             \lim_{n \to \infty} \left( (\gamma_n c_i^{\bmH})_* \rmm_{[\rmH^{\circ}]}  \right)^{\psi} 
             \\
             &= 
             \normm{I/\sim_1}^{-1}  \sum_{[i] \in I/\sim_1} \left(
             (c_i^{\bmH})_*\rmm_{[\rmL^{\circ}]}
             \right)^{\psi}
             \\
             &=   \normm{I/\sim_{\bmL}}^{-1} 
             \sum_{[i] \in I/\sim_{\bmL}}
             \left(
              (c_i^{\bmH})_*\rmm_{[\rmL^{\circ}]}
             \right)^{\psi}
             =  \rmm^{\psi}_{[\rmL^{\rmH}]}.
        \end{aligned}
    \end{equation*}
\end{proof}

Define an equivalence relation $\sim$ on $\INT:= \INT_{\Gamma}(\bmH,\bmG)$ by
\begin{equation*}
    \bmL_1 \sim \bmL_2 \iff
    \rmL_1^{\rmH} \Gamma = g\rmL_2^{\rmH} \Gamma,\;\exists\, g\in \rmG.
 \end{equation*}

\begin{lem}
    We can decompose $X^{\meas, \psi}_{\bmH} = \bigsqcup_{[\bmL]\in \INT/\sim} \alpha^{\psi}_{\rmG} \left(  \rmm^{\psi}_{[\rmL^{\rmH}]}
    \right)$.
\end{lem}

\begin{proof}
    It suffices to show that every $\nu \in X^{\meas, \psi}_{\bmH}$ is equal to $\alpha^{\psi}_g \left( \rmm^{\psi}_{[\rmL^{\rmH}]} \right)$ for some $\bmL\in\INT_{\Gamma}(\bmH,\bmG)$ and $g\in \rmG$.
    Indeed by Theorem \ref{theorem_nondivergence}, every sequence $(g_n)$ in $\rmG$, after passing to a subsequence, can be written as $g_n= \delta_n \gamma_n h_n$ with $(\delta_n)$ converging to some $\delta_{\infty}$, $(\gamma_n) \subset \Gamma^{\sma}$ with $\gamma_n \bmH \gamma_n^{-1}$ converging to some $\bmL \in \INT$ and $(h_n)\subset \rmH$. By Lemma \ref{lemma_limit_fat_homogeneous_measures_infinite}, $\lim \alpha^{\psi}_{g_n}  \left( \rmm^{\psi}_{[\rmH]} \right) =
    \alpha^{\psi}_{\delta_{\infty}} \left( \rmm^{\psi}_{[\rmL^{\rmH}]} \right)$. So we are done.
\end{proof}

Define $\Psi_{\bmH}^{\psi}(g.o):= \alpha_g^{\psi} \left( \rmm_{[\rmH]}^{\psi} \right)$.
Let $\rmX:= \wtrmX_{\rmH}^{\IMTP,\cor}$ and $\rmD:= \wtrmD_{\rmH}^{\IMTP,\cor}$ as before.
Similar to Theorem \ref{theorem_almost_algebraic_compactification_dominate_measure_compactification_finite_volume}, we have:

\begin{thm}\label{theorem_almost_algebraic_compactification_dominate_measure_compactification_infinite_volume}
    Assume that $\bmG$, $\bmH$ are connected, reductive,  $\bmZ_{\bmG}(\bmH)^{\circ} \subset \bmH$ and 
    $\psi$ is as above Equa.(\ref{definition_prob_psi}) where $\Prob^{\psi}(\rmG/\Gamma)$ is defined.
    Let $(x_n)$ be a sequence in $\rmX \setminus \rmD$ converging to some $x_{\infty} \in \rmX$. Then $\left( \Psi^{\psi}_{\bmH}(x_n) \right)$ is convergent in $ \Prob^{\psi}(\rmG/\Gamma)$. Consequently, $\Psi^{\psi}_{\bmH}$ extends to a continuous map $\overline{\Psi}^{\psi}_{\bmH} : \rmX \to \Prob^{\psi}(\rmG/\Gamma) $.
\end{thm}

Note that conditions  \hyperref[condition_F_1]{$(\rmF 1)$}, \hyperref[condition_F_2]{$(\rmF 2)$}, \hyperref[condition_C_1]{$(\rmC 1)$} and \hyperref[condition_N_1]{$(\rmN 1)$} hold under the assumptions made here.

\subsection{Compactifications of $U$}
Let $\bmU$ be a variety over $\Q$ equipped with a transitive $\bmG$-action and  assume that there exists $o\in \bmU(\Q)$ such that the stabilizer of $o$ in $\bmG$ is $\bmH$.
Let $\bmH_x$ denote the stabilizer of  $x\in \bmU(\Q)$. For different rational points, $\bmH_x$'s may not be isomorphic to each other. So there seems no canonical choice of $\bmH_x$ to identify $\bmU$ with $\bmG/\bmH_x$.

\subsubsection{Conjugacy between intermediate subgroups}

While $\bmH_x$ and $\bmH_y$ may not be isomorphic for different $x,y \in \bmU(\Q)$, the set of intermediate groups are related.

\begin{lem}\label{lemma_conjugacy_intermediate_groups}
Let $x,y$ be two rational points on $\bmU(\Q)$. Find $f \in \bmG(\overline{\Q})$ such that $f y= x$.
     Assume $\bmL$ is a $\Q$-subgroup normalized by $\bmH_y$. Then $f\bmL f^{-1}$ is defined over $\Q$. Therefore, $ \bmL \to f\bmL f^{-1}$ induces bijections between 
     $\INT_{\Q}(\bmH_y,\bmG )$ and $\INT_{\Q}(\bmH_x,\bmG)$; 
     $\INT^{\obs}_{\Q}(\bmH_y,\bmG )$ and $\INT^{\obs}_{\Q}(\bmH_x,\bmG)$; $\scrP_{\bmH_y}$ and $\scrP_{\bmH_x}$.
\end{lem}
\begin{proof}
     For $\sigma \in \Gal_{\Q}$, we have 
     \begin{equation*}
         fy=x= \sigma(x)= \sigma (fy) = \sigma(f)y
         \implies f^{-1}\sigma(f) \in \bmH_y.
     \end{equation*}
     Let $h_{\sigma}:= f^{-1}\sigma(f) $,
     then $\sigma(f \bmL f^{-1}) =  f h_{\sigma }\bmL h_{\sigma }^{-1} f^{-1} = f \bmL f^{-1} $ for all $\sigma \in \Gal_{\Q}$. So $f \bmL f^{-1}$ is defined over $\Q$.
     The rest of the claim follows from this.
\end{proof}

\subsubsection{Convergent intermediate subgroups}

However, it is not clear whether conjugation by $f$ maps $\INT_{\Gamma}(\bmH_y,\bmG)$ to $\INT_{\Gamma}(\bmH_x,\bmG)$.
Consequently, though Lemma \ref{lemma_sufficient_condition_B_1_one_orbit} or its analogue for $\bmX^{\IMTP}$ is great for counting a single $\Gamma$-orbit, it needs to be modified for counting integral points, which may consist of different $\Gamma$-orbits. 
Define
\begin{equation*}
    \INT_{\bmG(\overline{\Q})}^{\cvg}(\bmH,\bmG) = \left\{ \bmL \in \INT_{\Q}(\bmH,\bmG) \;\middle\vert\;
    \bmL= g \bmL' g^{-1} ,\;\exists\, g \in \bmG(\overline{\Q}),\;\bmL'\in \INT_{\Gamma}(\bmH,\bmG)
    \right\}.
\end{equation*}
\begin{lem}\label{lemma_conjugacy_intermediate_groups_convergence}
    Assumptions same as in Lemma \ref{lemma_conjugacy_intermediate_groups}. Then 
    $ \bmL \to f\bmL f^{-1}$ induces a bijection between 
    $ \INT_{\bmG(\overline{\Q})}^{\cvg}(\bmH_x,\bmG)$ and $ \INT_{\bmG(\overline{\Q})}^{\cvg}(\bmH_y,\bmG)$.
\end{lem}
\begin{proof}
    This follows from Lemma \ref{lemma_existence_convergence_to_full} and Lemma \ref{lemma_conjugacy_intermediate_groups} above.
\end{proof}
We formulate a condition that guarantees the finiteness of this set:
\begin{itemize}
    \item[(N2)]\label{condition_N_2} $\bmN_{\bmG}(\bmH)^{\circ}\subset \bmN_{\bmG}(\bmL)$ for every $\bmL \in \INT_{\Q}^{\obs}(\bmH,\bmG)$.
\end{itemize}
\begin{lem}
     Assume that one of the following is true:
     \begin{itemize}
         \item[(1)]  $\bmG$ and $\bmH$ are reductive and $\bmZ_{\bmG}(\bmH)^{\circ} \subset \bmH \cdot \bmZ(\bmG)$;
         \item[(2)] $\bmH$ contains a maximal unipotent subgroup of $\bmG$.
     \end{itemize}
     Then \hyperref[condition_N_2]{$(\rmN 2)$}  holds.
\end{lem}

In particular the example $\SL_n \times \bmD $ modulo the diagonal embedding of $\bmD$ mentioned in Section \ref{subsection_definition_compactification_intermediate_groups} also satisfies  \hyperref[condition_N_2]{$(\rmN 2)$}.

\begin{proof}
    Part 1 is easy as $\bmH\cdot \bmZ_{\bmG}(\bmH)^{\circ}= \bmN_{\bmG}(\bmH)^{\circ}$ (this is essentially a combination of the case of $\bmH$ being a torus, see \cite[Corollary 3.2.9]{Spr98}, and the case of $\bmH$ being semisimple, see \cite[Theorem 2.1.4]{PlatonovRapinchuk2023Vol1}). Part 2 follows from the structure of subgroups containing some maximal unipotent subgroup.
\end{proof}
\begin{lem}
    Assume \hyperref[condition_C_1]{$(\rmC 1)$}, \hyperref[condition_F_1]{$(\rmF 1)$} and \hyperref[condition_N_2]{$(\rmN 2)$} hold, then $\INT^{\cvg}_{\bmG(\overline{\Q})}(\bmH,\bmG)$ is finite.
\end{lem}
\begin{proof}

    By \hyperref[condition_C_1]{$(\rmC 1)$}, for each $\bmL\in \INT_{\Gamma}(\bmH,\bmG)$, there are finitely many elements $(f_i)_{i \in \scrJ_{\bmL}}$ such that 
    \begin{equation*}
        \bmZ(\bmH,\bmL)(\overline{\Q})
        = \bigsqcup_{i\in \scrJ_{\bmL}} \bmL(\overline{\Q}) f_i \bmN_{\bmG}(\bmH)^{\circ}(\overline{\Q}).
    \end{equation*}
    For $g\in  \bmZ(\bmH,\bmL)(\overline{\Q})$ such that $g^{-1}{\bmL}g$ is defined over $\Q$, write $g= l_g f_i n_g$ according to this decomposition, 
    \begin{equation*}
        g^{-1}\bmL g
        =n_g^{-1}(f_i^{-1}\bmL f_i) n_g = f_i^{-1}\bmL f_i
    \end{equation*}
   by  \hyperref[condition_N_2]{$(\rmN 2)$} applied to $n_{g}^{-1} f_i^{-1}\bmL f_i n_g \in \INT_{\Q}^{\obs}(\bmH,\bmG) $.
    Therefore,
    \begin{equation*}
        \INT^{\cvg}_{\bmG(\overline{\Q})}(
        \bmH, \bmG) \subset 
        \left\{ f_i^{-1} \bmL f_i \;\middle\vert\;
        \bmL \in \INT,\;i\in \scrJ_{\bmL},\;f_i^{-1}\bmL f_i \text{ is defined over }\Q
        \right\}
    \end{equation*}
    is finite by \hyperref[condition_F_1]{$(\rmF 1)$}.
\end{proof}

\subsubsection{Refined compactifications and condition  \hyperref[condition_B_1]{(B1)}}\label{subsection_compactification_IMTP'}
Similar to Definition \ref{defiComptfyX_1} and \ref{defiComptfyX_2}, define (abbreviate $\INT^{\cvg}_{\bmG(\overline{\Q})}:= \INT^{\cvg}_{\bmG(\overline{\Q}) }(\bmH,\bmG)$) 
\begin{equation}\label{equation_compactificatino_IMTP'}
    \bmX^{\IMTP'}_{\bmH} \subset \prod_{\bmL \in \INT^{\cvg}_{\bmG(\overline{\Q})}  \sqcup \scrV(\bmH)} \bmP(\bmV_{\bmL} \oplus \Q)
\end{equation}
and 
\begin{equation}\label{equation_compactificatino_IMTP'}
    \bmX^{\INTP'}_{\bmH} \subset  \bmP(\bmV_{\bmH}\oplus \Q) \times \prod_{\bmL \in \INT^{\cvg}_{\bmG(\overline{\Q})} } \bmP(\bmV_{\frakl} \oplus \Q) \times \prod_{\rho \in \scrV(\bmH)} \bmP(\bmV_{\rho} \oplus \Q)
\end{equation}
which admit natural morphisms onto $\bmX^{\IMTP}_{\bmH}$ or $\bmX^{\INTP}$.
In defining $\bmB^{\IMTP'}_{\bmH}$ (similar to Definition \ref{definition_B_INTP}), we replace the index $\INT_{\Gamma}(\bmH,\bmG)$ by the larger $ \INT^{\cvg}_{\bmG(\overline{\Q})}(\bmH,\bmG)$.
For simplicity refer $
    \left( 
       \wtbmX^{\IMTP'}_{\bmH}, \wtbmD^{\IMTP'}_{\bmH},\wtbmB^{\IMTP'}_{\bmH}
    \right)
$ as the \textit{decorated log smooth $\IMTP'$ compactification} of $\bmU$. \textit{Decorated log smooth $\INTP'$ compactification} is similarly defined.
We have the following:
\begin{lem}\label{lemma_existence_convergence_boundary}
Consider one of the following two situations:
\begin{itemize}
     \item[(1)]  Assume \hyperref[condition_F_1]{$(\rmF 1)$}, \hyperref[condition_F_2]{$(\rmF 2)$}, \hyperref[condition_C_1]{$(\rmC 1)$}, \hyperref[condition_N_1]{$(\rmN 1)$} and \hyperref[condition_N_2]{$(\rmN 2)$} hold. Let $(\bmX,\bmD,\bmB)$ be the decorated log smooth $\IMTP'$ compactification;
     \item[(2)] Assume that the unipotent radical of $\bmG$ lies in the center of $\bmG$. Assume \hyperref[condition_F_1]{$(\rmF 1)$}, \hyperref[condition_F_2]{$(\rmF 2)$},  \hyperref[condition_C_1]{$(\rmC 1)$}, \hyperref[condition_N_2]{$(\rmN 2)$} hold and that the modular character of $\bmL$ is trivial on $\bmH$ for all $\bmL\in \INT^{\cvg}_{\bmG(\overline{\Q})}(\bmH,\bmG)$. Let $(\bmX,\bmD,\bmB)$ be the decorated log smooth $\INTP'$ compactification.
\end{itemize}
    Write $\rmX:= \overline{\bmX\setminus \bmD (\R)}$. Let $x\in \bmU(\Q)$ and $(g_n)\subset \rmG$ be such that $(g_n.x)$ converges to some $\bmv \in \rmX$. Then
   \begin{equation*}
        \lim_{n\to \infty} \left[ (g_n)_* \rmm_{[\rmH^{\circ}_x]} \right]  = \left[ \rmm_{[\rmG]} \right] \iff
        \bmv \in \bmB(\R).
    \end{equation*}
    In particular, condition \hyperref[condition_B_1]{$(\rmB 1)$} holds for $(\bmX,\bmD,\bmB)$.
\end{lem}

\begin{proof}
The proof is similar to that of Lemma \ref{lemma_sufficient_condition_B_1_one_orbit} and is omitted.
\end{proof}

\subsection{Wrap-up}

By combining efforts made so far, we have arrived at a few theorems.

\begin{thm}\label{theorem_existence_good_height_1}
Assume $\bmG \neq \bmH$.
     Consider one of the following two situations:
\begin{itemize}
     \item[(1)]  Assume \hyperref[condition_F_1]{$(\rmF 1)$}, \hyperref[condition_F_2]{$(\rmF 2)$}, \hyperref[condition_C_1]{$(\rmC 1)$}, \hyperref[condition_N_1]{$(\rmN 1)$} and \hyperref[condition_N_2]{$(\rmN 2)$} hold. Let $(\bmX,\bmD,\bmB)$ be the decorated log smooth $\IMTP'$ compactification;
     \item[(2)] Assume that the unipotent radical of $\bmG$ commutes with $\bmG$. Assume \hyperref[condition_F_1]{$(\rmF 1)$}, \hyperref[condition_F_2]{$(\rmF 2)$}, \hyperref[condition_C_1]{$(\rmC 1)$}, \hyperref[condition_N_2]{$(\rmN 2)$} hold and that the modular character of $\bmL$ is trivial on $\bmH$ for all $\bmL\in \INT^{\cvg}_{\bmG(\overline{\Q})}(\bmH, \bmG)$. Let $(\bmX,\bmD,\bmB)$ be the decorated log smooth $\INTP'$ compactification.
\end{itemize} 
     Also, assume equivalent conditions in Lemma \ref{lemma_existence_convergence_to_full} hold for $\bmH_x \subset \bmG$ for all $x\in \bmU(\Q)$.
     Then condition \hyperref[condition_B_4]{$(\rmB 4)$} holds for $(\bmX,\bmD,\bmB)$. Therefore,
      there exists an effective divisor $\bmL$ supported on $\bmD$ such that $\Ht_{\bmL}$ is good if $\Pic(\bmU)$ is torsion and is good with weights in general.
\end{thm}
\begin{proof}
     By assumption, there exists $(\gamma_n) \subset \Gamma$ that is unbounded modulo the normalizer of $\bmH_x$. 
    So condition \hyperref[condition_B_4]{$(\rmB 4)$} holds by applying Theorem \ref{theorem_anticanonical_general_arithmetic}.  This plus Lemma \ref{lemma_exists_good_heights} 
    show the existence of height satisfying  \hyperref[condition_BH_1]{$(\rmB \rmH 1)$}.
    \hyperref[condition_B_1]{$(\rmB 1)$} holds by Lemma \ref{lemma_existence_convergence_boundary}.
    Lemma \ref{lemma_condition_good_height} then concludes the proof.
\end{proof}

\begin{thm}\label{theorem_log_anti_canonical_height_general_height}
    Assume that $\bmG$ and $\bmH$ are connected, reductive and $\bmZ_{\bmG}(\bmH)^{\circ} \subset \bmH$. Also, assume that $\bmG \neq \bmH$ and the projection of $\bmH$ to the compact $\Q$-factor of $\bmG$ is surjective.
    Then \hyperref[condition_B_2]{$(\rmB 2)$} and \hyperref[condition_K_1]{$(\rmK 1)$} hold for $(\bmX,\bmD,\bmB)$, the decorated log smooth $\IMTP'$ compactification. In particular, the log anti-canonical height is good if $\Pic(\bmU)$ is torsion and is good with weights in general. Every other height is ok if $\rmm_{[\rmH]}$ is finite and is ok with weights otherwise.
\end{thm}

Note that under these assumptions, \hyperref[condition_F_1]{$(\rmF 1)$}, \hyperref[condition_F_2]{$(\rmF 2)$}, \hyperref[condition_C_1]{$(\rmC 1)$}, \hyperref[condition_N_1]{$(\rmN 1)$}, \hyperref[condition_N_2]{$(\rmN 2)$} and equivalent conditions in Lemma \ref{lemma_existence_convergence_to_full} hold for $\bmH_x$ for all $x\in \bmU(\Q)$.

\begin{proof}
    Condition \hyperref[condition_K_1]{$(\rmK 1)$} follows from Theorem \ref{theorem_anticanonical_affine_reductive}.
    By Theorem \ref{theorem_nondivergence}, every limit measure $\mu$ is nonzero and hence homogeneous.
    By Lemma \ref{lemma_existence_convergence_to_full}, there exists $(\gamma_n) \subset \Gamma$ such that $\left[ (\gamma_n)_* \rmm_{[\rmm_{\rmH_x^{\circ}}]} \right]$ converges to $\left[ \rmm_{[\rmG]} \right] $ for every $x \in \bmU(\Q)$. Hence \hyperref[condition_B_2]{$(\rmB 2)$} holds.
    Lemma \ref{lemma_log_anti_canonical_good} gives condition \hyperref[condition_BH_1]{$(\rmB\rmH 1)$} and Lemma \ref{lemma_existence_convergence_boundary} gives 
    \hyperref[condition_B_1]{$(\rmB 1)$}. So we are done by invoking Lemma \ref{lemma_condition_good_height}.

    When $\rmm_{[\rmH]}$ is finite, $\bmH$ and $\bmG$ have no nontrivial $\Q$-characters under our assumption.
    Also, thanks to Theorem \ref{theorem_almost_algebraic_compactification_dominate_measure_compactification_finite_volume} and \ref{theorem_nondivergence}, condition \hyperref[condition_S_1]{$(\rmS 1)$} and \hyperref[condition_D_1]{$(\rmD 1)$} hold, which imply \hyperref[condition_H_1]{$(\rmH 1)$} by Lemma \ref{lemma_ok_heights}. Then one can conclude with Lemma \ref{lemma_H_1_imply_ok}.
    
    In the infinite-volume case, apply Theorem \ref{theorem_almost_algebraic_compactification_dominate_measure_compactification_infinite_volume}   and \ref{theorem_nondivergence} to get \hyperref[condition_S_1]{$(\rmS 1)$} and \hyperref[condition_D_2]{$(\rmD 2)$}.
    Then we get \hyperref[condition_H_2]{$(\rmH 2)$} by Lemma \ref{lemma_ok_heights_weights}.
    Finally one concludes with Lemma \ref{lemma_H_2_imply_ok_weights}.
\end{proof}

\subsection{Lift equidistributions}

We would like to extend Theorem \ref{theorem_existence_good_height_1} to a more general setting. This is necessary to prove Theorem \ref{theorem_Hardy_Littlewood} in full generality.

\subsubsection{Standing assumptions}\label{subsection_assumption_lift_equidistributions}

Take $\bmG$ to be a connected linear algebraic group over $\Q$,
$\bmU$ to be a homogeneous variety under $\bmG$ with a point $o\in \bmU(\Q)$ and $\bmH$ to be the stabilizer of $o$ in $\bmG$.
Assume 
\begin{assumption}\label{assumption_unipotent_radical_commute}
    $\bmR_{\bmu}(\bmG)$ commutes with the reductive part of $\bmG$. Namely, $\bmG = \bmG^{\red} \times \bmR_{\bmu}(\bmG)$.
\end{assumption}

Thus $\bmG$ can be written as $(\bmG^{\cpt}\cdot \bmG^{\nc} \cdot \bmZ(\bmG)) \times \bmR_{\bmu}(\bmG)$. 
Under these assumptions, here are a few natural quotient morphisms:
\begin{itemize}
    \item $p^{\stau}: \bmG \to \bmG/ \bmG^{\nc}$; $p_1: \bmG \to \bmG_1:= \bmG/\bmZ(\bmG)^{\spl}$,
    $p_{12}: \bmG_1 \to \bmG_2:= \bmG/ (\bmG^{\cpt}\cdot \bmZ(\bmG) )\times \bmR_{\bmu}(\bmG) $, and
     $p_{2}: \bmG \to \bmG_2$ is the composition $p_{12} \circ p_1$.
\end{itemize}
Regarding $\bmH$, we assume that
\begin{assumption}\label{assumption_H_is_not_small}
    $\bmH$ is a connected observable $\Q$-subgroup of $\bmG$,
    $p^{\stau}$ is surjective restricted to $\bmH$,
   $p_1(\bmH)$ is observable and  $\scrP^{\max}_{\bmH}$ is finite.
\end{assumption}
Note that under this assumption, $\bmG/\bmH$ is a homogeneous space under $\bmG^{\nc}$. However, we get some advantage by considering the larger automorphism subgroup $\bmG$.
Also, every parabolic $\Q$-subgroup of $\bmG$ must contain $ (\bmG^{\cpt}\cdot \bmZ(\bmG) )\times \bmR_{\bmu}(\bmG) $, thus $\scrP^{\max}_{\bmH}$ being finite is equivalent to $\scrP^{\max}_{\bmH_1}$ being finite, which is again equivalent to $\scrP^{\max}_{\bmH_2}$ being finite.

\begin{rmk}
    By Lemma \ref{lemma_existence_convergence_to_full}, $p^{\stau}$ being surjective is a necessary assumption for the existence of $(\gamma_n)\subset \Gamma$ such that $\left[ (\gamma_n)_* \rmm_{[\rmH^{\circ}]} \right] $ converges to $\left[ \rmm_{[\rmG]} \right]$ since maximal proper subgroups of a unipotent group is normal. Of course, when $\bmG^{\red}$ acts nontrivially on $\bmR_{\bmu}(\bmG)$, this is no longer necessary.
\end{rmk}

Recall that the stabilizer of a general point $x\in \bmU(\Q)$ in $\bmG$ is denoted as $\bmH_x$. Its image in $\bmG_i$ is denoted as $\bmH_{x,i}$ for $i=1,2$. 
We have assumed $\bmH_1:=p_1(\bmH)$ to be observable, thus $\bmH_2:=p_2(\bmH)$ is also observable (see Lemma \ref{lemma_projection_observable_compact_fibration}) and $\rmm_{[\rmH^{\circ}_{x,2}]}
$ makes sense.

\begin{assumption}\label{assumption}
     There exists (and we fix such) a $\bmG_2$-pair $(\bmX,\bmB)$ over $\Q$ with $\bmX$ being a  $\bmG_2$-equivariant compactification of $\bmG_2/\bmH_2$ such that
     \begin{itemize}
         \item[1.] for every $x\in (\bmG_2/\bmH_2)(\Q)$, the analytic closure of $\rmG_2.x $ intersects with $\bmB(\R)$;
         \item[2.] for every $x\in (\bmG_2/\bmH_2)(\Q)$ and $(g_n)\subset \rmG_2$ such that every limit point of $(g_n.x)$ is contained in $\bmB(\R)$, we have
          \begin{equation*}
        \lim_{n\to \infty} \left[ (g_n)_* \rmm_{[\rmH^{\circ}_{2,x}]}
        \right] = \left[ \rmm_{[\rmG_2]}  \right]
    \end{equation*}
    where $\bmH_{2,x}$ denotes the stabilizer of $x$ in $\bmG_2$.
     \end{itemize}
\end{assumption}
The compactification of $\bmG_2/\bmH_2$ as in the assumption can be easily upgraded to one of $\bmU\cong \bmG/\bmH$ with the same conclusion. Namely, for every $x\in \bmU(\Q)$ and every $(g_n)\subset \rmG$ with $\lim g_n.x \in \bmB(\R)$, one has $\lim \left[ p_2(g_n)_* \rmm_{[\rmH^{\circ}_{x,2}]}
        \right] = \left[ \rmm_{[\rmG_2]}  \right]$ and such a sequence $(g_n)$ does exist.
        Note that if $\phi$ is the quotient morphism $\bmG/\bmH \to \bmG_2 /\bmH_2$ and $x\in \bmG/\bmH(\Q)$, then $\bmH_{x,2} = \bmH_{2, \phi(x)}$.
Finally,
\begin{itemize}
  \item let  $\Gamma$ be an arithmetic subgroup of $\bmG$ and $\Gamma_i:= p_i(\Gamma)$.
\end{itemize}
Note that $\Gamma_i$ is an arithmetic subgroup of $\bmG_i$ for $i=1,2$ (see \cite[Theorem 4.7]{PlatonovRapinchuk2023Vol1}).

\subsubsection{Main theorem and outline of the proof}

\begin{thm}\label{theorem_existence_good_boundary}
    There exist a smooth $\bmG$-pair $(\bmX,\bmD)$ over $\Q$ with $\bmU:=\bmX \setminus \bmD$ equivariantly isomorphic to $\bmG/\bmH$, and $\bmB \subset \bmD$,  a union of irreducible components of $\bmD$ over $\Q$, such that \hyperref[condition_B_1]{$(\rmB 1)$} and \hyperref[condition_B_4]{$(\rmB 4)$} hold.
\end{thm}

\begin{proof}[Proof of Theorem \ref{theorem_existence_good_boundary}]
We explain that Theorem \ref{theorem_existence_good_boundary} follows from Lemma 
\ref{lemma_strong_convergence_equivalence} and \ref{lemma_existence_good_boundary}  below.

We fix some $x\in \bmU(\Q)$.

Let $(\bmX_1, \bmD_1,\bmB_1)$ be the one obtained by Lemma \ref{lemma_existence_good_boundary}.
    $(\bmX_3,\bmB_3)$ is constructed by
    firstly applying Theorem \ref{theorem_resolution_of_singularity} to $(\bmX_1,\bmB_1)$ to get
    $(\bmX_2,\bmD_2,\bmB_2)$ ($\bmD_2$ is the inverse image of $\bmD_1$).
    Apply    Theorem \ref{theorem_resolution_of_singularity} again to $(\bmX_2,\bmD_2)$ to get $(\bmX_3,\bmD_3,\bmB_3)$ where $\bmB_3$ is the inverse image of $\bmB_2$.
    Thus $(\bmX_3,\bmD_3)$ is a smooth $\bmG$-pair over $\Q$ and $\bmB_3$ is a union of irreducible components of $\bmD_3$ that is over $\Q$.

    Now take a sequence $(g_n)\subset \rmG$  such that $(g_n .x)$ converges to a point in $\bmB_3(\R)$.
    By Lemma \ref{lemma_existence_good_boundary},  we find $(h_n) \subset \rmH_{x}^{\circ}$,  a bounded sequence $(b_n)$ in $\rmG$ and $(\gamma_n) \subset \Gamma$ 
    such that $g_nh_n = b_n \gamma_n$ and
    \begin{equation*}
         \norm{\Ad(p_1(\gamma_n))\bmv_{\bmP}} \to \infty ,\quad \forall \, \bmP \in \scrP^{\max}_{\bmH_{x,1}}.
    \end{equation*}
 Moreover, Lemma \ref{lemma_existence_good_boundary} asserts that
    \begin{equation*}
        \lim_{n \to \infty} \left[ 
          p_1(\gamma_n)_* \rmm_{[\rmH_{x,1}^{\circ}]} 
        \right]   = \left[  \rmm_{[\rmG_1]} \right].
    \end{equation*}
    Hence  $
    \left(
         p_1(\gamma_n) \bmH_{x,1} p_1(\gamma_n)^{-1}
    \right)
    $  converges to $\bmG_1$.
    Passing to a subsequence if necessary, assume that $
    \left(
         p_1(\gamma_n) \bmH_{x,1} p_1(\gamma_n)^{-1}
    \right)
    $ strongly converges to some $\bmL$. Then by Lemma \ref{lemma_observablility_criterion}, $\bmL$ is observable  in  $\bmG_1$ and hence $\bmL=\bmG_1$. Invoking Lemma \ref{lemma_strong_convergence_equivalence}, we have that $
    \left(
         \gamma_n \bmH_x \gamma_n^{-1}
    \right)
    $  strongly converges to $\bmG$, showing that 
    \[
        \lim_{n \to \infty} \left[ 
          (\gamma_n)_* \rmm_{[\rmH_{x}^{\circ}]} 
        \right]   =
         \lim_{n \to \infty}  \left[ (g_n)_* \rmm_{[\rmH_{x}^{\circ}]} 
        \right]   =
         \left[  \rmm_{[\rmG]} \right].
    \]
    This verifies condition \hyperref[condition_B_1]{$(\rmB 1)$}.
        
     By Lemma \ref{lemma_existence_good_boundary}, we can find $(\gamma_n)\subset \Gamma$ with $\lim \gamma_n. x \in \bmB_3(\R)$. Thus, some irreducible component $B$ of $\bmB_3(\R)$ is contained in $\overline{\rmG.x}$. So we have some sequence $(g_n') \subset \rmG$ such that $(g_n'.x)$ converges to $B^{\circ}$ (the complement in $B$ of other boundary components).  By \hyperref[condition_B_1]{$(\rmB 1)$}, $\lim \left[ (g'_n)_* \rmm_{[\rmH_{x}^{\circ}]} 
        \right]   =
         \left[  \rmm_{[\rmG]} \right]$. In particular, $g_n'.x =b_n \gamma'_n.x$ for some bounded $(b_n)$ in $\rmG$ and $(\gamma_n') \subset \Gamma$. Thus  \hyperref[condition_B_4]{$(\rmB 4)$} follows from Theorem 
 \ref{theorem_anticanonical_general_arithmetic}   since $(\gamma_n')$ can not be bounded modulo the normalizer of $\bmH$.
\end{proof}

\subsubsection{Strongly convergence and observability}

Let us introduce  a simple criterion of observability.

\begin{lem}\label{lemma_observablility_criterion}
Assume that $\bmA$  is a connected linear algebraic group over $\Q$ with no nontrivial $\Q$-characters and $\bmC$ is a connected observable $\Q$-subgroup. 
Let $(\gamma_n)$ be a sequence in some fixed arithmetic subgroup of $\bmA$ and $\bmE$ be a connected $\Q$-subgroup of $\bmA$ such that 
\begin{itemize}
    \item[(1)] 
         $\inf_{\bmP \in \scrP^{\max}_{\bmC}}\norm{
            \Ad(\gamma_n) \bmv_{\bmP}
            } \to \infty$;
     \item[(2)] $(\gamma_{n} \bmC \gamma_{n}^{-1})$ strongly converges to $\bmE$.
    \end{itemize}
Then $\bmE$ is observable in $\bmA$.
\end{lem}

\begin{proof}
If $\bmE$ is contained some maximal proper parabolic $\Q$-subgroup $\bmP$, then $\bmC$ is contained in $\bmP_n:= \gamma_n^{-1}\bmP \gamma_n$ for all $n$. Thus $\norm{ \Ad(\gamma_n) \bmv_{\bmP_n }} \asymp \norm{ \Ad (\gamma_n \gamma_n^{-1} ) \bmv_{\bmP}}$ is bounded, contradicting our first assumption.

So $\bmE$, and hence $\bmE\bmR_{\bmu}(\bmA)$ are not contained in any parabolic $\Q$-subgroup. Thus $\bmE\bmR_{\bmu}(\bmA)/\bmR_{\bmu}(\bmA) \subset \bmA/\bmR_{\bmu}(\bmA)$ is reductive and hence observable. Equivalently,  $\bmE\bmR_{\bmu}(\bmA)$ is observable in $\bmA$. But $\bmE$ is observable in $\bmE\bmR_{\bmu}(\bmA)$ since $\bmE\bmR_{\bmu}(\bmA)/\bmE$ is a homogeneous space under a unipotent group and hence affine. So we are done.
\end{proof}

\subsubsection{Lift strong convergence}

Strongly convergence, unlike convergence, is liftable.

\begin{lem}\label{lemma_strong_convergence_equivalence}
    Let $(\gamma_n)$ be a sequence in $\Gamma$ and $(\bmL_n)$ be a sequence of connected observable $\Q$-subgroups of $\bmG$. 
    Assume that $p^{\stau}$ is surjective restricted to every $\bmL_n$.
    The followings are equivalent:
    \begin{itemize}
        \item[(1)] $(\bmL_n)$  strongly converges to $\bmG$;
        \item[(2)] $(p_1(\bmL_n))$  strongly converges to $\bmG_1$;
        \item[(3)] $(p_2(\bmL_n))$  strongly converges to $\bmG_2$.
    \end{itemize}
\end{lem}

\begin{proof}
    The nontrivial implication is $3. \implies 1.$

    By passing to a subsequence if necessary, we assume that $(\bmL_n)$ strongly converges to $\bmL$. Thus $\bmG_2=p_2(\bmL)$.      
    Write $\bmL= (  \bmL^{\nc} \cdot \bmL^{\cpt}  ) \cdot \bmR(\bmL)$, $\bmR(\bmL)$ being the radical of $\bmL$. Then the image of $\bmR(\bmL)$ under $p_2$ is contained in the radical of $p_2(\bmL)$, which is trivial. Similarly, $\bmL^{\cpt}$ is contained in the kernel of $p_2$. Thus, the subgroup $\bmL^{\nc}$ of $\bmG^{\nc}$ must be the full $\bmG^{\nc}$. In particular, $\bmL$ contains the kernel of $p^{\stau}$. But $p^{\stau}$ is surjective restricted to $\bmL_n$ and hence to $\bmL$, implying that $\bmL$ is equal to $\bmG$.
    
\end{proof}

\begin{rmk}
The analogous statement is wrong replacing ``strongly converges'' by ``converges''. For instance, let $\bmG:=\SL_2 \times \bmG_m$, $\bmG_m$ be embedded in $\SL_2$ as the diagonal torus, and $\bmH$ be the diagonal embedding of $\bmG_m$ in $\bmG$. So here $p_1=p_2:\bmG \to \SL_2=\bmG_2$ are the natural projection. Let $\bmN$ be the upper triangular unipotent subgroup of $\SL_2$ and let $\bmN':=\bmN \times \{\id\}$ in $\bmG$.
Take $\gamma_n:=
\left( \left[
\begin{array}{cc}
  1   & n \\
   0  & 1
\end{array} \right], \id \right).
$
Then one can verify that $\bmH \cdot \bmN'$ is an observable subgroup of $\bmG$ and  
$\left(
\gamma_n \bmH \gamma_n^{-1}
\right) $ converges to $\bmH\cdot \bmN'$. However, $ p_2(\bmH \bmN')$, being the upper triangular Borel subgroup,  is not observable in $\SL_2$ and $\left(
p_2(\gamma_n \bmH \gamma_n^{-1})
\right) $ actually converges to $\SL_2 $.
\end{rmk}

\subsubsection{Lift equidistribution through compact fibrations}

\begin{lem}\label{lemma_projection_observable_compact_fibration}
Let $\bmA$ be a linear algebraic group over $\Q$, $\bmB$ be a $\Q$-anisotropic normal subgroup of $\bmA$ and $\pi: \bmA \to \bmA/\bmB=:\overline{\bmA}$ be the natural quotient. Let $\bmC$ be a $\Q$-subgroup of $\bmA$. If $\bmC$ is observable in $\bmA$, then $\pi(\bmC) $ is observable in $\overline{\bmA}$.
\end{lem}

\begin{proof}
Let $\Gamma$ be an arithmetic subgroup of $\bmA$
and $\overline{\pi}$ be the induced quotient map $\bmA(\R)/\Gamma \to \overline{\bmA}(\R)/ \overline{\Gamma} $ where $\overline{\Gamma} := \pi(\Gamma)$ is an arithmetic subgroup of $\overline{\bmA}$.
By \cite[Corollary 7]{Weiss98}, it is sufficient to show that $\pi(\bmC)(\R) \overline{\Gamma}/\overline{\Gamma} $, or equivalently, $\pi(\bmC(\R))\overline{\Gamma}/\overline{\Gamma}  $ is closed in $\overline{\bmA}(\R)/ \overline{\Gamma}$.
This follows since $\overline{\pi}$ is a proper continuous map by the $\Q$-anisotropic assumption.

\end{proof}
As a consequence, the equivalence between the convergence when applying $p_{12}$ does hold.
\begin{lem}\label{lemma_equivalence_convergence_p_1p_2}
     Let $(g_n)$ be a sequence in $\rmG$ and $x\in \bmU(\Q)$.
    If $\left[ p_2(g_n)_* \rmm_{ [ \rmH^{\circ}_{x,2} ]} \right]$ converges to $\left[ \rmm_{[\rmG_2]} \right]$, then $\left[ p_1(g_n)_* \rmm_{[ \rmH^{\circ}_{x,1} ] } \right]$ converges to $\left[ \rmm_{[\rmG_1]} \right]$.
\end{lem}

\begin{proof}
    As the induced map $\overline{p}_{12}:  \rmG_1/\Gamma_1 \to \rmG_2/\Gamma_2$ is proper, 
    $p_1(g_n)\rmH_{x,1}^{\circ}\Gamma_1/\Gamma_1$ intersects with some bounded subset for all $n$.
    We may assume $p_1(g_n) = \lambda_n \in \Gamma_1$.
    So $( p_{12}( \lambda_n \bmH_{x,1} \lambda_n^{-1}) ) $ converges to  $\bmG_2$.
    After passing to a subsequence assume that $( \lambda_n \bmH_{x,1} \lambda_n^{-1})$ converges to $\bmL$, an observable $\Q$-subgroup of $\bmG_1$.  By Lemma \ref{lemma_projection_observable_compact_fibration}, $p_{12} (\bmL)$ is observable and thus has to be equal to $\bmG_2$.  
  But $\bmL \bmG_1^{\nc} = \bmG_1^{\nc}$  by assumption.
  Hence $\bmL$ has to be $\bmG_1$ by arguments similar to those in the proof of Lemma \ref{lemma_strong_convergence_equivalence}. And we conclude by invoking Theorem \ref{theorem_equidistribution_EMS}.
\end{proof}

\subsubsection{The general case}

For $x \in \bmU(\Q)$, fix a nonempty  open bounded subset $\Omega_x \subset \rmH^{\circ}_x\Gamma/\Gamma$.
Up to now, it only remains to prove the following:

\begin{lem}\label{lemma_existence_good_boundary}
    There exist a $\bmG$-pair $(\bmX,\bmD)$ over $\Q$ and a closed $\bmG$-invariant $\Q$-subvariety $\bmB$ of $\bmD$
    with $\bmX \setminus \bmD$ being $\bmG$-equivariantly isomorphic to $\bmU$ such that  the following holds.
      For every $x \in \bmU(\Q)$, the analytic closure of $\rmG.x$ intersects with $\bmB(\R)$. Moreover,
        for every $x \in \bmU(\Q)$ and $(g_n) \subset \rmG$ such that every limit point of $( g_n. x )$ is contained in $\bmB(\R)$, we have
        \begin{itemize}
            \item[(1)] 
            $
                \lim_{n \to \infty} \left[ (p_1(g_{n}) )_* \rmm_{[
                \rmH^{\circ}_{x,1}
                ]} \right] =
               \left[  \rmm_{[\rmG_1]}
              \right];
           $
            \item[(2)]  there exists a sequence $(h_n)$ in $\rmH^{\circ}_{x}$ such that
                 \begin{itemize}
                         \item[(2.1)]  $(g_nh_n \Omega_x)$ intersects with some bounded 
                          subset of $\rmG/\Gamma$ for all  $n$;   
                          \item[(2.2)] for every  $\bmP \in \scrP^{\max}_{\bmH_{x,1}}$,  $\norm{
                                              \Ad (p_1(g_nh_n) ) \bmv_{\bmP}
                                            } \to \infty$.
                \end{itemize}
        \end{itemize}
\end{lem}

\subsubsection{ Splitting of parabolic characters }

Let
$\bmM_{\bmH}:= \left( 
\bmH \cap  \left( \bmG^{\semisimple} \bmZ(\bmG)^{\an}  \bmR_{\bmu}(\bmG) \right)
\right)^{\circ}$.
Then $\bmH/\bmM_{\bmH}$ is a $\Q$-split torus.

 Let
 $\bmM_{\bmH,1}:= p_1(\bmM_{\bmH})$,
then $\bmH/\bmM_{\bmH}$ naturally surjects onto $\bmH_1/\bmM_{\bmH,1}$. Therefore, $\bmH_1/\bmM_{\bmH,1}$ is also a $\Q$-split torus.
So there exists a $\Q$-split subtorus $\wtbmS^{\bmH}_{1}$  of $\bmH_1$ that maps to $\bmH_1/\bmM_{\bmH,1}$ surjectively with finite kernel.
Choose another $\Q$-split subtorus $\bmS_{\bmH}^1$ of $\bmM_{\bmH,1}$
such that the restriction of $p^{\spl}_{\bmH_1}$ to $\wtbmS_{\bmH_1} : =  \bmS_{\bmH}^{1} \cdot \wtbmS^{\bmH}_{1}$ is surjective onto $\bmS_{\bmH_1}$ with finite kernel. By definition, 
\begin{equation}\label{equation_split_S_H_1}
    \Lie( \wtbmS_{\bmH_1}  )   = \Lie( \bmS_{\bmH}^{1} )\oplus \Lie( \wtbmS^{\bmH}_{1} ).
\end{equation}
Note that the image of $\bmM_{\bmH,1}$ under $p^{\spl}_{\bmH_1}$ is the same as that of $\bmS^1_{\bmH}$: certainly the image of the former contains $p^{\spl}_{\bmH_1}(\bmS^1_{\bmH})$, but it can not be strictly larger as it has to be (almost) disjoint from the image of $ \wtbmS^{\bmH}_{1}$.

For a  $\Q$-eigenvector $\bmv$ of $\bmH_{1}$, let $l_{\bmv}$ (resp. $\alpha_{\bmv}$) be the associated linear functional (resp. character) on $  \Lie( \wtbmS_{\bmH_1}  )   $ (resp. $\bmH$).

Let 
\begin{equation}
    \scrP^{+}_{ \bmH_{1} }  :=   \left\{
         \otimes^{\bmn_{\bmP}}\bmv_{\bmP}   \;\middle\vert\;
         \bmn_{\bullet} \in \Map^*( \scrP^{\max}_{\bmH_{1} } , \Z_{\geq 0})
    \right\}
\end{equation}
where $ \Map^*( \scrP^{\max}_{\bmH_{1} } , \Z_{\geq 0}) $ denotes the collection of set-theoretic maps from $\scrP^{\max}_{\bmH_1}$ to $\Z_{\geq 0}$ excluding the zero map. Note that $ \scrP^{\max}_{ \bmH_{1} } $ can be naturally viewed as a subset of $ \scrP^{+}_{ \bmH_{1} } $.

For a subset $\scrP $ of $ \scrP^{+}_{\bmH_{1}} $, let the corresponding $\calP:= \{ l_{\bmv} ,\; \bmv \in \scrP \} $.
We say that a subset of a $\Q$-vector space is \textit{nondegenerate} iff $\Q_{\geq 0}$-linear combinations of this subset equal to the whole space.

\begin{lem}\label{lemma_characters_parabolic_nondegenerate}
Let $\bmA$ be a connected linear algebraic group over $\Q$ with no nontrivial $\Q$-characters.
Let $\bmB$ be an observable $\Q$-subgroup of $\bmA$. Define $\scrP^{+}_{\bmB}$ in the same way as $ \scrP^{+}_{\bmH_{1}} $ above. And $\calP$ corresponding to some $\scrP \subset \scrP^{+}_{\bmB}$  is viewed as a subset of $\Lie( \bmS_{\bmB} )^{\vee}$ here.
Then $\calP_{\bmB}$, corresponding to $\scrP^{\max}_{\bmB}$, is nondegenerate.
\end{lem}

Before presenting the proof, we firstly deduce that

\begin{lem}\label{lemma_split_parabolic_characters}
    There exists a finite subset $\scrP \subset \scrP^{+}_{\bmH_{1} }$ such that the corresponding $\calP$ is nondegenerate and for each $\bmv \in \scrP$, either $\Lie( \bmS_{\bmH}^{1} )  \subset \ker( l_{\bmv})$ or $\Lie( \wtbmS^{\bmH}_{1} ) \subset \ker( l_{\bmv})$.
\end{lem}

\begin{proof}
    Indeed, it is easy to find a finite subset of $  \Lie( \wtbmS_{\bmH_1}  ) ^{\vee}$ satisfying the required properties. By Lemma \ref{lemma_characters_parabolic_nondegenerate}, replacing this subset by suitable $\Z^+$-multiples one gets a subset of $\calP^{+}_{\bmH_{1} }$ without losing any required properties.
\end{proof}

Henceforth, we fix such a choice of $\scrP_{o} = \scrP_{o}^1 \sqcup \scrP_{o}^{2}$ where
\begin{equation}\label{equation_decompose_characters}
\begin{aligned}
     &\scrP_{o}^1:= \left\{
          \bmv\in \scrP_{o} ,\;\Lie(\bmS_{\bmH}^1 )\subset \ker(l_{\bmv})
    \right\},\;\\
     &\scrP_{o}^{2}:= \left\{
          \bmv\in \scrP_{o} ,\;\Lie( \wtbmS^{\bmH}_{1} )\subset \ker(l_{\bmv})
    \right\}.
\end{aligned}
\end{equation}
Replacing elements $\bmv \in \scrP_{o}$  by $\bmv^{\otimes 2}$ if necessary, we assume that each $\alpha_{\bmv}$ factors through $p^{\spl}_{\bmH_1}\circ p_1$.
Note that for $\bmv \in \scrP_{o}^1$, the character $\alpha_{\bmv}$ vanishes on $\bmM_{\bmH}$ (since it factors through $p^{\spl}_{\bmH_1}(\bmM_{\bmH, 1})$, which is the same as $p^{\spl}_{\bmH_1}(\bmS^1_{\bmH})$) and factors through $\bmS_{\bmG}$.

\subsubsection{Proof of Lemma \ref{lemma_characters_parabolic_nondegenerate} }
For the purpose of the proof, we consider $\scrP_{\bmB}$, including all the proper parabolic $\Q$-subgroups containing $\bmB$, not just the maximal ones. Note that the character associated to a proper parabolic $\Q$-subgroup can be expressed as a $\Q_{\geq }$-linear combination of the characters of \textit{maximal} parabolic $\Q$-subgroups containing it.

By \cite[Theorem 7.3]{Grosshans97}, we find a nonzero highest weight $\Q$-vector $\bmv$ (let $\bmP_{[\bmv]}$ denote the  stabilizer of the line spanned by $\bmv$, a parabolic $\Q$-subgroup)  such that $\bmv$ is fixed by $\bmB$ (so $\bmB$ is contained in $\bmP_{[\bmv]}$) and $\bmR_{\bmu}(\bmB)$ is contained in $\bmR_{\bmu}(\bmP_{[\bmv]})$.  Let $\chi$ be the $\Q$-character of $\bmP_{[\bmv]}$ attached to $\bmv$, then $\chi$ can be written as a $\Q_{> 0}$(not just $\Q_{\geq 0}$ !)-linear combination of characters of maximal parabolic $\Q$-subgroups containing $\bmP_{[\bmv]}$.

Let $\bmB^{\red}$ be a  maximal reductive $\Q$-subgroup of $\bmB$ and
 $\wtbmS_{\bmB}$ be the $\Q$-split part of the central torus of  $\bmB^{\red}$.
Find a maximal reductive $\Q$-subgroup $\bmL_{\bmP_{[\bmv]}}$ of $\bmP_{[\bmv]}$ containing $\bmB^{\red}$. 
Then $\bmB^{\red}\subset \bmL_{\bmP_{[\bmv]}}^{\chi} := \bmL_{\bmP_{[\bmv]}} \cap \ker \chi $. 
Let $\wtbmS_{\bmP_{[\bmv]}}$ be the $\Q$-split part of the central torus of $\bmL_{\bmP_{[\bmv]}}$ and $\bmM_{\bmP_{[\bmv]}}:= ({}^{\circ} \bmL_{\bmP_{[\bmv]}} )^{\circ}$. Then  $\bmL_{\bmP_{[\bmv]}}$  is an almost direct product $\bmM_{\bmP_{[\bmv]}}\cdot \wtbmS_{\bmP_{[\bmv]}}$.

Let 
$\wtbmS^{\chi}_{\bmP_{[\bmv]}} := \wtbmS_{\bmP_{[\bmv]}} \cap \ker \chi$ and  
$\wtbmS_{\bmB}':= \wtbmS_{\bmB} \cdot \wtbmS_{\bmP_{[\bmv]}}^{\chi}$, 
contained in the centralizer of $\bmB^{\red}$. 
Then $\wtbmS_{\bmB}' $ is an almost direct product of  two subtori: 
\[
\bmS_1 := \wtbmS_{\bmP_{[\bmv]}}^{\chi} \text{ and } \bmS_2 := \left( \bmM_{\bmP_{[\bmv]}} \cap \wtbmS_{\bmB}' \right)^{\circ}. 
\]

The proof below actually proves the conclusion with $\bmB$ replaced by $\bmB':= \wtbmS_{\bmB}' \cdot \bmB \cdot \bmR_{\bmu}(\bmP_{[\bmv]})$, which is a stronger claim.
That is to say, we are going to show that the set of linear functionals on $\Lie(\wtbmS_{\bmB}')$ coming from $\scrP_{\bmB'}$ is nondegenerate.
For this purpose, we will divide $\scrP_{\bmB'}$ into two disjoint types.

The first type $\scrP_1$ consists of parabolic $\Q$-subgroups containing $\bmP_{[\bmv]}$.
The other type $\scrP_2$ is constructed as follows. Take a $\Q$-cocharacter $\bma_t : \bmG_m \to \bmS_2$, let 
\begin{equation*}
     \bmQ^1_{\bma_t}  := \left\{
        x\in \bmL_{\bmP_{[\bmv]}} \;\middle\vert\;
        \lim_{t\to 0} \bma_t x \bma_t^{-1} \text{ exists }
     \right\} ,\quad
     \bmQ_{\bma_t}  :=    \bmQ^1_{\bma_t} \ltimes \bmR_{\bmu}(\bmP_{[\bmv]}).
\end{equation*}
Then $   \bmQ^1_{\bma_t}  $ is a parabolic $\Q$-subgroup of $\bmL_{\bmP_{[\bmv]}}$ and $\bmQ_{\bma_t} $ is a parabolic $\Q$-subgroup of $\bmA$. Define $\scrP_2$  to be the collection of all $\bmQ_{\bma_t} $'s constructed this way.
By definition, every element from $\scrP_1 \cup \scrP_2 $ contains $\bmB'$.
Let $\calP_1$ and $\calP_2$ be the associated linear functionals (denoted as $l'_{\bmP}$ as $\bmP$ varies in $\scrP_1\cup \scrP_2$) on $\Lie (\wtbmS_{\bmB}')$. It suffices to show that $\calP_1\cup \calP_2\subset \Lie (\wtbmS_{\bmB}')^{\vee}$ is nondegenerate.

First we explain that $\calP_1$ is nondegenerate when restricted to $\Lie( \bmS_1  )$.
The cone in $\Lie(\wtbmS_{\bmP_{[\bmv]}})^{\vee}$ spanned  by those associated to $\scrP_1$ contains $\chi$ in its interior. Thus the orthogonal projection of this cone to the subspace orthogonal to $\chi$ must be full. Equivalently,  $\calP_1$ is nondegenerate restricted to $\Lie( \bmS_1  )= \Lie(\wtbmS_{\bmP_{[\bmv]}}^{\chi})$.

 
Next we show that $\calP_2$ is nondegenerate restricted to $\Lie( \bmS_2  )$. 
Take $\bmQ=\bmQ_{\bma_t} \in \scrP_2$.
Let $l''_{\bmQ}$ be the linear functional on $\Lie(\bmS_2)$ associated with $\bmQ^1_{\bma_t}$. As the determinant character of the adjoint action of $\bmS_2$ on $\bmR_{\bmu}(\bmP_{[\bmv]})$ is trivial, $l''_{\bmQ}$ coincides with $l'_{\bmQ}$ restricted to $\Lie(\bmS_2)$.
Since $\bmL_{\bmP_{[\bmv]}}$ is reductive, it admits an involution whose restriction to $\bmS_2$ brings each element to its inverse.
Consequently, the $\Q_{\geq 0}$-cone spanned by $\{l''_{\bmQ}\}_{\bmQ\in \scrP_2}$ is actually a $\Q$-linear subspace. So if it is not nondegenerate, then there exists a $\Q$-cocharacter $\bmb_t : \bmG_m \to \bmS_2$ such that $l''_{\bmQ}$ vanishes on the image of $\diff \bmb_t$ for all $\bmQ_{\bma_t}=\bmQ \in \scrP_2$. In particular, this is true when $\bma_t = \bmb_t$, which is a contradiction.

The general case follows since $\Lie(\bmS_2)$ is contained in the common kernel of elements from
$\calP_1$.

\subsubsection{Definition of the compactification}

For $\bmv \in \scrP_{o}^1$, let $V_{\bmv}$ be the $\bmG_1$-representation where $\bmv$ lives.
As is already noted, the character $\alpha_{\bmv} : \bmH \to \bmG_m$ defined by
\begin{equation*}
    p_1(h) .\bmv = \alpha_{\bmv} (h) \bmv,\;\forall \,h\in \bmH
\end{equation*}
factors through some $\beta_{\bmv}: \bmS_{\bmG}\to \bmG_m$.
Define a new $\Q$-linear action $\rho_{\bmv}$ of $\bmG$ on $V_{\bmv}$ twisting the given one by $\beta_{\bmv}^{-1}$:
\begin{equation*}
    \rho_{\bmv}(g)\cdot w :=
    \beta_{\bmv}(p^{\spl}(g))^{-1} \cdot p_1(g). w,\;\forall\, g\in \bmG, \,w\in V_{\bmv}.
\end{equation*}
Thus for every $h \in \bmH$,
\begin{equation*}
    \rho_{\bmv}(h)\bmv= \beta_{\bmv}(p^{\spl}(h))^{-1} \cdot \alpha_{\bmv} (h) \bmv =\bmv,
\end{equation*}
implying that $\bmv$ is fixed by $\bmH$.

 Let $(\bmX_1,\bmB_1)$ be as offered by Assumption \ref{assumption} (by the remark after the assumption, we assume this is a $\bmG$-pair). Recall that $\bmG/\bmH \cong \bmU$ with identity coset sent to $o\in \bmU(\Q)$.
For the rest of this subsection,
    let $\bmX$ be the Zariski closure of 
    \begin{equation*}
    \begin{aligned}
        \bmU &\hookrightarrow \bmX_{1}  \times 
        \bigoplus_{\bmv\in \scrP_{ o }^1} V_{\bmv}   \subset \bmX_{1} \times 
        \prod_{
        \bmv\in \scrP_{o}^1 } \bmP(V_{\bmv} \oplus \Q)
        \\
        g.o &\mapsto ( g.o , \oplus \rho_{\bmv}(g) \bmv),
    \end{aligned}
    \end{equation*}
        Also,  let $\bmD$ be the complement of $\bmU$ in $\bmX$ and
        \begin{equation*}
            \bmB:=\left\{
             \left(
             x, \prod [ \bmx_{\bmv} : t_{\bmv} ]
             \right) \;\middle\vert\;
             x\in \bmB_1;\; \bmt_{\bmv}=0 ,\,\forall \, \bmv \in \scrP_{ o }^1
            \right\}.
        \end{equation*}

   \subsubsection{Proof of Lemma \ref{lemma_existence_good_boundary}, I, equidistribution}

   Take $x \in \bmU(\Q)$ and $(g_n)\subset \rmG$ with $\lim g_n.x \in \bmB(\R)$.
   We firstly prove the ``Moreover, ...'' part of Lemma \ref{lemma_existence_good_boundary}.
   
   Find $g_x \in \bmG(\overline{\Q})$ such that $x = g_x. o$. 
   So $(g_n g_x.o)$ converges to some point in $\bmB(\R)$.
   By our assumption, we may and do require that $g_n$'s and $g_x$ are contained in $ \bmG^{\semisimple }\bmZ(\bmG)^{\an} \times \bmR_{\bmu}(\bmG)  $.
   Item ($1$) of Lemma \ref{lemma_existence_good_boundary} holds for free by Lemma \ref{lemma_equivalence_convergence_p_1p_2}.   
   In particular, replacing $g_n$ by $g_nh_n^1$ for some $(h_n^1) \subset \rmH_{x}^{\circ}$,  there are $(b_{n,1})$ bounded in $\rmG_1$, $(\gamma_{n,1})\subset \Gamma_1$ such that $p_1(g_n) = b_{n,1}\gamma_{n,1}$ and 
   $ \lim_{n \to \infty} \left[ (\gamma_{n,1})_* \rmm_{[
                \rmH^{\circ}_{x,1}
                ]} \right] =
               \left[  \rmm_{[\rmG_1]}
              \right]$.
   
   Now turn to item ($2$). Namely, we need to find $(h_n) \subset \bmH_x(\R)^{\circ}$ such that 
   \begin{itemize}
        \item[1.] $(g_n h_n \Omega_x)$ is nondivergent in $\rmG/\Gamma$;
        \item[2.] for every $\bmP \in \scrP^{\max}_{\bmH_{x,1}} $, $ \norm{ 
        \Ad(p_1(g_n h_n) ) \bmv_{\bmP}
        }   \to \infty$.
   \end{itemize}

   Thanks to Lemma \ref{lemma_conjugacy_intermediate_groups},  for every $\Q$-subgroup $\bmF$ of $\bmG_1$ that is normalized by $\bmH_1$, $p_1(g_x) \bmF p_1(g_x)^{-1}$ is defined over $\Q$.
In particular, conjugating by $p_1(g_x)$ induces a bijection between $\scrP^{\max}_{\bmH_1}$ and $\scrP^{\max}_{\bmH_{x,1}}$.
The proof of Lemma \ref{lemma_conjugacy_intermediate_groups} also shows that conjugating by $p_1(g_x)$ induces a bijection between
 $\frakX_{\Q}^*(\bmH_1)$ and $\frakX_{\Q}^{*}(\bmH_{x,1})$.
 
For $i=1,2$, let $\overline{\calP}_o^{i}$ be the characters of  $\bmH_1$ associated with  $\scrP_o^i$ and $\overline{\calP}_x^{i}$ (resp.  $\scrP^i_x$) be the subset of $\frakX_{\Q}^{*}(\bmH_{x,1})$ (resp. $\scrP^+_{\bmH_{x,1}}$) induced from conjugating by $p_1(g_x)$. 
By the definition of $\bmB$, noting that $\beta_{\bmv}(g_ng_x)=1$, we have
\[
    \norm{p_1(g_n g_x) . \bmv} \to + \infty ,\quad \forall \, \bmv \in \scrP^1_{o}.
\]
Then 
\begin{equation}\label{equation_diverge_first_parabolic_characters}
     \norm{p_1(g_n) .\bmv} \to +\infty ,\quad \forall\, \bmv \in \scrP^1_{x}.
\end{equation}

We choose a $\Q$-split subtorus $\wtbmS_{\bmH_{x,1}} \subset \bmH_{x,1}$ such that $p^{\spl}_{\bmH_{x,1}}$ restricted to $\wtbmS_{\bmH_{x,1}}$  is surjective onto $\bmS_{\bmH_{x,1}}$ with finite kernel\footnote{It is not clear that $p_1(g_x) \wtbmS_{\bmH_1} p_1(g_x)^{-1}$ is defined over $\Q$. So this can not be used as the definition of $\wtbmS_{\bmH_{x,1}}$.}.
For $i=1,2$, let $\calP_{x}^i$ collect the linear functionals on $\Lie( \bmS_{\bmH_{x,1}} )$ induced from  $\overline{\calP}_x^{i}$. Let $ \bmS^1_{\bmH_x} $ be  the identity component of the common kernel of $\overline{\calP}_x^{1}$ in  $\wtbmS_{\bmH_{x,1}}$ and $\wtbmS^{\bmH_x}_1 $ be that of $\overline{\calP}_x^{2}$.
Then
$\wtbmS_{\bmH_{x,1}}   =  \bmS^1_{\bmH_x}  \cdot \wtbmS^{\bmH_x}_1 $ is an almost direct product of two $\Q$-subtori. This is similar to what we did when $x=o$ in Equa.(\ref{equation_split_S_H_1}). The difference is that the conclusion of Lemma \ref{lemma_split_parabolic_characters} becomes our definition now.  We need the following lemma:
    \begin{lem}\label{lemma_orbit_on_v_P_diverge_or_bounded}
       Let $\bmA$ be a linear algebraic group over $\Q$ and $(\gamma_n)$ be contained in a fixed arithmetic subgroup of $\bmA$. 
       Let $\bmC$ be a connected $\Q$-subgroup and assume that $\scrP^{\max}_{\bmC}$ is finite.
       Let $ \wtbmS_{\bmC}$ be a $\Q$-split subtorus of $\bmC$ that surjects onto $\bmS_{\bmC}$ under $p^{\spl}_{\bmC}$  with a finite kernel.
        After passing to a subsequence, one of the following holds:
        \begin{itemize}
            \item[(1)] there exists $(z_n) \subset \Lie(   \wtbmS_{\bmC} (\R) ) $ such that 
            \begin{equation*}
                \norm{
                  \Ad(  \gamma_n \exp( z_n )  ) \bmv_{\bmP}
                } \to \infty,\;\forall \, \bmP \in
                \scrP^{\max}_{  \bmC   };
            \end{equation*}
            \item[(2)] there exists $\bmv \in \scrP^{+}_{  \bmC  }$ fixed by $ \bmC $ such that $(\norm{  \gamma_n. \bmv})$ is bounded.
        \end{itemize}
    \end{lem}
    
    \begin{proof}
    The proof is essentially contained in Section \ref{subsection_unbounded_polytopes}. Indeed, by \cite[Lemma 3.4]{zhangrunlinCompositio2021}, if the first alternative were not true, then there are positive integers $(a_i)$ and $(\bmP_i)\subset \scrP^{\max}_{\bmC}$ such that $\bmv:= \otimes \bmv_{\bmP_i}^{\otimes a_i}$ is fixed by $\bmC$ and $(\gamma_n. \bmv)$ is bounded.
    \end{proof}

   By our assumption, 
   Lemma \ref{lemma_orbit_on_v_P_diverge_or_bounded} is applicable to $\bmA=\bmG_1$, $\bmC=\bmH_{x,1}$ and $\gamma_n=\gamma_{n,1}$.  As the second alternative violates the equidistribution towards $\rmm_{[\rmG_1]}$, there exist  
   $(s_n)\subset \Lie(  \bmS^1_{\bmH_x} (\R))$ and $(t_n)\subset \Lie(\wtbmS^{\bmH_x}_1(\R))$ such that 
   \begin{equation*}
       \norm{
       \Ad (
       p_1 ({g_n} )  \exp(s_n+ t_n)
       ) \bmv_{\bmP}
       } \to \infty,\;\forall\,\bmP \in \scrP^{\max}_{\bmH_{x,1}}.
   \end{equation*}
   Therefore,
   \begin{equation*}
       \norm{
       p_1 ( {g_n} ) \exp(s_n+ t_n)
       . \bmv
       } \to \infty, \; \forall \, \bmv \in \scrP^+_{\bmH_{x,1}}.
   \end{equation*}
   Now observe that one can replace $(t_n)$ by zero.
   Indeed, for every $\bmv \in \scrP^2_x$,
   \begin{equation*}
       \norm{ 
       p_1 ( {g_n} )  \exp(s_n).
        \bmv 
       }
       =
       \norm{ 
       p_1 ( {g_n} )  \exp( s_n + t_n).
        \bmv 
       } \to \infty.
   \end{equation*}
   On the other hand, by Equa.(\ref{equation_diverge_first_parabolic_characters}),  \begin{equation*}
       \norm{ 
       p_1 ( {g_n} )  \exp( s_n).
        \bmv 
       }
       =
       \norm{ 
       p_1 ( {g_n} ) .
        \bmv 
       } \to \infty
   \end{equation*}
   for every $\bmv \in \scrP^1_x $.
   Find $(h_n) \subset \bmM_{\bmH_x}(\R)^{\circ}$ ($\bmM_{\bmH_x}:= \left( \bmH_x \cap {}^{\circ}\bmG \right)^{\circ}$)  such that $\exp(s_n) = p_1(h_n)$, then 
   \begin{equation}\label{equation_parabolic_diverge}
   \norm{
     \Ad(p_1(g_nh_n)) \bmv_{
     \bmP}
     }
     \to \infty,\; \forall \, \bmP \in \scrP^{\max}_{\bmH_{x,1}}
   .
   \end{equation}
   
To complete the proof, it only remains to verify that $(g_n h_n \Omega_x) $ is nondivergent in $\rmG/\Gamma$.
    Since $p^{\spl}( g_n h_n \Omega_x)= p^{\spl}(\Omega_x)
    $ is bounded in $\bmS_{\bmG}(\R)$, it suffices to show that 
    $\left( p_1( g_n h_n ) \Omega_{1,x} \right) $ is nondivergent for every nonempty open bounded subset $\Omega_{1,x}$ of $\rmH_{x,1}^{\circ} \Gamma_1/\Gamma_1$, which is true by  Equa.(\ref{equation_parabolic_diverge}) and Theorem \ref{theorem_nondivergence_bounded_piece}.

   \subsubsection{Proof of Lemma \ref{lemma_existence_good_boundary}, II, nonempty intersection}

   Here we prove the remaining assertion of Lemma \ref{lemma_existence_good_boundary}, namely, 
    \begin{claim}\label{lemma_B_2_nonempty}
        For $x\in \bmU(\Q)$,
        the closure of $\rmG.x$ under the analytic topology intersects with $\bmB(\R)$ nontrivially.
    \end{claim}

     \begin{proof}
    By Assumption \ref{assumption}, choose  a sequence $(g_n)$ in $\rmG$ such that every limit point of $(g_n.x)$ is contained in $\bmB_1(\R)$.
        We will build a new sequence $(g_n') \subset \rmG$ such that $(g_n'.x)$ converges to some point in $\bmB(\R)$.
      
       By Lemma \ref{lemma_equivalence_convergence_p_1p_2},
       \begin{equation}\label{equation_limit_in_B_1}
           \lim_{n\to \infty} \left[
           p_1(g_n)_*\rmm_{[\rmH^{\circ}_{x,1}]} \right] = \left[  \rmm_{[\rmG_1]}  \right].
       \end{equation}
        There exist a bounded sequence $(\delta_n)\subset \rmG_1$, a sequence 
         $(\gamma_n)\subset    \Gamma_1$ and $(h_n) \subset \rmH_{x, 1}^{\circ} $ such that 
               $p_1(g_n)= \delta_n \gamma_n h_n$.
               
Apply Lemma \ref{lemma_orbit_on_v_P_diverge_or_bounded} with $\bmC = \bmH_{x,1}$, $\bmA=\bmG_1$. 
The second alternative can not hold because of  Equa.(\ref{equation_limit_in_B_1}).
Therefore, there exists a sequence $(z_n)$ of $\Lie(\wtbmS_{\bmH_{x,1}}(\R))$ such that
       \begin{equation}\label{equation_divergence_parabolic_vectors}
           \norm{
           \Ad( p_1(g_n) \exp(z_n) ) \bmv_{\bmP}
           } \to \infty,\; \forall \, \bmP \in \scrP^{\max}_{\bmH_{x,1}}.
       \end{equation}
       Recall $g_x\in {}^{\circ}\bmG (\overline{\Q})$ and $x = g_x.o$.
       Take $g_n' \in {}^{\circ}\rmG$ such that $p_1(g_n') = p_1(g_n) \exp(z_n)$.
       We shall show that the sequence $(g_n'. x ) = ( g_n'  g_x.o )$ converges to some point in $\bmB(\R)$.
       For this, we only need to check that
       \begin{equation*}
        \norm{
          \rho_{{\bmv}} (g_n'  g_x) \bmv
        } =
        \norm{
         p_1 (g_n'  g_x )  \bmv
        } \to \infty, \quad \forall \,  \bmv \in \scrP^1_{o}.
       \end{equation*}
       As vectors in $\scrP_{o}^{1} \subset \scrP^+_{\bmH_1}$ are obtained by tensoring powers of $\bmv_{\bmP}$'s for some 
       $\bmP\in \scrP^{\max}_{\bmH_{1}}$, the above holds if
       \begin{equation*}
       \norm{
        \Ad( p_1 (g_n'  g_x ) )  \bmv_{\bmP}
        } \to \infty, \quad \forall \, \bmP  \in \scrP^{\max}_{\bmH_{1}}.
       \end{equation*}
      By Lemma \ref{lemma_conjugacy_intermediate_groups},  $\bmP_x:= p_1(g_x) \bmP p_1(g_x)^{-1}$ belongs to $\scrP^{\max}_{\bmH_{x,1}}$ for every $\bmP  \in \scrP^{\max}_{\bmH_{1}}$. So there is a nonzero complex number $c_{\bmP}$ such that $\Ad( p_1(g_x) )\bmv_{\bmP} = c_{\bmP} \bmv_{\bmP_x}$ and
\begin{equation*}
       \norm{
        \Ad( p_1 (g_n'  g_x ) )  \bmv_{\bmP}
        } =  \normm{c_{\bmP}}  \cdot \norm{
        \Ad( p_1(g_n ) \exp(z_n) )  \bmv_{\bmP_x}
        } ,
\end{equation*}
which diverges to $+\infty$ by Equa.(\ref{equation_divergence_parabolic_vectors}).
       
\end{proof}

\subsection{Proof of Theorem \ref{theorem_Hardy_Littlewood}}\label{subsection_proof_hardy_littlewood}

Let us prove something more general:
\begin{thm}\label{theorem_existence_good_height_counting}
    Let $\bmG$ be a connected semisimple linear algebraic group over $\Q$ without compact factors and $\bmH$ be a connected reductive $\Q$-subgroup without nontrivial $\Q$-characters. Then there exist a smooth $\bmG$-pair $(\bmX,\bmD)$ with $\bmU:= \bmX \setminus \bmD$ equivariantly isomorphic to $\bmG/\bmH$ and an effective divisor $\bmL$ supported on $\bmD$ such that for any smooth metric on $\calO_{\bmX}(\bmL)$, the associated height function $\Ht:=\Ht_{\bmL}$ satisfies:
     \begin{equation}\label{equation_orbit_count_pic_torsion}
        \lim_{R\to \infty}
        \frac{
        \# \Gamma.x \cap B^{\Ht}_R
        }{ \omega_{\infty}(B^{\Ht}_R \cap \rmG.x) }
        = \frac{ \normm{\rmm_{[\rmH_x]}} }{
        \normm{\rmm_{[\rmG]}}
        } ,\quad \forall \, x\in \bmU(\Q),\; \forall\, \text{arithmetic subgroup }\Gamma.
    \end{equation}
\end{thm}
Recall that $\rmm_{\rmG}, \rmm_{\rmH_{x}}$ and $\omega_{\infty}$ are assumed to be compatible.

\begin{rmk}
This result, combined with \cite[Theorem 4.3]{Wei_Xu_2016},  gives generalizations of Theorem \ref{theorem_Hardy_Littlewood}. 
\end{rmk}

\begin{proof}
Let $o\in \bmU(\Q)$ be the image of identity coset of the isomorphism $\bmG/\bmH \cong \bmU$.

    Write $\bmZ_{\bmG}(\bmH)^{\circ} $ as an almost direct product:
\begin{equation*}
    \bmZ_{\bmG}(\bmH)^{\circ} = \bmZ(\bmH)^{\circ} \cdot \bmM_{\bmZ} \cdot \bmT_{\bmZ} =
     \bmZ(\bmH)^{\circ} \cdot \bmM_{\bmZ} \cdot \bmT_{\bmZ}^{\an} \cdot \bmT_{\bmZ}^{\spl}
\end{equation*}
where $\bmM_{\bmZ} $ is semisimple, $\bmT^{\an}_{\bmZ}$ is a $\Q$-anisotropic torus and $\bmT^{\spl}_{\bmZ}$ is a $\Q$-split torus.
Let $\wtbmG:= \bmG \times \bmM_{\bmZ} \times \bmT_{\bmZ}^{\an} \times \bmT_{\bmZ}^{\spl}$. It acts on $\bmU \cong \bmG/\bmH$ by $(a, b_1,b_2,b_3).g\bmH := a g b^{-1}\bmH $ with $b:=b_1b_2b_3$.
Write $\bmF:=  \bmM_{\bmZ} \times \bmT_{\bmZ}^{\an} \times \bmT_{\bmZ}^{\spl}$ for simplicity.
Let $\iota_o: \bmF \to \bmG$ be the natural product map. And let $\Delta_{\iota_o}: \bmF \to \wtbmG $ be defined by $f \mapsto (\iota_o(f),f)$.
For a general $x\in \bmU(\Q)$, fix $g_x\in \bmG(\overline{\Q})$ with $x= g_x.o$.
Let $\iota_x (f): = g_x \iota_o(f) g_x^{-1}$ and $\Delta_{\iota_x}$ be defined similarly. Similar to the proof of Lemma \ref{lemma_conjugacy_intermediate_groups}, one shows that $\iota_x$ is actually defined over $\Q$.
Then, the stabilizer of $x$ in $\wtbmG$ is $\bmH_x \times \Delta_{\iota_x}(\bmF)$.
Let $\wtGamma:= \Gamma \times \bmM_{\bmZ}\cap \Gamma \times  \bmT_{\bmZ}^{\an}\cap \Gamma$.

We now apply discussions in previous sections to the pair $(\wtbmG,\wtbmH)$ in the place of $(\bmG,\bmH)$.
For this, we suppose that $(\wtbmG,\wtbmH)$ satisfies the assumption listed in Section \ref{subsection_assumption_lift_equidistributions}, verification of which will come later. Theorem \ref{theorem_existence_good_boundary} is applicable. 
Note that $\Pic(\bmU)$ is torsion under our assumption.
Combined with Lemma \ref{lemma_condition_good_height} and \ref{lemma_exists_good_heights}, this shows that for every $x$ and $\Gamma$,
\begin{equation*}
        \lim_{R\to \infty}
        \frac{
        \# \wtGamma.x \cap B^{\Ht}_R
        }{ \omega_{\infty}(B^{\Ht}_R \cap \widetilde{\rmG}.x ) 
        }
        = \frac{ \normm{\rmm_{[\widetilde{\rmH}^{\bmone}_x]}} }{
        \normm{\rmm_{[\wtrmG^{\bmone}]} }
        } ,\quad \forall \, x\in \bmU(\Q),\; \forall\, \text{arithmetic subgroup }\Gamma.
\end{equation*}
Note that:
\begin{itemize}
    \item $\wtGamma.x  = \Gamma.x$, $\widetilde{\rmG}.x = \rmG.x$;
    \item $\wtbmG^1 = \bmG \times \bmM_{\bmZ} \times \bmT^{\an}_{\bmZ}$ and $\wtbmH^1_x = \bmH_x \times \Delta_{\iota_x}(\bmM_{\bmZ} \times \bmT^{\an}_{\bmZ})$ and hence
    \begin{equation*}
        \normm{\rmm_{[\widetilde{\rmH}^{\bmone}_x]}}
        = \normm{ \rmm_{[\rmH_x]} } 
        \cdot \normm{ \rmm_{[\rmM_{\bmZ}]} } \cdot \normm{ \rmm_{[\rmT_{\bmZ}^{\an}]} } \;\text{and}\;
        \normm{  \rmm_{ [\wtrmG^{\bmone}] }  }
        = \normm{ \rmm_{[\rmG]} } 
        \cdot \normm{ \rmm_{[\rmM_{\bmZ}]} }  \cdot \normm{ \rmm_{[\rmT_{\bmZ}^{\an}]} }.
    \end{equation*}
\end{itemize}
Thus Equa.(\ref{equation_orbit_count_pic_torsion}) follows.

Now it only remains to verify the assumptions.
As $\wtbmG$ is reductive, Assumption \ref{assumption_unipotent_radical_commute} holds.
As for Assumption \ref{assumption_H_is_not_small},
\begin{itemize}
   \item[(1)] $\wtbmH$ being reductive implies the observability of $\wtbmH$ and $p_1(\wtbmH)$;
   \item[(2)] $\bmG= \bmG^{\nc}$  by assumption, so $p^{\stau}$, the quotient of $\wtbmG$ by $\wtbmG^{\nc}$, factors through that by $\bmG$. Thus $p^{\stau}\vert_{\wtbmH}$ is surjective by the definition of $\wtbmH$;
   \item[(3)] It is direct to verify that the identity component of the centralizer of $p_2(\wtbmH)$ is contained in itself and hence $\scrP_{p_2(\wtbmH)}$ is finite, which is equivalent to $\scrP^{\max}_{\wtbmH}$ being finite. 
\end{itemize}
Finally, let us take care of Assumption \ref{assumption}. By Theorem \ref{theorem_existence_good_height_1} and Lemma \ref{lemma_existence_convergence_boundary}, this reduces to the fact that the identity component of the centralizer of $p_2(\wtbmH)$ in $p_2(\wtbmG)$ is contained in $p_2(\wtbmH)$.

\end{proof}

\section{Example I, Representation of a binary quadratic form by a quaternary form}
\label{section_example_1}

Let $(Q_1,\Q^2)$ and $(Q_2,\Q^4)$ be quadratic froms of rank $2$ and  rank $4$ respectively. Let 
$\bmU_{Q_1,Q_2}$ be the $\Q$-variety defined by
\begin{equation*}
    \bmU_{Q_1,Q_2}(\Q) :=\left\{
      \phi\in \Hom_{\Q}(\Q^2,\Q^4)
       ,\;
      \phi^*Q_2=Q_1
    \right\}.
\end{equation*}

We assume that
\begin{itemize}
    \item  $(Q_1,\Q^2)$ is $\Q$-anisotropic, that is, $Q_1(\bmx)=0$ has only zero solution in $\Q^2$;
    \item  $(Q_1, \R^2)$ has signature $(1,1)$ and $(Q_2,\R^4)$ has signature $(2,2)$;
    \item $\bmU_{Q_1,Q_2}(\Q) \neq \emptyset$.
\end{itemize}

Thus we can find $(Q_1',\Q^2)$ such that $(Q_2,\Q^4) \cong (Q_1,\Q^2) \oplus (Q_1',\Q^2)$. We further assume that
\begin{itemize}
    \item $\SO_{Q_1'}$ is not isomorphic to $\SO_{Q_1}$ as an algebraic group over $\Q$ .
\end{itemize}
Equivalently, $Q_1'$ is not isomorphic to $Q_1$ over $\Q$ up to a scalar.

For $i=1,2$, let $ M_{Q_i} \in \Mat_2(\Q)$ or $\Mat_4(\Q)$  be symmetric matrices representing $Q_i$: $Q_i(v,w)= v^{\Tr}  M_{Q_i} w$ if $v,w \in \Q^2$ or $\Q^4$  are written as column vectors.
Then $\bmU_{Q_1,Q_2}$ is naturally embedded into $\Mat_{4,2}$ as 
\begin{equation*}
    \left\{
    M\in \Mat_{4,2}(\Q) \;\middle\vert\;
     M^{\Tr} M_{Q_2} M = M_{Q_1}
     \right\}.
\end{equation*}

Let $\calU_{Q_1,Q_2}$ be the integral model by taking closure of $\bmU_{Q_1,Q_2}$ in $\calMat_{4,2}$ which is naturally an $4\times 2=8$-dimensional affine space over $\Z$. Let $\omega$ be an invariant gauge form on $\bmU_{Q_1,Q_2}$ over $\Q$.
For each prime $p$,
 let $\normm{\omega}_{p}$ be the Haar measure on $\bmU_{Q_1,Q_2}(\Q_p)$ induced from $\omega$.
For a matrix $M$ with real coefficients, let $\norm{M}$ be the Euclidean norm of $M$.

\begin{thm}\label{theorem_main_example_1}
    Write $\calU:=\calU_{Q_1,Q_2}$ for simplicity. Then as $R$ tends to infinity,
    \begin{equation*}
        \# \left\{
             M\in \calU(\Z) ,\; \norm{M}\leq R
        \right\} 
        \sim 
        \prod_{p\in \Val_f} \omega_p
        \left(
        \calU(\Z_p)  \right)
        \cdot 
        \omega_{\infty}
        \left(
         \left\{  M\in \calU(\R) ,\; \norm{M}\leq R  \right\}
        \right).
    \end{equation*}
    Namely, $\bmU_{Q_1,Q_2}$ with this embedding is strongly Hardy-Littlewood.
    Also, there exists some constant $c>0$ such that 
    \begin{equation*}
         \# \left\{
             M\in \calU(\Z) ,\; \norm{M}\leq R
        \right\}  \sim c R^2 \log R.
    \end{equation*}
\end{thm}

By \cite[Theorem 4.3, Lemma 4.1]{Wei_Xu_2016} (One should use replace $\SO_{Q_2}$ by its simply connected cover when applying this result and verify that the stabilizer subgroup remains connected. Details are omitted.), for the first part, it suffices to verify certain orbital counting statement which then is implied by the equidistribution towards the full Haar measure of some sequence of measures. 
At this point, let us note that $\bmU_{Q_1,Q_2}$ is homogeneous under the action of $\bmG:=\SO_{Q_2}\times \SO_{Q_1}$.
The stabilizer of $\phi_0 \in \bmU_{Q_1,Q_2}(\Q)$ is
of the form $\bmH_{\phi_0}= \bmH_1 \times \bmH_2$ where $\bmH_1 $ is the subgroup of $\SO_{Q_2}$  preserving the $Q_2$-orthogonal complement of the image of $\phi_0$ and $\bmH_2$ is suitable diagonal embedding of $\SO_{Q_1}$.
Under the assumptions imposed above, the only connected intermediate $\Q$-group is the maximal torus $\bmT_{\phi_0}$ containing $\bmH_{\phi_0}$.
By equidistribution theorems, the Haar measure supported on $\bmH_{\phi_0}(\R)\Gamma/\Gamma$, where $\Gamma$ is an(y) arithmetic lattice of $\bmG $, when pushed by $(g_n)$, equidistributes towards the full Haar measure unless $(g_n)$ is bounded modulo $\bmT_{\phi_0}$ which one  can show is generically not true.
Our main purpose here is to explain how this is done using the approach taken in the current paper.

Let $\bmX^{\Mat}$ denote the Zariski closure of $\bmU_{Q_1,Q_2}$ inside $\bmP(\Mat_{4,2}\oplus \Q) \cong \bmP^8_{\Q}$. 
We will describe below an explicit resolution of singularity of this space.
To begin with, it is covered by another ``incidence compactification'' which we now describe.

For a matrix $M$, let $\la M \ra$ be the linear subspace spanned by the column vectors of $M$. Then $\bmX^{\INC}$ is the Zariski closure of 
\begin{equation*}
    \left\{
      ([M:1], \la M \ra ) ,\;
      M \in \bmU_{Q_1,Q_2}
    \right\}
    \text{ in }  \bmX^{\Mat} \times \Gr(2,4).
\end{equation*}

Let $\bmD^{\INC}$ be the complement of $\bmU_{Q_1,Q_2}$ in $\bmX^{\INC}$.
The pair $(\bmX^{\INC}, \bmD^{\INC})$ is still not smooth.
Luckily, it is equipped with a natural stratification, components of which one can try blowing up.

For simplicity, we often work with complex-valued points and a typical point in $\bmX^{\INC}$ is written as $([M:\lambda],\la N \ra)$ where $N$ is some $4$-by-$2$ matrix.

At the level of complex points, one sees that 
the functions $\rk(\la N \ra ,Q_2)$ (rank of a quadratic form restricted to a subspace), $\rk(M)$ (rank of a matrix) and $\rk( \la M^{\Tr} \ra ,Q_1)$ are $\bmG(\C)$-invariant. Note that $f\leq c$, for $f$ being one of these functions and $c$ an integer, is indeed (locally) defined by zeros of certain polynomials and defines a closed $\bmG$-invariant subvariety over $\Q$.

\begin{defi}
For a tuple of integers $(a,b,c)$, 
    \begin{itemize}
        \item let $\bmS_{a,b,c}$ be the subvariety of $\bmX^{\INC}$  whose complex-valued points consist of $([M:\lambda],\la N \ra)$  such that 
        \begin{equation*}
            \rk(\la N \ra ,Q_2) = a, \;
            \rk(M)=b, \;
            \rk( \la M^{\Tr} \ra ,Q_1)=c;
        \end{equation*}
        \item let $\bmD_{a,b,c}$ be the Zariski closure of $\bmS_{a,b,c}$, whose complex-valued points consist of  $([M:\lambda],\la N \ra)$  such that 
        \begin{equation*}
            \rk(\la N \ra ,Q_2) \leq a, \;
            \rk(M) \leq b, \;
            \rk( \la M^{\Tr} \ra ,Q_1) \leq c.
        \end{equation*}
    \end{itemize}
\end{defi}

When $c=0$, exactly one of the column vectors of $M$ must be zero, consequently each nonempty $\bmS_{a,b,0}$ splits into two disjoint $\bmS_{a,b,0^{+}}$ and $\bmS_{a,b,0^{-}}$. Similarly, $\bmD_{a,b,0}$ is also a disjoint union of $\bmD_{a,b,0^{+}}$ and $\bmD_{a,b,0^{-}}$. 
We will be working inside affine open subvarieties where only one of them appears, it is often unnecessary to distinguish them. 

The values of these functions are subject to certain constraints. Some easy-to-observe ones are
\begin{equation*}
\begin{aligned}
      &\rk(\la N \ra ,Q_2) \in \{0,1,2 \} ;\;
       \rk(M)\in \{1,2\};\;
           \rk( \la M^{\Tr} \ra ,Q_1) \in \{0,1,2\}  ;\;\\
         &\rk( \la M^{\Tr} \ra ,Q_1)\leq \rk(M).
\end{aligned}
\end{equation*}
Working harder, one finds that 
\begin{lem}
    The space $\bmX^{\INC}$ decomposes into nonempty locally closed subvarieties:
    \begin{equation*}
        \bmX^{\INC}
        = \bmS_{2,2,2}\bigsqcup \bmS_{2,1,0}\bigsqcup \bmS_{1,1,1} \bigsqcup \bmS_{1,1,0} \bigsqcup \bmS_{0,2,2}\bigsqcup \bmS_{0,1,1}\bigsqcup \bmS_{0,1,0}.
    \end{equation*}
    Moreover, $\bmX^{\INC}$ has the following stratified structure (for $\star = -$ or $+$):
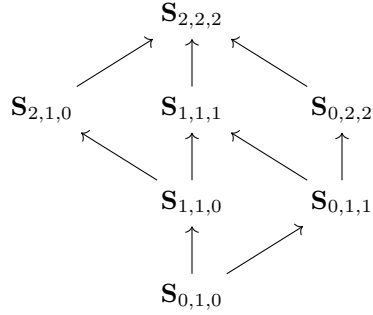
\[
    \begin{tikzcd}
                            & \bmS_{2,2,2} &  \\
     \bmS_{2,1,0^{\star}} \arrow[ru]  & \bmS_{1,1,1} \arrow[u] &   \bmS_{0,2,2} \arrow[lu]
     \\
                            & \bmS_{1,1,0^{\star}} \arrow[u] \arrow[lu] &  \bmS_{0,1,1} \arrow[lu] \arrow[u] \\
                            & \bmS_{0,1,0^{\star}} \arrow[u] \arrow[ru] &  
    \end{tikzcd}
\]
where $A \to B $ means ``$A$ is contained in the closure of $B$''.
\end{lem}

Now, a smooth pair can be obtained from the following three-step blowup process:
\begin{equation*}
    \begin{tikzcd}
        \bmX^3 := \Bl_{ (\bmE_2 \cap \bmD_{022}^{++})
        \sqcup \bmD_{111}^{++} \sqcup (
        \bmE_2 \cap \bmD_{210}^{++}
        )
        } \bmX^2
        \arrow[d,"\pi_3"]
        \\
        \bmX^2 := \Bl_{ 
        \bmD_{011}^+ \sqcup \bmD_{110}^+
        } \bmX^1
        \arrow[d,"\pi_2"]
        \\
        \bmX^1 := \Bl_{
        \bmD_{010}
        } \bmX^{\INC}
        \arrow[d,"\pi_1"]
        \\
        \bmX^{\INC} 
    \end{tikzcd}
\end{equation*}
where $\bmE_i$ denotes the exceptional divisor of $\pi_i$ and $\bmD^+$ the birational transform of a closed subvariety $\bmD$, which is well defined as long as $\bmD\setminus\bmC$ ($\bmC:=\text{center of the blowup}$) is dense in $\bmD$.

Let $\bmD^i$ be the complement of $\bmU_{Q_1,Q_2}$ in $\bmX^i$.
Let $\bmE_2^1$ (resp. $\bmE_2^2$) be the inverse image of $\bmD_{011}^+$ (resp. $\bmD_{110}^+$) under $\pi_2$. 
Let $\bmE^0_3$ (resp. $\bmE^1_3$, $\bmE^2_3$) be the inverse image of $\bmE_2^1 \cap \bmD_{022}^{++}$ (resp. $\bmD_{111}^{++}$, $\bmE^2_2 \cap  \bmD_{210}^{++}$).
If one takes the splitting of $\bmD_{a,b,0^{\pm}}$ into consideration, then
\[
\bmE_1^{++}= \bmE_{1^{-}}^{++}\sqcup\bmE_{1^+}^{++},\;
\bmD_{210}^{++} = \bmD_{210^{-}}^{++}\sqcup \bmD_{210^{+}}^{++} ,\;
\bmE_3^2 = \bmE_{3^-}^2 \sqcup \bmE_{3^+}^2 ,\;
(\bmE_{2}^2)^+=  (\bmE_{2^{-}}^2)^+\sqcup (\bmE_{2^{+}}^2)^+.
\]

Let $(\calL_{0}, s_{0})$ be the metrized line bundle and the global section over $\bmX^{\Mat}$ that underlies the definition of the matrix norm $\norm{\cdot }$ and let $(\wtcalL_{0}, \wts_{0})$ be their pull-backs over $\bmX^3$.

\begin{lem}\label{lemma_example_1_divisors_on_smooth_pair}
    $(\bmX^3,\bmD^3)$ is a smooth pair. Moreover, the irreducible components of $\bmD^3$ are 
    \begin{equation*}
        ( \bmE_{1^{\pm}}^{++}, \bmD_{210^{\pm}}^{++}, \bm\bmE_{3^{\pm}}^2, (\bmE_{2^{\pm}}^2)^+
        , \bmE_3^1, (\bmE_2^1)^+, \bmE_3^0, \bmD_{022}^+ ),
    \end{equation*}
    which in this order are labeled as $(\bmD^3_{1^{\pm}},...,\bmD^3_{4^{\pm}},\bmD^3_5,...,\bmD^3_8)$ (total number is $12$). Let $\bmD^3_i:= \bmD^3_{i^-}\sqcup \bmD^3_{i^+}$.
    The intersection pattern is described as follows:
   For $\star = -$ or $+$,  $\bmD_{1^{\star}}^3$ intersects with $(\bmD^3_{2^{\star}},...,\bmD^3_{4^{\star}},\bmD^3_5,...,\bmD^3_8)$ transversally. 
   Those with a $-$ in the subscript do not intersect with those with a $+$ in the subscript.
   For $i\neq 1$, $\bmD^3_i$ only intersects with $\bmD_1^3$ and $\bmD^3_{i \pm 1}$ (so $\bmD^3_{8}$ only intersects with $\bmD^3_1$ and $\bmD^3_7$).
    Moreover,
    \begin{equation*}
        -\divisor( \omega_{\bmU} )
        = 
        2\bmD^3_{1}+ 1\bmD^3_{2}+ 
        3\bmD^3_{3}+ 5\bmD^3_{4} + 3\bmD^3_{5} 
        +7\bmD^3_{6} +5\bmD^3_{7}+ 3\bmD^3_{8}
        =: \sum d_{\alpha} \bmD^3_{\alpha}
    \end{equation*}
    and 
    \begin{equation*}
        \divisor( \wts_0 )=
         1 \bmD^3_{1}+ 1\bmD^3_{2}+ 
        2 \bmD^3_{3} + 3 \bmD^3_{4} + 1 \bmD^3_{5} 
        +3 \bmD^3_{6} + 2 \bmD^3_{7} + 1 \bmD^3_{8}
        =: \sum \lambda_{\alpha} \bmD^3_{\alpha}.
    \end{equation*}
    Consequently,
    \begin{equation*}
        \left( 
        \frac{d_{\alpha}-1}{\lambda_{\alpha}}
        \right)_{\alpha=1}^{8} =
        \left(
            1, 0 ,1, 4/3, 2,2,2,2
        \right).
    \end{equation*}
\end{lem}

The proof of Lemma \ref{lemma_example_1_divisors_on_smooth_pair}, which can be found in the appendix,  relies on explicit construction of blowups and explicit calculation of the divisor of the invariant gauge form. 

Note that $\rmG:= \SO_{Q_2}(\R)$ is connected in analytic topology. Let $\bmH_x$ denote the stabilizer subgroup of $\bmG$ of some point $x\in \bmU(\Q)$.
Thanks to our $\R$-split assumption,
\begin{lem}
    For every $\alpha \in \{1,...,8\}$, $\bmD^3_{\alpha}(\R)\neq \emptyset$ and
    $\bmU(\R)$ consists of a single $\rmG$-orbit. 
    Consequently, for every $\alpha =1,...,8$ and $x\in \bmU(\Q)$, the analytic closure of $\rmG.x $ intersects with $\bmD_{\alpha}(\R)$.
\end{lem}

\begin{proof}
    That $\bmD_{\alpha}^3(\R)\neq \emptyset$ follows from the construction (see appendix for details). It is easy to see that the analytic closure of $\bmU(\R)$ intersects with every $\bmD_{\alpha}^3(\R)$. However, $\rmG.x $ is equal to $\bmU(\R)$ thanks to the following ``exact'' sequence
    \begin{equation*}
        \rmH_x \to \rmG \to \bmU(\R) \to H^1(\Gal(\C/\R),\bmH_x)
    \end{equation*}
    and the triviality of $H^1(\Gal(\C/\R),\bmH_x)$.
\end{proof}

Let $\Gamma$ be some arithmetic lattice.

\begin{lem}\label{lemma_example_1_measure_equidistribution}
    For every $x\in \bmU(\Q)$ and every sequence $(g_n)$ of $\rmG$, the following two are equivalent:
    \begin{itemize}
        \item[1.] every limit point of $(g_n.x)$  is contained in $\bigcup_{\alpha \neq 2} \bmD^3_{\alpha}(
        \R
        )$;
        \item[2.] under the weak-$*$ topology,  $\lim_{n\to \infty} (g_n)_*\rmm^{\bmone}_{[\rmH_x]} = \rmm^{\bmone}_{[\rmG]}$.
    \end{itemize}
\end{lem}

\begin{proof}
    By our assumption, the only intermediate connected closed $\Q$-subgroup between $\bmH_x$ and $\bmG$ is the centralizer of $\bmH_x$, a maximal torus in $\bmG$. And a sequence $(g_n)$ is bounded modulo this maximal torus iff the limit of $(g_n.x)$ is contained in $\bmD_2^3(\R )^{\circ} \cup \bmU(\R)$. The conclusion thus follows from Theorem \ref{theorem_equidistribution_EMS}.
\end{proof}

\begin{proof}[Proof of Theorem \ref{theorem_main_example_1}]
    For the first part,
    by \cite[Theorem 4.3]{Wei_Xu_2016} and the fact that $\Pic(\bmH)$ and hence $\mathrm{Br}(\bmU)/\mathrm{Br}(\Q)$ is trivial, it is sufficient to show that for every $x\in \bmU(\Q)$, and every arithmetic lattice $\Gamma$ of $\bmG$,
    \begin{equation*}
        \frac{
        \# \left\{
           M \in \Gamma.x \;\middle\vert\;
           \norm{M} \leq R
        \right\}
        }{
         \omega_{\infty} 
         \left\{
           M \in \rmG.x \;\middle\vert\;
           \norm{M} \leq R
        \right\}
        }
        \sim \frac{
        \normm{\rmm_{[\rmH_x]}}
        }{
        \normm{\rmm_{[\rmG]}}
        }
    \end{equation*}
    where     $\omega_{\infty}$ is identified with the Haar measure $\rmm_{\rmG/\rmH_x}$ which is assumed to be compatible with the Haar measures $\rmm_{\rmG}$ and $\rmm_{\rmH_x}$.

    By Lemma \ref{lemma_example_1_divisors_on_smooth_pair} and Theorem \ref{theorem_equidistribution_Chamber-Loir_Tschinkel}, the limit
    \begin{equation*}
    \nu:= \lim_{R \to \infty} \frac{\rmm_{\rmG/\rmH_x} \cdot \bmone_{B_{R,x}}
    }
    {\rmm_{\rmG/\rmH_x}(B_{R,x})}
\end{equation*}
exists and is supported on 
$\cup_{i=5}^{7} \bmD^3_i (\R) \cap\bmD^3_{i+1} (\R) $ where $B_{R,x}:= \{ M \in \rmG.x ,\; \norm{M}\leq R\}$.
It only remains to invoke Lemma \ref{lemma_condition_good_height}.

The second part follows from Theorem \ref{theorem_equidistribution_Chamber-Loir_Tschinkel} and the calculation of $(\frac{d_{\alpha}-1 }{ \lambda_{\alpha}})$ as in Lemma \ref{lemma_example_1_divisors_on_smooth_pair}.
\end{proof}

\section{Example II, $(n,1)$-splitting of $\Z^{n+1}$}\label{section_example_2}

Let 
\begin{equation*}
        \Lambda(n,1):=\left\{ (v,M)\;\middle\vert\;
        v\in \Z^n,\, M\in \Prim^{n}(\Z^{n+1}),\, \Z^{n+1}=\Z.v\oplus M
        \right\}
    \end{equation*}
    where $\Prim^{n}(\Z^{n+1})$ denotes the collection of the rank-$n$ subgroups $M$ of $\Z^{n+1}$ such that $\Z^{n+1}/M$ is torsion-free.

\begin{thm}\label{theorem_main_example_2}
    For two positive integers $\lambda_1,\lambda_2$, we have 
    \begin{equation*}
        \#\{(v,M)\in\Lambda(n,1), \;
        \norm{v} ^{\lambda_1} \norm{M}^{\lambda_2} \leq R \} \sim 
        \begin{cases}
            c_{\lambda_1,\lambda_2} \cdot R^{\frac{n}{\lambda_1}} \log (R) ,\quad
            & \lambda_1=\lambda_2
            \\
            c_{\lambda_1,\lambda_2} \cdot R^{\frac{n}{\min\{ \lambda_1,\lambda_2 \} } },\quad
            &  \lambda_1 \neq \lambda_2
        \end{cases}
    \end{equation*}
    for some $ c_{\lambda_1,\lambda_2} > 0$.
\end{thm}

Here $\norm{\star}$ denotes the covolume of $\star$ in the $\R$-subspace spanned by $\star$ with respect to the standard Euclidean norm.

\begin{rmk}
    As one sees in the proof, in the case $ {\lambda_1 \neq \lambda_2}$, ``focusing'' happens. 
\end{rmk}

Let $\bmG:=\SL_{n+1}$, which naturally acts on $\Q^{n+1}$ and its dual $(\Q^{n+1})^{\vee}$. Let  $(\bme_1,...,\bme_{n+1})$ be the standard basis of $\Q^{n+1}$ and $(\bme^{\vee}_1,...,\bme^{\vee}_{n+1})$ be its dual basis so that  $(\Q^{n+1})^{\vee}$ is identified with $\Q^{n+1}$ under this basis. Note that $\Q^{n+1}$ naturally embeds in $\bmP^{n+1}_{\Q}$ by $\bmx \mapsto [\bmx:1]$.
Let $y_0:=([\bme_1:1],[\bme_1^{\vee}:1])$ and $\bmH$ be its stabilizer in $\bmG$.

Let $\bmX$ be the Zariski closure of $\bmU:=\bmG.y_0$ in $\bmP^{n+1}_{\Q} \times \bmP^{n+1}_{\Q} $ and $\bmD$ be the complement of $\bmU$ in $\bmX$.
Write $([\bmx:s],[\bmalpha :t])$  for a point in $\bmX$. Define the following subvarieties of $\bmX$:
\begin{itemize}
    \item[1.] $\bmD_1$ is defined by $s=0$ and $\bmD_2$ by $t=0$;
    \item[2.] $\bmS_1$ (resp. $\bmS_2$) is defined by $s=0$ but $t\neq 0$ (resp. $t=0$ but $s\neq 0$);
    \item[3.] $\bmS_1^*$ (resp. $\bmS_1^0$) is the subvariety of $\bmS_1$ such that $\bmalpha \neq \bmzero$ (resp. $\bmalpha = \bmzero$).
    $\bmS_2^*$ and $\bmS_2^0$  are similarly defined.
\end{itemize}

\begin{lem}\label{lemma_example_2_smooth_pair}
    The pair $(\bmX,\bmD)$ is smooth. Moreover, $\bmD_1$ and $\bmD_2$ are the irreducible components of $\bmD$ whose intersection is nonempty and irreducible.
\end{lem}

As for real points, we have

\begin{lem}
    With respect to the analytic topology, $\bmU(\R)$ is connected and dense in $\bmX(\R)$.
\end{lem}

For convenience we fix some representatives in each $\rmG$-orbit as follows
\begin{itemize}
    \item $y_0:=([\bme_1:1],[\bme_1^{\vee}:1])$ as above is a representative of $\bmU(\R)$;
    \item $y^*_1:= ([\bme_2:0],[\bme_1^{\vee}:1]) \in \bmS_1^{*}(\R)$ and 
    $y^*_2:= ([\bme_1:1],[\bme_2^{\vee}:0]) \in \bmS_2^*(\R)$;
    \item $y^0_1:= ([\bme_2:0],[\bmzero:1]) \in \bmS_1^{0}(\R)$ and 
    $y^0_2:= ([\bmzero:1],[\bme_2^{\vee}:0]) \in \bmS_2^{0}(\R)$;
    \item $y_{12}:= ([\bme_1:0], [\bme_2^{\vee}:0] ) \in \bmD_1(\R)\cap \bmD_2(\R) $.
\end{itemize}

\begin{lem}
    Under the $\rmG$-action, we have
    \begin{equation*}
    \begin{aligned}
        \bmX(\R) 
        &= \bmU(\R) \bigsqcup \bmS_1^*(\R) \bigsqcup \bmS_1^0(\R) \bigsqcup \bmS_2^*(\R)  \bigsqcup \bmS_2^0(\R)  \bigsqcup
        \left( \bmD_1(\R) \cap \bmD_2(\R)  \right)
        \\
        &=
        \rmG.y_0 \bigsqcup
        \rmG.y_1^*  \bigsqcup
        \rmG.y_1^0  \bigsqcup
        \rmG.y_2^*  \bigsqcup
        \rmG.y_2^0  \bigsqcup
        \rmG.y_{12}  .
    \end{aligned}
    \end{equation*}
\end{lem}

Let $\Gamma:=\SL_{n+1}(\Z)$, $\bmM_1$ be the stabilizer of $\bme_1$ in $\bmG$ and $\bmM_2$ be the stabilizer of $\bme_1^{\vee}$.
We define a map $\Phi$ from $\bmX(\R)$ to $\Prob(\rmG/\Gamma)\cup \{\bmzero\}$ by
\begin{equation*}
\Phi(x) =
\begin{cases}
     g_* \rmm_{[\rmH]}^{\bmone}  \quad &  x= g.y_0 \in \bmU(\R)  \\
     g_* \rmm_{[\rmM_2]}^{\bmone}  \quad & x = g.y^*_1 \in \bmS_1^*(\R)\\
     g_* \rmm_{[\rmM_1]}^{\bmone} \quad & x = g.y^*_2 \in \bmS_2^*(\R) \\
     \bmzero \quad & x \in \bmS_1^{0}(\R) \cup \bmS_2^{0}(\R)\\
     \rmm^{\bmone}_{[\rmG]}   \quad & x\in \bmD_1(\R) \cap \bmD_2(\R).
\end{cases}
\end{equation*}

\begin{lem}\label{lemma_example_2_map_to_measures}
    The map $\Phi$ is well-defined (independent of the choice of $g$) and continuous.
\end{lem}

Let $\omega_{\bmU}$ be the $\bmG$-invariant gauge form on $\bmU$.

\begin{lem}\label{lemma_example_2_divisor_anti_canonical}
  The anti-canonical divisor is
   $ -\divisor (\omega_{\bmU}) = (n+1) \bmD_1 + (n+1) \bmD_2$.
\end{lem}

See Section \ref{section_example_2_divisor_invariant_gauge} for the proof.

Finally, note the metric line bundle underlying the height function.
\begin{lem}\label{lemma_example_2_height_divisor}
    There exist smooth metrics $\norm{\cdot}$ on $\calO_{\bmX}(\bmD_i)$ for $i=1,2$ such that for $(v,M)\in \Lambda(n,1)$,
    \begin{equation*}
        \norm{\bmone_{\bmD_1} }_{(v,M)} =
        \left( \sum v_i^2 \right)^{-1/2},\quad
         \norm{\bmone_{\bmD_2} }_{(v,M)} =
        \left( \sum \alpha_i^2 \right)^{-1/2}
    \end{equation*}
    where $v=\sum v_i \bme_i$ and $M = \sum \alpha_i \bme^{\vee}_i$.
\end{lem}

See Section
\ref{section_example_2_height_and_metrics} for details on this construction.

\begin{proof}[Proof of Theorem \ref{theorem_main_example_2}]
From Lemma \ref{lemma_example_2_height_divisor}, there exists a smooth metric $\norm{\cdot}$ on $\calO_{\bmX}(\lambda_1\bmD_1 + \lambda_2 \bmD_2)$ such that 
\[
\norm{\bmone_{\bmD_1}^{\otimes \lambda_1}  \otimes \bmone_{\bmD_2}^{\otimes \lambda_2} }_{(v,M)} =   \norm{v} ^{-\lambda_1} \norm{M}^{-\lambda_2}  .
\]
By Lemma \ref{lemma_example_2_map_to_measures} and \ref{lemma_ok_heights_possible_diverges}, there exists some $c>0$ such that 
    \begin{equation*}
         \#\{(v,M)\in \Lambda(n,1), \;
        \norm{v} ^{\lambda_1} \norm{M}^{\lambda_2} \leq R \} \sim
         c \cdot \rmm_{\rmG/\rmH}(B_{R,y_0}).
    \end{equation*}
    Thus the conclusion follows from Theorem \ref{theorem_equidistribution_Chamber-Loir_Tschinkel} and Lemma
    \ref{lemma_example_2_divisor_anti_canonical}.
\end{proof}

\section{Example III, the space of triangles}\label{section_example_3}

In this section we consider the space of three linearly independent lines on $\Q^3$. In fact, there exists a $\Q$-variety $\bmM_3$ such that 
\begin{equation*}
    \bmM_3(\Q) = \left\{
     (\bml_1 ,\bml_2, \bml_3) \in (\bmP^2(\Q))^3 
     \;\middle\vert\;
     (\bml_1 ,\bml_2, \bml_3) \text{ are linearly independent}
    \right\}.
\end{equation*}

One can show that $\bmM_3$ is dense in $\bmP_{\Q}^2 \times \bmP_{\Q}^2 \times \bmP_{\Q}^2$. Let $\bmD^0$ be its complement. $\bmP_{\Q}^2$ has a natural integral model over $\Z$, denoted as $\bmP^{2}_{\Z}$.
 Let $\calD^0$ be the closure of $\bmD^0$ in $\bmP^{2}_{\Z}\times \bmP^{2}_{\Z} \times \bmP^{2}_{\Z}$. Let $\calM_3$ be the complement of $\calD^0$ in  $\bmP^{2}_{\Z}\times \bmP^{2}_{\Z} \times \bmP^{2}_{\Z}$.
Then one can verify that, viewed as a subset of $\bmM_3(\Q)$ as above, 
\begin{equation*}
\begin{aligned}
        \calM_3(\Z)
     &=   \left\{
     (\bml_1 ,\bml_2, \bml_3) \in \bmM_3(\Q)
     \;\middle\vert\;
     \exists \, \bmv_i \in \Z^3 \text{ with }\bml_i = \Q .\bmv_i,\; \det(\bmv_1, \bmv_2, \bmv_3) = \pm 1
    \right\}
    \\
      &=\left\{
     (\bml_1 ,\bml_2, \bml_3) \in \bmM_3(\Q)
     \;\middle\vert\;
     \Z^3 = (\bml_1 \cap \Z^3 ) \oplus (\bml_2 \cap \Z^3 ) \oplus 
     (\bml_3 \cap \Z^3 ) 
    \right\}.
\end{aligned}
\end{equation*}
To measure the complexity of an element in $\calM_3(\Z)$, we introduce the following functions:

For a discrete $\Z$-submodule $\Lambda$ of $\R^3$, let 
\begin{equation*}
    \norm{\Lambda}:= \Vol( \Lambda \otimes_{\Z} \R /\Lambda)
\end{equation*}
where $\Vol$ is induced from the standard Euclidean metric on $\R^3$. By convention, $\norm{ \{ \bmzero \} }:=1$.
Let $(\Lambda_1,\Lambda_2,\Lambda_3)$ be three rank one free $\Z$-submodules of $\R^3$ that form an $\R$-linear basis of $\R^3$.
For $I \subset \{1,2,3\}$, let $\Lambda_I:= \oplus_{i\in I} \Lambda_i$. By default, $\Lambda_{\emptyset}:= \{ \bmzero \}$.
For $I,J \subset \{1,2,3\}$, 
\begin{equation*}
    d_{IJ}( (\Lambda_i)_{i=1}^3 ) := \frac{
    \norm{ \Lambda_I } \norm{ \Lambda_J }
    }{
    \norm{ \Lambda_{I\cap J} } \norm{ \Lambda_{I\cup J} }
    }.
\end{equation*}
For $\bmx = (\bml_1,\bml_2,\bml_3) \in \bmM_3(\R) $, take some nonzero $\bmv_i \in \bml_i$, and let $\Lambda_i:= \Z. \bmv_i$.
For $ l=1,2$, define
\begin{equation*}
    \begin{aligned}
        \Ht_{l} (\bmx) := \prod_{I,J\subset \{1,2,3\},\,
        |I|=|J|=l
        }  d_{IJ}( (\Lambda_i)_{i=1}^3 ).
    \end{aligned}
\end{equation*}
One can check that this definition is independent of the choice of $\bmv_i$'s.

Now we define weights for points in $\calM_3(\Z)$.
Take some $\bmx = (\bml_1,\bml_2,\bml_3) \in \calM_3(\Z)$ and let $\Lambda_i:= \bml_i \cap \Z^3$.
Consider, for $\eta>0$,
\begin{equation*}
    \Omega_{\bmx,\eta}:=  \left\{
     \bmt=(t_1,t_2,t_3) \in \R^3 \;\middle\vert\;
     t_1+t_2+t_3 =0,\;
     \sum_{i\in I} t_i \geq  -\ln \norm{\Lambda_I} + \ln{\eta},\,\forall\,
     I \subset \{1,2,3\}
\right\}
\end{equation*}
where summation over an empty set is set to be $0$.
The weight function is just 
\begin{equation*}
    \bmw_{\bmx}= \min \{ \Vol ( \Omega_{\bmx,\eta} ) ^{-1} , 1  \}.
\end{equation*}

\begin{thm}\label{theorem_main_example_3}
    Let $\kappa_1,\kappa_2 >0$ and 
    $\Ht(\bmx):= \Ht_{1}(\bmx)^{\kappa_1} \Ht_{2}(\bmx)^{ \kappa_2 }$ for $\bmx \in \bmM_3(\R)$.
    Then for $\eta>0$ small enough and some constant $c_{\kappa_1,\kappa_2}>0$,
    \begin{equation*}
        \sum_{
         \left\{
          \bmx \in \calM_3(\Z) \;\middle\vert\;
          \Ht(x) \leq R
        \right\} }  \bmw_x
        \sim 
        \begin{cases}
            c_{\kappa_1,\kappa_2} \cdot R^{\frac{8}{3} \max\{\kappa_1^{-1},\kappa_2^{-1} \} } 
            \quad & \kappa_1 \neq \kappa_2
            \\
             c_{\kappa_1,\kappa_2} \cdot R^{\frac{8}{3} \kappa_1^{-1} } \cdot \log(R)
            \quad & \kappa_1 =\kappa_2.
        \end{cases}
    \end{equation*}
\end{thm}

\begin{rmk}
    Without weights, we expect an additional $\log(R)^2$ factor in the asymptotic. We hope to discuss this  in a future work. This counting problem has been studied in \cite{ShaZhe18} with respect to a different height function.
\end{rmk}

It will be clear that the counting problem is naturally related to the embedding of $\bmM_3$ into the following incidence variety:
\begin{equation*}
    \bmX^{\INC} = \left\{
      (\bml_{1}, \bml_{2}, \bml_{3}, \bml_{12}, \bml_{13},\bml_{23})
      \in  (\bmP^2)^3 \times (\bmGr_{2,3})^3 
      \;\middle\vert\;
      \bml_{I}\subset \bml_{J} ,\;\forall\, I \subset J
    \right\}.
\end{equation*}
The morphism $\bmM_3 \to \bmX^{\INC}$ is given by mapping $(\bml_{1}, \bml_{2}, \bml_{3})$ to $(\bml_{1}, \bml_{2}, \bml_{3}, \bml_{12}, \bml_{13},\bml_{23})$ where $\bml_{ij}$ is the unique plane spanned by $\bml_i$ and $\bml_j$. One can check that the morphism is an open embedding with dense image.
As points in $\bmX^{\INC}$ are tuple of points and lines in $\bmP^2$ satisfying incidence conditions, they describe ``triangles'' in $\bmP^2$. The space is thus sometimes referred to as ``space of triangles''. It is also related to certain counting problem studied by Schubert. See \cite{Semple_1954_triangle, Roberts_Speiser_1984_triangles, Roberts_1988_triangle_varieties} for further discussions.

The space $\bmX^{\INC}$ is naturally stratified.
Let
\begin{equation*}
\begin{aligned}
     &\bmD_{123}^{1} :=\left\{
       \bml_1  =  \bml_2  = \bml_3
    \right\},\; 
    \bmD_{123}^{2} :=\left\{
      \bml_{13}=\bml_{13}=  \bml_{23}
    \right\},\\
    &
    \bmD_{12,3} := \left\{
      \bml_1 = \bml_2,\, \bml_{13} =\bml_{23}
    \right\},\;
    \bmD_{13,2} := \left\{
      \bml_1 = \bml_3,\, \bml_{12} =\bml_{23}
    \right\},\;\\
    &
    \bmD_{23,1} := \left\{
      \bml_2 = \bml_3,\, \bml_{12} =\bml_{13}
    \right\}.
\end{aligned}
\end{equation*}
Also let the corresponding $\bmS^{\star}_{\bullet}$ denote points contained in $\bmD^{\star}_{\bullet}$ but not the other $\bmD$'s.

Let $\bmX^1:= \Bl_{\bmD_{123}^1} (\bmX^{\INC}) $ and $\bmD^1$ be the complement of $\bmM_3$ in $\bmX^1$.

\begin{lem}\label{lemma_example_3_smooth_pair}
   The pair $(\bmX^1,\bmD^1)$  is smooth. Irreducible components of $\bmD^1$ are 
   \begin{equation*}
       \left( \bmE, \bmD_{12,3}^+, \bmD_{13,2}^+, \bmD_{23 ,1}^+, (\bmD_{123}^2)^+  \right) 
   \end{equation*}
    where $\bmE$ denotes the exceptional divisor.
    Moreover,  the intersection pattern is given as follows:
    $ \bmD_{ij,k}^+$'s do not intersect with each other and for each $i,j,k$, the triple $\left( \bmD_{ij,k}^+, \bmE, (\bmD_{123}^2)^+ \right)$ intersects transversally.
\end{lem}

Details are provided in Section \ref{section_example_3_blow_up} (see also \cite{Roberts_Speiser_1984_triangles}).
We henceforth label $( \bmE, \bmD_{12,3}^+, \bmD_{13,2}^+, \bmD_{23 ,1}^+, (\bmD_{123}^2)^+ )$ as $(\bmD_1,...,\bmD_5)$.
To continue, note that $\bmG:=\SL_3$ acts on $\bmM_3$ and also $\bmX^{\INC}$ by linear transformations, making $\bmM_3$ into a homogeneous space. If  $o\in \bmM_3(\Z)$ is the base point $(\Z \bme_1, \Z \bme_2, \Z \bme_3)$, then the stabilizer $\bmH$ of $o$ in $\bmG$ is equal to the full diagonal torus in $\SL_3$.

Let $\omega_{\bmM_3}$ be the invariant gauge form, then
\begin{lem}\label{lemma_example_3_divisor_gauge_form}
    In $\bmX^1$, the anti-canonical divisor is
    $-\divisor(\omega_{\bmM_3})
    = 9 \bmD_1 + 6 \bmD_2 + 6 \bmD_3 + 6 \bmD_4 + 9 \bmD_5
    $.
\end{lem}

See Section \ref{section_example_3_pole_invariant_gauge} for the proof.

Let $\psi$ be a non-negative function on $\rmG/\Gamma$ (here $\rmG:=\SL_3(\R)$ and $\Gamma:= \SL_3(\Z)$) whose support is large enough. We construct a continuous map from $\bmX^{\INC}(\R)$ to $\Prob^{\psi}(\rmG/\Gamma)$.
Let 
\begin{itemize}
    \item $\bmx_1: = (
        \Z \bme_1, \Z \bme_1, \Z \bme_1, \Z \bme_1 \oplus \Z \bme_2,
        \Z \bme_1 \oplus \Z \bme_3, \Z \bme_1 \oplus \Z (\bme_2+\bme_3)
    )$  in $\bmS_{123}^1(\R)$;
    \item $\bmx_2: = (
        \Z \bme_1, \Z \bme_1, \Z \bme_3, 
        \Z \bme_1 \oplus \Z \bme_2,
        \Z \bme_1 \oplus \Z \bme_3, \Z \bme_1 \oplus \Z \bme_3
    )$ in $\bmS_{12,3}(\R)$;
    \item $\bmx_3: = (
        \Z \bme_1, \Z \bme_2, \Z \bme_1,
        \Z \bme_1 \oplus \Z \bme_2,
        \Z \bme_1 \oplus \Z \bme_3, \Z \bme_1 \oplus \Z \bme_2
    )$  in $\bmS_{13,2}(\R)$;
    \item $\bmx_4: = (
        \Z \bme_1, \Z \bme_2, \Z \bme_2,
        \Z \bme_1 \oplus \Z \bme_2,
        \Z \bme_1 \oplus \Z \bme_2, \Z \bme_2 \oplus \Z \bme_3
    )$  in $\bmS_{23,1}(\R)$;
    \item $\bmx_5: = (
        \Z \bme_1, \Z \bme_2, \Z (\bme_1+\bme_2), 
        \Z \bme_1 \oplus \Z \bme_2,
        \Z \bme_1 \oplus \Z \bme_2, 
        \Z \bme_1 \oplus \Z \bme_2
    )$ in $\bmS_{123}^{2}(\R)$.
\end{itemize}

\begin{lem}\label{lemma_example_3_strata_incidence_variety}
    We have 
    \begin{equation*}
        \bmX^{\INC}(\R) = \rmG.o \bigsqcup  \rmG.x_2 \bigsqcup
        \rmG.x_3  \bigsqcup \rmG.x_4 \bigsqcup (\bmD^1_{123}(\R) \cup \bmD^{2}_{123}(\R) ).
    \end{equation*}
\end{lem}

For a partition $ \{i,j\} \sqcup \{k\}$ of $\{1,2,3\}$,
let $\bmH_{ij}$ be the simultaneous stabilizer in $\bmG$ of the line spanned by $\bme_k$ and the plane spanned by $\bme_i, \bme_j$.
Define $\Phi_{\bmH}: \bmX^{\INC}(\R) \to  \Prob^{\psi}(\rmG/\Gamma)$ by
\begin{equation*}
    \Phi_{\bmH}(\bmx):=
    \begin{cases}
        \alpha_g^{\psi}(\rmm_{[\rmH]})  \quad & \text{if }\,  \bmx =g.o \in \bmM_3(\R)  \\
        \alpha_g^{\psi}(\rmm_{[\rmH_{12}]}) \quad & \text{if }\,  \bmx =g.x_2 \in \bmS_{12,3}(\R)  \\
        \alpha_g^{\psi}(\rmm_{[\rmH_{13}]})  \quad & \text{if }\,  \bmx =g.x_3 \in \bmS_{13,2}(\R) \\
        \alpha_g^{\psi}(\rmm_{[\rmH_{23}]})  \quad & \text{if } \, \bmx =g.x_4 \in \bmS_{23,1}(\R) \\
        \rmm^{\psi}_{[\rmG]}  \quad & \text{if }\,  \bmx \in \bmD^1_{123}(\R) \cup \bmD^{2}_{123}(\R) .
    \end{cases}
\end{equation*}

\begin{lem}\label{lemma_example_3_map_to_measure_compactification}
    The map  $\Phi_{\bmH}$ is well-defined and continuous.
\end{lem}

This is a corollary to discussions in Section \ref{section_example_3_measure_classification}.

For $\{i,j\} \subset \{1,2,3\}$, let 
\begin{equation*}
    \bmX^{\INC}_{ij}:=
    \left\{
     (\bml_i,\bml_j,\bml_{ij}) \in \bmP^2 \times \bmP^2    \times \bmGr_{2,3}
        \,\middle\vert\,
      \bml_i, \bml_j \subset \bml_{ij}
 \right\}
\end{equation*}
together with a natural morphism $\pi_{ij}:\bmX^{\INC} \to \bmX^{\INC}_{ij}$. Let $\bmM_{ij}$ denote the open subvariety where $\bml_i$ and $\bml_j$ are linearly independent and $\bmD_{ij}$ be its complement.

Similarly, for 
$\{\{i,j\}, \{i,k\} \}\subset \{\{1,2\},\{1,3\},\{2,3\}  \}$, we define $\bmX^{\INC}_{ij,ik}$, a subvariety of $ \bmP^2   \times \bmGr_{2,3}\times \bmGr_{2,3} $, and $\pi_{ij,ik}:\bmX^{\INC} \to \bmX^{\INC}_{ij,ik}$. Let $\bmM_{ij,ik}$ denote the open subvariety where $\bml_{ij}\neq \bml_{ik}$ and $\bmD_{ij,ik}$ its complement.

\begin{lem}\label{lemma_example_3_divisor_height}
    The pair $(\bmX^{\INC}_{ij}, \bmD_{ij} )$  is smooth.
    Moreover, there exists a smooth metric on $\calO_{\bmX^{\INC}_{ij}}(\bmD_{ij})$ such that if $\bmv_i$ and $\bmv_j$ are nonzero vectors on $\bml_i$ and $\bml_j$ respectively, then for $\bmx=(\bml_i,\bml_j,\bml_{ij})$ in $\bmM_{ij}(\R)$,
    \begin{equation*}
        \norm{ \bmone_{\bmD_{ij}} }_{\bmx} = \frac{
        \norm{ \bmv_i \wedge \bmv_j }
        }{ 
        \norm{\bmv_i} \norm{\bmv_j}
        }.
    \end{equation*}
    Similarly,  the pair $(\bmX^{\INC}_{ij,ik}, \bmD_{ij,ik} )$  is smooth. And $\calO_{\bmX^{\INC}_{ij,ik}}(\bmD_{ij,ik})$ can be equipped with a smooth metric such that 
    \begin{equation*}
        \norm{ \bmone_{ \bmD_{ij,ik} } }_{\bmx} = \frac{
        \norm{ \bmv_i} \norm{\bmv_i \wedge \bmv_j \wedge \bmv_k }
        }{ 
        \norm{\bmv_i \wedge \bmv_j} \norm{\bmv_i \wedge \bmv_k}
        }
    \end{equation*}
    where $\bml_{ij}$ is spanned by $\bmv_i,\bmv_j$ and $\bml_{ik}$ is spanned by $\bmv_i,\bmv_k$.
\end{lem}


Let $\bms_{ij}$ be the pull-back of the section $ \bmone_{\bmD_{ij}}$ to a section of the pull-back of $\calO_{\bmX^{\INC}_{ij}}(\bmD_{ij})$ to a line bundle over $\bmX^1$. Similarly define $\bms_{ij,ik}$.
\begin{lem}\label{lemma_example_3_divisor_pull_back}
    For $\{i,j\} \subset \{1,2,3\}$, $\divisor(\bms_{ij})= \bmE + 
    \bmD_{ij,k}^+$ where $k$ is such that $\{i,j,k\}=\{1,2,3\}$. And 
    $\divisor(\bms_{ij,ik})= (\bmD^2_{123})^{+}  + \bmD_{jk,i}^+$.
\end{lem}

\begin{proof}[Proof of Theorem \ref{theorem_main_example_3}]
Let $\bmB:= \bmE \cup (\bmD^{2}_{123})^+  = \bmD_1 \cup \bmD_5$.
By Lemma \ref{lemma_condition_good_height}, it suffices to check condition \hyperref[condition_B_1]{$(\rmB 1)$} and \hyperref[condition_BH_1]{$(\rmB \rmH 1)$}.
\hyperref[condition_B_1]{$(\rmB 1)$}  follows from Lemma \ref{lemma_example_3_map_to_measure_compactification}.
For any $(\kappa_1,\kappa_2)$,
\[
      \max \left\{  \frac{9-1}{\kappa_1}, \frac{9-1}{\kappa_2} \right\} > \frac{6-1}{\kappa_1+\kappa_2},
\]
hence the limit
\begin{equation*}
    \nu:= \lim_{R \to \infty} \frac{\rmm_{\rmG/\rmH} \cdot \bmone_{B_{R,o}}
    }
    {\rmm_{\rmG/\rmH}(B_{R,o})}
\end{equation*}
exists and is supported on $\bmB(\R)=\bmE(\R)\cup (\bmD_{123}^2)^+(\R)$ by Lemma \ref{lemma_example_3_divisor_gauge_form} and Theorem \ref{theorem_equidistribution_Chamber-Loir_Tschinkel}.
Consequently, \hyperref[condition_BH_1]{$(\rmB \rmH 1)$} holds.

\end{proof}

\begin{rmk}
     Let $\bmL:=\sum_{i=1}^5 \lambda_i \bmD_i$ be a divisor on $\bmX^1$ with $\lambda_i >0$. If we count with respect to a height function associated with $\bmL$, then  ``focusing'' happens exactly when 
      $\min\{\lambda_1,\lambda_2\} > \frac{4}{3} \min\{ \lambda_3,\lambda_4,\lambda_5 \}$. 
\end{rmk}

\appendix

\section{Details on the Examples I, II, III}

\subsection{Example I}

Here we fill in the various missing details from Section \ref{section_example_1}.

\subsubsection{Change of coordinates}\label{subsubsection_example_1_change_of_coordinates}

\begin{defi}
    \begin{itemize}
        \item Define quadratic forms $Q_1^{0}$ and $Q_2^{0}$ by $Q_1^0(x_1,x_2)=2x_1x_2$ and  $Q_2^0(x_1,x_2, x_3, x_4)=2x_1x_4+2x_2x_3$;
         \item Identify vectors with column vectors and  linear maps with matrices.
         Let $M_{Q_i^0}$ ($i=1,2$) be the matrix representations of $Q_i^0$ ($i=1,2$), we have
        \begin{equation*}
            M_{Q_1^0}=
            \left[\begin{array}{cc}
               0  & 1 \\
                1 & 0
            \end{array} \right]
             ,\quad
             M_{Q_2^0}=
            \left[\begin{array}{cccc}
               0  & 0 &0&1\\
                0 & 0&1&0\\
                0&1&0&0\\
                1&0&0 &0
            \end{array} \right];
        \end{equation*}
        \item Let $x_0$ be the point $([M_0:\lambda_0],\la N_0 \ra)$  with $\lambda_0=1$, 
        \[
        M_0 =\left[\begin{array}{cc}
           0 & 0 \\
           1 & 0 \\
           0 & 1 \\
           0 & 0
        \end{array}\right],\quad
        N_0 = \left[ 
        \begin{array}{cc}
           0 & 0 \\
           1 & 0 \\
           0 & 1 \\
           0 & 0
        \end{array}
        \right].
        \]
    \end{itemize}
\end{defi}

 Note that $(Q_1,\R^2)$ is isomorphic to $(Q_1^0,\R^2)$ and $(Q_2, \R^4)$ is isomorphic to $(Q_2^0,\R^4)$.
Therefore, over $\R$, we may and do replace $Q_2$, $Q_1$ by $Q_2^0$ and $Q_1^0$ in the definition of $\bmX^{\INC}$.

\subsubsection{Local coordinates: ideals of the closure}

Without loss of generality, we assume that $([M:\lambda], \la N \ra)$ lies in an open affine subvariety $\bmO$ where, since $\la M \ra$ is contained in $\la N \ra$, they take the form
\begin{equation}\label{equation_example_1_local_expression}
M= \left[ \begin{array}{cc}
      1 & \beta_1 \\
      \alpha_{2} & \beta_{2} \\
      x_1+\alpha_2y_1 & \beta_1x_1 +\beta_2 y_1 \\
      x_2 + \alpha_2 y_2 & \beta_1 x_2 + \beta_2 y_2
    \end{array} \right],\;
    N= \left[ \begin{array}{cc}
      1 & 0 \\
      0 & 1 \\
      x_1 & y_1 \\
      x_2 & y_2
    \end{array} \right].
\end{equation}
Hence we regard  $\bmO$ as a closed subvariety of the affine space 
\begin{equation*}
    \Spec \R[\alpha_2, \beta_1,\beta_2, \lambda, x_1, x_2, y_1, y_2].
\end{equation*}
We need to find the (prime) ideals that correspond to $\bmX^{\INC}$ and various strata.
By restricting to the open subvariety $\lambda\neq 0, \beta_2-\beta_1\alpha_2\neq 0$, we have the equations (write $M=[\bmalpha,\bmbeta]$)
\begin{equation*}
\begin{aligned}
     &\begin{cases}
        Q_2^0(\bmalpha, \bmalpha)=0 \\
        Q_2^0(\bmalpha,\bmbeta)=\lambda^2 \\
        Q_2^0(\bmbeta, \bmbeta)=0 
    \end{cases} 
    \iff
    \begin{cases}
        Q_2^0(\bmalpha, \bmalpha)=0\\
        Q_2^0(\bmalpha,\bmbeta-\beta_1\bmalpha)=\lambda^2\\
        Q_2^0(\bmbeta-\beta_1\bmalpha, \bmbeta- \beta_1\bmalpha)+2\beta_1\lambda^2=0
    \end{cases} \\
    \iff &
    \begin{cases}
        (x_2+\alpha_2 y_2)+ \alpha_2(x_1+\alpha_2y_1)=0\\
        (\beta_2-\beta_1\alpha_2)(y_2+x_1+2\alpha_2y_1) -\lambda^2 =0   \\
        (\beta_2-\beta_1\alpha_2 )
        \left[
          (\beta_2-\beta_1\alpha_2)y_1 + \beta_1 (y_2+ x_1 +2\alpha_2y_1)
        \right] =0 
    \end{cases} \\
    \iff &
    \begin{cases}
        (x_2+\alpha_2 y_2)+ \alpha_2(x_1+\alpha_2y_1)=0
        &(1)
        \\
        (\beta_2-\beta_1\alpha_2)(y_2+x_1+2\alpha_2y_1) -\lambda^2 =0
        &(2)
        \\
          (\beta_2-\beta_1\alpha_2)y_1 + \beta_1 (y_2+ x_1 +2\alpha_2y_1) =0 
          &(3)
    \end{cases} 
\end{aligned}
\end{equation*}
The last ``$\iff$'' is because $\beta_2-\beta_1\alpha_2$ is invertible in this open subvariety.
If $\calI_{\bmX^{\INC}}$ stands for the prime ideal corresponding to $\bmO$, then we have seen that 
\begin{equation*}
    \left\{
    \parbox{60mm}
    {
    \raggedright
    $(x_2+\alpha_2 y_2)+ \alpha_2(x_1+\alpha_2y_1)$,
    $(\beta_2-\beta_1\alpha_2)(y_2+x_1+2\alpha_2y_1) -\lambda^2$,
    $(\beta_2-\beta_1\alpha_2)y_1 + \beta_1 (y_2+ x_1 +2\alpha_2y_1)$
    }
    \right\} \subset \calI_{\bmX^{\INC}}
\end{equation*}
It will be shown that $\calI_{\bmX^{\INC}}$ is actually generated by these polynomials. This follows once we know that the ideal generated by the left hand side, denoted as $\calI_0$ for the moment, is prime. To see this, we first do a few change of variables.

By equations above, we eliminate $x_2$ and replace $\beta_2$ and $y_2$ by
\begin{equation}\label{equation_quadratic_new_local_coordinates}
    \beta_2':= \beta_2- \beta_1 \alpha_2,
    \;
    y_2':= y_2 + x_1 + 2 \alpha_2y_1.
\end{equation}

Then $\bmO$ is the closed subvariety of $\Spec\R[\lambda,\beta_1,\beta_2',y_1,y_2',\alpha_2,x_1]$ corresponding to the smallest prime ideal  containing
\begin{equation}
    \calI_0= \la
     \beta_2'y_2' -\lambda^2,
     \beta_2'y_1+ \beta_1 y_2'
    \ra.
\end{equation}

\begin{lem}\label{lemma_local_equation_closure}
    In $R:=\C[x,y,w,a,b]$, the ideal $I$ generated by $xy-w^2$ and $xa+yb$ is prime. Hence $\calI_0$ is a prime ideal and $\calI_0=\calI_{\bmX^{\INC}}$.
\end{lem}

\begin{proof}
    Once we know that $\calI_0$ is prime, the closed subscheme $X_{0}$ cut out by $\calI_0$ is irreducible. In particular, the closure of $X_0 \cap \{\lambda \neq 0\}$ is $X_0$.  But $\bmX^{\INC}$, by definition, is the closure of $\bmX^{\INC} \cap \{\lambda \neq 0\}=X_0 \cap \{\lambda \neq 0\} $, so $X_0=\bmX^{\INC}$ and $\calI_0 = \calI_{\bmX^{\INC}}$.

    To show $I$ is prime, we first prove that $x$ is not a zero divisor in $R/I$.
    Otherwise, there exist $\phi_0,\phi_1,\phi_2\in R$, $\phi_0 \notin I$ such that 
    \begin{equation}\label{equation_x_not_zero_divisor}
        x\cdot \phi_0 = \phi_1 (xy-w^2) + \phi_2 (xa+yb).
    \end{equation}
    Quotienting by $\la x \ra$, we have  
    \begin{equation*}
        \overline{\phi}_1w^2 = \overline{\phi}_2yb.
    \end{equation*}
    Since $\C[y,w,a,b]$ is a unique factorization domain, we have
    \begin{equation*}
        \phi_1 =x \phi_3 + yb \phi_4,\;
        \phi_2= x \phi_5 + w^2 \phi_4
    \end{equation*}
    for some $\phi_3,\phi_4,\phi_5 \in R$.
    Inserting them back to Equa.(\ref{equation_x_not_zero_divisor}), we have
    \begin{equation*}
        \begin{aligned}
              x\cdot \phi_0
             &= x \cdot  [  \phi_3(xy-w^2) + \phi_5 (xa+yb) 
             +\phi_4 ( y^2b + w^2a)] \\
             \implies
               \phi_0 &=  \phi_3(xy-w^2) + \phi_5 (xa+yb) 
             +\phi_4 ( y^2b + w^2a)  \\
             &= \phi_3(xy-w^2) + \phi_5 (xa+yb) 
             +\phi_4 ( -a (xy-w^2)  + y (xa+yb) ) \in I,
        \end{aligned}
    \end{equation*}
    a contradiction.
    So $x$ is not a zero divisor modulo $I$. By inverting $x$, we see that $\la I \ra = \la y-x^{-1}w^2 , a+ x^{-1}yb \ra$  in $\R[x,{x}^{-1},y,w, a,b]$. In particular, it is a prime ideal. So if $\phi\cdot \psi \in I$, then there exists $n$ such that $x^n \phi \in I$ or $x^n \psi \in I$. But $x$ is not a zero divisor modulo $I$, thus, $\phi\in I$ or $\psi \in I$. This shows that $I$ is prime.
     
\end{proof}

\subsubsection{Local coordinates: ideals of the stratum}

Now we give the explicit local equations of $\bmD_{a,b,c}$. 
Let $I_{a,b,c}$ be the ideal corresponding to $\bmD_{a,b,c}$ in the local affine subscheme given by $\Spec\R[\lambda,\beta_1,\beta_2',y_1,y_2',\alpha_2,x_1]$, same as last subsection.

\begin{lem}
  \begin{equation*}
  \begin{aligned}
            &I_{2,1,0^+}= \la \lambda, \beta_1, \beta_2' \ra,\;
      I_{1,1,1}= \la \lambda, \beta_2', y_2' \ra,\;
      I_{0,2,2}= \la \lambda, y_1, y_2' \ra,\;
      I_{1,1,0^+}= \la \lambda, \beta_1, \beta_2', y_2' \ra,\;\\
      &
      I_{0,1,1}= \la \lambda, y_1, y_2', \beta_2' \ra,\;
      I_{0,1,0^+}= \la \lambda, \beta_1, \beta_2', y_1, y_2' \ra.
  \end{aligned}
  \end{equation*}
\end{lem}

\begin{proof}
    Via column operations, transform $M$ into
    \begin{equation*}
        M'=\left[\begin{array}{cc}
          1   &  0 \\
          \alpha_2   & \beta_2'\\
          x_1+ \alpha_2y_1 & \beta_2' y_1\\
          x_2+ \alpha_2 y_2 & \beta_2' y_2
        \end{array}\right].
    \end{equation*}
    Thus $\rk(M)\leq 1$ iff all two minors of $M'$ vanishes, which happens exactly when $\beta_2'=0$.

    To detect $\rk(\la N \ra,Q_2^0)$, write $N=[\bmx,\bmy]$ (as in Equa.(\ref{equation_example_1_local_expression})) and
    \begin{equation*}
        B_{N}:= 
        \left[\begin{array}{cc}
          Q_2^0(\bmx ,\bmx)   &  Q_2^0(\bmx,\bmy) \\
          Q_2^0(\bmx,\bmy)  & Q_2^0(\bmy, \bmy)
        \end{array}\right]
        =
        \left[\begin{array}{cc}
          2x_1  &  y_2+x_1 \\
          y_2+x_1  & 2y_1
        \end{array}\right].
    \end{equation*}
    Thus the ideal corresponding to $\rk(\la N \ra,Q_2^0)\leq 1$ is the radical of $\la x_2y_1 - (y_2+x_1)^2 , \lambda \ra$, which, after eliminating $x_2$, is equal to 
    $
    \la (y_2+x_1 + 2\alpha_2 y_1)^2, \lambda \ra
    $. So its radical is $
    \la y_2'=y_2+x_1 + 2\alpha_2 y_1, \lambda \ra$.
     And the ideal corresponding to $\rk(\la N \ra,Q_2^0)=0$ is then $\la x_2, y_1, y_2+x_1, \lambda \ra$,  which, after eliminating $x_2$, is equal to $\la y_1,y_2',\lambda \ra$.

    To see $\rk([M^{\Tr}],Q_1^0)$, first we perform row operation on $M$ to get
    \begin{equation*}
        M''=\left[\begin{array}{cc}
          1   &  \beta_1 \\
          0  & \beta_2'\\
          0 & \beta_2' y_1\\
          0 & \beta_2' y_2
        \end{array}\right]
        =  \left[\begin{array}{c}
          \bmv_1 \\
          \bmv_2 \\
          \bmv_3 \\
          \bmv_4
        \end{array}\right].
    \end{equation*}
    Then we compute
    \begin{equation*}
    \begin{aligned}
        B_{M''}
        &=
        \left[\begin{array}{cccc}
          Q_1^0(\bmv_1,\bmv_1)  & Q_1^0(\bmv_1,\bmv_2)
          &  Q_1^0(\bmv_1,\bmv_3) &  Q_1^0(\bmv_1,\bmv_4) 
          \\
           Q_1^0(\bmv_1,\bmv_2)  & Q_1^0(\bmv_2, \bmv_2)  
           & Q_1^0(\bmv_2,\bmv_3)  & Q_1^0(\bmv_2,\bmv_4) 
           \\
          Q_1^0(\bmv_1,\bmv_3) & Q_1^0(\bmv_2,\bmv_3) 
          & Q_1^0(\bmv_3, \bmv_3) &  Q_1^0(\bmv_3,\bmv_4) 
          \\
          Q_1^0(\bmv_1,\bmv_4)  & Q_1^0(\bmv_2,\bmv_4) 
          & Q_1^0(\bmv_2,\bmv_3)  &Q_1^0(\bmv_4, \bmv_4)
        \end{array}\right]  
        \\
        &=  \left[\begin{array}{cccc}
           2\beta_1  & \beta_2'
          &  y_1\beta_2' &  y_2\beta_2'
          \\
           \beta_2'  &  0 &  0 &  0
           \\
          y_1\beta_2' & 0 & 0  &  0
          \\
          y_2\beta_2'  & 0 & 0  &  0
        \end{array}\right].
    \end{aligned}
    \end{equation*}
    Thus the ideal for $\rk( \la M^{\Tr} \ra ,Q_1^0)\leq 1$ is $\la \beta_2',\lambda \ra$ and the ideal for $\rk( \la M^{\Tr} \ra ,Q_1^0) =0$ is $\la \beta_1, \beta_2',\lambda \ra$.

    Now it is not hard to check the claim using these computations.
\end{proof}

\subsubsection{Explicit constructions of a resolution, step I}

By Lemma \ref{lemma_local_equation_closure}, $\bmX^{\INC}$ restricted to this affine open subvariety corresponds to
\begin{equation*}
    \la \beta_2' y_2'- \lambda^2 , \beta_2'y_1+ \beta_1 y_2' \ra \subset \R[\beta_1,\beta_2',\lambda,y_1,y_2']\times \R[\alpha_2,x_1].
\end{equation*}
 As the variables $\alpha_2,x_1$ play no role in the calculation, we will omit them for simplicity. So we work inside a 5-dimensional affine space.

First we blow up $\bmD_{010^{+}}$.

Then $\Bl_{\bmD_{010^{+}}}\bmX^{\INC}$ can be viewed as a subvariety of $\Bl_{\bmD_{010^{+}}}\bmA^5 \subset \bmA^5 \times \bmP^{4}$. If $(x_1:x_2:...:x_5)$ (resp. $(x_1,x_2,...,x_5)$) is a typical point of $\bmP^4$ (resp. $\bmA^5$), we let $\rmU_i$ (resp. $\bmA_i^5$) be the open subvariety of $\Bl_{\bmD_{010^{+}}}\bmA^5$ (resp. $\bmA^5$) with $x_i \neq 0$. 
Therefore, $\Bl_{\bmD_{010^{+}}}\bmX^{\INC} \cap \rmU_i$ is viewed as a subvariety of $\bmA^1 \times \bmA^4 $. Recall that we are going to ignore $\alpha_2,x_1$.

\subsubsection{$\rmU_1$}

On $\rmU_1:=  \{\wtlambda\neq 0 \}$,  we let $\wtbeta_1, \wtbeta_2, \wty_1,\wty_2$ be the coordinates of $\bmA^4$. 
So $\rmU_1$ is identified with
$
    \Spec\R[\lambda] \times \Spec\R[\wtbeta_1, \wtbeta_2, \wty_1,\wty_2]
$
and we have
\begin{itemize}
    \item[1.]  $\beta_1 =\lambda \wtbeta_1,\;
    \beta_2' =\lambda \wtbeta_2,\;
     y_1 =\lambda \wty_1,\;
    y_2' =\lambda \wty_1$;
    \item[2.] The ideal of $\Bl_{\bmD_{010^{+}}}\bmX^{\INC}$ on $\rmU_1$ is $
    \la \wtbeta_2 \wty_2 -1 ,   \wtbeta_2^2 \wty_1 + \wtbeta_1
    \ra
    $;
    \item[3.] The proper transforms of $\bmD_{210}, \bmD_{111}$ and $\bmD_{022}$ are contained in the complement;
    \item[4.] The exceptional divisor $\bmE_1$ is given by $\la \lambda \ra$.
\end{itemize}
According to 2, we can eliminate $\wtbeta_1 = -\wtbeta_2^2 \wty_1$ and view  $\Bl_{\bmD_{010^{+}}}\bmX^{\INC}$ on $\rmU_1$ as a closed subvariety of $\Spec\R[\lambda, \wtbeta_2, \wty_1,\wty_2]$ defined by the ideal $ \la \wtbeta_2 \wty_2 -1 \ra$. With this new coordinate,
\begin{itemize}
    \item[3'.] The proper transforms of $\bmD_{210}, \bmD_{111}, \bmD_{022}$ are contained in the complement;
    \item[4'.] The exceptional divisor $\bmE_1$ is given by $\la \lambda \ra$.
\end{itemize}
So it is clear that we have arrived at a smooth pair on $\rmU_1$.

\subsubsection{$\rmU_2$}
Similarly $\rmU_2= \{ \wtbeta_1 \neq 0 \}$ is  identified with
$
    \Spec\R[\beta_1] \times \Spec\R[\wtlambda, \wtbeta_2, \wty_1,\wty_2]
$.
And
\begin{itemize}
    \item[1.] $ \lambda =\beta_1  \wtlambda , \beta_2' = \beta_1 \wtbeta_2 , y_1 = \beta_1 \wty_1, y_2'= \beta_1 \wty_2
    $;
    \item[2.] The ideal of $\Bl_{\bmD_{010^{+}}}\bmX^{\INC}$ on $\rmU_2$ is $
    \la \wtbeta_2 \wty_2 - \wtlambda^2 ,   \wtbeta_2 \wty_1 + \wty_2
    \ra
    $;
    \item[3.] The proper transform of $\bmD_{210}$ is contained in the complement of $\rmU_2$;
    \item[4.] The proper transform $\bmD_{111}^+$ is $\la \wtlambda, \wtbeta_2,\wty_2 \ra$, $\bmD_{022}^+$ is $\la \wtlambda, \wty_1,\wty_2 \ra$ and $\bmD_{011}^+$ is $\la \wtlambda, \wty_1,\wty_2, \wtbeta_2\ra$;
    \item[5.] The exceptional divisor $\bmE_1$ is given by $\la \beta_1 \ra$.
\end{itemize}
Eliminate $\wty_2 = - \wtbeta_2 \wty_1$ by 2. So view   $\Bl_{\bmD_{010^{+}}}\bmX^{\INC}$ on $\rmU_2$ as the closed subvariety of $  \Spec\R[\beta_1,\wtlambda, \wtbeta_2, \wty_1]$ defined by $ \la \wtbeta_2^2 \wty_1 + \wtlambda^2 
    \ra$. The proper transforms of divisors (on $\rmU_2$) become
\begin{equation*}
        \bmD_{210}^+: \emptyset;  \;
        \bmD_{111}^+: \la \wtlambda, \wtbeta_2 \ra;  \;
        \bmD_{022}^+: \la \wtlambda, \wty_1  \ra;  \;
        \bmD_{011}^+: \la \wtlambda,  \wtbeta_2 , \wty_1  \ra;  \; 
        \bmE_1: \la \beta_1 \ra.
\end{equation*}

\subsubsection{$\rmU_3$}

$\rmU_3$ is identified with
$
    \Spec\R[\beta_2'] \times \Spec\R[
    \wtlambda,\wtbeta_1,\wty_1,\wty_2
    ]$.
And
\begin{itemize} 
    \item[1.] $ \lambda =\beta_1  \wtlambda , \beta_2' = \beta_1 \wtbeta_2 , y_1 = \beta_1 \wty_1, y_2'= \beta_1 \wty_2
    $;
    \item[2.] The ideal of $\Bl_{\bmD_{010^{+}}}\bmX^{\INC}$ is $
    \la \wty_2- \wtlambda^2, \wty_1+ \wtbeta_1 \wty_2
    \ra
    $;
    \item[3.] The proper transforms of $\bmD_{210}, \bmD_{111}$ are contained in the complement;
    \item[4.] The proper transform $\bmD_{022}^+$ is $\la \wtlambda \ra$;
    \item[5.] The exceptional divisor $\bmE_1$ is given by $\la \beta_2' \ra$.
\end{itemize}
Eliminating $\wty_2 = \wtlambda^2 $, $\wty_1= -\wtbeta_1 \wty_2 $ by 2, we see that $\Bl_{\bmD_{010^{+}}}\bmX^{\INC} \cap \rmU_3$ is isomorphic to the affine space $\Spec\R[\beta_2',   \wtlambda,\wtbeta_1] $, with
\begin{itemize} 
    \item[3'.] The proper transforms of $\bmD_{210}, \bmD_{111}$ are contained in the complement;
    \item[4'.] The proper transform $\bmD_{022}^+$ is $\la \wtlambda \ra$;
    \item[5'.] The exceptional divisor $\bmE_1$ is given by $\la \beta_2' \ra$.
\end{itemize}
So we have arrived at a smooth pair.

\subsubsection{$\rmU_4$}

 $\rmU_4$ is identified with
$
    \Spec\R[y_1] \times \Spec\R[\wtlambda, \wtbeta_1, \wtbeta_2,\wty_2].
$
It turns out that $\wtbeta_2=-\wtbeta_1 \wty_2$ can be eliminated and we  get
\begin{itemize}
    \item[1.] $\Bl_{ \bmD_{010^{+}} }\bmX^{\INC} \cap \rmU_4$ is a closed subvariety of
    $\Spec\R[y_1, \wtlambda, \wtbeta_1, \wty_2]$ defined by
     $
    \la \wty_2^2 \wtbeta_1 + \wtlambda^2 
    \ra
    $;
    \item[2.] The proper transform of $\bmD_{022}$ is contained in the complement;
    \item[3.] The proper transform $\bmD_{210}^+$ is $\la \wtlambda, \wtbeta_1 \ra$,
     $\bmD_{111}^+$ is $\la \wtlambda, \wty_2 \ra$ and
     $\bmD_{110}^+$ is $\la \wtlambda, \wtbeta_1,\wty_2 \ra$;
    \item[4.] The exceptional divisor $\bmE_1$ is given by $\la y_1 \ra$.
\end{itemize}

\subsubsection{$\rmU_5$}
 $\rmU_5$ is identified with
$
    \Spec\R[y_2'] \times \Spec\R[\wtlambda, \wtbeta_1, \wtbeta_2,\wty_1].
$
It turns out that both $\wtbeta_1= \wtlambda^2,\wtbeta_2=-\wtbeta_2 \wty_1$ can be eliminated and 
\begin{itemize}
    \item[1.]  $\Bl_{\bmD_{010^{+}}}\bmX^{\INC}$ on $\rmU_5$ is identified with the affine space $\Spec\R[y_2', \wtlambda,\wty_1]$;
    \item[2.] The proper transforms of $\bmD_{022}$ and $\bmD_{111}$ are contained in the complement of $\rmU_5$;
    \item[3.] The proper transform $\bmD_{210}^+$ is $\la \wtlambda\ra$.
    \item[4.] The exceptional divisor $\bmE_1$ is  $\la y_2' \ra$.
\end{itemize}
So we also get a smooth pair here.

\subsubsection{Summary}
On $\rmU_1,\rmU_3,\rmU_5$, our $\Bl_{\bmD_{0,1,0^+}}\bmX^{\INC}$ is already a smooth pair. Below is the intersection patterns of the proper transforms of various boundary components:
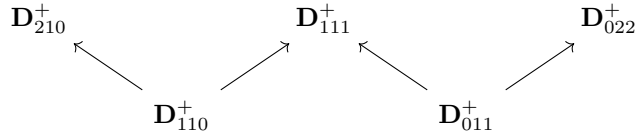
\[
    \begin{tikzcd}
     \bmD_{210}^+  &         &  \bmD_{111}^+ &&  \bmD_{022}^+\\
                      &  \bmD_{110}^+ \arrow[lu] \arrow[ru] &    &   \bmD_{011}^+ \arrow[ru] \arrow[lu]  &
    \end{tikzcd}
\]
where the arrow means ``is contained in'', or more precisely, `` is the intersection of ".
The exceptional divisor $\bmE_1$, which does not appear, intersects all of them transversally.

\subsubsection{Explicit constructions of a resolution, step II}
Call $\bmX^1:= \Bl_{\bmD_{0,1,0^+}}\bmX^{\INC}$.
As the second step, we blow up $\bmD^{+}_{011}\bigsqcup \bmD^{+}_{110}$ on $\rmU_2$ and $\rmU_4$.
Let $\bmE^1_2$ (resp. $\bmE^2_2$) the exceptional divisor corresponding to $\bmD^{+}_{011}$ (resp. $\bmD^{+}_{110}$).

\subsubsection{$\rmU_{21}$}
Over $\rmU_2$ we would have $\rmU_{21}\cup \rmU_{22} \cup \rmU_{23}$, with $\bmD_{011}^+$  being blown up.

$\rmU_{21}$ is identified with
$ \Spec \R[\beta_1, \wtlambda] \times \Spec \R[\wttbeta_2,\wtty_1]$.
We have
\begin{itemize}
    \item[1.] The ideal of $ \Bl_{\bmD^+_{011}}\bmX^{1} \cap \rmU_{21}$ is $\la \wttbeta_2^2 \wtlambda \wtty_1  +1 \ra$;
    \item[2.] The proper transforms of $\bmD_{210}^+$, $\bmD_{022}^+$ and $\bmD_{111}^+$ are contained in the complement;
    \item[3.] The proper transform $\bmE_1^+$ is $\la \beta_1 \ra$.
    \item[4.] The exceptional divisor $\bmE^1_2$ is in the complement.
\end{itemize}

\subsubsection{$\rmU_{22}$}
 $\rmU_{22}$ is identified with
$
 \Spec \R[\beta_1, \wtbeta_2] \times \Spec \R[\wttlambda,\wtty_1].
$
We have
\begin{itemize}
    \item[1.] The ideal of $\Bl_{\bmD^+_{011}}\bmX^{1} \cap \rmU_{22}$ is 
    $\la \wtbeta_2 \wtty_1 + \wttlambda^2 \ra$;
    \item[2.] The proper transforms of $\bmD_{210}^+$ and $\bmD_{111}^+$ are contained in the complement of $\rmU_{21}$;
    \item[3.] The proper transform $\bmD_{022}^{++}$ (of $\bmD_{022}^{+}$) is $\la \wttlambda, \wtty_1\ra$.
    \item[4.] The proper transform $\bmE_1^+$ is $\la \beta_1 \ra$ (strictly speaking, $\bmE_1^+$ should be the intersection of the zero of $\la \beta_1 \ra$ with $\Bl_{\bmD^+_{011}}\bmX^{1}$. Similar abuse of notation also appears below).
    The exceptional divisor $\bmE^1_2$ is $\la \wtbeta_2, \wttlambda \ra$.
\end{itemize}

\subsubsection{$\rmU_{23}$}
 $\rmU_{23}$ is identified with
$
 \Spec \R[\beta_1, \wty_1] \times \Spec \R[\wttlambda,\wttbeta_2]
$
and
\begin{itemize}
    \item[1.] The ideal of $\Bl_{\bmD^+_{011}}\bmX^{1} \cap \rmU_{23}$ is $\la \wttbeta^2_2 \wty_1 + \wttlambda^2 \ra$;
    \item[2.] The proper transforms of $\bmD_{210}^+$ and $\bmD_{022}^{+}$  are contained in the complement of $\rmU_{21}$;
    \item[3.] The proper transform $\bmD_{111}^{++}$ (of $\bmD_{111}^+$)  is $\la \wttlambda, \wttbeta_2\ra$;
    \item[4.] The proper transform $\bmE_1^+$ is $\la \beta_1 \ra$.
    The exceptional divisor $\bmE^1_2$ is $\la \wty_1, \wttlambda \ra$.
\end{itemize}

\subsubsection{$\rmU_{41}$}

Over $\rmU_4$ we would have $\rmU_{41}\cup \rmU_{42} \cup \rmU_{43}$, with $\bmD^+_{110}$ being blown up.

$\rmU_{41}$ is identified with
$
 \Spec \R[y_1, \wtlambda] \times \Spec \R[\wttbeta_1,\wtty_2]
$
and
\begin{itemize}
    \item[1.] The ideal of $\Bl_{\bmD^+_{110}}\bmX^{1} \cap \rmU_{41}$ is 
    $\la \wtlambda   \wttbeta_1  \wtty_2^2+1 \ra$;
    \item[2.] The proper transforms of $\bmD_{210}^+$, $\bmD_{022}^+$ and $\bmD_{111}^+$ are contained in the complement;
    \item[3.] The proper transform $\bmE_1^+$ is $\la y_1 \ra$;
    \item[4.] The exceptional divisor $\bmE^2_2$ is in the complement.
\end{itemize}

\subsubsection{$\rmU_{42}$}
 $\rmU_{42}$ is identified with
$
\Spec \R[y_1, \wtbeta_1] \times \Spec \R[\wttlambda,\wtty_2]
$
and
\begin{itemize}
    \item[1.] The ideal of $\Bl_{\bmD^+_{110}}\bmX^{1} \cap \rmU_{42}$ is
    $\la \wtty_2^2 \wtbeta_1 + \wttlambda^2 \ra$;
    \item[2.] The proper transforms of $\bmD_{210}^+$ and $\bmD_{022}^{+}$  are contained in the complement of $\rmU_{42}$;
    \item[3.] The proper transform $\bmD_{111}^{++}$ (of $\bmD_{111}^+$)  is $\la \wttlambda, \wtty_2\ra$;
    \item[4.] The proper transform $\bmE_1^+$ is $\la y_1 \ra$.
    The exceptional divisor $\bmE^2_2$ is $\la \wtbeta_1, \wttlambda \ra$.
\end{itemize}

\subsubsection{$\rmU_{43}$}
$\rmU_{43}$ is identified with
$
\Spec \R[y_1, \wty_2] \times \Spec \R[\wttlambda,\wttbeta_1]
$
and
\begin{itemize}
    \item[1.] The ideal of $\Bl_{\bmD^+_{110}}\bmX^{1} \cap \rmU_{43}$ is
    $\la \wttbeta_1 \wty_2 + \wttlambda^2 \ra$;
    \item[2.] The proper transforms of $\bmD_{111}^+$ and $\bmD_{022}^{+}$  are contained in the complement;
    \item[3.] The proper transform $\bmD_{210}^{++}$ of $\bmD_{210}^+$  is $\la \wttlambda, \wttbeta_1\ra$;
    \item[4.] The proper transform $\bmE_1^+$ is $\la y_1 \ra$.
    The exceptional divisor $\bmE^2_2$ is $\la \wty_2, \wttlambda \ra$.
\end{itemize}

\subsubsection{Summary}
On $\rmU_{21}$ and $ \rmU_{41}$, we already have a good pair. It remains to deal with $\rmU_{22},\rmU_{23},\rmU_{42}$ and $\rmU_{43}$.
The intersection pattern after this blowup is as follows:
\[
    \begin{tikzcd}
     \bmD_{210}^{++} \arrow[r, leftrightarrow] &  
     \bmE_2^2  \arrow[r, leftrightarrow]        & 
      \bmD_{111}^{++} \arrow[r, leftrightarrow] &
      \bmE_2^1  \arrow[r, leftrightarrow] &
        \bmD_{022}^{++}  
     \end{tikzcd}
\]
where $\leftrightarrow$ means ``intersects with''. As before, $\bmE_1^+$ intersects every divisor here transversally.

\subsubsection{Explicit constructions of a resolution, step III}
Call $\bmX^2:= \Bl_{\bmD_{011}^+\sqcup \bmD^+_{110}}\bmX^{1} $.
As the final step, we blow up 
$ \bmC:=(\bmE_2^1\cap \bmD_{022}^{++})\bigsqcup \bmD_{111}^{++} \bigsqcup
(\bmD_{210}^{++}\cap \bmE^2_2)$.
Let $\bmE^0_3$ (resp. $\bmE^1_3$, $\bmE^2_3$) be the total transform of $\bmE_2^1\cap \bmD_{022}^{++}$ (resp. $\bmD_{111}^{++}$, $\bmD_{210}^{++}\cap \bmE^2_2$).

\subsubsection{$\rmU_{221}$}
Over $\rmU_{22}$ we are going to have $\rmU_{221}\cup \rmU_{222}\cup \rmU_{223}$.

$\rmU_{221}$ is identified with
$
    \Spec\R[\beta_1,\wttlambda] \times \Spec\R[ 
    \wttbeta_2, \wttty_1
    ]
$
and
\begin{itemize}
    \item[1.] The ideal of $\Bl_{\bmC}\bmX^{2} \cap \rmU_{221}$ is
    $\la   \wttbeta_2 \wttty_1 +1 
    \ra$;
    \item[2.] The proper transforms of $\bmD^{++}_{210}$, $\bmD_{111}^{++}$, $\bmD_{022}^{++}$ and $\bmE_2^1$ are in the complement of $\rmU_{221}$;
    \item[3.] The proper transform $\bmE_1^{++}$ is $\la \beta_1 \ra$.
    The exceptional divisor $\bmE^0_3$ is $\la \wttlambda \ra$.
\end{itemize}

\subsubsection{$\rmU_{222}$}
$\rmU_{222}$ is identified with
 $   \Spec\R[\beta_1,\wtbeta_2] \times \Spec\R[ 
    \wtttlambda, \wttty_1
    ] $.
It turns out that one can eliminate $\wttty_1=-\wtttlambda^2$. After this is done, we have: 
\begin{itemize}
    \item[1.]  $\Bl_{\bmC}\bmX^{2} \cap \rmU_{222}$ is identified with the affine space $\Spec\R[\beta_1,\wtbeta_2,  \wtttlambda
    ] $;
    \item[2.] The proper transforms of $\bmD^{++}_{210}$, $\bmD_{111}^{++}$ and $\bmE_2^1$ are in the complement of $\rmU_{222}$.
    The proper transform $\bmD_{022}^{+++}$ is  $\la \wtttlambda \ra$;
    \item[3.] The proper transform $\bmE_1^{++}$ is $\la \beta_1 \ra$.
    The exceptional divisor $\bmE^0_3$ is $\la \wtbeta_2 \ra$.
\end{itemize}

\subsubsection{$\rmU_{223}$}
$\rmU_{223}$ is identified as
$
    \Spec\R[\beta_1,\wtty_1] \times \Spec\R[ 
    \wtttlambda, \wttbeta_2
    ]
$ and we can eliminate $ \wttbeta_2 = - \wtttlambda^2 $. Then we have
\begin{itemize}
    \item[1.] $\Bl_{\bmC}\bmX^{2} \cap \rmU_{223}$ is identified with the affine space $\Spec\R[\beta_1,\wtty_1,\wtttlambda]$;
    \item[2.] The proper transforms of $\bmD^{++}_{210}$, $\bmD_{111}^{++}$ and $\bmD_{022}^{++}$ are in the complement of $\rmU_{223}$;
    \item[3.] The proper transform $\bmE_1^{++}$ is $\la \beta_1 \ra$, $(\bmE_2^1)^{+}$ is $\la \wtttlambda \ra$ and
    the exceptional divisor $\bmE^0_3$ is $\la \wtty_1 \ra$.
\end{itemize}

\subsubsection{$\rmU_{231}$}
Over $\rmU_{23}$ we are going to have $\rmU_{231}\cup \rmU_{232}$.

$\rmU_{231}$ is identified as
$
    \Spec\R[\beta_1, \wty_1, \wttlambda] \times \Spec\R[ 
    \wtttbeta_2
    ]
$
and
\begin{itemize}
    \item[1.] The ideal of $\Bl_{\bmC}\bmX^{2} \cap \rmU_{231}$ is
    $\la   \wtttbeta_2^2 \wty_1 +1 
    \ra$;
    \item[2.] The proper transforms of $\bmD^{++}_{210}$, $\bmD_{111}^{++}$, $\bmD_{022}^{++}$ and $\bmE_2^1$ are in the complement of $\rmU_{221}$;
    \item[3.] The proper transform $\bmE_1^{++}$ is $\la \beta_1 \ra$.
    The exceptional divisor $\bmE^1_3$ is $\la \wttlambda \ra$.
\end{itemize}

\subsubsection{$\rmU_{232}$}
 $\rmU_{232}$ is identified as
$
    \Spec\R[\beta_1, \wty_1, \wttbeta_2] \times \Spec\R[ 
    \wtttlambda
    ]
$
and one eliminates $\wty_1=-\wtttlambda^2$. 
After that, we have:
\begin{itemize}
    \item[1.]  $\Bl_{\bmC}\bmX^{2} \cap \rmU_{232}$ is identified with the affine space $ \Spec\R[\beta_1,  \wttbeta_2,  \wtttlambda] $;
    \item[2.] The proper transforms of $\bmD^{++}_{210}$, $\bmD_{111}^{++}$ and $\bmD_{022}^{++}$ are in the complement of $\rmU_{232}$;
    \item[3.] The proper transform $\bmE_1^{++}$ is $\la \beta_1 \ra$, $(\bmE_2^1)^{+}$ is $\la \wtttlambda \ra$
    and the exceptional divisor $\bmE^1_3$ is $\la \wttbeta_2 \ra$.
\end{itemize}

\subsubsection{$\rmU_{421}$}
Over $\rmU_{42}$ we  have $\rmU_{421}\cup \rmU_{422}$ with the ideal $\la \wttlambda, \wtty_2 \ra$ of $\bmD_{111}^{++}$ being blown up.

$\rmU_{421}$ is identified as
$
    \Spec\R[y_1,\wtbeta_1,\wttlambda] \times \Spec\R[ 
    \wttty_2
    ]
$
and we have
\begin{itemize}
    \item[1.] The ideal of $\Bl_{\bmC}\bmX^{2}$ is
    $\la   \wttty_2^2 \wtbeta_1 +1 
    \ra$;
    \item[2.] The proper transforms of $\bmD^{++}_{210}$, $\bmD_{111}^{++}$, $\bmD_{022}^{++}$ and $\bmE_2^2$ are in the complement of $\rmU_{421}$;
    \item[3.] The proper transform $\bmE_1^{++}$ is $\la y_1 \ra$.
    The exceptional divisor $\bmE^1_3$ is $\la \wttlambda \ra$.
\end{itemize}

\subsubsection{$\rmU_{422}$}
$\rmU_{422}$ is identified as
$
     \Spec\R[y_1,\wtbeta_1,\wtty_2] \times \Spec\R[ 
    \wtttlambda
    ]
$
and we eliminate $\wtbeta_1=-\wtttlambda^2$. Then,
\begin{itemize}
    \item[1.] $\Bl_{\bmC}\bmX^{2} \cap \rmU_{422}$ is identified with the affine space $\Spec\R[y_1,\wtty_2,\wtttlambda]$;
    \item[2.] The proper transforms of $\bmD^{++}_{210}$, $\bmD_{111}^{++}$ and $\bmD_{022}^{++}$ are in the complement of $\rmU_{422}$;
    \item[3.] The proper transform $\bmE_1^{++}$ is $\la y_1 \ra$ and $(\bmE_2^2)^{+}$ is $\la \wtttlambda \ra$.
    The exceptional divisor $\bmE^1_3$ is $\la \wtty_2 \ra$.
\end{itemize}

\subsubsection{$\rmU_{431}$}
Over $\rmU_{43}$, we are going to have $\rmU_{431}\cup \rmU_{432}\cup \rmU_{433}$, with the ideal $\la \wttlambda , \wttbeta_1, \wty_2 \ra$ of $\bmD^{++}_{210} \cap \bmE_2^2$ being blown up.

 $\rmU_{431}$ is identified as
$
    \Spec\R[ y_1,\wttlambda] \times \Spec\R[ 
    \wtttbeta_1, \wtty_2
    ]
$
and we have
\begin{itemize}
    \item[1.] The ideal of $\Bl_{\bmC}\bmX^{2}$ is
    $\la   \wtty_2 \wtttbeta_1 +1 
    \ra$;
    \item[2.] The proper transforms of $\bmD^{++}_{210}$, $\bmD_{111}^{++}$, $\bmD_{022}^{++}$ and $(\bmE_2^2)^{+}$  are in the complement of $\rmU_{221}$;
    \item[3.] The proper transform $\bmE_1^{++}$ is $\la y_1 \ra$.
    The exceptional divisor $\bmE^2_3$ is $\la \wttlambda \ra$.
\end{itemize}

\subsubsection{$\rmU_{432}$}
$\rmU_{432}$ is identified as
$
    \Spec\R[\wty_1,\wttbeta_1] \times \Spec\R[ 
    \wtttlambda, \wtty_2
    ]
$.
One can eliminate $\wtty_2=-\wtttlambda^2$ here after which we have
\begin{itemize}
    \item[1.] $\Bl_{\bmC}\bmX^{2}$ is identified with the affine space $ \Spec\R[\wty_1,\wttbeta_1,  \wtttlambda  ]$;
    \item[2.] The proper transforms of $\bmD^{++}_{210}$, $\bmD_{111}^{++}$ and $\bmD_{022}^{++}$ are in the complement of $\rmU_{432}$;
    \item[3.] The proper transform $\bmE_1^{++}$ is $\la y_1 \ra$, $(\bmE_2^2)^{+}$ is $\la \wtttlambda \ra$.
    The exceptional divisor $\bmE^0_3$ is $\la \wttbeta_1 \ra$.
\end{itemize}

\subsubsection{$\rmU_{433}$}
 $\rmU_{433}$ is identified with
$
     \Spec\R[\wty_1,\wty_2] \times \Spec\R[ 
    \wtttlambda, \wtttbeta_1
    ]
$
and we eliminate  $\wtttbeta_1=-\wtttlambda^2$ here. 
After that,
\begin{itemize}
    \item[1.]  $\Bl_{\bmC}\bmX^{2}$ is identified with the affine space $\Spec\R[\wty_1,\wty_2, \wtttlambda]$;
    \item[2.] The proper transforms of $\bmD^{++}_{111}$, $\bmD_{022}^{++}$ and  $(\bmE_2^1)^{+}$  are in the complement of $\rmU_{433}$.
    The proper transform $\bmD_{210}^{+++}$ is  $\la \wtttlambda \ra$;
    \item[3.] The proper transform $\bmE_1^{++}$ is $\la y_1 \ra$.
    The exceptional divisor $\bmE^2_3$ is $\la \wty_2 \ra$.
\end{itemize}

Finally, we have obtained a smooth pair  everywhere. Let $\bmX^3:= \Bl_{\bmC}\bmX^{2}$. Since each center of the blowup is $\bmG$-invariant, all the blowup morphisms above are $\bmG$-equivariant.
Here is an illustration of the intersection pattern. As above, $\bmE_1^{++}$ intersects everything transversally.
\[
    \begin{tikzcd}
     \bmD_{210}^{+++} \arrow[r, leftrightarrow] &  
     \bmE_3^2  \arrow[r, leftrightarrow]        & 
      (\bmE_2^2)^+  \arrow[r, leftrightarrow] &
       \bmE_3^1  \arrow[r, leftrightarrow]        & 
       (\bmE_2^1)^+  \arrow[r, leftrightarrow] &
         \bmE_3^0  \arrow[r, leftrightarrow]        &
        \bmD_{022}^{+++}  
     \end{tikzcd}.
\]

\subsubsection{Explicit invariant gauge form}

First we present the invariant gauge form on $\SL_2(\R) \times \SL_2(\R)$. In coordinates, elements in 
$\SL_2(\R) \times \SL_2(\R)$ can be written as
\begin{equation}\label{equation_example_1_invariant_form_1}
    \left(\left[
    \begin{array}{cc}
       u_1  & u_2 \\
       u_3  & u_4
    \end{array}
    \right], \,
    \left[
    \begin{array}{cc}
       w_1  & w_2 \\
       w_3  & w_4
    \end{array}
    \right]\right), \text{ satisfying }
    u_1 u_4 - u_2 u_3=w_1w_4-w_2w_3=1.
\end{equation}
Then the invariant form is given by
\begin{equation*}
    \omega_0:= \frac{\diff u_1\wedge \diff u_2 \wedge \diff u_3}{u_1}
    \wedge \frac{\diff w_1\wedge \diff w_2 \wedge \diff w_3}{w_1}.
\end{equation*}
Let us put one more restriction:
\begin{equation}\label{equation_example_1_invariant_form_2}
    f= u_1 w_4 + u_4 w_1 - u_2 w_3 - u_3 w_2 =0,
\end{equation}
which is invariant under the $\SL_2(\R)\times \SL_2(\R)$ action:
$
    (g,h)\cdot (A,B) := (gAh^{-1},gBh^{-1}).
$
One can verify that
\begin{equation*}
    \omega_1:= \frac{\diff u_1 \wedge \diff u_2 \wedge \diff u_3 \wedge \diff w_1 \wedge \diff w_2}{(u_1w_2-u_2 w_1)u_1}
\end{equation*}
satisfies 
\begin{equation*}
    \omega_1 \wedge \diff f =  \omega_0 \text{ restricted to }f=0.
\end{equation*}
And hence $\omega_1$ is an invariant volume form on $\{f=0\}$.

Let $\bmu':=(u'_1, u'_2, u'_3,u'_4) := (u_1,-u_2,u_3,u_4)$ and 
$\bmw':=(w'_1, w'_2, w'_3,w'_4) := (w_1,-w_2,w_3,w_4)$, then Equa.(\ref{equation_example_1_invariant_form_1}, \ref{equation_example_1_invariant_form_2}) become
\begin{equation}\label{equation_u'w'_quadratic}
    Q_2^0(\bmu',\bmu')=Q_2^0(\bmw',\bmw')=2,\;
    Q_2^0(\bmu',\bmw')=0.
\end{equation}
Let us further set
\begin{equation*}
    \bmgamma := \frac{\bmu' + i\bmw'}{2},\;
    \bmtheta:= \frac{\bmu'-i\bmw'}{2}.
\end{equation*}
Then Equa.(\ref{equation_u'w'_quadratic}) becomes
\begin{equation*}
    Q_2^0(\bmgamma,\bmgamma)=Q_2^0(\bmtheta,\bmtheta)=0,\;
    Q_2^0(\bmgamma,\bmtheta)=1,
\end{equation*}
which is exactly the model of homogeneous space discussed above in Section \ref{subsubsection_example_1_change_of_coordinates}.
Using these new variables, 
\begin{equation*}
    \omega_1 = \frac{
    \diffgamma_1 \wedge\difftheta_1 \wedge \diffgamma_2 \wedge \difftheta_2 \wedge ( \diffgamma_3 + \difftheta_3)
    }{
    (\gamma_1 + \theta_1)(\gamma_1\theta_2 - \gamma_2 \theta_1)
    }.
\end{equation*}
Recall the coordinates as introduced in Equa.(\ref{equation_quadratic_new_local_coordinates}) and
\begin{equation*}
    [ 
    \left[
     \begin{array}{cc}
        \gamma_1  &  \theta_1 \\
        \gamma_2  &  \theta_2 \\
       \gamma_3  &  \theta_3  \\
       \gamma_4  &  \theta_4
     \end{array}
    \right] :1
    ] = \left[
     \begin{array}{cc}
       1  &  \beta_1 \\
        \alpha_2  &  \beta_2 \\
       x_1+ \alpha_2 y_1  & \beta_1x_1+ \beta_2 y_1 \\
        x_2+ \alpha_2 y_2   & \beta_1x_1+ \beta_2 y_2
     \end{array}
    \right] :  \lambda
    ].
\end{equation*}
We have
\begin{equation*}
    \begin{aligned}
        \begin{cases}
            \gamma_1 =\frac{1}{\lambda} \implies 
            \diffgamma_1 = - \lambda^{-1}\difflambda
            \\
            \theta_1 =  \frac{\beta_1}{\lambda} \implies
            \difftheta_1 = \pm
            \lambda^{-1}\diffbeta_1 \; \mod \, \difflambda
            \\
            \gamma_2 = \frac{\alpha_2}{\lambda} \implies
            \diffgamma_2 = \lambda^{-1} \diffalpha_2\;
            \mod \, \difflambda
            \\
            \theta_2=\frac{\beta_2}{\lambda} \implies
            \difftheta_2 =\lambda^{-1} \diffbeta_2' \mod\,
            \difflambda,\,\diffbeta_1,\,\diffalpha_2 , 
        \end{cases}
    \end{aligned}
\end{equation*}
and finally,
\begin{equation*}
\begin{aligned}
   \gamma_3 + \theta_3 &=\lambda^{-1} (x_1 + \alpha_2 y_1 +\beta_1 x_1 +\beta_2 y_1) \\
   &\implies 
   \diffgamma_3 +\difftheta_3 = \pm \frac{1+\beta_1}{\lambda}\diff x_1
   \; \mod \difflambda,\,\diffbeta_1,\,\diffalpha_2 ,\, \diffbeta_2' \\
   \gamma_1 +\theta_1 &= \lambda^{-1} (1+ \beta_1),\quad
   \gamma_1\theta_2 - \gamma_2 \theta_1 =\lambda^{-2}\beta_2'.
\end{aligned}  
\end{equation*}
In sum, we get the invariant gauge form on  $\bmX^{\INC}$:
\begin{equation*}
    \omega_1  =  \pm
    \frac{1}{\lambda^3 \beta_2'}\difflambda \wedge \diffbeta_1 \wedge\diffbeta_2' \wedge \diff x_1 \wedge \diffalpha_2.
\end{equation*}
Similar to the above, we can ignore the variables $x_1$ and $\alpha_2$ and focus on
\begin{equation}
    \omega_2:= \frac{1}{\lambda^3 \beta_2'}\difflambda \wedge \diffbeta_1 \wedge\diffbeta_2'.
\end{equation}

We will compute poles of $\omega_2$, or equivalently, $\omega_1$, along the boundaries in $\bmX^3$.

\subsubsection{The order of pole along $\bmE_1^{++}$ is 2}

This can be checked on $\rmU_1$ where
\begin{equation*}
    \beta_1 =\lambda \wtbeta_1,\;\beta_2' = \lambda \wtbeta_2',\;
    \wtbeta_1= -\wtbeta_2^{2} \wty_1.
\end{equation*}
In particular,
\begin{equation*}
    \diff\wtbeta_1 =- \wtbeta_2^2 \diff\wty_1 \mod \diff\wtbeta_2.
\end{equation*}
Consequently,
\begin{equation}
    \omega_2 = \frac{\wtbeta_2}{\lambda^2} \difflambda \wedge\diff\wty_1\wedge \diff\wtbeta_2.
\end{equation}
Since the local equation of $\bmE_1$ is given by $\lambda=0$, we get
$\Pole(\omega_2, \bmE_1)=2$. But the blowup morphisms are trivial over $\rmU_1$ hence $\Pole(\omega_2, \bmE^{++}_1)=2$.

\subsubsection{The order of pole along $\bmD_{022}^{+++}$ is 3}

This can be found on $\rmU_3$ where
$
    \lambda= \wtlambda\beta_2',\;
    \beta_1= \wtbeta_1\beta_2'.
$
Thus,
\begin{equation*}
    \omega_2 = \frac{1}{\wtlambda^3 \beta_2'^4} \diff(\wtlambda\beta_2')\wedge \diff(\wtbeta_1\beta_2')\wedge\diffbeta_2'
    = \frac{1}{\wtlambda^3 (\beta_2')^2} \diff\wtlambda \wedge
    \diff\wtbeta_1\wedge \diffbeta_2'.
\end{equation*}
Since the local equation of $\bmD^{+}_{022}$ is given by $\wtlambda=0$, we have $\Pole(\omega_2, \bmD^{+}_{022})=3$. But the blowup morphism is trivial over $\rmU_3$, so $\Pole(\omega_2, \bmD^{+++}_{022})=3$.

\subsubsection{The order of pole along $\bmD_{210}^{+++}$ is 1}

 This can be checked on $\rmU_5$ where one has the coordinates $y_2',\wty_1,\wtlambda$ and $\wtbeta_2=\wtlambda^2$, $\wtbeta_1= -\wtbeta_2 \wty_1=-\wtlambda^2 \wty_1$ are eliminated. Also, $\lambda= \wtlambda y_2'$, $\beta_1=\wtbeta_1 y_2'$ and $\beta_2 = \wtbeta_2 y_2'$.
 The local equation of $\bmD^+_{210}$ is $\wtlambda=0$.
Then one can compute 
\begin{equation*}
    \begin{aligned}
        \omega_2 &= \frac{1}{
        y_2'^4 \wtlambda^3 \wtbeta_2
        } \diff(\wtlambda y_2') \wedge \diff(
        \wtbeta_1 y_2'
        ) \wedge \diff( \wtbeta_2 y_2' )\\
        &= \pm \frac{1}{(y_2')^2 \wtlambda} \diff y_2'\wedge \diff\wty_1\wedge \diff\wtlambda.
    \end{aligned}
\end{equation*}
Now one sees directly that $\Pole(\omega_2,\bmD_{210}^{+++} )=\Pole(\omega_2,\bmD_{210}^{+}) =1$.

\subsubsection{The order of pole along $\bmE_{3}^{0}$ is 5, along $(\bmE_2^1)^+$ is 7}
These can be checked on $\rmU_{222}$.

First we have that on $\rmU_2$, the coordinates are given by
$\beta_1, \wtlambda, \wtbeta_2 ,\wty_1$ and $\wty_2=-\wtbeta_2 \wty_1$ is eliminated. The local equation of $\bmX^1$ is 
$
    \wtbeta_2^2 \wty_1 + \wtlambda^2 =0.
$
And one has the relations:
\begin{equation*}
    \lambda=\wtlambda \beta_1,\;
    \beta_2= \wtbeta_2 \beta_1,\;
    y_i =\wty_i \beta_1 ,\;i=1,2.
\end{equation*}
From this one can compute
\begin{equation*}
    \omega_2 = \frac{1}{
    \wtlambda^3 \beta_1^2 \wtbeta_2
    } \diff\wtlambda \wedge \diffbeta_1 \wedge \diff\wtbeta_2.
\end{equation*}

On $\rmU_{22}$ one has coordinates $\beta_1, \wtbeta_2, \wttlambda, \wtty_1$ with
$
    \wtlambda= \wtbeta_2 \wttlambda,\;
    \wty_1= \wtbeta_2 \wtty_1
$.
Local equation of $\bmX^2$ is 
$
    \wtbeta_2 \wtty_1 + \wttlambda^2 =0.
$
Then 
\begin{equation*}
    \omega_2 =
    \frac{1}{
    (\wttlambda)^3 \beta_1^2 \wtbeta_2^3
    } \diff\wttlambda \wedge \diffbeta_1 \wedge \diff\wtbeta_2.
\end{equation*}

Finally over $\rmU_{222}$, we have the coordinates $\beta_1, \wtty_1, \wtttlambda$  and $\wttbeta_2=-\wtttlambda^2$ is eliminated. And $\bmX^3 \cap \rmU_{222}$ is the full affine space here. One has $\wttlambda= \wtty_1 \wtttlambda$, $\wttbeta_2= \wtty_1 \wtbeta_2$.
The local equation of $\bmE^0_3$ is $\wtty_1=0$ and that of $(\bmE_2^1)^+$ is $\wtttlambda=0$.
Now compute
\begin{equation*}
    \begin{aligned}
        \omega_2 
        &= \frac{1}{
        (\wtttlambda )^3 \beta_1^2 (\wttbeta_2)^3 (\wtty_1)^6
        } \diff(\wtty_1 \wtttlambda) \wedge \diffbeta_1\wedge \diff(\wtty_1 \wtttlambda^2)\\
        &= \pm
        \frac{1}{
        (\wtttlambda )^7 \beta_1^2 (\wtty_1)^5
        }\diff\wtty_1 \wedge \diffbeta_1 \wedge \diff\wtttlambda.
    \end{aligned}
\end{equation*}
Hence $\Pole(\omega_2, \bmE^0_3)=5$ and $\Pole(\omega_2, (\bmE_2^1)^+)=7$.

\subsubsection{The order of pole along $\bmE_{3}^{1}$ is 3, along $(\bmE_2^2)^+$ is 5}

We shall check this on $\rmU_{422}$.

To cut the story short, on $\rmU_{422}$, we have coordinates $y_1,\wtty_2,\wtttlambda$ and $\bmX^3$ is the full affine space.
The local equation of $\bmE^1_3$ is $\wtty_2=0$ and that of $(\bmE_2^2)^+$ is $\wtttlambda=0$.
The relations with the old coordinates are 
\begin{equation*}
    \begin{aligned}
        &\lambda
        = y_1 \wtlambda =y_1\wtbeta_1 \wttlambda
        =-y_1 (\wtttlambda)^3 \wtty_2,\\
        &\beta_1
        = - y_1 (\wtttlambda)^2,\;\;
        \beta_2'  
        = -y_1 (\wtttlambda)^4 \wtty_2.
    \end{aligned}
\end{equation*}
Now we can compute
\begin{equation*}
    \begin{aligned}
        \omega_2 
        &= \pm \frac{
        \diff(y_1 (\wtttlambda)^3 \wtty_2) \wedge
        \diff(y_1 (\wtttlambda)^4 \wtty_2 ) \wedge
        \diff(y_1 (\wtttlambda)^2  )
        }{
        y_1^3 (\wtttlambda)^9  (\wtty_2)^3 y_1
        (\wtttlambda)^4 \wtty_2
        } \\
        & = \pm
         \frac{1}{y_1^4 (\wtttlambda)^{13} (\wtty_2)^4 } 
         \diff(y_1 (\wtttlambda)^3 \wtty_2 )\wedge
         \left( y_1 (\wtttlambda)^3 \wtty_2 \right)  \diff\wtttlambda
         \wedge \diff(y_1 (\wtttlambda)^2  )\\
        & = \pm
        \frac{1}{y_1^2 (\wtttlambda)^5 (\wtty_2)^3 } \diff\wtty_2 \wedge \diff\wtttlambda \wedge \diff y_1.
    \end{aligned}
\end{equation*}

Hence $\Pole(\omega_2, \bmE^1_3)=3$ and $\Pole(\omega_2, (\bmE_2^2)^+)=5$.

\subsubsection{The order of pole along $\bmE_{3}^{2}$ is 3}
This will be checked on $\rmU_{432}$.
We have here coordinates $y_1,\wttbeta_1,\wtttlambda$ and there are no further constrains on  $\bmX^3$.
Moreover, the local equation of $\bmE^2_3$ is $\wttbeta_1=0$.
The relations with the old coordinates are given by
\begin{equation*}
    \begin{aligned}
        &\lambda
        = y_1 \wtlambda = y_1\wty_2 \wttlambda 
        = y_1 \wttbeta_1^2 \wtty_2 \wtttlambda
        =-y_1 (\wttbeta_1)^2 (\wtttlambda)^3, \\
        &\beta_1= y_1 \wtbeta_1 = y_1 \wty_2 \wttbeta_1=
        -y_1 (\wttbeta_1)^2 (\wtttlambda)^2,\;\\
        &\beta_2'=y_1\wtbeta_2=- y_1 \wtbeta_1\wty_2
        =-y_1\wttbeta_1 \wty_2^2
        =- y_1 (\wttbeta_1)^3 (\wtty_2)^2 
        =- y_1 (\wttbeta_1)^3 (\wtttlambda)^4.
    \end{aligned}
\end{equation*}
Therefore,
\begin{equation*}
    \begin{aligned}
        \omega_2
        &= \pm \frac{
        \diff(
        y_1 (\wttbeta_1)^2 (\wtttlambda)^3
        )\wedge \diff(
        y_1 (\wttbeta_1)^2 (\wtttlambda)^2
        ) \wedge 
        \diff(
        y_1 (\wttbeta_1)^3 (\wtttlambda)^4
        )
        }{
        \left(
        y_1^3 (\wttbeta_1)^6 (\wtttlambda)^9
         \right) \cdot \left(
            y_1 (\wttbeta_1)^3 (\wtttlambda)^4
         \right)
        } \\
        &= \pm
        \frac{1}{  y_1^4 (\wttbeta_1)^9 (\wtttlambda)^{13}}
        \left( 
         y_1 (\wttbeta_1)^2 (\wtttlambda)^2
        \right)\diff\wtttlambda \wedge 
        \diff(
        y_1 (\wttbeta_1)^2 (\wtttlambda)^2
        ) \wedge 
        \diff(
        y_1 (\wttbeta_1)^3 (\wtttlambda)^4
        )
        \\
        &= \pm
        \frac{1}{  y_1^3 (\wttbeta_1)^7 (\wtttlambda)^{11}}
        \diff\wtttlambda\wedge
        (\wtttlambda)^2
        \diff(
        y_1 (\wttbeta_1)^2 
        ) \wedge
        (\wtttlambda)^4
        \diff(
        y_1 (\wttbeta_1)^3 
        ) \\
        &= \pm
        \frac{1}{  y_1^2 (\wttbeta_1)^3 (\wtttlambda)^{5}}
        \diff\wtttlambda\wedge
        \diff y_1  \wedge \diff\wttbeta_1,
    \end{aligned}
\end{equation*}
from which one sees that $\Pole(\omega_2, \bmE^2_3)=3$.



\subsection{Example II}

Here we provide details for the example presented in Section \ref{section_example_2}.
Recall that
$\bmX$ is the Zariski closure of $\bmU:=\bmG.y_0$ with $y_0:=([\bme_1:1],[\bme_1^{\vee}:1])$ in $\bmP^{n+1}_{\Q} \times \bmP^{n+1}_{\Q} $ and $\bmD$ is the complement of $\bmU$ in $\bmX$.
If a general (closed) point in $\bmX$ is denoted as $[\bmv:s],[\bmalpha : t]$, then
$\bmD_1$ is defined by $s=0$ and $\bmD_2$ by $t=0$.

\subsubsection{Local equations of stratum}

We define the following affine open subvarieties of $\bmX$:
\begin{equation*}
    \begin{aligned}
        & \bmO_1: \quad 
            \left\{
        ([1:v_2:...:v_{n+1}:s],[\alpha_1:1:\alpha_3:...:\alpha_{n+1}:t])
        \right\};\\
        & \bmO_2: \quad 
            \left\{
             ([v_1:v_2:...:v_{n+1}:1],[1:\alpha_2:...:\alpha_{n+1}:t])
           \right\}; \\
         & \bmO_3: \quad 
            \left\{
           ([1:v_2:...:v_{n+1}:s],[\alpha_1:\alpha_2:...:\alpha_{n+1}:1])
           \right\}.
    \end{aligned}
\end{equation*}
Using these local coordinates, one finds that\footnote{Compared to last subsection, it is direct to verify the ideals generated by $f_i$ is prime.}
\begin{lem}
    As a closed subvariety of $\Spec \Q[v_2,...,v_{n+1},\alpha_1,\alpha_3,...,\alpha_{n+1} ,s,t]$, the ideal of $\bmO_1$ is generated by 
    \begin{equation*}
        f_1:=\alpha_1 + v_2 +  \sum_{i = 3}^{n+1}  v_i \alpha_i -st.
    \end{equation*}
    The ideal of $\bmD_1 \cap \bmO_1$ (resp. $\bmD_2 \cap \bmO_1$) is  generated by $f_1,s$ (resp. $f_1,t$).
    \\
    As a closed subvariety of $\Spec \Q[v_1,...,v_{n+1},\alpha_2,...,\alpha_{n+1} ,t]$, the ideal of $\bmO_2$ is generated by 
    \begin{equation*}
        f_2:= v_1 +  \sum_{i= 2}^{n+1}   v_i \alpha_i -t.
    \end{equation*}
    The ideal of  $\bmD_2 \cap \bmO_2$ is  generated by $f_2,t$ and
    $\bmD_1 \cap \bmO_2$  is empty.
    \\
     As a closed subvariety of $\Spec \Q[v_2,...,v_{n+1},\alpha_1,...,\alpha_{n+1} ,s]$, the ideal of $\bmO_3$ is generated by 
    \begin{equation*}
        f_3:= \alpha_1 + \sum_{i = 2}^{n+1} v_i \alpha_i -s.
    \end{equation*}
    The ideal of  $\bmD_1 \cap \bmO_3$ is  generated by $f_2, s$ and
    $\bmD_2 \cap \bmO_3$  is empty.
\end{lem}

\subsubsection{The divisor of the invariant gauge form}
\label{section_example_2_divisor_invariant_gauge}

First we work in the affine space $\bmA^{n+1}\times \bmA^{n+1}$ with coordinates
$[v_1:...:v_{n+1}:1]\times [\alpha_1:...:\alpha_{n+1}:1]$ with equation $f=\sum v_i \alpha_i -1=0$.
In the ambient affine space, the differential form
\[
     \omega_1:= \diff v_1 \wedge ... \wedge \diff v_{n+1} \wedge \diffalpha_1 \wedge...\wedge \diffalpha_{n+1}
\]
is $\bmG$-invariant. Then, the unique (up to a scalar) solution $\omega_{\bmU}$ on $f=0$ to
\begin{equation*}
    \omega_1 = \omega_{\bmU} \wedge \diff f
\end{equation*}
is the invariant gauge form.
So when $\alpha_1\neq 0$, 
\[
       \omega_{\bmU}=   \alpha_1^{-1} \diff v_1 \wedge \diff v_3 \wedge... \wedge \diff v_{n+1} \wedge \diffalpha_1 \wedge...\wedge \diffalpha_{n+1}.
\]
Turning into coordinates on $\bmO_1$, we get
\begin{equation*}
\begin{aligned}
        \omega_{\bmU} =&
     (t^{-1})^{-1} \diff \frac{1}{s} \wedge \diff \frac{v_3}{s} \wedge... \wedge \diff \frac{v_{n+1}}{s} \wedge 
     \frac{\diffalpha_1}{t} \wedge \frac{1}{t}\wedge \frac{\alpha_3}{t} \wedge...\wedge \frac{\diffalpha_{n+1}}{t}
     \\
    =& \pm \frac{1}{s^{n+1}t^{n+1}} \diff s \wedge \diff v_3 \wedge...\wedge \diff v_{n+1} \wedge \diff t \wedge \diffalpha_1 \wedge \diffalpha_3 \wedge...\wedge \diffalpha_{n+1}.
    \end{aligned}
\end{equation*}
Now it is direct to see that 
\[
    \Pole( \omega_{\bmU}, \bmD_1) =  \Pole( \omega_{\bmU}, \bmD_2) = n+1. 
\]

\subsubsection{Measure compactification}\label{section_example_3_measure_classifications}

The lemma below follows from the equidistribution theorem of \cite{EskMozSha96} and the nondivergence theorem of \cite{DanMar91} (stated more precisely in \cite[Theorem 4.6]{DawGoroUll18}). Note that there are exactly two rational parabolic subgroups containing $\bmH$. And for a sequence $(g_n)$ in $\rmG$, to test whether adjoint orbits based at the two vectors representing the Lie algebras of these two parabolic subgroups go to $0$, is the same as to test whether $(g_n\bme_1)$ and $(g_n.\bme^{\vee}_1)$ go to $0$.

\begin{lem}\label{lemma_example_3_limit_measure_classifications}
    Let $(g_n)$ be a sequence in $\rmG$. By passing to a subsequence, we are in exactly one of the following situations:
    \begin{itemize}
        \item[1.] $(g_n.\bme_1)$ or $(g_n \bme_1^{\vee})$ tends to $\bmzero$, in which case, 
        \begin{equation*}
            \lim_{n \to \infty} (g_n)_* \rmm^{\bmone}_{[\rmH]} = \bmzero;
        \end{equation*}
        \item[2.] $(g_n)$ converges to $\delta\in \rmG$ modulo $\rmH$, in which case, 
        \begin{equation*}
            \lim_{n \to \infty} (g_n)_* \rmm^{\bmone}_{[\rmH]} = 
             \delta_* \rmm^{\bmone}_{[\rmH]}
            ;
        \end{equation*}
        \item[3.] $(g_n)$ is unbounded modulo $\rmH$, but converges to $\delta\in \rmG$ modulo $\rmM_1$, in which case, 
        \begin{equation*}
            \lim_{n \to \infty} (g_n)_* \rmm^{\bmone}_{[\rmH]} = 
             \delta_* \rmm^{\bmone}_{[\rmM_1]}
            ;
        \end{equation*}
        \item[4.] $(g_n)$ is unbounded modulo $\rmH$, but converges to $\delta\in \rmG$ modulo $\rmM_2$, in which case, 
        \begin{equation*}
            \lim_{n \to \infty} (g_n)_* \rmm^{\bmone}_{[\rmH]} = 
             \delta_* \rmm^{\bmone}_{[\rmM_2]}
            ;
        \end{equation*}
        \item[5.] $(g_n.\bme_1)$ and $(g_n. \bme_1^{\vee})$ both tend to infinity, in which case, 
        \begin{equation*}
            \lim_{n \to \infty} (g_n)_* \rmm^{\bmone}_{[\rmH]} = \rmm_{[\rmG]}^{\bmone}.
        \end{equation*}
    \end{itemize}
\end{lem}

To deduce Lemma \ref{lemma_example_2_map_to_measures} from Lemma \ref{lemma_example_3_limit_measure_classifications}, one observes that $(g_n)$ is bounded modulo $\rmM_1$ iff $(g_n.\bme_1)$ is bounded away from $\bmzero$ and infinity and similarly, it is bounded modulo $\rmM_2$ iff $(g_n.\bme^{\vee}_1)$ is bounded away from $\bmzero$ and infinity.

\subsubsection{The metric line bundles}\label{section_example_2_height_and_metrics}

Let $\rho$ be a non-negative smooth function on $\R^{n+1}$ such that 
\begin{equation*}
    \rho(v) = \begin{cases}
        0 & \text{ if }  \norm{v}\geq 1\\
        1 & \text{ if } \norm{v}\leq 0.5.
    \end{cases}
\end{equation*}
For a positive integer $\kappa$, we define a smooth metric $\norm{\cdot}_{2\kappa}$ on $\calO_{\bmX}(\bmD_1)$ by
\begin{equation*}
    \norm{\bm1_{\bmD_1}([\bmv:s],[\bmalpha:t])}_{2\kappa}
    :=
    \left(
    {\frac{s^{2\kappa}}
    {\sum  v_i^{2\kappa} + s^{2\kappa}
    \rho(\frac{v_1}{s},...,\frac{v_{n+1}}{s}
     )  }  
    }\right)^{\frac{1}{2\kappa}}.
\end{equation*}
Similarly,
define a smooth metric $\norm{\cdot}_{2\kappa}$ on $\calO_{\bmX}(\bmD_2)$ by
\begin{equation*}
    \norm{\bm1_{\bmD_2}([\bmv:s],[\bmalpha:t])}_{2\kappa}
    :=
    \left(
    {\frac{t^{2\kappa}  }
    {  \sum  \alpha_i^{2\kappa}  + t^{2\kappa}
    \rho(\frac{\alpha_1}{t},...,\frac{\alpha_{n+1}}{t}
     )  }  
    }\right)^{\frac{1}{2\kappa}}.
\end{equation*}
One can check they are indeed smooth metrics using local coordinates. Let us check the one on $\calO_{\bmX}(\bmD_1)$.
On $\bmO_1\cup \bmO_3$, we have
\begin{equation*}
\begin{aligned}
        &\norm{ \frac{1}{s}  \bm1_{\bmD_1}([1:v_2:...:v_{n+1}:s],[\bmalpha:t])}_{2\kappa}  \\
    =  &
      \left(\frac{1}{
        (1+v_2^{2\kappa}+...+v_{n+1}^{2\kappa})+ {s}^{2\kappa} \rho(
        \frac{1}{s},\frac{v_2}{s},...,\frac{v_{n+1}}{s}
        )
     }\right)^{\frac{1}{2\kappa}}.
\end{aligned}
\end{equation*} 
Note that  $\rho$ vanishes when $\normm{s} \leq 1$. On $\bmO_2$, we have
\begin{equation*}
\begin{aligned}
        \norm{\bm1_{\bmD_1}([v_1:...:v_{n+1}:1],[\bmalpha:t] )}_{2\kappa}
    =  
     \left(  \frac{1}{
        \sum v_i^{2\kappa}+ \rho(
        v_1,...,v_{n+1}
        )
    } \right)
  ^{\frac{1}{2\kappa}}.
\end{aligned}
\end{equation*}
Since $\calO_{\bmX}(\bmD_1)$ is generated by $\bmone_{\bmD_1}$ on $\bmO_2$ and generated by $\frac{1}{s} \bmone_{\bmD_1}$ on $\bmO_1\cup \bmO_3$, the above computation shows the smoothness of the metric on $\calO_{\bmX}(\bmD_1)$. Note that $\rho$ vanishes when evaluated on integral points of $\bmM_3$.

\subsection{Example III}

Here we provide details for the example presented in Section \ref{section_example_3}.
Recall that we are concerned with certain homogeneous variety $\bmM_3 \cong \bmG/\bmH$ with $\bmG=\SL_3$ and $\bmH$ being the full diagonal torus, which is compactified by
\begin{equation*}
    \bmX^{\INC} = \left\{
      (\bml_{1}, \bml_{2}, \bml_{3}, \bml_{12}, \bml_{13},\bml_{23})
      \in  (\bmP^2)^3 \times (\bmGr_{2,3})^3 
      \;\middle\vert\;
      \bml_{I}\subset \bml_{J} ,\;\forall\, I \subset J
    \right\}.
\end{equation*}
It is equipped with $\bmG$-invariant closed subvarieties:
\begin{equation*}
\begin{aligned}
     &\bmD_{123}^{1} :=\left\{
       \bml_1  =  \bml_2  = \bml_3
    \right\},\; 
    \bmD_{123}^{2} :=\left\{
      \bml_{13}=\bml_{13}=  \bml_{23}
    \right\},\\
    &
    \bmD_{12,3} := \left\{
      \bml_1 = \bml_2,\, \bml_{13} =\bml_{23}
    \right\},\;
    \bmD_{13,2} := \left\{
      \bml_1 = \bml_3,\, \bml_{12} =\bml_{23}
    \right\},\;\\
    &
    \bmD_{23,1} := \left\{
      \bml_2 = \bml_3,\, \bml_{12} =\bml_{13}
    \right\}.
\end{aligned}
\end{equation*}

\subsubsection{An open affine subvariety}

Using Plucker coordinates, we identify $\bmGr_{2,3}$ with $\bmP^2$. So $\bmX^{\INC}$ is a closed subvariety of $(\bmP^2)^6$. For our purpose, we may restrict to the open subvariety $\bmO$ where all the first coordinates are nonzero. That is, points of the form
\begin{equation*}
    \left \{([1:a_2:a_3],[1:b_2:b_3],[1:c_2:c_3],
    [1:\alpha_{13}:\alpha_{23}], [1:\beta_{13}:\beta_{23}], [1:\gamma_{13}:\gamma_{23}] ) \right\} 
\end{equation*}
that are contained in $\bmX^{\INC}$.
Thus, $\bmO$ is a closed subvariety of 
\[
\Spec\R[a_2,a_3, b_2,b_3, c_2,c_3,\alpha_{13},\alpha_{23},\beta_{13},\beta_{23},\gamma_{13},\gamma_{23}]\]
defined by the prime ideal generated by (these functions come from the incidence relations)
\begin{equation*}
\begin{aligned}
    &  f_1 := \alpha_{23}-a_2 \alpha_{13}+ a_3, \;
    f_2 := \beta_{23} - a_2 \beta_{13} + a_3,\;
    f_3:=\alpha_{23} - b_2 \alpha_{13} + b_3 ,\;\\
    &f_4:= \gamma_{23} - b_2 \gamma_{13} + b_3,\;
    f_5:= \beta_{23} -c_2 \beta_{13}+ c_3,\;
    f_6:= \gamma_{23} -c_2 \gamma_{13} + c_3.
\end{aligned}
\end{equation*}
That they indeed generate a prime ideal will be clear later.

\subsubsection{The invariant gauge form}

We compute the invariant gauge form on $\bmM_3$ under these coordinates.
Set
\begin{equation*}
\begin{aligned}
     &\eta_1:= \frac{1}{(\gamma_{23} - a_2 \gamma_{13}+a_3)^3} \diff a_2\wedge\diff a_3 \wedge \diffgamma_{13}\wedge \diffgamma_{23};\\
     &\eta_2:= \frac{1}{(\beta_{23} - b_2 \beta_{13}+b_3)^3} \diff b_2 \wedge \diff b_3 \wedge \diffbeta_{13}\wedge \diffbeta_{23};\\
     &\eta_3:= \frac{1}{(\alpha_{23} - c_2 \alpha_{13}+c_3)^3} \diff c_2\wedge\diff c_3 \wedge \diffalpha_{13}\wedge \diffalpha_{23}.
\end{aligned}
\end{equation*}

\begin{lem}
    For $i=1,2,3$, each $\eta_i$ is $\bmG$-invariant.
\end{lem}

\begin{proof}
    We only present the proof for $\eta_1$. The other two cases are similar.
    
    We write $\bm{a}$ (resp. $\bmgamma$) as a shorthand for $(a_2,a_3)$ (resp. $(\gamma_{13},\gamma_{23})$). Also set $\diff\bm{a}:= \diff a_2\wedge \diff a_3$ and $\diff\bmgamma:=\diffgamma_{13} \wedge \diffgamma_{23}$.
    Let $\phi_1(\bm{a},\bmgamma):= {(\gamma_{23} - a_2 \gamma_{13}+a_3)^3}$ and
    \begin{equation*}
        \begin{aligned}
            g_1:=\left[
            \begin{array}{ccc}
               1  & x &0\\
               0  & 1 & 0\\
              0 & 0 & 1
            \end{array}
            \right],\,
             g_2:=\left[
            \begin{array}{ccc}
               1  & 0 &x\\
               0  & 1 & 0\\
              0 & 0 & 1
            \end{array}
            \right],\,
             g_3:=\left[
            \begin{array}{ccc}
               1  & 0 &0\\
               0  & 1 & x\\
              0 & 0 & 1
            \end{array}
            \right],\,
             g_4:=\left[
            \begin{array}{ccc}
               1  & 0 &0\\
               x  & 1 & 0\\
              y & z & 1
            \end{array}
            \right].
        \end{aligned}
    \end{equation*}
    It suffices to verify that $\eta_1$ is invariant under $g_1, g_2, g_3, g_4$ for all $x,y,z$.
    We only verify the $g_2$-invariance and the other cases are similar.

    For functions $f= a_2, a_3, \bm{a}, \gamma_{13}, \gamma_{23}$ or $\bmgamma$, let $f^{g}:= g_2^{*} (f)$.
    Since
    \begin{equation*}
        \left[
            \begin{array}{ccc}
               1  & 0 &x\\
               0  & 1 & 0\\
              0 & 0 & 1
            \end{array}
            \right] 
            \left[
            \begin{array}{c}
               1  \\
               a_2  \\
              a_3
            \end{array}
            \right] =
            \left[
             \begin{array}{c}
               1+xa_3  \\
               a_2  \\
               a_3
            \end{array}
            \right],
    \end{equation*}
       we have
    \begin{equation*}
        a_2^g = \frac{a_2}{1+xa_3},\;
         a_3^g = \frac{a_3}{1+xa_3},\;
         \text{ and }
         \diff\bm{a}^g = (1+xa_3)^{-3} \diff\bm{a}.
    \end{equation*}
    Under the basis $(\bme_1\wedge \bme_2, \bme_1\wedge \bme_3,\bme_2\wedge \bme_3)$, we have
    \begin{equation*}
        \left( \wedge^2 \left[
            \begin{array}{ccc}
               1  & 0 &x\\
               0  & 1 & 0\\
              0 & 0 & 1
            \end{array}
            \right]   \right)
            \left[
            \begin{array}{c}
               1  \\
               \gamma_{13}  \\
              \gamma_{23}
            \end{array}
            \right]
            =
            \left[
            \begin{array}{ccc}
               1  & 0 &-x\\
               0  & 1 & 0\\
              0 & 0 & 1
            \end{array}
            \right] 
            \left[
            \begin{array}{c}
               1  \\
              \gamma_{13}  \\
              \gamma_{23}
            \end{array}
            \right]
            =
            \left[
            \begin{array}{c}
               1-x\gamma_{23}  \\
               \gamma_{13}  \\
              \gamma_{23}
            \end{array}
            \right].
    \end{equation*}
         Thus,
         \begin{equation*}
             \gamma_{13}^g= \frac{\gamma_{13}}{1-x\gamma_{23}},\;
             \gamma_{23}^g= \frac{\gamma_{23}}{1-x\gamma_{23}},\;
             \diff\bmgamma^g = (1-x\gamma_{23})^{-3}\diff\bmgamma.
         \end{equation*}
       On the other hand,
       \begin{equation*}
       \begin{aligned}
               \phi_1(\bm{a}^g,\bmgamma^g)
               &= 
               (1-x\gamma_{23})^{-1} \gamma_{23} - (1+xa_3)^{-1}a_2\cdot 
               (1-x\gamma_{23})^{-1}\gamma_{13} + (1+xa_3)^{-1}a_3\\
                &=(1-x\gamma_{23})^{-1}(1+xa_3)^{-1} \phi_1(\bm{a},\bmgamma).
       \end{aligned}
       \end{equation*}
     It follows that 
         \begin{equation*}
             \phi_1^{-3}(\bm{a}^g,\bmgamma^g)\diff\bm{a}^g \wedge \diff\bmgamma^g = \phi_1^{-3}(\bm{a},\bmgamma)\diff\bm{a} \wedge \diff\bmgamma.
         \end{equation*}
    That is to say, $\eta_1$ is $g_2$-invariant.
\end{proof}

It can be checked that $f_i$'s intersect transversally in an open set. Thus, the  solution $\omega$ to the equation
\begin{equation}\label{equation_space_of_triangle_restriction_Haar_form}
    \eta_1\wedge\eta_2\wedge\eta_3 =\omega \wedge \diff f_1 \wedge...\wedge \diff f_6
\end{equation}
is unique when restricted to the common zero set of $f_i$'s. As $\eta_i$'s and $\diff f_i$'s are $\bmG$-invariant. It follows that such an $\omega$ would also be $\bmG$-invariant.

\begin{lem}
    Let 
    \begin{equation*}
        \omega_{\bmM_3}:= 
        \frac{
         \diff a_2\wedge\diff a_3 \wedge\diff b_2 \wedge\diff c_2 \wedge \diffalpha_{13}\wedge \diffgamma_{13}
        }{
           (a_2-b_2)^{6} (\alpha_{13}-\gamma_{13})^{3}
        (\beta_{13}-\alpha_{13})^{6} (a_2-c_2)^{4}
        }.
    \end{equation*}
    Then $\omega:= \omega_{\bmM_3}$ solves Equa.(\ref{equation_space_of_triangle_restriction_Haar_form}) above.
\end{lem}

\begin{proof}
    Note that by $f_1=...=f_6=0$, we have
    \begin{equation*}
        \begin{aligned}
            \gamma_{23} -a_2\gamma_{13} +a_3
            &=b_2\gamma_{13} - b_3 -a_2\gamma_{13}+a_3 
            = b_2\gamma_{13} +\alpha_{23} -b_2 \alpha_{13}-a_2 \gamma_{13}+a_3\\
            &=(b_2-a_2)\gamma_{13}+(a_2-b_2)\alpha_{13}
            =(a_2-b_2)(\alpha_{13}-\gamma_{13});\\
            \beta_{23}-b_2\beta_{13} +b_3 
            &=a_2\beta_{13} -a_3 -b_2\beta_{13}+ b_2\alpha_{13} -\alpha_{23} \\
            &=(a_2-b_2)\beta_{13} -a_3 +b_2\alpha_{13} -a_2\alpha_{13} +a_3
            =(a_2-b_2)(\beta_{13}-\alpha_{13});
            \\
            \alpha_{23}-c_2\alpha_{13}+c_3
            &= 
            a_2 \alpha_{13} -a_3 -c_2 \alpha_{13} + c_2\beta_{13} -\beta_{23}\\
            &=
            (a_2-c_2)\alpha_{13} +c_2\beta_{13} -a_2\beta_{13} 
            =(a_2-c_2) (\alpha_{13}-\beta_{13}).
        \end{aligned}
    \end{equation*}
    And modulo (the kernel of) $\omega_{\bmM_3}\wedge \cdot $, we have
    \begin{equation*}
        \begin{aligned}
            &\diff f_1 \wedge \diff f_2 \wedge \diff f_3\wedge \diff f_4 \wedge \diff f_5 \wedge \diff f_6 \\
            \equiv  \, &
            \diffalpha_{23} \wedge (\diffbeta_{23}-a_2 \diffbeta_{13}) \wedge
            \diff b_3 \wedge \diffgamma_{23} \wedge (-c_2\diffbeta_{13}+\diffbeta_{23}) \wedge \diff c_3
            \\
            \equiv \,&
            (a_2-c_2)\diffalpha_{23} \wedge \diffbeta_{23} \wedge \diffbeta_{13} \wedge
            \diff b_3 \wedge \diffgamma_{23} \wedge \diff c_3.
        \end{aligned}
    \end{equation*}
    Putting these computations together yields the result.
\end{proof}

\subsubsection{New coordinates}

We perform a few change of variables and eliminate a few redundant variables to make the description of $\bmO$ clearer. The explicit formula of the group action in the new coordinates will be quite complicated but we are not concerned about it.

Using $f_1,f_2,f_3,f_4,f_5=0$, we can eliminate
\begin{equation*}
    \begin{cases}
        \alpha_{23} = a_2\alpha_{13} - a_3, \\
        \beta_{23}= a_2\beta_{13} - a_3 , \\
        b_{3} = b_2\alpha_{13}  - a_2\alpha_{13} + a_3, \\
        \gamma_{23} = b_2\gamma_{13} - b_2\alpha_{13}  + a_2\alpha_{13} - a_3 ,\\
        c_3 = c_2\beta_{13} - a_2\beta_{13} + a_3
    \end{cases}
\end{equation*}
And $f_6=0$ is equivalent to
\begin{equation*}
\begin{aligned}
    &( b_2\gamma_{13} - b_2\alpha_{13}  + a_2\alpha_{13} - a_3)
    -c_2\gamma_{13} + (c_2\beta_{13} - a_2\beta_{13} + a_3) = 0
    \\
    \iff \,&
     b_2(\gamma_{13} -\alpha_{13})  + a_2(\alpha_{13}
     - \beta_{13} ) +c_2 (\beta_{13}-\gamma_{13}) = 0 \\
     \iff \,&
     - (a_2-c_2) (\beta_{13}-\alpha_{13}) + (b_2-c_2) (\gamma_{13}-\alpha_{13}) =0.
\end{aligned}
\end{equation*}
Now we replace $\beta_{13}, \gamma_{13}, a_2'$ and $b_2'$ by
\begin{equation*}
    \beta_{13}' := \beta_{13} -\alpha_{13},\;
    \gamma_{13}' := \gamma_{13} - \alpha_{13},\;
    a_2' := a_2 - c_2,\; 
    b_2' := b_2 - c_2.
\end{equation*}
Hence the open subvariety $\bmO \subset \bmX^{\INC}$ is identified with the
zero set  of the prime ideal
\begin{equation*}
    \la -a_2'\beta_{13}' + \gamma_{13}' b_2' \ra
    \text{ in }
    \Spec\Q[ a_3, c_2, \alpha_{13} ] \times \Spec\Q[ a_2', b_2',\beta_{13}', \gamma_{13}'].
\end{equation*}
Thus $\bmO$ is singular exactly when $a_2'=\beta_{13}'=\gamma_{13}'=b_2'=0$.
The first three coordinates $a_3, c_2, \alpha_{13}$ have no effect on the discussions below and are often omitted for simplicity.

\subsubsection{Invariant gauge form in the new coordinates}
Under the new coordinates, we have
\begin{equation*}
    \omega_{\bmM_3} = \pm 
    (\diff a_3 \wedge \diff c_2 \wedge \diffalpha_{13} ) \wedge 
    \frac{
     \diff a_2' \wedge \diff b_2' \wedge \diffgamma_{13}'
    }{
    (a_2' - b_2')^{6} (\gamma_{13}')^{3} (\beta_{13}')^{6} (a_2')^{4}
    }.
\end{equation*}

\subsubsection{Boundary stratum in the new coordinates}

We list the ideal for $\bmO$ and various boundary stratum intersecting $\bmO$:
\begin{equation*}
\text{The ideal of }
    \begin{cases}
    \bmO \text{ is generated by } &  -a_2'\beta_{13}' + \gamma_{13}' b_2'
    \\
     \bmD_{123}^1 \text{ is generated by } &  a_2', b_2'
    \\
     \bmD_{12,3} \text{ is generated by } &  a_2'-b_2', \beta_{13}'-\gamma_{13}'
    \\
     \bmD_{13,2} \text{ is generated by } &  a_2', \gamma_{13}'
     \\
      \bmD_{23,1} \text{ is generated by } &  b_2', \beta_{13}'
      \\ 
      \bmD_{123}^2 \text{ is generated by } &   \beta_{13}', \gamma_{13}'
    \end{cases}.
\end{equation*}

\subsubsection{Blowup}\label{section_example_3_blow_up}
We blow up the ideal $\la a_2', b_2' \ra$ of $\bmD_{123}^1$ and verify in local coordinates that the obtained pair $(\bmX^1,\bmD^1)$ is, as has been pointed out in \cite{Roberts_Speiser_1984_triangles},  smooth.  Also we calculate the invariant gauge form in this smooth pair.

Naturally $\Bl_{\bmD_{123}^1} (\bmO) $ is the union of two affine open subvarieties $\bmU_1 \cup \bmU_2$. 
For simplicity, coordinates $a_3,c_2,\alpha_{13}$ are omitted below.

\subsubsection{Blowup, $\bmU_1$}
Here we have coordinates $\beta_{13}',\gamma_{13}',a_2',\wtb_2$ with
 $ b_2'= a_2 ' \wtb_2$. And the equation defining $\Bl_{\bmD_{123}^1} (\bmO) $  is $\beta_{13}'=\gamma_{13}' \wtb_2$. Thus we can eliminate $\beta_{13}'$ and treat $\bmU_1$ as $\Spec \Q[\gamma_{13}',a_2',\wtb_2]$.

 One can compute that (let $\bmE$ denote the exceptional divisor)
\begin{itemize}
     \item[1.] $\omega_{\bmM_3} = \pm {
       (a_2')^{-9} (1-\wtb_2)^{-6} (\gamma_{13}')^{-9} (\wtb_2)^{-6}
       }\cdot 
       {\diff a_2'\wedge \diff\wtb_2 \wedge \diffgamma'_{13}}$;
     \item[2.] $\bmD_{12,3}^+ =\{1-\wtb_2=0 \},\, \bmD_{13,2}^{+} \text{ is in the complement},\,\bmD_{23,1}^+ =\{\wtb_2=0\}$;
     \item[3.] $\bmE= \{ a_2' =0\},\, (\bmD_{123}^{2})^+ = \{ \gamma_{13}'=0 \}$.
\end{itemize}

\subsubsection{Blowup, $\bmU_2$}
Here we have coordinates $\beta_{13}',\gamma_{13}',b_2',\wta_2$ with
 $ a_2'= b_2 ' \wta_2$. And the equation defining $\Bl_{\bmD_{123}^1} (\bmO) $ in $\bmU_1$ is $\gamma_{13}'=\beta_{13}' \wta_2$. Thus we can eliminate $\gamma_{13}'$ and treat $\bmU_2$ as $\Spec \Q[\beta_{13}',b_2',\wta_2]$.

 One can compute that (let $\bmE$ denote the exceptional divisor)
\begin{itemize}
     \item[1.] $ \omega_{\bmM_3} = \pm
       {
       (b_2')^{-9} (1-\wta_2)^{-6} (\beta_{13}')^{-9} (\wta_2)^{-6}
       }\cdot 
       {\diff b_2'\wedge \diff\wta_2 \wedge \diffbeta'_{13}} $;
     \item[2.] $ \bmD_{12,3}^+ = \{1 - \wta_2=0 \}  ,\, 
       \bmD_{13,2}^{+} = \{\wta_2=0\}  ,\,
       \bmD_{23,1}^+  \text{ is in the complement}  $;
     \item[3.] $  \bmE= \{ b_2' =0\},\, (\bmD_{123}^{2})^+ = \{ \beta_{13}'=0 \}  $.
\end{itemize}

\subsubsection{The divisor of the invariant gauge form}\label{section_example_3_pole_invariant_gauge}
Labeling the boundary divisors on $\bmX^2$ as 
\begin{equation*}
    \bmD_1 := \bmE,\; 
    \bmD_2 := \bmD_{12,3}^+,\;
    \bmD_3 := \bmD_{13,2}^{+},\;
    \bmD_4 := \bmD_{23,1}^+,\;
    \bmD_5:= (\bmD_{123}^{2})^+.
\end{equation*}
From the description above, 
one gets
\begin{equation*}
   - \divisor (\omega_{\bmM_3}) =
   9\bmD_1 + 6\bmD_2 + 6 \bmD_3 + 6\bmD_4 + 9\bmD_5.
\end{equation*}


\subsubsection{Mapping to the measure compactification space}\label{section_example_3_measure_classification}

Let $\Gamma \leq \bmG(\Q)$ be an arithmetic subgroup. 
We have the following, implied by \cite[Theorem 2.4]{ShaZhe18}.

\begin{thm}
    Let $(g_n)$ be a sequence in $\rmG$. 
    \begin{itemize}
        \item[1.] If $(g_n)$ is unbounded modulo $\rmH_{ij}$ for every $\{i,j\}\subset\{1,2,3\}$, then 
        \begin{equation*}
            \lim_{n \to \infty} \left[(g_n)_*\rmm_{[\rmH]}  \right]
         =
        \left[\rmm_{[\bmG]}   \right];
        \end{equation*}
        \item[2.] If $(g_n)$ is unbounded modulo $\rmH$ but convergent to $\delta \in \rmG$ modulo $\rmH_{ij}$ for some (necessarily unique) $\{i,j\}$, then 
        \begin{equation*}
            \lim_{n \to \infty} \left[(g_n)_*\rmm_{[\rmH]} \right] =
        \left[\delta_*\rmm_{[\rmH_{ij}]}  \right].
        \end{equation*}
     \end{itemize}
\end{thm}
As a corollary, we have
\begin{coro}
Let $X^{\meas}$ be the closure of $\{ \alpha^{\psi}_{\rmG}(\rmm^{\psi}_{[\rmH]})  \}$ in $\Prob^{\psi}(\rmG/\Gamma)$. Then 
\begin{equation*}
    X^{\meas} =  \alpha^{\psi}_{\rmG}(\rmm^{\psi}_{[\rmH]})
    \bigsqcup 
    \left(
      \alpha^{\psi}_{\rmG}(\rmm^{\psi}_{[\rmH_{12}]})  \cup  
       \alpha^{\psi}_{\rmG}(\rmm^{\psi}_{[\rmH_{13}]}) 
     \cup  
      \alpha^{\psi}_{\rmG}(\rmm^{\psi}_{[\rmH_{23}]}) 
    \right) \bigsqcup  \left\{ \rmm^{\psi}_{[\bmG]} \right\}.
\end{equation*}
\end{coro}



\section{Divergence of Translates of Homogeneous Closed Subsets}

Let $\bmG$ be a connected linear algebraic group over $\Q$, $\bmH$ be a connected observable $\Q$-subgroup, and $\Gamma$ be an arithmetic subgroup of $\bmG$. Fix a maximal reductive connected  $\Q$-subgroup $\bmG^{\red}$ of $\bmG$ and hence $\bmG=\bmG^{\red} \ltimes \bmR_{\bmu}(\bmG)$. Also fix a Cartan involution and hence a maximal compact subgroup $\rmK$ of $\bmG^{\red}(\R)$.
Depending on this choice, each parabolic $\Q$-subgroup $\bmP$ is associated with a $\Q$-split subtorus $\bmA_{\bmP,\rmK}$ that is isomorphic to $\bmS_{\bmP}$ under $p^{\spl}_{\bmP}$.
 Unlike the main body of the paper, in this appendix, the Roman letter $\rmL$ is used to denote $\bmL(\R)^{\circ} $ for an algebraic group over $\R$.

Let $\Phi^{\red}(\bmA_{\bmP,\rmK},\bmP)$ be the nontrivial characters of $\bmA_{\bmP,\rmK}$ appearing in $\frakr_{\bmu}(\frakp)/\frakr_{\bmu}(\frakg)  \subset \frakg/\frakr_{\bmu}(\frakg)$. Let $\Delta^{\red}(\bmA_{\bmP,\rmK},\bmP)$ be the subset of cardinality $\dim \bmA_{\bmP,\rmK}$ whose $\Z_{\geq 0}$-span equals to $\Phi^{\red}(\bmA_{\bmP,\rmK},\bmP)$. When $\bmG$ is already reductive, the 
superscript ``${\red}$'' is dropped.

\begin{thm}\label{theorem_divergence_observable}
 Given a sequence $(g_n) \subset \rmG$, after passing to a subsequence, there exist a sequence $(h_n)\subset \rmH$, $(\gamma_n)\subset \Gamma$ and a parabolic $\Q$-subgroup $\bmP$ such that the following holds. Write $g_n h_n \gamma_n^{-1} = k_n a_n p_n$ using horospherical coordinates of  $(\bmP,\rmK)$.
 Then
 \begin{itemize}
     \item[(1)]  $(p_n)$ is bounded;
     \item[(2)]  $\alpha(a_n) \to 0$ for every $\alpha\in \Delta^{\red}(\bmA_{\bmP,\rmK},\bmP)$;
     \item[(3)]  if $(a_n)$ is unbounded, there exist a $\Q$-representation $\bmV$ of $\bmG$ factoring through $\bmG/\bmR_{\bmu}(\bmG)$ and $\bmv\in \bmV(\Q)$ such that the line spanned by $\bmv$ is preserved by $\bmP$, $\bmv$ is fixed by  $\gamma_n \bmH \gamma_n^{-1}$ for all $n$ and $\lim_{n\to \infty} a_n. \bmv = \bmzero$;
     \item[(4)] $(\gamma_n \bmH \gamma_n^{-1})$ strongly converges to some observable subgroup of $\bmG$.
 \end{itemize}
 \end{thm}
 
 Compared to \cite[Theorem 5.2]{zhangrunlinDcds2022}, 3 and 4 are new. Note that whereas one quotients $\Gamma$ from the left in \cite{zhangrunlinDcds2022}, here we choose to quotient from the right. Also, when one has $\alpha(a_n) \to \infty$ in \cite{zhangrunlinDcds2022}, we have $\alpha(a_n)\to 0 $ here.
 We shall follow the logic of \cite{zhangrunlinDcds2022} by first reducing to the case of $\bmG=\SL_N$. And this special case will be handled with the help of ``canonical polygons''. The proof is independent from the work of \cite{DawGoroUllLi19} and can actually be used to give an alternative treatment of their main results. The connection between canonical polygons and nondivergence property of unipotent flows has been noted in \cite{deSaxce_2023}.
 
 For simplicity, by saying the $(\bmP,\rmK)$ coordinate of some element, we mean its \textbf{horospherical coordinate} associated to $(\bmP,\rmK)$ (see \cite[Lemma 2.1]{zhangrunlinDcds2022}, \cite[Chapter 9]{BorJi06} or \cite[Proposition 1.5]{BorSer73}).
\subsection{Reduction to the reductive case}
 
 We explain how the general case follows from the case when $\bmG$ is reductive.
 Indeed, by projecting everything to $p^{\red}: \bmG \to \bmG/\bmR_{\bmu}(\bmG)$, after passing to a subsequence, we find 
 $(h_n)\subset \rmH$, $(\lambda_n)\subset \Gamma$ and a parabolic $\Q$-subgroup $\bmP$ such that if $g_n h_n \lambda_n^{-1} = k_n a_n q_n$  is the $(\bmP,\rmK)$ coordinate of $g_n h_n \lambda_n^{-1}$, then
  \begin{itemize}
     \item[1.]  $q_n = b_n u_n$ for some $(b_n)$ bounded in $\rmP$  and $(u_n)\subset \bmR_{\bmu}(\bmG)(\R)$.
     \item[2.]  $\alpha(a_n) \to 0$ for every $\alpha\in \Delta^{\red}(\bmA_{\bmP,\rmK},\bmP)$;
     \item[3.] there exist a $\Q$-representation $\bmV$ of $\bmG$ factoring through $p^{\red}$ and $\bmv\in \bmV(\Q)$ such that the line spanned by $\bmv$  
                    is preserved by $\bmP$, $\bmv$ is fixed by  $\gamma_n \bmH \gamma_n^{-1}$ for all $n$ and $\lim_{n\to \infty} a_n. \bmv = \bmzero$;         
      \item[4.] $\left( p^{\red}\left( \lambda_n \bmH \lambda_n^{-1} \right) \right)$ strongly converges to some 
                     observable subgroup of $\bmG/\bmR_{\bmu}(\bmG)$.
 \end{itemize}
 
Write $u_n =b_n' \lambda_n'$ for some bounded sequence $(b_n')\subset \bmR_{\bmu}(\bmG)(\R)$ and $( \lambda_n' )\subset \bmR_{\bmu}(\bmG)(\R) \cap \Gamma$. Let $\gamma_n:= \lambda_n' \lambda_n$. It only remains to verify item 4. in the theorem. Since $(\lambda'_n)$ are contained in the kernel of $p^{\red}$, we still have $\left( p^{\red}\left( \gamma_n \bmH \gamma_n^{-1} \right) \right)$ strongly converges to some observable subgroup $\bmF$ of $\bmG/\bmR_{\bmu}(\bmG)$. After passing to a subsequence assume $(\gamma_n \bmH \gamma_n^{-1})$ strongly converges to $\bmL$. Then $p^{\red}(\bmL)= \bmF$. Note that $(p^{\red})^{-1}(\bmF)$  is observable in $\bmG$ by assumption. On the other hand, $(p^{\red})^{-1} (\bmF) = \bmL\cdot \bmR_{\bmu}(\bmG)$ and hence $(p^{\red})^{-1} (\bmF) / \bmL \cong \bmR_{\bmu}(\bmG)/\bmU$ for some $\Q$-subgroup $\bmU$ of $\bmR_{\bmu}(\bmG)$. Therefore, $ \bmL $ is observable in $(p^{\red})^{-1} (\bmF)$ and hence in $\bmG$. This completes the proof.

\subsection{Reduction to the $SL_N$ case}

Let $\bmG$ be a connected reductive linear algebraic group over $\Q$. Without loss of generality, we assume that $\bmG$ is a $\Q$-subgroup of some $\SL_N$ that is invariant under taking transpose. We let $\rmK:= \bmG(\R)\cap \SO_N(\R)$ be the maximal compact subgroup of $\bmG(\R)$ that is transpose invariant.

The concept of Siegel sets can be generalized to reductive groups as follows. Write $\bmG= \bmG^{\semisimple} \cdot \bmZ(\bmG)^{\an}\cdot \bmZ(\bmG)^{\spl}$. Then every Siegel set $\frakS$ of $\bmG$ is of the form $\frakS' \cdot F \cdot \bmZ(\bmG)^{\spl}(\R)^{\circ}$ where $\frakS \subset \bmG^{\semisimple}(\R)$ is a Siegel set for $\bmG^{\semisimple}$ and $F\subset \bmZ(\bmG)^{\an}(\R)$ is a compact subset that is left invariant under the maximal compact subtorus.
For every $\Q$-minimal parabolic subgroup $\bmP$ of $\bmG$, there exist finitely many $c_1,...,c_l\in \bmG(\Q)$ and a Siegel set $\frakS$ associated with $(\bmP,\rmK)$ such that 
\[
\bmG(\R) = \bigcup_{i=1}^l \frakS \cdot c_i \cdot \Gamma.
\]

The main argument of \cite[Section 5]{zhangrunlinDcds2022} carries through. We sketch the proof below.
Without loss of generality, assume $\Gamma= \SL_N(\Z) \cap \bmG(\R)$.

By the $\SL_N$ case, after passing to a subsequence, there are $(h_n)\subset \rmH$, $(\gamma_n)\subset \SL_N(\Z)$ and a parabolic $\Q$-subgroup $\bmQ\subset \SL_N$ such that if the $(\bmQ,\rmK')$ coordinate (for simplicity we write $\rmK':=\SO_N(\R)$) of $g_nh_n\gamma_n^{-1}$ is $k_na_nq_n$, then
\begin{itemize}
    \item[1.] $(q_n)$ is bounded;
    \item[2.] $\alpha(a_n) \to 0$ for every $\alpha \in \Delta(\bmA_{\bmQ,\rmK'},\bmQ)$;
    \item[3.] $\gamma_n\bmH \gamma_n^{-1} \subset \bmQ$ for all $n$. And there exist a $\SL_N$-representation and a $\Q$-vector $\bmv$ fixed by $\gamma_n \bmH \gamma_n^{-1}$ and the line spanned by $\bmv$ is preserved by $\bmQ$ with a nontrivial character;
    \item[4.] $(\gamma_n \bmH \gamma_n^{-1}) $ strongly converges to an observable subgroup $\bmL$ of $\SL_N$.
\end{itemize}
Note that $a_n.\bmv \to \bmzero$ follows automatically.

\subsubsection{Step 1}
Assume $\fraksl_N(\Z) \cap \Lie( \bmR_{\bmu}(\bmQ)) = \oplus \Z w_i^{\bmQ}$ and let $\bmw:= \oplus w_i^{\bmQ}$. Then 
\[ 
\lim_{n\to \infty} g_nh_n\gamma_n^{-1} .\bmw =\bmzero.
\]
 By geometric invariant theory, there exist $n_0\in \Z^+$ and $\bma_t \in \frakX^{\Q}_{*}(\bmG)$ such that 
\[
\lim_{t \to \infty} \bma_t \gamma_{n_0}^{-1} .\bmw =0.
\]
Let (similar notations are also used for cocharacters other than $\bma_t$)
\begin{equation*}
    \begin{aligned}
        \bmQ_{\bma_t} &:= \left\{
        x \in \SL_N \midd \lim_{t \to \infty} \bma_t x \bma_t^{-1} \text{ exists}
        \right\};\\
        \bmP_{\bma_t}&:= \left\{
          x \in \bmG \midd \lim_{t \to \infty} \bma_t x \bma_t^{-1} \text{ exists}
        \right\} .    
    \end{aligned}
\end{equation*}
Then $\bmR_{\bmu}(\bmQ_{\bma_t}) $ contains $\bmR_{\bmu}(\gamma_{n_0}^{-1}\bmQ \gamma_{n_0})$, or equivalently, $\bmQ_{\bma_t} \subset \gamma_{n_0}^{-1}\bmQ \gamma_{n_0}$.  Replacing $\bmQ$ by $\gamma_{n_0}^{-1}\bmQ \gamma_{n_0}$ and $\gamma_n$ by $\gamma_{n_0}^{-1}\gamma_n$, we assume without loss of generality that $\gamma_{n_0}$ is the identity element.

\subsubsection{Step 2}

Let $\bmL_{\bmP_{\bma_t},\rmK}$ (resp. $\bmL_{\bmQ_{\bma_t},\rmK'}$) be the maximal reductive subgroup of $\bmP_{\bma_t}$ (resp. $\bmQ_{\bma_t}$) that is invariant under taking transpose. 
Let $\bmA_{\bmP_{\bma_t},\rmK}$ be the central torus of $\bmL_{\bmP_{\bma_t},\rmK}$, which is conjugate to a maximal $\Q$-split torus $\bmS_{\max}$ in $\bmP_{\bma_t}$.
We have $\bmL_{\bmP_{\bma_t},\rmK} \subset  \bmL_{\bmQ_{\bma_t},\rmK'} $ and $\bmR_{\bmu}(\bmP_{\bma_t}) \subset \bmR_{\bmu}(\bmQ_{\bma_t}) $.
There exists $\bmb_t \in \frakX_{*}^{\Q}(\bmS_{\max})$ such that 
\begin{itemize}
    \item[1.] the centralizer of $\{\bmb_t\}$ in $\SL_N$ is equal to the centralizer of $\bmS_{\max}$ in $\SL_N$. In particular $\bmP_{\bmb_t}$ is a $\Q$-minimal parabolic subgroup of $\bmG$;
    \item[2.] $\bmQ_{\bmb_t} \subset \bmQ_{\bma_t}$. In particular, $\bmP_{\bmb_t} \subset \bmP_{\bma_t}$.
\end{itemize}
Thus $\bmA_{\bmP_{\bmb_t},\rmK}$ is contained in $\bmA_{\bmQ_{\bmb_t},\rmK'}$.
Consequently, the  $( \bmP_{\bmb_t} , \rmK)$ coordinate of $g\in \bmG(\R)$  is the same as its $(\bmQ_{\bmb_t} ,\rmK' )$ coordinate.

\subsubsection{Step 3}

Because $\bmP_{\bmb_t}$ is $\Q$-minimal, by passing to a subsequence, there exist $c_1\in \bmG(\Q)$, $(\lambda_n)\subset \Gamma$ and a Siegel set $\frakS(\bmP_{\bmb_t})$ associated to $( \bmP_{\bmb_t} , \rmK)$ such that 
\begin{equation*}
    g_n h_n \in \frakS(\bmP_{\bmb_t}) c_1 \lambda_n ,\quad
     \forall\, n\in \Z^+.
\end{equation*}
In particular, if   $g_nh_n\lambda_n^{-1}c_1^{-1}= k_n'a_n' p_n'$ under $( \bmP_{\bmb_t} , \rmK)$ coordinate, then $(p_n')$ is bounded and for some $t_0 >0$,  for every $n$ and $\alpha \in \Delta(\bmA_{\bmP_{\bmb_t},\rmK},\bmP_{\bmb_t})$, one has $\alpha(a_n') < t_0$.

Let $\Phi(\bmA_{\bmQ_{\bmb_t},\rmK'}, \fraksl_N)$ be the collection of nontrivial weights of $\bmA_{\bmQ_{\bmb_t},\rmK'}$ appearing in the adjoint action on $\fraksl_N$.
Passing to a subsequence, assume that for every $\beta\in \Phi(\bmA_{\bmQ_{\bmb_t},\rmK'}, \fraksl_N)$, the sequence $\left( \beta (a_n') \right)$ is bounded (away from $0$ and $+\infty$), or converges to $0$, or diverges to $+\infty$.
Define a parabolic $\Q$-subgroup by
\[
\bmQ_{(a_n')}:= \left\{
 x\in \SL_N \midd \left( a_n' x a_n'^{-1} \right) \text{ is bounded}
\right\}.
\]
One notes that $\bmP_{\bmb_t}$ is contained in $\bmQ_{(a_n')}$.

Since $(a_n')$ is a sequence in $\rmA_{\bmP_{\bmb_t},\rmK}$, which is conjugate to $\rmS_{\max}$ in $\bmP_{\bmb_t}$, we can
choose $\bmb'_t \in \frakX_*^{\Q}(\bmS_{\max})$ such that
\[
 \bmQ_{\bmb_t'}  \subset \bmQ_{(a_n')}  ,\;
\bmP_{\bmb'_t} = \bmP_{\bmb_t},\;
\bmZ_{\SL_N}(\{\bmb_t'\}) = \bmZ_{\SL_N}(\bmS_{\max}).
\]
Then the $(\bmP_{\bmb_t}=\bmP_{\bmb'_t} ,\rmK )$ coordinate of $g_nh_n\lambda_n^{-1}c_1^{-1}$  is the same as its $(\bmQ_{\bmb'_t},\rmK')$ coordinate. From the definition of $\bmQ_{(a'_n)}$, one deduces that for any weight $\alpha \in \Phi(\bmA_{\bmP_{\bmb'_t},\rmK},\bmQ_{(a'_n)}   )$ and in particular $\alpha \in \Phi(\bmA_{\bmP_{\bmb'_t},\rmK},\bmQ_{\bmb_t'}   )$, we have $\alpha(a_n')< t_1$ for some $t_1 >0$ and for all $n$. 
We conclude from here that $g_nh_n\lambda_n^{-1}c_1^{-1}$ is contained in a Siegel set $\frakS(\bmQ_{\bmb_t'})$ associated to $(\bmQ_{\bmb_t'},\rmK')$.

\subsubsection{Step 4}

Though $\bmP_{\bmb_t}= \bmP_{\bmb'_t}$, it is not clear whether $\bmQ_{\bmb_t}$ and $\bmQ_{\bmb'_t}$ contain a common Borel subgroup.
Find $w_1 \in \SL_N(\Q)$ such that $\bmQ':= w_1 \bmQ_{\bmb_t'}w_1^{-1} \cap \bmQ_{\bmb_t}$ is a parabolic subgroup. Thus $g_nh_n\lambda_n^{-1}c_1^{-1}w_1^{-1}$ is contained in some 
Siegel set of $(w_1 \bmQ_{\bmb'_t} w_1^{-1}, \rmK') $ and hence some 
Siegel set $\frakS(\bmQ')$ associated to $(\bmQ',\rmK')$.
On the other hand, $(g_nh_n\gamma_n^{-1})$ is known to be contained in some Siegel set attached to $(\bmQ,\rmK')$ and hence in some Siegel set $\frakS(\bmQ')'$ attached to $(\bmQ' ,\rmK')$. 
Replacing by a larger Siegel set if necessary, assume $\frakS(\bmQ')=\frakS(\bmQ')'$. Hence
\[
g_nh_n \gamma_n^{-1} \in \frakS(\bmQ') \cap \frakS(\bmQ') w_1 c_1 \lambda_n \gamma_n^{-1}.
\]
In particular, this intersection is nonempty. But for every $q\in \SL_N(\Q)$, the set $\{\gamma \in \SL_N(\Z),\; \frakS(\bmQ') \cap \frakS(\bmQ')q \gamma\}$ is finite. Therefore, after passing to a subsequence, there exists $\gamma_1'\in \SL_N(\Z) $ such that $\lambda_n\gamma_n^{-1}=\gamma_1'$ for all $n$.

\subsubsection{Step 5}
Let $\bmQ'':= \gamma_1' \bmQ \gamma_1'^{-1}$, which contains $\lambda_n \bmH \lambda_n^{-1}$ for all $n$. As in Step 1 and 2, find some maximal $\Q$-split torus $\bmS_{\max}'$ of $\bmG$ and  $ \bmb_t'' \in \frakX_{*}^{\Q}(\bmS_{\max}')$ such that 
\begin{itemize}
    \item[1.] $\bmZ_{\SL_N}(\{\bmb_t''  \}) = \bmZ_{\SL_N}(\bmS_{\max}')$. In particular, $\bmP_{\bmb''_t}$ is a $\Q$-minimal parabolic subgroup;
    \item[2.] $\bmQ_{\bmb_t''}\subset \bmQ''$.
\end{itemize}
By assumption, $g_nh_n\lambda_n^{-1}=g_nh_n \gamma_n^{-1}\gamma_1'^{-1}$ belongs to some Siegel set associated to $(\bmQ'' ,\rmK' )$ and hence to $(\bmQ_{\bmb_t''} ,\rmK' )$. Let 
\[
     g_nh_n\lambda_n^{-1}= k_n'' a_n'' p_n''
\]
be the $(\bmP_{\bmb''_t }  , \rmK)$ coordinate, which is the same as the $(\bmQ_{\bmb''_t}, \rmK')$ coordinate.
We find that 
\begin{itemize}
    \item[1.] $(p_n'')$ is bounded;
    \item[2.] there exists $t_2>0$ such that for every $n$ and $\alpha \in \Delta(\bmA_{\bmP_{\bmb''_t},\rmK},\bmP_{\bmb''_t})$, we have $\alpha(a_n'') < t_2$;
    \item[3.] $\bmQ''=\bmQ_{(g_nh_n\lambda_n^{-1})} = \bmQ_{(a_n'')}$.
\end{itemize}
 Choose $\bmc_t \in \frakX_*^{\Q}(\bmS_{\max}')$ such that $\bmQ_{(a_n'')}= \bmQ_{\bmc_t}$.
After passing to a subsequence, assume that for every $\alpha \in \Delta(\bmA_{\bmP_{\bmb''_t},\rmK},\bmP_{\bmb''_t})$, either $(\alpha(a_n'') )$ converges to $0 $ or is bounded away from $0$.
Let $I''$ be those $\alpha$ such that $(\alpha(a_n'') )$ is bounded away from $0$. Then 
$(\bmP_{\bmb_t''})_{I''} = \bmP_{\bmc_t}
= \bmQ_{\bmc_t} \cap \bmG = \bmQ'' \cap \bmG$.
 In particular, $\lambda_n \bmH \lambda_n^{-1} $ is contained in $\bmP_{\bmc_t}$ for all $n$. Also, if $g_nh_n \lambda_n^{-1}= k_n'''a_n''' p_n'''$ is the $(\bmP_{\bmc_t} ,\rmK )$ coordinate, then $( p_n''')$ is bounded and $\alpha(a_n''') \to 0$ for all $\alpha \in \Delta(\bmA_{\bmP_{\bmc_t},\rmK}, \bmP_{\bmc_t})$.
 
 Finally, by assumption, there exist certain $\SL_N$-representation and a $\Q$ vector $\bmv''$ fixed by $\lambda_n \bmH \lambda_n^{-1}$ and the line spanned by $\bmv''$ is preserved by $\bmQ''$ with a nontrivial character, which must be nontrivial on $\bmc_t$ since $\bmQ''= \bmQ_{\bmc_t}$. In particular, $\bmv''$ is preserved by $\bmP_{\bmc_t}$ with a nontrivial character. So we are done.

\subsection{The case of $SL_N$}

Fix an observable $\Q$-subgroup $\bmH$ of $\SL_N$ and a sequence $(g_n)\subset \SL_N(\R)$.  We prove Theorem \ref{theorem_divergence_observable} in this case using the idea of canonical polygons (see \cite{Casselman_2004, Grayson_semistability_1984}). For simplicity we abbreviate $\bmA_{\bmP,\rmK'}$ as $\bmA_{\bmP}$ below and continue call $\rmK':= \SO_N(\R)$.

\subsubsection{Lattice flags}

A subgroup $\Delta$ of $\Z^N$ is said to be primitive iff $(\Delta \otimes \Q) \cap \Z^N = \Delta$.
Let $\Prim(\Z^N)$ be the collection of all primitive subgroups of $\Z^N$.
For a $\Q$-character $\alpha$ of $\bmH$, let $\Prim^{\bmH}_{\alpha}(\Z^N)$ denote those $\Delta = \Z \bmv_1\oplus ... \oplus \Z \bmv_k$ with 
\[
h.  \bmv_1\wedge ...\wedge \bmv_k = \alpha(h)  \bmv_1\wedge ...\wedge \bmv_k,
\quad \forall \, h\in \rmH.
\]
Let $\Prim^{\bmH}(\Z^N)$ be the union of these $\Prim^{\bmH}_{\alpha}(\Z^N)$'s.

A totally ordered subset of $\Prim(\Z^N)$ (without loss of generality, assumed to contain $\{\bmzero\}$ and $\Z^N$) is referred to as a \textbf{lattice flag}, written as 
\begin{equation*}
    \calF= \left\{ \{\bmzero \} =  \Delta_0 \subsetneq \Delta_1 \subsetneq ...\subsetneq \Delta_l =\Z^N \right\}.
\end{equation*}
The number $l-1$ is called the rank of $\calF$.
For $i=1,..., l$, we define
\begin{equation}\label{equation_quotient_lattice}
    \overline{\Delta_i}:=
    \left( \Delta_i/\Delta_{i-1} \right) ^{'}:= 
    \left( \frac{
    \norm{\Delta_i}
    }{
    \norm{\Delta_{i-1}}
    } \right)^{ -
    \frac{1}{\rank(\Delta_i)-\rank(\Delta_{i-1})}
    } \cdot \Delta_i/\Delta_{i-1},
\end{equation}
which has covolume one in the quotient Euclidean space $\Delta_i\otimes \R/\Delta_{i-1}\otimes \R$. 
We also define a parabolic $\Q$-subgroup by
\[
\bmP_{\calF}:= \left\{ g\in \SL_N(\C) \midd
g \text{ preserves } \Delta_i\otimes \C ,\; \forall \, i
\right\}.
\]
Given a subset $\calI $ of $\{0,1,...,N\}$ containing $\{0,N\}$, ordered as 
$
\calI = \left\{
    0=i_0 < i_1 <...< i_l=N
\right\}
$,
let $\calF_{\calI}$ be the \textbf{standard flag} associated to $\calI$ defined by
\[
\calF_{\calI} := \left\{ \{\bmzero \}  \subsetneq \oplus_{i=1}^{i_1} \Z \bme_{i} \subsetneq ... \subsetneq \oplus_{i=1}^{i_{l-1}} \Z \bme_{i} \subsetneq \Z^N \right\}.
\]
For simplicity we write $\bmP_{\calI}:= \bmP_{\calF_{\calI}}$.

\subsubsection{Quotient lattices}\label{subsubsection_quotient_lattices}

Fix some $\calI = \left\{
    0=i_0 < i_1 <...< i_l=N
\right\}$ for this subsection.

Define, for $k=1,...,l$,
\begin{equation*}
    \begin{aligned}
        \Lambda_k^{\Std} := \bigoplus_{i=1}^{i_k} \Z.\bme_i,\quad
        V_k:= \Lambda_k^{\Std} \otimes_{\Z} \R ,\quad
        \olV_k:=V_k/V_{k-1}, \quad
        \olV_k(\Z):= \Lambda_k^{\Std}/\Lambda_{k-1}^{\Std},
    \end{aligned}
\end{equation*}
and $\pi_k: \rmP_{\calI}  := \bmP_{\calI}(\R)^{\circ} \to \SL(\olV_k)$ by
\begin{equation*}
    p \mapsto c_p \cdot \left( p\vert_{V_k} \;(\mathrm{ mod } \; V_{k-1}) \right)
\end{equation*}
where $c_p$ is the unique positive real number such that the right hand side has determinant one. $\pi_k$ descends to a continuous map 
\begin{equation*}
    \olpi_k: \rmP_{\calI}/\rmP_{\calI} \cap \SL_N(\Z) \to \SL(\olV_k)/\SL(\olV_k(
    \Z)).
\end{equation*}
that is equivariant with respect to $\pi_k$.
The natural map 
\begin{equation*}
    p \mapsto \left\{ \{\bmzero\} = p \Lambda_0^{\Std} \subsetneq  p \Lambda_1^{\Std} 
    \subsetneq ... \subsetneq p \Lambda_{l}^{\Std} = \Z^N
    \right\}
\end{equation*}
induces a bijection between 
\begin{equation*}
     \rmP_{\calI}/\rmP_{\calI} \cap \SL_N(\Z)  \cong
         \mathrm{LF}_{\calI}:= \left\{
        (\Delta_{\star}) \text{ is a lattice flag } \midd
        (\Delta_i)_{\R} = V_i,\,\forall\,i=1,...,l
       \right\} .
\end{equation*}
If we additionally identity 
\begin{equation*}
    \SL( \olV_k)/\SL( \olV_k(
    \Z)) \cong 
    \left\{
       \text{ unimodular lattices in } \olV_k\,
     \right\}
\end{equation*}
by $g\mapsto g \olV_k(\Z)$,
then $\olpi_k$ is nothing but
\begin{equation*}
    (\Delta_{\star}) \mapsto \olDelta_k
\end{equation*}
defined in last subsection (see Equa.(\ref{equation_quotient_lattice})).

Let $\olp^{\spl}: \rmP_{\calI}/\rmP_{\calI}\cap \SL_N(\Z) \to \rmS_{\bmP_{\calI}}$. It is direct to check the following
\begin{lem}\label{lemma_flag_lattice}
    For each $\calI$, the map
     $\rmP_{\calI}/\rmP_{\calI}\cap \SL_N(\Z) \to
     \rmS_{\bmP_{\calI}} \times \prod_{k=1}^l  \SL(V_k)/\SL(V_k(
    \Z))$ defined by $\Phi_{\calI}:= \olp^{\spl} \times \prod_{k} \olpi_k$ is a proper continuous map.
\end{lem}

For $k=1,...,l$, define $j_k:= i_k-i_{k-1}$ and
 $f_k, \varphi_{k},\alpha_k: \rmA_{\bmP_{\calI}} \to \R^+$ by
\begin{equation*}
   f_k(a):= \norm{a \Lambda^{\Std}_k},
   \quad
   \varphi_k (a) := \left( \frac{ \norm{a \Lambda^{\Std}_k}
    }{
    \norm{a \Lambda^{\Std}_{k-1} }
    } \right) ^{\frac{1}{j_k} },
    \quad
    \alpha_k(a) := \varphi_{k}(a) / \varphi_{k+1}(a).
\end{equation*}
One can check that
\begin{equation}\label{equation_root_lattice}
    \Delta(\rmA_{\bmP_{\calI}} , \bmP_{\calI})
    = \left\{ \alpha_k \midd k=1,...,l-1 \right\}.
\end{equation}

\subsubsection{Canonical flags}
In this subsection we fix an element $g\in \SL_N(\R)$.
Define
\begin{equation*}
\begin{aligned}
    \Plot^{\bmH}(g)  &:= 
    \left\{
         \left( \rank(\Delta), \log (\norm{g\Delta}) \right) \midd
         \Delta \in \Prim^{\bmH}(\Z^N)
    \right\};\\
    \Poly^{\bmH}(g) &:=
    \text{closure of the convex hull of }\Plot^{\bmH}(g).
\end{aligned}
\end{equation*}
$\Poly^{\bmH}(g)$ is a polygon of finitely many sides. Let
\begin{equation*}
    \begin{aligned}
        \calE^{\bmH}(g):= 
        \left\{
          \Delta\in \Prim^{\bmH}(\Z^N) \midd
            \left( \rank \Delta, \log(\norm{g \Delta}) \right) \text{ is an extreme point of }\Poly^{\bmH}(g)
        \right\}.
    \end{aligned}
\end{equation*}
\begin{lem}
    For every $g\in \SL_N(\R)$, $\calE^{\bmH}(g)$ is a lattice flag containing $\{\bmzero\}$ and $\Z^N$.
\end{lem}
A proof can be found in \cite[Section 4]{Casselman_2004}. Let $r_g:= \rank(\calE^{\bmH}(g))$ and write
\begin{equation*}
     \calE^{\bmH}(g) = 
     \left\{ \{\bmzero \} =  \Delta_0(g) \subsetneq \Delta_1(g) \subsetneq ...\subsetneq \Delta_{r_g+1}(g) =\Z^N \right\},
\end{equation*}
called the \textbf{canonical flag}. Define
\begin{equation*}
    \begin{aligned}
        \type(g) &:= 
        \left\{
            0=\rank(\Delta_0(g)) < \rank(\Delta_1(g)) <...< \rank(\Delta_{r_{g}+1}(g)) = N
        \right\};\\
        \slope_k(g)  &:= \frac{
        \log \norm{(g \Delta_k(g)} ) - \log( \norm{g \Delta_{k-1}(g)} )
        }{
        \rank(\Delta_k(g) ) -  \rank(\Delta_{k-1}(g) )
         },\;
         k =1,...,r_g+1;
         \\
         d_k(g)& := \slope_{k+1}(g) - \slope_{k}(g),\;
         k = 1,...,r_g.
    \end{aligned}
\end{equation*}

\subsubsection{Variation of canonical flags}

Fixing some $g\in \SL_N(\R)$, we would like to search for a better one among  $\calE^{\bmH}(gh)$ as $h$ varies in $\rmH$ by (locally) optimizing certain quantities. This part is an extra key step required to treat general connected observable subgroups instead of just those without nontrivial $\Q$-characters.

Find a $\Q$-split subtorus $\bmS^{\bmH}$ such that $\bmH$ is an almost direct product of $({}^{\circ}\bmH)^{\circ}$ and $\bmS^{\bmH}$. Thus the Lie algebra $\frakh= {}^{\circ}\frakh \oplus \fraks^{\frakh}$. Write $\fraks^{\frakh}_{\R}$ for the $\R$-span of $\fraks^{\frakh}$. For a $\Q$-character $\alpha$ of $\bmH$, we regard $\diff \alpha$ as a linear functional on $\fraks^{\frakh}_{\R}$ since it vanishes on ${}^{\circ}\frakh_{\R}$.
Let $\Phi_{\bmH}^{\Std}$ collect all $\diffalpha$ such that $\Prim^{\bmH}_{\alpha}(\Z^N)$ is nonzero.

For $\lambda \in \R$, let $\CVX_{\lambda}$ be the collection of all convex closed polygons $\calP$ such that
\begin{itemize}
    \item[1.] $[0,N]\times \R_{\geq 0} \subset \calP \subset \{(x,y) ,\; y\geq \lambda, \, 0\leq x \leq N\}$;
    \item[2.] $(x,y)\in \Extre(\calP)\implies x \in \{0,1,...,N\}$ and $\{(0,0),(N,0)\} \subset \Extre(\calP)$.
\end{itemize}
Equip $\cup \CVX_{\lambda}$ with the Chabauty topology.
Note that $\Poly^{\bmH}(g)$ belongs to $\CVX_{\lambda}$ for some $\lambda$.
Given some $\calP \in \CVX_{\lambda}$ for some $\lambda \in \R$, we let 
\begin{itemize}
    \item[1.] $y_i(\calP):= \inf\left\{ y \midd (i,y)\in \calP  \right\}$ for $i\in\{0,...,N\}$. So $y_0=y_N=0$;
    \item[2.] $s_i(\calP):= y_i(\calP)-y_{i-1}(\calP)$ for $i\in \{1,...,N\}$;
    \item[3.] $d_i(\calP) = s_{i+1}(\calP) -s_{i}(\calP) = y_{i+1}(\calP)+y_{i-1}(\calP)-2y_i(\calP)$ for $i\in \{1,...,N-1\}$.
\end{itemize}
Thus $d_i(g) = d_{\rank(\Delta_i(g))}(\Poly^{\bmH}(g)) $.
Also let $\tau_{\calP}$ be a permutation on $\{1,...,N-1\}$ such that
\begin{itemize}
    \item[1.] $\left\{d_i(\calP) \midd i=1,...,N-1 \right\} = \left\{ d_{\tau_{\calP}(1)}(\calP) \geq d_{\tau_{\calP}(2)}(\calP)  \geq .... \geq d_{\tau_{\calP}(N-1)}(\calP) \right\} $;
    \item[2.] if $d_{\tau_{\calP}(i)}(\calP) = d_{\tau_{\calP}(i+1)}(\calP) $ then $\tau_{\calP}(i) < \tau_{\calP}(i+1)$.
\end{itemize}
Define a partial order $\prec$ on $\cup \CVX_{\lambda}$ by $\calP \prec \calP'$ iff $\calP = \calP'$ or 
$y_{\tau_{\calP}(m)} (\calP) <  y_{\tau_{\calP'}(m)}(\calP') $ with $m:= \min\left\{ i\midd y_{\tau_{\calP}(i)} (\calP) \neq y_{\tau_{\calP'}(i)}(\calP')  \right\}$. The following explains why it is possible to find a local maxima.
\begin{lem}\label{lemma_limit_polygon}
    Let  $g\in \SL_N(\R)$, $\lambda\in \R$ and $\calP \in \CVX_{\lambda}$.
    If there exists a sequence $(h_n)\subset \rmH$ such that 
    \[
    \calP= \lim_{n \to \infty} \Poly^{\bmH}(g h_n),
    \]
    then there exists $h\in \rmH$ such that $\calP = \Poly^{\bmH}(gh)$.
\end{lem}

\begin{proof}
    Let $V_{\fraks}:= \bigcap_{\diff\alpha \in \Phi_{\bmH}^{\Std}} \ker (\diff\alpha)$
    and choose some complementary $\R$-linear subspace $W_{\fraks}$ such that $\fraks^{\frakh}_{\R} = V_{\fraks} \oplus W_{\fraks}$.
    Without loss of generality, assume the existence of $(w_n)\subset W_{\fraks}$ such that 
    \[
    \calP= \lim_{n \to \infty} \Poly^{\bmH}(g \exp(w_n)).
    \]
    So we can find $\lambda' \in \R$ such that $\Poly^{\bmH}(g \exp(w_n)) \in \CVX_{\lambda'}$ for all $n$. In particular,  there exists $C_1 \in \R$ such that $\diffalpha(w_n)>C_1$ for all $n$ and $\diffalpha \in \Phi_{\bmH}^{\Std}$.
    
    Passing to a subsequence, assume that for every $\diff\alpha\in \Phi_{\bmH}^{\Std}$, either $\diff \alpha(w_n) \to +\infty$ or remains bounded.
    Let $\Phi_{\bdd}\subset \Phi_{\bmH}^{\Std}$ correspond to the bounded ones and $\Phi_{\infty}$ be its complement.
    Let $W':= \bigcap_{\diff\alpha \in \Phi_{\bdd}} \ker \diff\alpha \cap W_{\fraks}$ and take another subspace $W''$ such that $W_{\fraks} = W'\oplus W''$.
    Write $w_n= w_n'+ w_n''$. Then $(w_n'')$ is contained in some bounded subset $B\subset W''$ and $\diff\alpha(w_n') \to +\infty$ for every $\diff\alpha \in \Phi_{\infty}$.

    Find $u_0\in W'$ such that for every $v\in B$ one has
    \begin{equation}\label{equation_appendix_divergence_unbounded_characters}
        \inf_{\diff\alpha \in \Phi_{\infty}}
        \inf_{\Delta \in \Prim^{\bmH}_{\alpha}(\Z^N)} 
        \log \left(
           \norm{  g\exp(u_0+ v) \Delta }
        \right) >0.
    \end{equation}
    Define $w_n^{\new}:= u_0 + w_n''$.
    Since $\norm{g\exp(w_n^{\new}) \Delta } =  \norm{g\exp(w_n) \Delta }$ for every $\Delta \in \Prim^{\bmH}_{\alpha}(\Z^N)$ and $\diffalpha\in \Phi_{\bdd}$, we have (by Equa.(\ref{equation_appendix_divergence_unbounded_characters}))
    \[
    \Poly^{\bmH} (g \exp(w_n) ) =\Poly^{\bmH} ( g\exp(w_n^{\new} ) )  .
    \]
    Being bounded, we can select a subsequence such that $(w_n^{\new})$ converges to some $w_{\infty}$. Then
    $
    \calP= \Poly^{\bmH}(g \exp(w_{\infty})).
    $
\end{proof}
Thanks to Lemma \ref{lemma_limit_polygon}, for every $g\in \SL_N(\R)$, there exists (and we fix such an) $h_g\in \rmH$ such that a local maxima with respect to $\prec$ is taken at $\Poly^{\bmH}(gh_g)$ among 
$\{ \Poly^{\bmH}(gh),\; h\in \rmH \}$. 
\begin{lem}\label{lemma_singular_canonical_weights_levels}
    Take some $g\in \SL_N(\R)$.
    For $C\in \R$, let $\scrI_C:= \{i=1,...,r_g ,\; d_i(gh_g)>C\}$.
    For $\Delta \in \Prim^{\bmH}(\Z^N)$, we let $\alpha_{\Delta}$ be the unique $\alpha \in \frakX^*_{\Q}(\bmH)$ such that $\Delta\in \Prim_{\alpha}^{\bmH}(\Z^N)$ and regard its differential $\diff\alpha_{\Delta}$ as a linear functional on $\fraks^{\frakh}_{\R}$.
    For every $C\in \R$ such that $\scrI_C$ is nonempty, there exists a set of positive integers $(a_i)_{i\in \scrI_C}$  such that $\sum_{i\in \scrI_C} a_i \diff\alpha_{\Delta_i(gh_g)} =0$.
\end{lem}
\begin{proof}
    If the conclusion were wrong for some $C\in \R$, then $\bmzero$ would not belong to the relative interior of the convex hull of $\left\{ \diff\alpha_{\Delta_i(gh_g)} ,\; i \in \scrI_C\right\}$, which is nonempty.
   Therefore,
   there exists $w_0 \in \fraks^{\frakh}_{\R}$ such that $\diff\alpha_{\Delta_i(gh_g)} (w_0) \geq 0 $ for all $i \in \scrI_C$ and $\diff\alpha_{\Delta_i(gh_g)} (w_0) > 0 $ for some $i\in \scrI_C$.    
    For $\ep>0$ small enough, the polygon $\Poly^{\bmH}(gh_g) \prec \Poly^{\bmH}(gh_g \exp(\ep w_0)) $ strictly. This is a contradiction against the choice of $h_g$.
\end{proof}

\subsubsection{Weights}\label{subsubsection_appendix_divergence_weights}

By passing to a subsequence, assume without loss of generality that for all $n\in \Z^+$, $r_{g_nh_{g_n}}=r$ for some $r$, $\type(g_nh_{g_n}) = \type(g_1 h_{g_1})$ and $\alpha_{\Delta_k(g_nh_{g_n})}$ remains constant for each $k$.
So there exist $(\lambda_n) \subset \SL_N(\Z)$ and 
$\calI=\left\{ 0=i_0 <i_1 <...< i_{r+1}=N \right\} = \type(g_nh_{g_n})$ with $\calE^{\bmH}
(g_nh_{g_n})  = \lambda_n.\calF_{\calI}$ for every $n\in \Z^+$.
Hence
\begin{equation*}
    \calF_{\calI} = \lambda_n^{-1} \calE^{\bmH}
(g_nh_{g_n})  =\calE^{ \lambda_n^{-1} \bmH \lambda_n}
(g_nh_{g_n}\lambda_n).
\end{equation*}
Also note that $\lambda_n^{-1} \bmH \lambda_n \subset \bmP_{\calI}$ for all $n$.

Now that $\calI$ is fixed, it is safe to adopt the notation from Section \ref{subsubsection_quotient_lattices} with $l=r+1$.
Using the adjoint action, define $\Prim^{\bmH}(\fraksl_N(\Z))$ to be the primitive subgroups of $\fraksl_N(\Z)$ whose $\C$-spans are $\Ad(\bmH)$-invariant.
For $\alpha \in \frakX_{\Q}^*(\bmH) $, define $\Prim_{\alpha}^{\bmH}(\fraksl_N(\Z))$ similarly to $\Prim_{\alpha}^{\bmH}(\Z^N)$ and let $\calV_{\alpha,n}^{\Ad}(\Z) := \lambda_n^{-1} \Prim_{\alpha}^{\bmH}(\fraksl_N(\Z))$. To make notations uniform, define 
\begin{equation*}
    \begin{aligned}
         &
         \calV_{\alpha,n}^{\Std}(\Z) :=  \lambda_n^{-1} \Prim_{\alpha}^{\bmH}(\Z^N),\quad
         \calV_{\alpha,n}(\Z):= \calV_{\alpha,n}^{\Ad}(\Z)\bigsqcup \calV_{\alpha,n}^{\Std}(\Z), \quad
        \calV_{\alpha}(\Z):= \lambda_n. \calV_{\alpha,n}(\Z);\\
         &
         \Phi^{\Ad}_{\bmH}:= \left\{ \diffalpha : \fraks^{\frakh}_{\R} \to \R \midd
      \calV_{\alpha,n}^{\Ad}(\Z) \neq \{\{\bmzero\} \}
    \right\},\quad
    \Phi_{\bmH}:=\Phi^{\Std}_{\bmH} \bigcup \Phi^{\Ad}_{\bmH}.
    \end{aligned}
\end{equation*}

\subsubsection{Nondivergence}
Using $(\bmP_{\calI},\rmK')$ coordinates, write 
$
    g_n h_{g_n} \lambda_n  = k_n a_n p_n.
$
 \begin{lem}\label{lemma_Siegel_set_with_respect_to_H}
     For every $n$ and $\alpha\in \Delta(\rmA_{\bmP_{\calI}}, \bmP_{\calI})$, one has $\alpha(a_n) \leq  1$. Moreover, for every nonempty bounded open subset $\calO_{\rmH}$, there exists a compact subset $\scrC$ of $\rmP_{\calI}/\rmP_{\calI}\cap \SL_N(\Z)$ such that  $ p_n \lambda_n^{-1}[\calO_{\rmH}]\lambda_n \cap \scrC \neq \emptyset$ for all $n$.
 \end{lem}

 Here $[\calO_{\rmH}]$ denotes the image of $\calO_{\rmH}$ in $\rmP_{\calE^{\bmH}(g_nh_{g_n})}/\rmP_{\calE^{\bmH}(g_nh_{g_n})} \cap \SL_N(\Z)$ and $\lambda_n^{-1}[\calO_{\rmH}]\lambda_n$ denotes the further image under the map from $\rmP_{\calE^{\bmH}(g_nh_{g_n})}/\rmP_{\calE^{\bmH}(g_nh_{g_n})} \cap \SL_N(\Z) $ to $\rmP_{\calI}/\rmP_{\calI}\cap \SL_N(\Z)$ induced by $x\mapsto \lambda_n^{-1} x \lambda_n$.

 \begin{proof}[Proof of the first part of Lemma \ref{lemma_Siegel_set_with_respect_to_H}]
    For $\alpha \in \Delta(\rmA_{\bmP_{\calI}},\bmP_{\calI})$, find $k=1,...,r+1$ such that $\alpha (a_n) = \alpha_k(a_n)$ by Equa.(\ref{equation_root_lattice}). Then
    \begin{equation*}
    \begin{aligned}
         \alpha(a_n) =& 
        \varphi_k(a_n) \varphi_{k+1}(a_n)^{-1}
         =
         \left(
           \frac{
           \norm{a_n \Lambda^{\Std}_{k+1} }
           }{
           \norm{a_n \Lambda^{\Std}_{k} }
           }
         \right)^{-j_{k+1}^{-1}}
         \cdot 
         \left(
           \frac{
           \norm{a_n \Lambda^{\Std}_{k} }
           }{
           \norm{a_n \Lambda^{\Std}_{k-1} }
           }
         \right)^{j_{k}^{-1}} 
         \\
         =&
         \left(
           \frac{
           \norm{g_nh_{g_n}\lambda_n \Lambda^{\Std}_{k+1} }
           }{
           \norm{g_nh_{g_n}\lambda_n \Lambda^{\Std}_{k} }
           }
         \right)^{-j_{k+1}^{-1}}
         \cdot 
         \left(
           \frac{
           \norm{g_nh_{g_n}\lambda_n \Lambda^{\Std}_{k} }
           }{
           \norm{g_nh_{g_n}\lambda_n \Lambda^{\Std}_{k-1} }
           }
         \right)^{j_{k}^{-1}} 
         \\
         =&
         \exp \Bigg(
             - \frac{
             \ln \norm{ 
             g_n h_{g_n} \Delta_{k+1}(g_nh_{g_n})
             }
             -   \ln \norm{ 
             g_n h_{g_n} \Delta_{k}(g_nh_{g_n})
             }
             }{
             j_{k+1}
             }  \\
             &+
             \frac{
             \ln \norm{ 
             g_n h_{g_n} \Delta_{k}(g_nh_{g_n})
             }
             -   \ln \norm{
             g_n h_{g_n} \Delta_{k-1}(g_nh_{g_n})
             }
             }{
             j_{k}
             }
         \Bigg) \\
         =& \exp(-d_{k}(g_nh_{g_n})) \leq 1.
    \end{aligned}
    \end{equation*}
   The last inequality follows from the convexity of the polygons.
 \end{proof}

To prove the second part, recall the following  nondivergence criterion proved in \cite[Theorem 4.1]{zhangrunlinDcds2022} based on Kleinbock-Margulis \cite{KleMar98} and some geometry of numbers.
\begin{thm}\label{theorem_nondivergence_criterion_SL_N}
    Fix some nonempty bounded open subset $\calO_{\rmH}$ of $\rmH$ and a smooth probability measure $\rmm_{\calO}$ on $[\calO_{\rmH}]$.
    Also fix $k\in \{1,...,r\}$ and some $\ep\in (0,1)$.
    If there exists $\delta>0$ such that for every $n\in \Z^+$, one has 
    \begin{equation*}
        \norm{  \pi_k(p_n)\Delta} \geq \delta,\quad \forall\,
        \Delta \in \Prim^{\pi_k(\lambda_n^{-1} \bmH \lambda_n)} (\olV_k(\Z)),
    \end{equation*}
    then there exists a compact subset $\scrC_k \subset \SL( \olV_k)/\SL( \olV_k(\Z))$ such that 
    \begin{equation*}
        \rmm_{\calO}\left\{
           x\in [\calO_{\rmH}] \midd
           \olpi_k(p_n \lambda_n x \lambda_n^{-1}) \notin \scrC_k
        \right\}  < \ep.
    \end{equation*}
\end{thm}

\begin{proof}[Proof of the second part of Lemma \ref{lemma_Siegel_set_with_respect_to_H}]
    By rigidity of diagonalizable groups (see \cite[3.2.8]{Spr98}), 
    \[
    \left\{
      \olp^{\spl}( p_n\lambda_n^{-1} [\calO_{\rmH}] \lambda_n )=
       \olp^{\spl}( \lambda_n^{-1} [\calO_{\rmH}] \lambda_n )  \midd n\in \Z^+
    \right\}
    \]
    consists of finitely many bounded subsets of $\rmS_{\bmP_{\calI}}$ and thus remains bounded.
   In light of Lemma \ref{lemma_flag_lattice} and Theorem \ref{theorem_nondivergence_criterion_SL_N}, it suffices  to show that
    \begin{equation*}
        \norm{\pi_k(p_n)\Delta} \geq \delta,\quad \forall\,
        \Delta \in \Prim^{\pi_k(\lambda_n^{-1} \bmH \lambda_n)} (\olV_k(\Z)),
    \end{equation*}
    for some $\delta>0$ independent of $k,n$.
    Fix such a $\Delta$, let $\wtDelta \in \Prim^{\lambda_n^{-1} \bmH \lambda_n}(\Z^N)$ be the unique element such that
    \begin{equation*}
        \Lambda^{\Std}_{k-1} \subset \wtDelta \subset \Lambda_k^{\Std},\quad
        \Delta = 
        {\wtDelta/ \Lambda^{\Std}_{k-1}}.
    \end{equation*}
     By assumption we have
     \begin{equation*}
         \left(
         \rank(\wtDelta), \log \norm{
          g_n h_{g_n}\lambda_n \wtDelta
         }
         \right)  \in \Plot^{\lambda_n^{-1}\bmH\lambda_n}(g_nh_{g_n}\lambda_n) = \Plot^{\bmH}(g_nh_{g_n})
     \end{equation*}
     Therefore,
     \begin{equation*}
         \begin{aligned}
             \frac{
             \log\norm{ \pi_k(p_n) \Delta} 
             }{
             \rank(\wtDelta) - i_{k-1}
             } 
             =&\frac{
             \log\norm{p_n \wtDelta} -\log \norm{ p_n \Lambda_{k-1}^{\Std}}
             }{
             \rank(\wtDelta) - i_{k-1}
             }  \\
             =&
             \frac{
             \log\norm{k_na_np_n \wtDelta} -\log \norm{ k_na_np_n \Lambda_{k-1}^{\Std}}
             }{
             \rank(\wtDelta) - i_{k-1}
             } - \log(\varphi_{k}(a_n) )  \\
             =&
             \frac{
             \log\norm{k_na_np_n \wtDelta} -\log \norm{ g_nh_{g_n} \Delta_{k-1}(g_nh_{g_n})    }
             }{
             \rank(\wtDelta) - i_{k-1}
             } - \log(\varphi_{k}(a_n) )  
             \\
             \geq &
             \frac{
             \log\norm{ g_nh_{g_n} \Delta_{k}(g_nh_{g_n})    }
             -\log \norm{ g_nh_{g_n} \Delta_{k-1}(g_nh_{g_n})    }
             }{
             j_k
             } - \log(\varphi_{k}(a_n) ) \\
             =&
             \frac{
             \log\norm{ p_n \Lambda_{k}^{\Std}} -\log \norm{ p_n \Lambda_{k-1}^{\Std}}
             }{
             j_k
             } =0.
         \end{aligned}
     \end{equation*}
     That is to say, $\norm{ \pi_k(p_n) \Delta}  \geq 1$ and we are done.
\end{proof}

\subsubsection{Deep in the polytopes}

Let $\calV^*_{\alpha,n}(\Z)  :=  \calV_{\alpha,n}(\Z) \setminus \{\{ \bmzero \} \}$ and 
\begin{equation*}
    M_{\alpha,n}:= \inf_{\Delta \in \calV^*_{\alpha,n}(\Z)} \log \norm{p_n\Delta},\quad \diffalpha\in \Phi_{\bmH}.
\end{equation*}
By Lemma \ref{lemma_Siegel_set_with_respect_to_H}, $M_{\alpha,n}>C_2$ for some $C_2\in \R$ and for all $\diffalpha\in \Phi_{\bmH}$ and $n$.
Passing to a subsequence, assume that for each $\diffalpha \in \Phi_{\bmH}$, either $(M_{\alpha,n})$ remains bounded or diverges to $+\infty$.
Let $\Phi(\bdd)$ correspond to those bounded ones and $\Phi_{\infty}$ collect the rest. Let $M_n(\infty):= \inf_{\diff \alpha\in \Phi_{\infty}} M_{\alpha, n}$, which diverges to $+\infty$. Passing to a further subsequence, assume $M_n(\infty)$ is positive for each $n$.
Let 
\begin{equation*}
    \begin{aligned}
        \Phi_0:= \left\{
            \diffalpha \in \Phi(\bdd) \midd
            -\diffalpha \in \Cone\left(\Phi(\bdd)\right)
        \;\right\},\;
        \Phi_1:=\Phi(\bdd) \setminus \Phi_0
    \end{aligned}
\end{equation*}
where $\Cone(-)$ refers to the $\R_{\geq 0}$-linear span of a subset.
Note that $\Cone(\Phi_0)$ is a $\R$-linear subspace and the convex hull of $\Phi_1$ is compact and disjoint from $\Cone(\Phi_0)$.
So there exists $w_0 \in \fraks^{\frakh}_{\R}$ such that 
\begin{equation*}
    \diffalpha(w_0)\, 
    \begin{cases}
       \; =0 \quad & \forall \, \diffalpha \in \Phi_0 \\
        \; >0 \quad & \forall \, \diffalpha \in \Phi_1
    \end{cases}.
\end{equation*}
Also, fix $C_3>0$ such that 
\begin{equation*}
    \normm{\diffalpha(w_0)} < C_3 ,\quad
    \forall\, \diffalpha \in \Phi_{\bmH}.
\end{equation*}
Define a sequence of positive numbers $(\kappa_n)$ by
\begin{equation*}
    \kappa_n:= \frac{M_{n}(\infty)}{2C_3}.
\end{equation*}
Using $(\bmP_{\calI},\rmK')$ coordinates, write
\begin{equation*}
    g_nh'_{g_n} \lambda_n=
    k'_n a_n'p_n',\quad \text{where }\; h_{g_n}':= h_{g_n}\exp(\kappa_n w_0).
\end{equation*}
Then\footnote{Strictly speaking $k_n=k_n'$ is only true modulo $\rmK'\cap \rmM_{\bmP,\rmK'}$. As this does not affect the proof, we choose to ignore it.}
\begin{equation*}
    k_n=k_n',\quad a_n p_n \cdot ( \lambda_n^{-1} \exp(\kappa_n w_0) \lambda_n) = a_n'p_n'.
\end{equation*}
Moreover, by Lemma \ref{lemma_singular_canonical_weights_levels}, the character associated to $\Delta_{k}(g_nh_{g_n})=\lambda_n \Lambda^{\Std}_k$ belongs to $\Phi_0$ for each $k=1,...,r$. We have
\begin{equation}\label{equation_appendix_divergence_deep_polytope_w_0_fix_parabolic}
   \lambda_n^{-1} \exp(\kappa_n w_0) \lambda_n. \Lambda^{\Std}_k
   =\Lambda^{\Std}_k.
\end{equation}
Hence $\lambda_n^{-1} \exp(\kappa_n w_0) \lambda_n \in {}^{\circ}\rmP_{\calI}$ and
\begin{equation*}
    a_n'=a_n ,\quad
    p_n' = p_n \cdot ( \lambda_n^{-1} \exp(\kappa_n w_0) \lambda_n).
\end{equation*}
For $\diffalpha \in \Phi_{\bmH}$ and $\Delta \in \calV_{\alpha,n}(\Z)$, we have
\begin{equation}\label{equation_appendix_divergence_deep_polytope_1}
    \log\norm{  p_n' \Delta}
    =\kappa_n\diffalpha(w_0) + \log \norm{ p_n\Delta}.
\end{equation}
We have already seen by Equa.(\ref{equation_appendix_divergence_deep_polytope_w_0_fix_parabolic}) that:
\begin{equation*}
    \log\norm{g_nh_{g_n}'\lambda_n  \Lambda^{\Std}_k }
    = \log\norm{g_nh_{g_n}\lambda_n  \Lambda^{\Std}_k }.
\end{equation*}
Next we claim that 
\begin{claim}
Let $\Delta \in \lambda_n^{-1} \Prim^{\bmH}(\Z^N) $.
The point $\left(\rank(\Delta), \log \norm{g_nh_{g_n}'\lambda_n \Delta } \right)$ lies above the piecewise-linear line connecting  $\left( \rank(\Lambda^{\Std}_{k}), \log \norm{g_nh_{g_n}'\lambda_n  \Lambda_k^{\Std} } \right)$ as $k$ goes from $0$ to $r+1$.
\end{claim}
\begin{proof}[Proof of the claim]
      First we assume that $\Delta$ is compatible with $\calF_{\calI}$, namely $\Lambda^{\Std}_{k-1} \subset \Delta \subset\Lambda^{\Std}_{k}$ for some $k$, and prove
\begin{equation}\label{equation_appendix_divergence_deep_polytope_2}
    \frac{
       \log \norm{g_nh_{g_n}'\lambda_n \Delta }
       -  \log\norm{g_nh_{g_n}'\lambda_n  \Lambda_{k-1}^{\Std} }
     }{ \rank(\Delta) - i_{k-1}}
     \geq \frac{
       \log\norm{g_nh_{g_n}'\lambda_n\Lambda_{k}^{\Std} }
       -  \log\norm{g_nh_{g_n}'\lambda_n  \Lambda_{k-1}^{\Std} }
     }{ i_k - i_{k-1}}.
\end{equation}
Since  $\Delta$ is assumed to be compatible with $\calF_{\calI}$, this is equivalent to
\begin{equation*}
    \frac{
       \log\norm{p_n' \Delta }
       -  \log\norm{p_n'  \Lambda_{k-1}^{\Std} }
     }{ \rank(\Delta) - i_{k-1}}
     \geq \frac{
       \log\norm{p_n' \Lambda_{k}^{\Std} }
       -  \log\norm{p_n' \Lambda_{k-1}^{\Std} }
     }{ i_k - i_{k-1}},
\end{equation*}
which, thanks to Equa.(\ref{equation_appendix_divergence_deep_polytope_1}), is again equivalent to
\begin{equation}\label{equation_appendix_divergence_deep_polytope_3}
    { \kappa_n \diffalpha(w_0) +
       \log\norm{p_n \Delta }
     }
     \geq 0.
\end{equation}
For $\diff\alpha \in \Phi(\bdd)$,
\begin{equation*}
    { \kappa_n \diffalpha(w_0) +
       \log\norm{p_n \Delta }
     }
     \geq \log\norm{p_n \Delta } \geq 0.
\end{equation*}
For $\diff\alpha \in \Phi_{\infty}$,
\begin{equation*}
    { \kappa_n \diffalpha(w_0) +
       \log\norm{p_n \Delta }
     }
     \geq \log\norm{p_n \Delta } - C_3 \kappa_n \geq 0.
\end{equation*}
The truth of  Equa.(\ref{equation_appendix_divergence_deep_polytope_3}) and hence of Equa.(\ref{equation_appendix_divergence_deep_polytope_2}) is verified.

Now we drop the compatibility assumption on $\Delta$.
Fix some $m< l$ such that
 $\Delta$ is contained in  $\Lambda^{\Std}_l$ but not $\Lambda^{\Std}_{l-1}$;
$\Delta$ contains  $\Lambda^{\Std}_m$ but not $\Lambda^{\Std}_{m+1}$.
Let $\Delta_j := \Delta \cap \Lambda_{m+j} $ as $j$ ranges over $\{0,1,....,l-m\}$ and $i'_j:= \rank (\Delta_j) $. So $\Delta_0= \Lambda_m$ and $\Delta_{l-m}= \Delta$. 
We are going to explain that  for $j=0,1,...,l-m-1$,
\begin{equation}\label{equation_appendix_divergence_deep_polytope_4}
    \frac{
       \log \norm{g_nh_{g_n}'\lambda_n \Delta_{j+1} }
       -  \log\norm{g_nh_{g_n}'\lambda_n \Delta_{j} }
     }{ i'_{j+1} - i'_j }
     \geq \frac{
       \log\norm{g_nh_{g_n}'\lambda_n\Lambda_{m+j+1}^{\Std} }
       -  \log\norm{g_nh_{g_n}'\lambda_n  \Lambda_{m+j}^{\Std} }
     }{ i_{m+j+1} - i_{m+j} },
\end{equation}
which is sufficient to conclude the proof since $ i'_{j+1} - i'_j  \leq  i_{m+j+1} - i_{m+j}$.
As the argument has nothing to do with the explicit form of $g_nh_{g_n}' \lambda_n$, we abbreviate $\norm{g_nh_{g_n}' \lambda_n . \Lambda } =: \norm{\Lambda}_g$ for every discrete subgroup $\Lambda$.
For two primitive subgroups $\Lambda_1,\Lambda_2 $ of $\Z^N$, let 
\[
   \Lambda_1 +' \Lambda_2 :=  (\Lambda_1 +  \Lambda_2) \otimes  \Q \cap \Z^N,
\]
which is again primitive.
We have that (see \cite[Corollary 4.2]{Casselman_2004})
\begin{equation*}
\begin{aligned}
    & \frac{
    \norm{ \Delta_{j+1} }_g
    }{   \norm{\Delta_j}_g
    }  \geq 
    \frac{
    \norm{ \Delta+' \Lambda^{\Std}_{m+j}  }_g
    }{ 
    \norm{ \Lambda^{\Std}_{m+j}}_g
     } \\
    \implies &
    \frac{
       \log \norm{\Delta_{j+1} }_g
       -  \log\norm{ \Delta_{j} }_g
     }{ i'_{j+1} - i'_{j}}
     \geq \frac{
       \log\norm{ \Delta+' \Lambda^{\Std}_{m+j}  }_g
       -  \log \norm{ \Lambda^{\Std}_{m+j}}_g
     }{  \rank ( \Delta+' \Lambda^{\Std}_{m+j} ) - i_{m+j} }.
\end{aligned}
\end{equation*}
Equa.(\ref{equation_appendix_divergence_deep_polytope_4}) then follows by
applying Equa.(\ref{equation_appendix_divergence_deep_polytope_2}) to $\Delta+' \Lambda^{\Std}_{m+j}$, which is indeed compatible with $\calF_{\calI}$.

\end{proof}
By the claim and Equa.(\ref{equation_appendix_divergence_deep_polytope_w_0_fix_parabolic}),
\begin{equation*}
    \begin{aligned}
        &\Poly^{\lambda_n^{-1} \bmH \lambda_n}({g_n h_{g_n}\lambda_n})
        =  \Poly^{\lambda_n^{-1} \bmH \lambda_n}({g_n h'_{g_n}\lambda_n}) \\
        &
        \calE^{\lambda_n^{-1} \bmH \lambda_n}({g_n h_{g_n}\lambda_n})
        =  \calE^{\lambda_n^{-1} \bmH \lambda_n}({g_n h'_{g_n}\lambda_n}) =
        \left\{
          \{0 \} \subsetneq \Lambda_1^{\Std} \subsetneq ... \subsetneq \Lambda_{r+1}^{\Std} = \Z^N
        \right\}.
    \end{aligned}
\end{equation*}
Also for $\diff \alpha \in \Phi_1$ and $\Delta \in \calV_{\alpha,n}(\Z)$,
\begin{equation*} 
    \log\norm{p_n'\Delta} \geq \kappa_n \diffalpha(w_0) + C_2
\end{equation*}
diverges to $+\infty$.  
Similarly, $ \log\norm{p_n'\Delta}$ with $\Delta \in \calV_{\alpha,n}(\Z)$ has a common lower bound for other $\diff\alpha \in \Phi_{\bmH}$.
Repeating the proof of the second part of Lemma \ref{lemma_Siegel_set_with_respect_to_H}, we have that
\begin{equation*}
    \left( p_n' \lambda_n^{-1}[\calO_{\rmH}]\lambda_n \right) \text{ is nondivergent in }\rmP_{\calI}/\rmP_{\calI} \cap \SL_N(\Z).
\end{equation*}
As a consequence, we find some
$(\lambda_n')\subset \rmP_{\calI}\cap \SL_N(\Z)$, bounded sequence $(b_n)$ in $\rmP_{\calI}$ and another
bounded sequence $(\omega_n)$ in $\calO_{\rmH}$ such that 
\begin{equation*}
    p_n' \lambda_n^{-1} \omega_n \lambda_n = b_n \lambda_n'.
\end{equation*}
Letting $h_{g_n}'':= h_{g_n}'\omega_n$ and $\gamma_n:=\lambda'_n \lambda_n^{-1}$, we have
\begin{equation*}
    g_nh_{g_n}'' \gamma_n^{-1}
    =
    k_n a_n b_n, \quad 
    \gamma_n \bmH \gamma_n^{-1} \subset \bmP_{\calI} \text{ fixes }\bmv^{\calI}.
\end{equation*}
Rewrite $k_na_nb_n=k_na_n''b_n''$ using ($\bmP_{\calI},\rmK'$) coordinates, then $(a''_n)$ (resp. $(b''_n)$) is bounded away from $(a_n)$ (resp. $(b_n)$).
In particular, $(b_n'')$ is bounded and there exists $C_4 >0$ such that $\alpha(a_n'')< C_4 $ for all 
$\alpha \in \Delta(   \rmA_{ \bmP_{\calI} }, \bmP_{\calI})$. 
For $\diffalpha \in \Phi_{\bmH}$, let $M'_{\alpha,n}:= \inf\left\{ 
    \log\norm{p_n'\Delta} \midd \Delta \in \calV^*_{\alpha,n}(\Z)   \right\}$.
By what we know about $(M_{\alpha,n})$ and Equa.(\ref{equation_appendix_divergence_deep_polytope_1}),
$ (M'_{\alpha,n}) $ diverges to $+\infty$ for $\alpha \in \Phi_{\infty} \cup \Phi_1$ and $(M'_{\alpha,n}=M_{\alpha,n})$ remains bounded for $\alpha \in \Phi_0$.
In particular, if $(M'_{\alpha,n})$ is bounded, then $\sum_{i=0}^m a_i \diffalpha_i  = 0$
for some $(a_i)_{i=0}^m\subset \R^+$  and $(\diffalpha_i)_{i=0}^m \subset \Phi_{\bmH}$ with $\alpha_0=\alpha$ and $(M'_{\alpha_i,n})$ bounded for each $i$. As $\alpha_i$'s are all $\Q$-characters, we may take $(a_i)$'s to be positive integers.
Let 
\begin{equation*}
    M''_{\alpha,n}:= \inf\left\{ 
    \log\norm{\gamma_n\Delta} \midd \Delta \in \calV^*_{\alpha}(\Z)   \right\}.
\end{equation*}
Since 
\begin{equation*}
    \log \norm{\gamma_n \Delta }=  \log \norm{b_n^{-1} p_n' \lambda_n^{-1} \omega_n \Delta} \text{ is bounded from }
    \log \norm{ p_n'\lambda_n^{-1} \Delta},
\end{equation*}
the divergence/boundedness of $(M'_{\alpha,n})$ is the same as that of $(M''_{\alpha,n})$.
\begin{lem}
    Assume that $\left( \gamma_n \bmH \gamma_n^{-1} \right)$ has a subsequence converging to $\bmL$. Then $\bmL$ is observable in $\SL_N$.
\end{lem}
\begin{proof}
    It follows from the proof of \cite[Lemma 4.10]{zhangrunlinCompositio2021}. Let us briefly recall how.
    Denote by $\bmO$ the observable hull of $\bmL$ in $\SL_N$. It suffices to show that $\bmL$ is normalized by $\bmO$.
    Assume $\Delta_{\bmL} \in \Prim(\fraksl_N(\Z))$ spans the Lie algebra of $\bmL$. By abuse of notation, also view $\Delta_{\bmL}$ as a vector in the appropriate wedged vector space. Then we only need to show that the line spanned by $\Delta_{\bmL}$ is preserved by $\bmO$.
    Passing to a subsequence, assume for some fixed $\diffalpha\in \Phi_{\bmH}$, $\gamma_n^{-1} \Delta_{\bmL} \in \calV_{\alpha}(\Z)$ for all $n\in \Z^+$. But it is a tautology that  
    \[
     \left( \log \norm{\gamma_n.(\gamma_n^{-1}\Delta_{\bmL})} \right)
     \text{ is bounded.}
    \]
    So we can find $(q_i)\subset \Z^+$, $(\diffalpha_i)\in \Phi_{\bmH}$ and $(\Delta_i ) \subset \gamma_n.\calV_{\alpha_i}(\Z)$ with $\Delta_0 = \Delta_{\bmL}$ and $\sum q_i \diffalpha_i = 0$. Therefore,
    \begin{equation*}
        \begin{aligned}
            & \bigotimes \Delta_{i}^{\otimes q_i} \text{ is fixed by }\gamma_n \bmH \gamma_n^{-1} \text{ for every }n \\
            \implies &
           \bigotimes \Delta_{i}^{\otimes q_i} \text{ is fixed by }\bmO
           \\
            \implies &
            \text{ the line spanned by }\Delta_{0}=\Delta_{\bmL} \text{ is fixed by }\bmO.
        \end{aligned}
    \end{equation*}
    So we are done.
\end{proof}

Finally, passing to a subsequence, assume that for each $\alpha \in \Delta( \rmA_{\bmP_{\calI}}, \bmP_{\calI})$, either $\left( \alpha(a_n'') \right)$ remains bounded away from $0$ and $\infty$ or converges to $0$. Let $J$ collect the bounded ones. 
For simplicity let $\bmP_J:= (\bmP_{\calI})_J$.
If we write 
\begin{equation*}
    g_n h_{g_n}'' \gamma_n^{-1} = k_n^J a_n^J b_n^J
\end{equation*}
using $(\bmP_J ,\rmK')$ coordinates, then $(b_n^J)$ is bounded and $\alpha(a_n^J) \to 0$ for all $\alpha \in \Delta(\rmA_{\bmP_J},\bmP_J)$.

To conclude proof, we note that by Equa.(\ref{equation_root_lattice}), 
\begin{equation*}
    J \cong \left\{ 
     k=1,...,l-1 \midd \alpha_k(a_n)  \nrightarrow 0
    \right\}
    =\left\{ 
       k=1,...,l-1 \midd -\log \alpha_k(a_n)  \nrightarrow +\infty
    \right\}.
\end{equation*}
But 
\begin{equation*}
    \begin{aligned}
        &  -\log \alpha_k(a_n)  = \log \varphi_{k+1}(a_n) - \log\varphi_k(a_n)  
        \\
        =\,&
        \frac{
        \log \norm{a_n \Lambda_{k+1}^{\Std} } - \log \norm{a_n \Lambda_{k}^{\Std} }
        }{
        i_{k+1} - i_{k}
        } - 
        \frac{
        \log \norm{a_n \Lambda_{k}^{\Std} } - \log \norm{a_n \Lambda_{k-1}^{\Std} }
        }{
        i_{k} - i_{k-1}
        } 
        \\
        =\,&
        \frac{
        \log \norm{g_n h_{g_n}\lambda_n \Lambda_{k+1}^{\Std} } - \log \norm{g_n h_{g_n}\lambda_n \Lambda_{k}^{\Std} }
        }{
        i_{k+1} - i_{k}
        } - 
        \frac{
        \log \norm{g_n h_{g_n}\lambda_n \Lambda_{k}^{\Std} } - \log \norm{g_n h_{g_n}\lambda_n \Lambda_{k-1}^{\Std} }
        }{
        i_{k} - i_{k-1}
        } 
        \\
        =\,& d_k(g_nh_{g_n}).
    \end{aligned}
\end{equation*}
Choose $C_5 \in \R$ such that 
\[
    d_k(g_n h_{g_n}) \leq C_5,\quad \forall \, n,\; \forall \, k  \in J.
\]
Passing to a subsequence we may assume that $ d_k(g_n h_{g_n}) > C_5$ for all $n$ and $k\notin J$.
By Lemma \ref{lemma_singular_canonical_weights_levels} applied to $C=C_5$, we find $(a_k)_{k\in \{1,...,l-1\}\setminus J} \subset \Z^+$ such that $\sum_{k\notin J} a_k \diff\alpha_{\Delta_k(g_nh_{g_n})} = 0$ (recall we had assumed that $\alpha_{\Delta_k(g_nh_{g_n})}$ is independent of $n$ for each $k$).
Therefore, $\sum_{k \notin J} a_k \diff f_k  =0$.  By abuse of notation, treat $\Lambda_k^{\Std}$ as a vector for every $k$.
If $\bmv:= \otimes_{k \notin J} (\Lambda_k^{\Std})^{\otimes a_k}$, then $\bmv$ is fixed by $\gamma_n \bmH \gamma_n^{-1}$ for all $n$.
Note that $\Lambda_k^{\Std}$ for $k\notin J$ and in particular $\bmv$ are nontrivial eigenvectors of $\bmP_{J}$.
The proof of Theorem \ref{theorem_divergence_observable} in the $\SL_N$ case is thus complete with $\bmP= \bmP_J$.

\section{Indices of Various Conditions}

\begin{itemize}
     \item[($\rmB 1$)]\label{condition_B_1_appendix} For every $x\in \bmU(\Q)$ and $(g_n)\subset \rmG$ such that $\lim g_n.x \in \bmB(\R)$, one has 
        \begin{equation*}
            \lim_{n \to \infty} \left[
            (g_n)_* \rmm_{[\rmH^{\circ}_x]}
              \right] = \left[
            \rmm_{[\rmG]}   \right].
        \end{equation*}
\end{itemize}

\begin{itemize}
    \item[{$(\rmB \rmH 1)$}]  For every $x\in \bmU(\Q)$, $\nu_{\infty,x}:=\lim_{R\to \infty} \nu_{R,x}$ exists and $\supp(\nu_{\infty,x}) \subset \bmB(\R)$.
\end{itemize}

\begin{itemize}
     \item[($\rmB 2$)] For every $x\in \bmU(\Q)$ and $y \in \overline{\rmG.x}$, one has $\overline{\rmG.y}\cap \bmB(\R)\neq \emptyset$.
\end{itemize}

\begin{itemize}
    \item[($\rmK 1$)]\label{condition_K_1_appendix} For every $\alpha \in \scrA_{\R,\Q}^{\an}$, $d_{\alpha} -1 >0$.
\end{itemize}

\begin{itemize}
     \item[$(\rmB 4)$]\label{condition_B_4_appendix} For every $x\in \bmU(\Q)$, $d_{\alpha} -1 >0$ for some $\alpha \in \scrB_{\R,x}^{\an}$. In particular, $\scrB_{\R,x}^{\an}$ is nonempty.
\end{itemize}

\begin{itemize}
    \item[$(\rmH 1)$]\label{condition_H_1_appendix} For every $x\in \bmU(\Q)$ and arithmetic subgroup $\Gamma \subset \bmG(\Q)\cap \rmG$, $\bmH_x^{\circ}$ has no nontrivial $\Q$-characters and $\lim_{R\to \infty} \mu_{R,x} $ exists in $\Prob(\rmG/\Gamma)$.
\end{itemize}

\begin{itemize}
    \item[$(\rmD 1)$] \label{condition_D_1_appendix}
    For every $x\in \bmU(\Q)$ and arithmetic subgroup $\Gamma \subset \bmG(\Q)\cap \rmG$, $\bmH_x^{\circ}$ has no nontrivial $\Q$-characters and $\Psi_x$
    extends continuously to $\overline{\rmG.x}  \to  \Prob(\rmG/\Gamma) \cup \{ \bmzero \}$ where the closure is taken in $\rmX^{\cor}$, the manifold with corners associated with $(\bmX,\bmD)$ (see Section \ref{subsection_manifolds_corners}).
    \item[$(\rmS 1)$]\label{condition_S_1_appendix} For every $x\in \bmU(\Q)$ and arithmetic subgroup $\Gamma$ in $\bmG(\Q)\cap \rmG$,  there exists a bounded subset $B$ of $\rmG/\Gamma$ such that $ g\rmH_x^{\circ}\Gamma/\Gamma$ intersects with $B$ for every $g\in \rmG$.
\end{itemize}

\begin{itemize}
    \item[$(\rmD \rmS 1)$]\label{condition_DS_1_appendix}
    There exists a closed subset $D$ of $\bmB(\R)$ such that for every $(g_n)$ and every $x\in \bmU(\Q)$ with $\lim g_n.x $ not in $D$,
    there exists a bounded subset $B$ of $\rmG/\Gamma$ such that $g_n\rmH_x^{\circ}\Gamma/\Gamma$ intersects with $B$ for every $n$.
    Moreover, for every $F\subset \scrC^{\an}_{\R,\Q}$, if $\Leb_F$ denotes a smooth measure on $D_F$, then $\Leb_F(B)=0$.
\end{itemize}

\begin{itemize}
    \item[$(\rmH 2)$]\label{condition_H_2_appendix} For every $x\in \bmU(\Q)$ and arithmetic subgroup $\Gamma \subset \bmG(\Q)\cap \rmG$,  there exists a non-negative compactly supported continuous function $\psi$ on $\rmG/\Gamma$ such that  $\la \psi, g_*\rmm_{[\rmH_x^{\circ}]} \ra \neq 0$
    for all $g \in \rmG$
    and for every such $\psi$, $\lim_{R\to \infty} \mu^{\psi}_{R,x} $ exists in $\Prob^{\psi}(\rmG/\Gamma)$.
\end{itemize}

\begin{itemize}
    \item[$(\rmD 2)$]\label{condition_D_2_appendix} For every $x\in \bmU(\Q)$ and arithmetic subgroup $\Gamma \subset \bmG(\Q)\cap \rmG$, there exists $\psi$ satisfying the paragraph above  and $\Psi^{\psi}_x$
    extends continuously to $\overline{\rmG.x}  \to  \Prob^{\psi}(\rmG/\Gamma)$ where the closure is taken in $\rmX^{\cor}$, the manifold with corners associated with $(\bmX,\bmD)$.
\end{itemize}

\begin{itemize}
    \item[(C1)]\label{condition_C_1_appendix} $\bmZ(\bmH,\bmL)$ decomposes into finitely many orbits under the action of $\bmL \times \bmN_{\bmG}(\bmH)$ for every $\bmL \in \INT(\bmH,\bmG)$;
    \item[(N1)]\label{condition_N_1_appendix} $\Gamma \cap \bmH$ is a finite index subgroup of $\Gamma \cap \bmN_{\bmG}(\bmH)$ for an(y) arithmetic subgroup $\Gamma$.
\end{itemize}

\begin{itemize}
    \item[(F1)]\label{condition_F_1_appendix} $\INT_{\Gamma}(\bmH,\bmG)$ is finite.
\end{itemize}

\begin{itemize}
    \item[(F2)]\label{condition_F_2_appendix} there are only finitely many parabolic $\Q$-subgroups containing $\bmH$, i.e., $\scrP_{\bmH}$ is finite.
\end{itemize}

\begin{itemize}
    \item[(N2)]\label{condition_N_2_appendix} $\bmN_{\bmG}(\bmH)^{\circ}\subset \bmN_{\bmG}(\bmL)$ for every $\bmL \in \INT(\bmH,\bmG)$.
\end{itemize}

\section*{Ackonwledgement}
We are grateful to useful discussions with Jinpeng An, Yue Gao, Zhizhong Huang, Osama Khalil, Fei Xu and Pengyu Yang. We also thank Osama Khalil, Hee Oh, Fei Xu for comments and questions on an earlier version of this paper.
The author is supported by NSFC, grant \#12201013.

\bibliographystyle{amsalpha}
\bibliography{ref}

\end{document}